# The integration theory of linear ordinary differential equations

**Nezbailo T.G.**

**Notations and equalities used in the text:**

- $N$ - a set of natural numbers. - C – an arbitrary constant.

- $I$ - an imaginary unit.

- $\Gamma(n) = (n-1)!$, where $n$ is a natural number is called a complete gamma-function. $\Psi(x) = \frac{d}{dx}\ln(\Gamma(x))$

- $C_j^i = C(i,j)$ - a binomial coefficient. For the $j, i, i \le j, 0 \le i$ integers it is defined by the following formula: $C(i,j) = \frac{j!}{i!\,(j-i)!}$. For any values of $i, j$ it is defined by the following formulas: $C(i,j) = \frac{\Gamma(i+1)}{\Gamma(j+1)\,\Gamma(j-i+1)}$ or with the use of limiting relations following from this formula.

- $[f(x)]_n$ - ***n***-multiple integral of a ***f(x)*** function to the variable ***x;*** $n = 0, 1, 2$ ....

For instance, $[f(x)]_0 = f(x)$, $[f(x)]_1 = \int f(x)\,dx$, $[f(z)]_2 = \iint f(z)\,dz\,dz$, $[f(s)]_3 = \iiint f(s)\,ds\,ds\,ds$ etc.

**Introduction. Definition of the Problem.** *An indefinitely small value is always greater than nothing*

Fundamentals of the **differential and integral calculus** theory were independently outlined by **I.Newton and G.Leibniz** in **1666 - 1702** [1], and, since then, were substantially extended by many prominent mathematicians like **J.Bernoulli, L.Euler, J.Lagrange, A.Cauchy, N.I. Lobachevsky** and others. The basic fundamental scientific achievement of **differential and integral calculus** was to define **the operation of differentiation:** $\frac{d}{dx}[\,]$ for an arbitrary analytic function $f(x)$, as a mathematical operation described by the following rule:

$$\frac{d}{dx}f(x) = \lim_{\delta \to 0}\frac{f(x+\delta) - f(x)}{\delta} \qquad (1.1)$$

Opposite to the operation of **differentiation:** $\frac{d}{dx}[\,]$ is the operation of **integration:** $\int[\,]\,dx$, which defines the impact on the analytic function $f(x)$, (indefinite integral of $f(x)$) according to the following equation:

$$\int f(x)\,dx = F(x) + C \qquad (1.2)$$

where $F(x)$ - is a function, which is called **primitive function** of $f(x)$ and satisfies the following equality:

$$\frac{d}{dx}F(x) = f(x) \qquad (1.3)$$

However, (1.2) **does not define the algorithm for the operation of integration:** $\int[\,]\,dx$, since it does not indicate how, that is, with the help of what conventional mathematical operations (and in what sequence) we can arrive from a well-known analytic function $f(x)$ to its primitive function: $F(x)$. Thus, the modern **theory of differential and integral calculus <u>still has not solved a fundamental problem:</u> <u>there is no defined mathematical algorithm for the operation of integration.</u>**

Calculation of indefinite integral for a **f(x)** function, that is, definition of $\int f(x)\,dx$, is understood as a set of mathematical algorithms (methods),which use artificial methods, like the substitution of integration variable, or the use of integration formula in parts, etc. to transform the expression $\int f(x)\,dx$ in a form belonging to a known set of integrals (1.2), which satisfy the condition (1.3). If this can be done, the $\int f(x)\,dx$ integral is deemed solved; if not, the problem remains unsolved. Because of the absence of such universal general algorithm, the calculation of indefinite integral depends upon the ability (intuition) to select (guess) the desired substitution, or a known method for transformation. Obviously, there are lots (actually, an endless number) of integrals without known calculation methods, such like: $\int \sin(\sin(x))\,dx$ etc.

At the same time, as soon as the operations of differentiation and integration are opposite to one another, that is, defined by the equalities

$$\frac{d}{dx}\left(\int f(x)\,dx\right) = \int \frac{d}{dx} f(x)\,dx = f(x) \qquad (1.4)$$

this means that **these operations are somehow related, and the objective of this publication is to define this relation.**

## 2. The *n*-th derivatives.

**2.1.Definition and calculation of the derivative of *n*-th order.** Let's assume (for the clarity) that all the functions, which are described below, are rather smooth, that is ***n*** times differentiable and integrable.

**Definition. Derivative of the *n*-th order of an analytic function $f(x)$ is a $G(n, x)$ function containing an $n$ parameter (which is a natural number) and defined by the following formula:**

$$G(n, x) = \frac{d^n}{dx^n} f(x) \qquad (2.1)$$

The $G(n, x)$ function forms the following sequence

$G(1, x), G(2, x), G(3, x)$ ........

and has the following obvious **properties:**

**a) recurrence:**

$$G(n, x) = \frac{\partial}{\partial x} G(n-1, x) \qquad (2.2)$$

**Therefore:** *If the parametric function* $G(n, x)$ is an $n$ *-th derivative, then the following formula applies as well:*

$$G(n+k, x) = \frac{\partial^k}{\partial x^k} G(n, x) \qquad (2.3)$$

where $k$ is a natural number.

**b) Inversion of the order.**

$$\frac{\partial^n}{\partial x^n} G(k, x) = \frac{\partial^k}{\partial x^k} G(n, x) = G(n+k, x) \qquad (2.4)$$

where $n, k$ **are natural numbers.**

**It is obvious that a given parametric function** $S(n, x)$ shall be an ***n***-th derivative only if it satisfies the conditions (2.2) - (2.4). Since the conditions (2.3) and (2.4) follow from the condition (2.2), this means that (2.2) is **a necessary and sufficient** condition to ensure that a given parametric function $S(n, x)$ is an $n$-th derivative of a certain function. **For example,** ***the parametric function*** $S(i, x) = x^i$ *is not an* $i$ *-th derivative, because it does not satisfy the condition* (2.2). Indeed, since $S(i+1, x) = x^{i+1}$ and $\frac{\partial}{\partial x} S(i, x) = \frac{\partial}{\partial x} x^i = i\, x^{i-1}$, *and the condition* (2.2) *is not satisfied, then the initial parametric function* $S(i, x) = x^i$ *is not an* $i$ *-th derivative. At the same time,* ***the parametric function*** $S(i, x) = k^i e^{kx}$ *is an* $i$ *-th derivative, because it satisfies the condition* (2.2): ( $S(i+1, x) = k^{i+1} e^{kx}$ , $\frac{\partial}{\partial x} S(i, x) = k^{i+1} e^{kx}$ )

Since the operation of differentiation transforms the $G(n, x)$ function into $G(n+1, x)$ , the question is whether there are such functions that are not changed in such cases, that is, are invariant to the operation of differentiation.

Obviously, these functions must satisfy both conditions:

$$G(n, x) = G(n+1, x) \qquad G(n, x) = \frac{\partial}{\partial x} G(n, x)$$

By solving the recurrence and differential equation, we identify that there is only one function of such type:

$$G(n, x) = C e^x \tag{2.5}$$

where $C$ is an arbitrary constant, which does not depend on the $\boldsymbol{n}$ parameter**.**

**2.2. Methods for the definition of *n*-th derivatives.** It is obvious that the $G(n, x)$ function is, in a general case, fundamentally different from $f(x)$ , because it contains a new parameter: $n$ . Thus, the $f(x)$ function gives no means to make judgments on the $G(n, x)$ , function, and the definition of $G(n, x)$ presents, in a general case, a sophisticated heuristic problem, that is, it is required, on the basis of whatever number of sequence members: $G(1, x), G(2, x), G(3, x)$ ........ to establish the form of $G(n, x)$ function.

**Example No. 1: Calculate the *n*-th derivative of the $f(x) = e^{kx}$ function, k-const.**

**Solution:** Since:

$$e^{kx}, \quad \frac{d}{dx} f(x) = k\, e^{kx}, \ \frac{d^2}{dx^2} f(x) = k^2 e^{kx}, \ \frac{d^3}{dx^3} f(x) = k^3 e^{kx} \tag{2.6}$$

then, by induction, the derivative of ***n***-th order takes the following form:

$$\frac{\partial^n}{\partial x^n} e^{kx} = k^n e^{kx} \tag{2.7}$$

Indeed, if we sequentially substitute $n = 0, n = 1, n = 2, n = 3$ ..... from (2.7), we arrive to the sequence described by (2.6). Let's verify the condition (2.2).

$$k^{i+1} e^{kx} = \frac{\partial}{\partial x} (k^i e^{kx}) = k^{i+1} e^{kx}$$

Thus, (2.7) actually defines the ***n***-th derivative of the $f(x) = e^{kx}$ function.

**Example No. 2. Calculate the *n*-th derivative of the $f(x) = x^m$ function,** where ***m*** is a natural number.

**Solution:** In this case, the sequence of derivatives takes the following form:

$$x^m, \quad \frac{d}{dx} f(x) = m\, x^{m-1}, \ \frac{d^2}{dx^2} f(x) = m\,(m-1)\, x^{m-2}, \ \frac{d^3}{dx^3} f(x) = m\, m(-1)\,(m-2)\, x^{m-3}$$

Obviously, the ***n***-th member of this sequence is defined, with $n \le m$, by the following formula:

$$\frac{\partial^n}{\partial x^n} x^m = A_n(m)\, x^{m-n}, A_n(m) = \frac{m!}{(m-n)!} \qquad (2.8)$$

The validity of this formula is verified as described in the example No. 1.

In some cases, the ***n***-th derivative can be defined only for the $G(0, x), G(1, x), G(2, x), G(3, x)$ ........sequence limited "at the bottom" (or "from the top").

**Example No. 3. Calculate the *n*-th derivative of the following function:** $f(x) = \ln(x)$ .

Since:

$\ln(x)$ , $\frac{d}{dx} f(x) = \frac{1}{x}, \frac{d^2}{dx^2} f(x) = -\frac{1}{x^2}, \frac{d^3}{dx^3} f(x) = \frac{2}{x^3}$ and so on.

In this case:

$$\frac{d^n}{dx^n} \ln(x) = \frac{(-1)^{n-1}(n-1)!}{x^n} \qquad (2.9)$$

As we can see, in this case it is easy to define a general parametric function only for a limited sequence of $\frac{d}{dx} f(x) = \frac{1}{x}, \frac{d^2}{dx^2} f(x) = -\frac{1}{x^2}, \frac{d^3}{dx^3} f(x) = \frac{2}{x^3}$ etc. This formula is valid only for the natural values of ***n*** parameter. Indeed, it is used in this particular representation. At the same time (as presented below) for the calculation of ***n***-th derivative of $f(x)$ function, it is essential that the resulting $G(n, x)$ function should satisfy the condition defined by the following **Definition:**

**If there is a limit for the following expression at the** $n = 0$ **point**

$$\lim_{n \to 0} G(n, x) = f(x) \qquad (2.10)$$

**and for other points** $n = k$

$$\lim_{n \to k} G(n, x) = G(k, x)$$

**where** $G(k, x)$ **is a certain function, then the *n*-th derivative of** $G(n, x)$ **is called normal with** $n = k$ **.**
**Otherwise it is called a singular *n*-th derivative with** $n = k$ **.**

For example, the ***n***-th derivative (2.7) is called normal, because it satisfies the condition (2.10):

$$\left[\frac{\partial^n}{\partial x^n} e^{kx}\right]_{n=0} = \left[k^n e^{kx}\right]_{n=0} = e^{kx}$$

Also, the ***n***-th derivative (2.8) satisfies the condition (2.10):

$$\left[\frac{\partial^n}{\partial x^n} x^m\right]_{n=0} = \left[A_n(m)\, x^{m-n}\right]_{n=0} = x^m$$

(2.9) fails to satisfy the condition (2.10), so it's called a singular ***n***-th derivative. Nevertheless, there is the ***n***-th derivative of the $f(x) = \ln(x)$ **function,** which satisfies the condition (2.10).

**Statement No. 1.** The normal ***n***-th derivative of the $f(x) = \ln(x)$ function is defined by the following formula:

$$\frac{d^n}{dx^n}\ln(x) = -\frac{(-\ln(x)+\gamma+\Psi(n)+\pi\cot(\pi n))\,x^{-n}}{\Gamma(-n+1)} \qquad (2.11)$$

The proof is as follows:

$$\lim_{n\to 0}\left(-\frac{(-\ln(x)+\gamma+\Psi(n)+\pi\cot(\pi n))\,x^{-n}}{\Gamma(-n+1)}\right) = \ln(x)$$

Condition (2.10) is satisfied. In a similar manner (by way of direct calculations), we can make sure that the following equalities hold true:

$$\frac{d}{dx}f(x) = \frac{1}{x},\ \frac{d^2}{dx^2}f(x) = -\frac{1}{x^2},\ \frac{d^3}{dx^3}f(x) = \frac{2}{x^3}$$

Let's prove that (2.11) is satisfied for any values of the ***n*** parameter.

By substituting $n+1$ for the ***n*** index in (2.11), we obtain:

$$\frac{d^{n+1}}{dx^{n+1}}\ln(x) = -\frac{(-\ln(x)+\gamma+\Psi(n+1)+\pi\cot(\pi(n+1)))\,x^{-n-1}}{\Gamma(-n)}$$

Since

$$\frac{d}{dx}\ln(x) = \frac{1}{x}$$

we arrive to the following formula:

$$\frac{d^n}{dx^n}\left(\frac{1}{x}\right) = -\frac{(-\ln(x)+\gamma+\Psi(n+1)+\pi\cot(\pi(n+1)))\,x^{-n-1}}{\Gamma(-n)}$$

As soon as

$$\lim_{n\to 0}\left(\frac{d^n}{dx^n}\left(\frac{1}{x}\right)\right) = \frac{1}{x}$$

and

$$\lim_{n\to 0}\left(-\frac{(-\ln(x)+\gamma+\Psi(n+1)+\pi\cot(\pi(n+1)))\,x^{-n-1}}{\Gamma(-n)}\right) = \frac{1}{x}$$

this proves the validity of (2.11).

It must be noted that the calculation of ***n***-th order derivative to the $f(x)$ function in some cases gives results, which look "paradoxical" at a first glance. For example, let's calculate the ***n***-th derivative of the $f(x) = C$ function, where ***C*** is an arbitrary constant. As soon as $f(x) = C$ **is a special case of a more general function:** $f(x) = Cx^m$, $x \neq 0$, with $m = 0$, then, according to (2.8), $\frac{\partial^n}{\partial x^n}(Cx^m) = C\left(\frac{\partial^n}{\partial x^n}x^m\right) = \frac{Cm!\,x^{m-n}}{(m-n)!}$. Assuming in this equation $m = 0$, we eventually obtain the desired formula:

$$\frac{\partial^n}{\partial x^n}C = \frac{Cx^{-n}}{\Gamma(-n+1)} \qquad x \neq 0 \qquad (2.12)$$

Indeed, with $n=0$ we receive $diff(C, x\$0) = \frac{C x^0}{\Gamma(1)}$. As soon as $diff(C, x\$0) = C$; $x^0 = 1$, ($x \neq 0$); $\Gamma(1) = 1$, we obtain the following identity: $C = C$, so this formula holds true with $n=0$. For all natural $n$: $\Gamma(-n+1) = \infty$, so in this case we always have: $\frac{\partial^n}{\partial x^n} C := \frac{C x^{-n}}{\infty} = 0$ ; $x \neq 0$. Thus, (2.12) satisfies the existing requirement (on the vanishing of the derivative of any constant), but it is defined by an expression, which is not equal to zero.

**2.2. Derivatives of *n*-th order of a product of two functions and complex functions.**

**A)** In the case where the $f(x)$ function can be represented as a product of two or more functions, such as $f(x) = u(x)\,v(x)$, we can use Leibniz formula [1] to calculate the derivative of the ***n***-th order to the $f(x)$ function:

$$\frac{d^n}{dx^n}(u(x)\,v(x)) = \sum_{i=0}^{n} C_n^i \left(\frac{d^i}{dx^i} u(x)\right)\left(\frac{d^{n-i}}{dx^{n-i}} v(x)\right) \qquad (2.13)$$

**Example: Let's calculate the *n*-th derivative of the following function:**

$$f(x) = x\,e^x \qquad (2.14)$$

**Solution:** In accordance with (2.13)

$$\frac{d^n}{dx^n}(x\,e^x) = \sum_{i=0}^{n} C_n^i \left(\frac{d^i}{dx^i} x\right)\left(\frac{d^{n-i}}{dx^{n-i}} e^x\right).$$

*As soon as, according to* (2.8) $\frac{d^i}{dx^i} x = \frac{x^{1-i}}{(1-i)!}$, $\frac{d^{n-i}}{dx^{n-i}} e^x = e^x$ *then*

$\frac{d^n}{dx^n}(x\,e^x) = e^x \left(\sum_{i=0}^{n} \frac{C_n^i\, x^{1-i}}{(1-i)!}\right)$ *Since* $\sum_{i=0}^{n} \frac{C_n^i\, x^{1-i}}{(1-i)!} = x + n$, *the desired formula looks as follows:*

$$\frac{d^n}{dx^n}(x\,e^x) = e^x\,(x+n) \qquad (2.15)$$

Obviously, the obtained function satisfies the condition (2.10), so it is a normal derivative.

**B)** In the case where the initial $f(x)$ function is a complex function, that is, can be represented as follows:

$$f(x) = W(R),\ R = R(x) \qquad (2.16)$$

the formula for its derivative of ***n***-th order takes the form [2]:

$$\frac{\partial^n}{\partial x^n} W(R) = \sum_{k=1}^{n} \left(\sum_{s=0}^{k-1} \frac{(-1)^s\, C(k,s)\, R^s \left(\frac{\partial^n}{\partial x^n} R^{k-s}\right)\left(\frac{d^k}{dR^k} W(R)\right)}{k!}\right) \qquad (2.17)$$

It is defined only for the natural values of ***n*.**

**3. Calculation formula for the indefinite integral.**

**3.1. The sum of *n*-th order derivatives.**

Let's assume that the ***n***-th derivative of $G(x,n)$ to $f(x)$ is normal with $n=0$ **.** Let's calculate the $S_n$ sum

$$S_n = \sum_{i=0}^{n} G(i,x) \qquad (3.1)$$

for the functional sequence $G(0,x), G(1,x), G(2,x), G(3,x)$ ..... where, obviously, $G(0,x) = f(x)$ .

For this purpose we can differentiate the left and right side of equality (3.1) and obtain the following formula:

$$\frac{\partial}{\partial x} S_n = \frac{\partial}{\partial x}\left(\sum_{i=0}^{n} G(i,x)\right) = \sum_{i=0}^{n}\left(\frac{\partial}{\partial x} G(i,x)\right) = \sum_{i=0}^{n} G(i+1,x) \qquad (3.2)$$

Subtracting equality (3.1) from equality (3.2), with consideration of (2.2), we obtain:

$$\frac{\partial}{\partial x} S_n - S_n = \sum_{i=0}^{n} G(i+1,x) - \left(\sum_{i=0}^{n} G(i,x)\right) = G(n+1,x) - G(0,x) \qquad (3.1.3)$$

Since $G(n,x)$ is a normal ***n***-th derivative, then equality (3.3) takes the following form:

$$\frac{\partial}{\partial x} S_n - S_n = G(n+1,x) - f(x) \qquad (3.4)$$

Thus, we obtain a linear differential equation of first order to the $S_n$ function.

Equation (3.4) has the following general solution:

$$S_n = e^x\left(\int (G(n+1,x) - f(x))\, e^{-x}\, dx\right) + Ce^x \qquad \textbf{C}\text{-const.} \qquad (3.5)$$

The arbitrary constant C is determined from the following condition:

$$[S_n]_{n=0} = f(x) \; .$$

In this case, the equation (3.5) takes the following form:

$$\sum_{i=0}^{n} G(i,x) = e^x\left(\int (G(n+1,x) - f(x))\, e^{-x}\, dx\right) \qquad (3.6)$$

Consequently, if the parametric function $G(i,x)$ is normal at the $n=0$ point of $i$ -th derivative to the $f(x)$ , then (3.6) determines the sum of this series.

**Conclusion: If $G(n,x)$ is a parametric function, for which the following condition is satisfied:**

$$\lim_{n\to\infty} G(n+1,x) = 0$$

**then the following formula holds true:**

$$\sum_{i=0}^{\infty}\left(\frac{d^i}{dx^i} f(x)\right) = -e^x\left(\int f(x)\, e^{-x}\, dx\right) \qquad (3.7)$$

**3.2. The Fundamental Theorem.**

**Let's prove the theorem, hereinafter referred to as "The fundamental theorem":**

**Suppose that for a given $n$ -times differentiable $f(x)$ function the following sequence:**

$$f(x),\ \frac{d}{dx}f(x),\ \frac{d^2}{dx^2}f(x),\ \frac{d^3}{dx^3}f(x)\ ....$$

**allows to define its *n*-th derivative:**

$$\frac{d^n}{dx^n}f(x) = G(n, x)$$ **(where $n$ - is a natural number),**
**which is normal.**

**If we expand the domain of the *n* parameter to the set of actual numbers, there is a limit expression:**

$$\lim_{n \to -1} G(n, x)$$

**and the following formula holds true:**

$$\int f(x)\,dx = \lim_{n \to -1} G(n, x) + const$$

**or, equivalently:**

$$\int f(x)\,dx = \lim_{n \to -1}\left(\frac{d^n}{dx^n}f(x)\right) + const \tag{3.8}$$

**The proof is as follows:** *Let's expand the domain of **n** parameter to the set of real numbers. Due to the fact that, by definition, **n**-th derivative is normal, we obtain:*

$$\lim_{n \to -1} G(n+1, x) = f(x)$$

*Therefore,* (3.6), *with consideration of this equation, takes the following form:*

$$\lim_{n \to -1}\left(\sum_{i=0}^{n} G(i, x)\right) = 0 \tag{3.9}$$

*As follows from the above:*

$$\lim_{n \to -1}\left(\sum_{i=1}^{n} G(i, x)\right) + f(x) = 0$$

*And, thus:*

$$f(x) = -\left(\lim_{n \to -1}\left(\sum_{i=1}^{n} G(i, x)\right)\right)$$

*Let's replace the summation index $i$ by $i+1$ in the right side of this equality.*

$$f(x) = -\left(\lim_{n \to -1}\left(\sum_{i=0}^{n-1} G(i+1, x)\right)\right)$$

*As soon as, according to* (2.2)

$$G(i+1, x) = \frac{\partial}{\partial x} G(i, x)$$

*the previous equation takes the following form:*

$$f(x) = -\left(\frac{\partial}{\partial x}\left(\lim_{n \to -1}\left(\sum_{i=0}^{n-1} G(i, x)\right)\right)\right)$$

*The right side of this equality shall eventually remain unchanged, if we add and subtract the* $G(n, x)$ *summand:*

$$f(x) = -\left(\frac{\partial}{\partial x}\left(\lim_{n \to -1}\left(\sum_{i=0}^{n} G(i, x) - G(n, x)\right)\right)\right)$$

*In this case this equality takes the following form:*

$$f(x) = -\left(\frac{\partial}{\partial x}\left(\lim_{n \to -1}\left(\sum_{i=0}^{n} G(i, x)\right)\right)\right) + \frac{\partial}{\partial x}\left(\lim_{n \to -1} G(n, x)\right)$$

*As soon as, according to* (3.9)

$$\lim_{n \to -1}\left(\sum_{i=0}^{n} G(i, x)\right) = 0$$

*the previous equation takes the following form:*

$$f(x) = \frac{\partial}{\partial x}\left(\lim_{n \to -1} G(n, x)\right) \qquad (3.10)$$

Since, according to the conditions of the Theorem, the limit expression $\lim_{n \to -1} G(n, x)$ exists, that is, there is a certain function, we obtain **the desired formula (3.8)** by integrating the equation (3.10).

**The Theorem is proved.**

Thus, **(3.8) defines the desired algorithm for the operation of integration, that is, the definition of** $\int f(x)\,dx$ **expression.** Let's note, that an important condition is the normality of the defined ***n***-th derivative to the initial function, because otherwise (3.8) can not be used.

**Conclusion: If** $G(n, x)$ is a normal derivative with $n = -1$, then (3.8) can be used in a simplified form**:**

$$\int f(x)\,dx = \left\{\frac{d^n}{dx^n} f(x)\right\}_{n=-1} + C \qquad (3.11)$$

**Example No. 1:** ***Calculate the integral of C constant:***

$$\int C\,dx$$

**Solution.** As soon as, according to (2,12)

$$\frac{\partial^n}{\partial x^n} C = \frac{C x^{-n}}{(-n)!}, x \neq 0$$

then, since this ***n***-th derivative is normal, the required integral equals:

$$\int C\,dx = \lim_{n \to -1} \frac{C x^{-n}}{(-n)!} = C x + const$$

which matches the tabulated value.

**Example No. 2: Let's calculate the following integral**

$\int \sin(a\,x)\,dx$ , **where** $a$ is a constant.

**Solution:** *As soon as*

$$f(x) = \sin(a\,x), \quad G(n,x) = a^n \sin\left(a\,x + \frac{n\,\pi}{2}\right)$$

*then, since* $G(0,x) = \sin(a\,x)$ **, *according to*** *(3.11) we obtain*

$$\int \sin(a\,x)\,dx = \left[a^n \sin\left(a\,x + \frac{n\,\pi}{2}\right)\right]\Big|_{n=-1} = -\frac{\cos(a\,x)}{a}$$

*Thus,*

$$\int \sin(a\,x)\,dx = -\frac{\cos(a\,x)}{a} + const$$

*which matches the tabulated value.*

**Example No. 4: Let's calculate the following integral:**

$\int \ln(x)\,dx$ ,

**Solution:** *As soon as, according to* (2.11)

$$f(x) = \ln(x), \quad G(n,x) = \frac{d^n}{dx^n}\ln(x) = -\frac{(-\ln(x) + \gamma + \Psi(n) + \pi\cot(\pi\,n))\,x^{-n}}{\Gamma(-n+1)}$$

*then, since*

$$G(0,x) = \lim_{n\to 0}\left(-\frac{(-\ln(x) + \gamma + \Psi(n) + \pi\cot(\pi\,n))\,x^{-n}}{\Gamma(-n+1)}\right) = \ln(x)$$ *according to* (3.8) *we obtain*

$$\int \ln(x)\,dx = \lim_{n\to -1}\left(-\frac{(-\ln(x) + \gamma + \Psi(n) + \pi\cot(\pi\,n))\,x^{-n}}{\Gamma(-n+1)}\right) = x\,(\ln(x) - 1)$$

*Therefore,*

$$\int \ln(x)\,dx = x\,(\ln(x) - 1) + const$$

*which matches the tabulated value.*

**Example No. 4.** *Let's calculate the following integral:*

$$\int \frac{1}{(a\,x+b)\,(c\,x+m)}\,dx$$

*where* $a, b, c, m$ *are constants.*

**Solution:** *By the calculation of the* ***n****-th derivative of the following function*:

$$\frac{1}{(a x+b)(c x+m)}$$

*we obtain the following formula:*

$$\frac{\partial^n}{\partial x^n}\left(\frac{1}{(a x+b)(c x+m)}\right) = \frac{n!\,(-1)^n (c x+m)^{-n-1}\left(a^{n+1}(a x+b)^{-n}(c x+m)^{n+1}-c^{n+1} a x-c^{n+1} b\right)}{(a x+b)(a m-b c)}$$

*Since*

$$\left[\frac{\partial^n}{\partial x^n}\left(\frac{1}{(a x+b)(c x+m)}\right)\right]_{n=0} = \frac{1}{(a x+b)(c x+m)}$$

*according to* (3.8) *we obtain:*

$$\lim_{n\to -1}\frac{n!\,(-1)^n (c x+m)^{-n-1}\left(a^{n+1}(a x+b)^{-n}(c x+m)^{n+1}-c^{n+1}(a x+b)\right)}{(a x+b)(a m-b c)} = \frac{-\ln(c x+m)+\ln(a x+b)-\ln(a)+\ln(c)}{a m-b c}$$

*Thus, eventually we obtain:*

$$\int\frac{1}{(a x+b)(c x+m)}\,dx = \frac{-\ln(c x+m)+\ln(a x+b)-\ln(a)+\ln(c)}{a m-b c} + const$$

*Tabulated value of the integral is as follows:*

$$\int\frac{1}{(a x+b)(c x+m)}\,dx = \frac{-\ln(c x+m)+\ln(a x+b)}{a m-b c} + const$$

*The difference is a* $\frac{-\ln(a)+\ln(c)}{a m-b c}$ *constant, which is neutralized by the choice of an arbitrary constant.*

**Note:** *A simple substitution of* $n=-1$ *to the following expression:*

$$\frac{n!\,(-1)^n (c x+m)^{-n-1}\left(a^{n+1}(a x+b)^{-n}(c x+m)^{n+1}-c^{n+1}(a x+b)\right)}{(a x+b)(a m-b c)}$$

*does not give the desired result (this leads to the uncertainty of* $\{\infty\cdot 0\}$ *type), so it is essential to have a limit in* (3.8).

### 3.3. Other formulas, derived from the Fundamental Theorem.

**Let's prove the Theorem: If** $G(n,x)$ **is a normal derivative at the** $n=0,1,2..s$ **points, and the** $\lim_{n\to -k} G(n,x),\ k=1..s$ **functions are differentiable and integrable, then the** $s$ **-multiple integral of a** $f(x)$ **function is defined by the following formula:**

$$[f(x)]_s = \lim_{n\to -s}\left(\frac{d^n}{dx^n}f(x)\right) + \sum_{i=0}^{s-1}\frac{c_i x^i}{i!} \qquad (3.12)$$

where $c_i$ are arbitrary constants.

**The proof is as follows:** *Indeed,*

$$[f(x)]_s = \left[\lim_{n\to -1}\left(\frac{d^n}{dx^n}f(x)\right)\right]_{s-1} = \lim_{n\to -1}\left(\frac{d^{n-s+1}}{dx^{n-s+1}}f(x)\right)$$

*Substituting* $k+s+1$ *for the* $n$ *index, we get the following:*

$$[f(x)]_s = \lim_{k\to -s}\left(\frac{d^k}{dx^k}f(x)\right)$$

*From the other side, because the* $G(-k,x), k=1..s$ *functions* ***are differentiable and integrable, then***

$$\iint f(x)\,dx\,dx := \iint (G(-1,x)+c_1)\,dx\,dx = \lim_{n\to -s} G(n,x) + \sum_{i=0}^{s-1}\frac{c_i x^i}{i!}$$

*The Theorem is proved.*

**Conclusion: If** $G(n,x)$ **is a normal derivative at the** $n=-s$ **points, and the** $G(-k,x), k=1..s$ **functions are differentiable and integrable, then the** $s$ **-multiple integral of a** $f(x)$ **function is also defined by the following formula:**

$$[f(x)]_s = \left\{\frac{d^n}{dx^n}f(x)\right\}_{n=-s} + \sum_{i=0}^{s-1}\frac{c_i x^i}{i!} \tag{3.13}$$

*where* $c_i$ *are arbitrary constants.*

**Example:** *Let's calculate the following integrals:*

$$\int e^{-x^2}\,dx \text{ и } \iint e^{-x^2}\,dx\,dx$$

**Solution:** *As soon as*

$$\frac{d^n}{dx^n}e^{-x^2} = -\frac{x^{-n}\,\mathrm{hypergeom}\left(\left[\frac{1}{2},1\right],\left[1-\frac{1}{2}n,-\frac{1}{2}n+\frac{1}{2}\right],-x^2\right)}{n\,\Gamma(-n)}$$

*and, according to the Fundamental Theorem and (3.12),*

$$\int e^{-x^2}\,dx = \lim_{n\to -1}\left(\frac{d^n}{dx^n}e^{-x^2}\right) \qquad \iint e^{-x^2}\,dx\,dx = \lim_{n\to -2}\left(\frac{d^n}{dx^n}e^{-x^2}\right)$$

*upon the calculation of*

$$\lim_{n\to -k}\left(\frac{d^n}{dx^n}e^{-x^2}\right) = \lim_{n\to -k}\left(-\frac{x^{-n}\,\mathrm{hypergeom}\left(\left[\frac{1}{2},1\right],\left[1-\frac{1}{2}n,-\frac{1}{2}n+\frac{1}{2}\right],-x^2\right)}{n\,\Gamma(-n)}\right), \quad k=1,2$$

*we eventually obtain:*

$$\int e^{-x^2}\,dx = x\,\mathrm{hypergeom}\left(\left[\frac{1}{2}\right],\left[\frac{3}{2}\right],-x^2\right); \quad \iint e^{-x^2}\,dx\,dx = \frac{1}{2}x^2\,\mathrm{hypergeom}\left(\left[\frac{1}{2},1\right],\left[\frac{3}{2},2\right],-x^2\right)$$

## 4. Calculation formula for the definite integral.

Conventionally, definite integrals are defined according to formula [1]:

$$\int_a^b f(x)\,dx = F(b) - F(a) \qquad (4.1)$$

where $F(x)$ is a primitive function of $f(x)$ .

Since, in accordance with the Fundamental Theorem, the primitive function $F(x)$ is defined as:

$$F(x) = \lim_{n \to -1} \left( \frac{d^n}{dx^n} f(x) \right)$$

we obtain the calculation formula for the definite integral:

$$\int_a^b f(x)\,dx = \left[ \lim_{n \to -1} \left( \frac{d^n}{dx^n} f(x) \right) \right]_{[x=b]} - \left[ \lim_{n \to -1} \left( \frac{d^n}{dx^n} f(x) \right) \right]_{[x=a]} \qquad (4.2)$$

In the case where the limiting operation can be replaced by substitution, this formula looks as follows:

$$\int_a^b f(x)\,dx = \left[ \frac{d^n}{dx^n} f(x) \right]_{[n=-1,\, x=b]} - \left[ \frac{d^n}{dx^n} f(x) \right]_{[n=-1,\, x=a]} \qquad (4.3)$$

**Example No. 1: Let's calculate the following definite integral:**

$$\int_0^1 \sin(x)\,dx$$

**Solution:** Since

$$\frac{d^n}{dx^n} \sin(x) = \sin\left( x + \frac{1 \pi n}{2} \right)$$

*according to* (4.3) *we obtain:*

$$\int_0^1 f(x)\,dx = \left[ \sin\left( x + \frac{1 \pi n}{2} \right) \right]_{[n=-1,\, x=1]} - \left[ \sin\left( x + \frac{1 \pi n}{2} \right) \right]_{[n=-1,\, x=0]} = -\sin\left( -\frac{1 \pi}{2} \right) + \sin\left( 1 - \frac{1 \pi}{2} \right)$$

*Thus, eventually we obtain:*

$$\int_0^1 f(x)\,dx = 1 - \cos(1)$$

*which corresponds to the results of calculations according to* (4.1).

For the integration of a complex function of $W(R(x))$ type we obtain the following formula:

$$\int W(R(x))\,dx = \sum_{k=1}^{\infty} \left( \sum_{s=0}^{k-1} \frac{(-1)^s\, C(k,s)\, R(x)^s \left( \int R(x)^{k-s}\,dx \right) \left( \frac{d^k}{dR^k} W(R) \right)}{k!} \right) \qquad (4.4)$$

**Example No. 1.** *Let's calculate the following integral:*

$$\int \sin(\sin(x))\,dx$$

**Solution.** Let's introduce the following notation: $R = \sin(x)$ . Then, according to (4.4) we obtain:

$$\int \sin(\sin(x))\,dx = \sum_{k=1}^{\infty} \frac{1 \left( \frac{d^k}{dR^k} \sin(R) \right) \left( \sum_{s=0}^{k-1} (-1)^s C(k,s) \sin(x)^s \left( \int \sin(x)^{k-s}\,dx \right) \right)}{k!}$$

*As soon as*

$$\frac{d^k}{dR^k} \sin(R) = \sin\left( R + \frac{k\pi}{2} \right)$$

we obtain the desired formula:

$$\int \sin(\sin(x))\,dx = \sum_{k=1}^{\infty} \left( \frac{1}{k!} \right) .\sin\left( R + \frac{k\pi}{2} \right) \left( \sum_{s=0}^{k-1} (-1)^s C(k,s) \sin(x)^s \left( \int \sin(x)^{k-s}\,dx \right) \right)$$

Which, with a high precision, gives correct result in the [-1, 1] interval.

## The theory of linear second-order ODE.

**Notations:- DE - differential equation; - LODE - linear ordinary differential equation; - ODE - ordinary differential equation- N - natural number; - C(n,i)- binomial coefficient – $C_n^i$ In line with such form of record, in many cases we shall use an equivalent RECORD of binomial coefficient:** $C_n^i = \begin{bmatrix} i \\ n \end{bmatrix}$. **Unless specifically stated otherwise, are arbitrary constants.** $[[0,0,0]]$ **under summation sign means dots - ... - (because it is impossible to enter dots under summation sign in the software program). - It is assumed that all LDE coefficients are functions differentiable to the required degree.**

**1).** $[f(x)]_k$ **-** $k$ **- multiple integral of a** $f(x)$ **function.**

$$[f(x)]_3 = \int\int\int f(x)\, dx\, dx\, dx$$

and so on. **2). Special notation for a sum of indexes from $i_{k_1}$ to $i_{k_2}$**

**For example:**

$= i_0 \; : \; = i_3 \; : \; = i_0 + i_1 \; : \; =$

**The notations $(0, 0) = 0$ and, generally, $(k, k) = k$, can be applied as well.**

**3). Special notation for a sequence of nested sums, for example:**

$$\sum_{s_{3,1} = [k_3, k_1]}^{[p_3, p_1]} G(s_1, s_2, s_3) = \sum_{s_3 = k_3}^{p_3} \left( \sum_{s_2 = k_2}^{p_2} \left( \sum_{s_1 = k_1}^{p_1} G(s_1, s_2, s_3) \right) \right)$$

**here from the lower index $s_3 = k_1$ to its upper limit $s_3 = p_1$ for the outer sum and $s_1 = k_3$ to $s_1 = p_3$ of the first nested sum. Example:**

$$\sum_{z_{0,3} = [-1, 2]}^{[0, 3]} G(z_0, z_1, z_2, z_3) = \sum_{z_0 = -1}^{0} \left( \sum_{z_1 = 0}^{1} \left( \sum_{z_2 = 1}^{2} \left( \sum_{z_3 = 2}^{3} G(z_0, z_1, z_2, z_3) \right) \right) \right)$$

**Here it is assumed that the lower boundary defined by values from -1 to 2 is divided by the value of the range, in this case - 4 (3,2,1,0), that is, the change from -1 to 2 occurs with an increment of 1. The same applies to the upper boundary. In case where the change from one value to another for any boundary is not obvious, the values are written down, for example as follows:**

$$\sum_{z_{0,3} = [-3, 2, 2]}^{[1, k+1, k+3]} G(z_0, z_1, z_2, z_3) = \sum_{z_0 = -3}^{1} \left( \sum_{z_1 = 2}^{k+1} \left( \sum_{z_2 = 2}^{k+2} \left( \sum_{z_3 = 2}^{k+3} G(z_0, z_1, z_2, z_3) \right) \right) \right)$$

Below there is an example where the value of the upper or lower boundary does not change:

$$Sum(G(z_0, z_1, z_2, z_3), z_{3,1} = 0..[4, 1]) = \sum_{z_3 = 0}^{4} \left( \sum_{z_2 = 0}^{1} \left( \sum_{z_1 = 0}^{1} G(z_1, z_2, z_3) \right) \right)$$

$$\sum_{z_{3,1} = 0}^{[4, 1]} G(z_0, z_1, z_2, z_3) = \sum_{z_3 = 0}^{4} \left( \sum_{z_2 = 0}^{1} \left( \sum_{z_1 = 0}^{1} G(z_1, z_2, z_3) \right) \right)$$

**4). Notation for a sequence in specially nested functions according to the following rule:**

$$[b_{1+(0,s)}(x)]_{i_{0,l}+k_{0,l}}\, b_{k_{0,l}}(x) = \left[\left[\left[[b_{1+i_0+i_1+i_2+i_3..i_s}(x)]_{i_0+k_0}\, b_{k_0}(x)\right]_{i_1+k_1} b_{k_1}(x)\right]_{i_2+k_2} b_{k_2}(x)\right]_{i_3+k_3} b_{k_3}(x)\,..$$

For instance:

$$[b_{1+(0,0)}(x)]_{i_{0,0}+k_{0,0}}\, b_{k_{0,0}}(x) = [b_{1+i_0}(x)]_{i_0+k_0}\, b_{k_0}(x)$$

$$[b_{1+(0,1)}(x)]_{i_{0,1}+k_{0,1}}\, b_{k_{0,1}}(x) = \left[[b_{1+i_0+i_1}(x)]_{i_0+k_0}\, b_{k_0}(x)\right]_{i_1+k_1} b_{k_1}(x)$$

$$[b_{1+(0,4)}(x)]_{i_{0,3}+k_{0,3}}\, b_{k_{0,3}}(x) = \left[\left[\left[[b_{1+i_0+i_1+i_2+i_3+i_4}(x)]_{i_0+k_0}\, b_{k_0}(x)\right]_{i_1+k_1} b_{k_1}(x)\right]_{i_2+k_2} b_{k_2}(x)\right]_{i_3+k_3} b_{k_3}(x)$$

and so on.

# Introduction

The history of differential calculus, started from the mathematical introduction of variables by Descartes and finalized by Newton and Leibniz in about 1667, has evloved to a comprehensive theory of ODE. The fundamental problem of the ODE theory [Collegiate Dictionary. Math. Scientific ed. "Great Russian Encyclopedia", Moscow, 1998] is to develop well-known and invent new ways of solving boundary value problems for ODEs.

Obviously, the solution of any of the known boundary value problems for the linear ODEs of *m*-order requires knowledge of all *m* massive of partial solutions, and, thus, the fundamental problem of the theory of ordinary differential equations is the integration of a given ODE. The best-known methods are, among others, replacement of the dependent and independent variables, special substitutions, group analysis, operational method, etc. [1,2,3]. However, so far there is no universal method to find partial solutions of ODEs. For the linear second-order ODEs the most common way of finding partial solutions is the method of power series. The solution of a given linear ODE, in the case where its coefficients are holomorphic functions in the neighborhood of $x_0$, is as follows:

(1.0)

where $a_k$ are constants defined by the substitution to the initial ODE. In fact, the method of power series is a way to find partial solutions in the neighborhood of a given $x_0$ point with the use of Taylor series; so this approach is rather an approximate method, and cannot be called universal. Moreover, it is difficult to acknowledge that different linear ODEs, which describe the most sophisticated dynamic processes, have a solution represented essentially in identical form (1.0). Moreover, the method of power series is applicable only under certain conditions; in particular, the coefficients of the initial LODE should be presented in the neighborhood of a given point through a convergent series and, at the same time, must not have (for example, for the linear second-order ODEs) poles above first or second degree (pp. 184-185, A.A. Esipov, L.I. Sazonov, V.I. Yudovich. Differential Equations. Moscow, "Higher School Publishers", 2001)**.** It is also obvious that finding solutions in the form of (1.1) becomes much more difficult with the increase of the order of ODEs. This obstacle significantly limits practical application.

**The objective of this publication is to build an universal integration theory for linear non-homogeneous ODEs of arbitrary order.**

It is also required to consider the form of representation of the desired solution. As known, in a general case such representation is possible either in exact form, that is, as a certain function, or in approximate form, with the use of analytic functions. It is important to understand that exact solutions of linear ODEs are, essentially, analytic functions.

Example No. 1. A linear second-order ODE

has the following general solution:

which is exact. However, as soon as

$$\mathbf{e}^x = \sum_{i=0}^{\infty} \frac{x^i}{i!}, \ \mathbf{e}^{(-x)} = \sum_{i=0}^{\infty} \frac{(-1)^i x^i}{i!}$$

this solution has an equivalent analytic representation:

$$\mathrm{Y}(x) = C_1 \left( \sum_{i=0}^{\infty} \frac{x^i}{i!} \right) + C_2 \left( \sum_{i=0}^{\infty} \frac{(-1)^i x^i}{i!} \right)$$

that is, an analytic function.

Example No. 2. A linear second-order ODE

has a general solution defined by Airy functions:

$$\mathrm{Y}(x) = C_1 \, \mathrm{AiryAi}(x) + C_2 \, \mathrm{AiryBi}(x)$$

which are hyper-geometric functions, that is, infinite converging sequences [1]. Conventionally, this solution is deemed exact, but essentially it is an analytic function.
As follows from the above, in a general case partial solutions of linear ODEs with variable coefficients should be found in the category of analytic functions.
Basing upon this postulate, let's introduce the definition of linear ODE integration algorithm.

**Definition 1.1.** The general integration algorithm for homogenous Linear Ordinary Differential Equations:

$$\left(\frac{\partial^m}{\partial x^m} y\right) + \left(\sum_{p=1}^{m} a_p(x) \left(\frac{\partial^{m-p}}{\partial x^{m-p}} y\right)\right) = 0 \qquad (1.1)$$

where $y = \{y_1(x), y_2(x) .. y_m(x)\}$ is the desired function, $a_p(x)$ are specified functions (differentiable and integrable to the required degree). Let's assume the following form of the desired solution:

$y = \mathrm{F}(x, \alpha)$ , $\mathrm{F}(x, \alpha) = \{F_i(x, \alpha), i = 1 .. m\}$ , where

is an analytic function,

$$\alpha = \sum_{i=0}^{\infty} \xi_i(a_1(x), a_2(x) .. a_m(x))$$

$\mathrm{F}(x, \alpha)$ is the required relation.

As follows from this definition:

- coefficients of the initial ODE need not (unless specifically required) to be represented through a convergent series.
- partial solutions of are defined as the following functions:

$$y_i(x) = F_i(x, \alpha), i = 1, 2, 3 .. m \qquad (1.2)$$

to be constructed on the basis of mathematical conclusions.
- $\alpha_k$ components in this general formula are the following series:

$$\alpha_k = \sum_{j=0}^{\infty} \xi_{k,j}(a_1(x), a_2(x) .. a_m(x)), k = 1, 2, 3 .. m \qquad (1.3)$$

from $\xi_{k,j}(a_1(x), a_2(x) .. a_m(x))$ functions, which are defined in some formulas through the initial coefficients: .

**Thus, the problem of finding a general algorithm for the integration of linear ODEs of (1.1) type consists of two successive tasks:**

**- to define the form to represent $F_i(x, \alpha)$ functions for the required partial solutions (1.2).**

**- to define the $\xi_{k,j}(a_1(x), a_2(x) .. a_m(x))$ functions.**

It is obvious that a general algorithm for the integration of linear ODEs with variable coefficients cannot be build in view of the existing concepts and definitions used in the modern theory of ODE, just because, if it was the case, this problem would be solved long time ago (for over 350 years). Therefore, the overall structure of the new theory of LODE integration requires the introduction of new mathematical ideas. One of them is the concept of Leibniz operator.

### 1.1. Leibniz operator

**Definition 1.2: Leibniz operator for index $i$ of $u$ , $v$ functions and $n$ parameter of $\mathrm{f}(i)$ is the expression $L[\mathrm{f}(i)]$ , defined by the following equation:**

$$L[\mathrm{f}(i)] = \sum_{i=0}^{n} \mathrm{f}(i)\, \mathrm{C}(n, i) \left(\frac{\partial^{i}}{\partial x^{i}} u\right) \left(\frac{\partial^{n-i}}{\partial x^{n-i}} v\right) \tag{1.4}$$

The task is to find the sum of the following series:

$$\sum_{i=0}^{n} \mathrm{f}(i)\, \mathrm{C}(n, i) \left(\frac{\partial^{i}}{\partial x^{i}} u\right) \left(\frac{\partial^{n-i}}{\partial x^{n-i}} v\right)$$

For example, if $\mathrm{f}(i) = 1$ ; then (1.4) takes the following form:

$$L[1] = \sum_{i=0}^{n} \mathrm{C}(n, i) \left(\frac{\partial^{i}}{\partial x^{i}} u\right) \left(\frac{\partial^{n-i}}{\partial x^{n-i}} v\right)$$

Since, according to Leibniz formula:

$$\sum_{i=0}^{n} \mathrm{C}(n, i) \left(\frac{\partial^{i}}{\partial x^{i}} u\right) \left(\frac{\partial^{n-i}}{\partial x^{n-i}} v\right) = \frac{\partial^{n}}{\partial x^{n}} (u\, v) \tag{1.5}$$

then, eventually, we get the desired equality:

(1.6)

**Proposition 1.1: Impact of Leibniz operator for index $i$ of $u$ , $v$ functions and $n$ parameter on the $\mathrm{f}(i) = i$ function is defined by the following formula:**

$$L[i] = n \left(\frac{\partial^{n-1}}{\partial x^{n-1}} \left(\left(\frac{\partial}{\partial x} u\right) v\right)\right) \tag{1.7}$$

**Proof**: In accordance with (1.4):

$$\sum_{i=0}^{n} i\, \mathrm{C}(n, i) \left(\frac{\partial^{i}}{\partial x^{i}} u\right) \left(\frac{\partial^{n-i}}{\partial x^{n-i}} v\right) = \sum_{i=0}^{n} \left(\left(\frac{n!\, i}{i!\,(n-i)!}\right) \cdot \left(\frac{\partial^{i}}{\partial x^{i}} u\right)\right) \left(\frac{\partial^{n-i}}{\partial x^{n-i}} v\right) =$$
$$\sum_{i=0}^{n} \left(\left(\frac{n\,(n-1)!}{(i-1)!\,(n-1-(i-1))!}\right) \cdot \left(\frac{\partial^{i}}{\partial x^{i}} u\right)\right) \left(\frac{\partial^{n-i}}{\partial x^{n-i}} v\right) =$$
$$n \left(\sum_{i=0}^{n} \mathrm{C}(n-1, i-1) \left(\frac{\partial^{i}}{\partial x^{i}} u\right) \left(\frac{\partial^{n-i}}{\partial x^{n-i}} v\right)\right)$$

Let's establish the value of the expression:

$$\sum_{i=0}^{n} \mathrm{C}(n-1, i-1) \left(\frac{\partial^{i}}{\partial x^{i}} u\right) \left(\frac{\partial^{n-i}}{\partial x^{n-i}} v\right)$$

If we replace the summation index $i$ by $i+1$ in the obtained value, we establish that

$$\sum_{i=0}^{n} C(n-1,i-1)\left(\frac{\partial^i}{\partial x^i}u\right)\left(\frac{\partial^{n-i}}{\partial x^{n-i}}v\right) = \sum_{i=0}^{n-1} C(n-1,i)\left(\frac{\partial^{i+1}}{\partial x^{i+1}}u\right)\left(\frac{\partial^{n-1-i}}{\partial x^{n-1-i}}v\right) =$$

$$\sum_{i=0}^{n-1} C(n-1,i)\,\mathrm{Diff}\left(\frac{\partial}{\partial x}u, x \$ i\right)\left(\frac{\partial^{n-1-i}}{\partial x^{n-1-i}}v\right)$$

In accordance with Leibniz formula (1.5), where $n$ parameter is replaced by $n-1$, and $u$ function is replaced by $\frac{\partial}{\partial x}u$, we obtain the following:

$$\sum_{i=0}^{n-1} C(n-1,i)\,\mathrm{Diff}\left(\frac{\partial}{\partial x}u, x \$ i\right)\left(\frac{\partial^{n-1-i}}{\partial x^{n-1-i}}v\right) = \frac{\partial^{n-1}}{\partial x^{n-1}}\left(\left(\frac{\partial}{\partial x}u\right)v\right)$$

Thus, returning back, we obtain (1.7).
**Proposition 1.1 is proved.**

Similarly we can prove
**Proposition 1.2.: Impact of Leibniz operator for index $i$ of $u$, $v$ functions and $n$ parameter on the $f(i)=i^2$ function is defined by the following formula:**

$$L[i^2] = n\left(\frac{\partial^{n-1}}{\partial x^{n-1}}\left(\left(\frac{\partial}{\partial x}u\right)v\right)\right) + n(n-1)\left(\frac{\partial^{n-2}}{\partial x^{n-2}}\left(\left(\frac{\partial^2}{\partial x^2}u\right)v\right)\right)$$

Now let's prove a general concept, that shall be much needed in the future:

**Proposition 1.3: Let's assume that $L[\ ]$ is Leibniz operator for index $i$ of $u$, $v$ functions and $n$ parameter. Then the following equation applies for any $k$ in the range of $[0,N]$:**

$$L[i^k] = \sum_{i=0}^{k}\left(\sum_{s=0}^{i}\frac{s^k(-1)^{(s+i)}}{\Gamma(i-s+1)\Gamma(s+1)}\right)C(n,i)\,i!\left(\frac{\partial^{n-i}}{\partial x^{n-i}}\left(\left(\frac{\partial^i}{\partial x^i}u\right)v\right)\right) \tag{1.8}$$

**Proof**: In accordance with the definition of Leibniz operator (1.4), (1.8) can be presented as follows:

$$\sum_{i=0}^{n} i^k C(n,i)\left(\frac{\partial^i}{\partial x^i}u\right)\left(\frac{\partial^{n-i}}{\partial x^{n-i}}v\right) =$$
$$\sum_{i=0}^{k}\left(\sum_{s=0}^{i}\frac{s^k(-1)^{(s+i)}}{\Gamma(i-s+1)\Gamma(s+1)}\right)C(n,i)\,i!\left(\frac{\partial^{n-i}}{\partial x^{n-i}}\left(\left(\frac{\partial^i}{\partial x^i}u\right)v\right)\right)$$

Let's check whether this equality is true for . Indeed, assuming that $k=0$ we get:

$$\sum_{i=0}^{n} C(n,i)\left(\frac{\partial^i}{\partial x^i}u\right)\left(\frac{\partial^{n-i}}{\partial x^{n-i}}v\right) =$$
$$\sum_{i=0}^{0}\left(\sum_{s=0}^{i}\frac{1(-1)^{(s+i)}}{\Gamma(i-s+1)\Gamma(s+1)}\right)C(n,i)\,i!\left(\frac{\partial^{n-i}}{\partial x^{n-i}}\left(\left(\frac{\partial^i}{\partial x^i}u\right)v\right)\right)$$

And, following that:

$$\sum_{i=0}^{n} C(n,i)\left(\frac{\partial^i}{\partial x^i}u\right)\left(\frac{\partial^{n-i}}{\partial x^{n-i}}v\right)=\frac{\partial^n}{\partial x^n}(u\,v)$$

However, this equality is Leibniz formula; hence (1.8) includes the well-known Leibniz formula as a special case.

Assuming that $k=1$ in (1.8), we obtain:

$$L[i]:=\sum_{i=0}^{1}\left(\sum_{s=0}^{i}\frac{s\,(-1)^{(s+i)}}{\Gamma(i-s+1)\,\Gamma(s+1)}\right)C(n,i)\,i!\left(\frac{\partial^{n-i}}{\partial x^{n-i}}\left(\left(\frac{\partial^i}{\partial x^i}u\right)v\right)\right)=$$
$$\left(\sum_{s=0}^{0}\frac{s\,(-1)^{s}}{\Gamma(1-s)\,\Gamma(s+1)}\right)C(n,0)\left(\frac{\partial^n}{\partial x^n}(u\,v)\right)$$
$$+\left(\sum_{s=0}^{1}\frac{s\,(-1)^{(s+1)}}{\Gamma(2-s)\,\Gamma(s+1)}\right)C(n,1)\left(\frac{\partial^{n-1}}{\partial x^{n-1}}\left(\left(\frac{\partial}{\partial x}u\right)v\right)\right)$$

As soon as in this formula

$$\left(\sum_{s=0}^{0}\frac{s\,(-1)^{s}}{\Gamma(1-s)\,\Gamma(s+1)}\right)C(n,0)\left(\frac{\partial^n}{\partial x^n}(u\,v)\right)=0$$

$$\left(\sum_{s=0}^{1}\frac{s\,(-1)^{(s+1)}}{\Gamma(2-s)\,\Gamma(s+1)}\right)C(n,1)\left(\frac{\partial^{n-1}}{\partial x^{n-1}}\left(\left(\frac{\partial}{\partial x}u\right)v\right)\right)=n\left(\frac{\partial^{n-1}}{\partial x^{n-1}}\left(\left(\frac{\partial}{\partial x}u\right)v\right)\right)$$

then:

$$L[i]=n\left(\frac{\partial^{n-1}}{\partial x^{n-1}}\left(\left(\frac{\partial}{\partial x}u\right)v\right)\right)$$

Since this formula coincides with (1.7), the condition is satisfied.
Let's assume that (1.8) is valid for all $k=p$ values.

$$\sum_{i=0}^{n} i^p\,C(n,i)\left(\frac{\partial^i}{\partial x^i}u\right)\left(\frac{\partial^{n-i}}{\partial x^{n-i}}v\right)=$$
$$\sum_{i=0}^{p}\left(\sum_{s=0}^{i}\frac{s^p\,(-1)^{(s+i)}}{\Gamma(i-s+1)\,\Gamma(s+1)}\right)C(n,i)\,i!\left(\frac{\partial^{n-i}}{\partial x^{n-i}}\left(\left(\frac{\partial^i}{\partial x^i}u\right)v\right)\right) \tag{1.9}$$

Let's prove its validity for $k=p+1$.

$$\sum_{i=0}^{n} i^{(p+1)}\,C(n,i)\left(\frac{\partial^i}{\partial x^i}u\right)\left(\frac{\partial^{n-i}}{\partial x^{n-i}}v\right)=$$
$$\sum_{i=0}^{p+1}\left(\sum_{s=0}^{i}\frac{s^{(p+1)}\,(-1)^{(s+i)}}{\Gamma(i-s+1)\,\Gamma(s+1)}\right)C(n,i)\,i!\left(\frac{\partial^{n-i}}{\partial x^{n-i}}\left(\left(\frac{\partial^i}{\partial x^i}u\right)v\right)\right) \tag{1.10}$$

Indeed

$$\sum_{i=0}^{n} i^{(p+1)} \mathrm{C}(n,i)\left(\frac{\partial^i}{\partial x^i}u\right)\left(\frac{\partial^{n-i}}{\partial x^{n-i}}v\right) = \sum_{i=0}^{n} i^p\, i\, \mathrm{C}(n,i)\left(\frac{\partial^i}{\partial x^i}u\right)\left(\frac{\partial^{n-i}}{\partial x^{n-i}}v\right) =$$

$$\sum_{i=0}^{n}\left(\left(\frac{i^p\, n!\, i}{i!\,(n-i)!}\right)\cdot\left(\frac{\partial^i}{\partial x^i}u\right)\right)\left(\frac{\partial^{n-i}}{\partial x^{n-i}}v\right) = \sum_{i=0}^{n}\left((i^p\, n\, \mathrm{C}(n-1,i-1))\,.\left(\frac{\partial^i}{\partial x^i}u\right)\right)\left(\frac{\partial^{n-i}}{\partial x^{n-i}}v\right)$$

If we change the summation index $i$ to $i+1$ in the resulting value we receive:

$$\sum_{i=0}^{n} i^{(p+1)} \mathrm{C}(n,i)\left(\frac{\partial^i}{\partial x^i}u\right)\left(\frac{\partial^{n-i}}{\partial x^{n-i}}v\right) = n\left(\sum_{i=0}^{n-1} i^p\, \mathrm{C}(n-1,i)\, \mathrm{Diff}\left(\frac{\partial}{\partial x}u, x \$ i\right)\left(\frac{\partial^{n-1-i}}{\partial x^{n-1-i}}v\right)\right)$$

Since (1.9), by definition, is valid, we get the following:

$$\sum_{i=0}^{n-1} i^p\, \mathrm{C}(n-1,i)\, \mathrm{Diff}\left(\frac{\partial}{\partial x}u, x \$ i\right)\left(\frac{\partial^{n-1-i}}{\partial x^{n-1-i}}v\right) =$$

$$\sum_{i=0}^{p}\left(\sum_{s=0}^{i}\frac{s^p\,(-1)^{(s+i)}}{\Gamma(i-s+1)\,\Gamma(s+1)}\right)\mathrm{C}(n-1,i)\, i!\left(\frac{\partial^{n-1-i}}{\partial x^{n-1-i}}\left(\mathrm{Diff}\left(\frac{\partial}{\partial x}u, x \$ i\right)v\right)\right)$$

Thus,

$$\sum_{i=0}^{n} i^{(p+1)} \mathrm{C}(n,i)\left(\frac{\partial^i}{\partial x^i}u\right)\left(\frac{\partial^{n-i}}{\partial x^{n-i}}v\right) =$$

$$n\left(\sum_{i=0}^{p}\left(\sum_{s=0}^{i}\frac{s^p\,(-1)^{(s+i)}}{\Gamma(i-s+1)\,\Gamma(s+1)}\right)\mathrm{C}(n-1,i)\, i!\left(\frac{\partial^{n-1-i}}{\partial x^{n-1-i}}\left(\mathrm{Diff}\left(\frac{\partial}{\partial x}u, x \$ i\right)v\right)\right)\right)$$

If we change the summation index $i$ to $i+1$ in the right side, we obtain (1.10).
**The proposition is proved.**
Using the Proposition 1.3. proved above, let's prove the following

**Proposition 1.4. Let's assume that $L[\ ]$ is Leibniz operator for index $i$ of $u$, $v$ functions and $n$ parameter. Then the following applies for any $k$ in the range of $[0,N]$ :**

$$L[\mathrm{C}(i-z,k)] = \left(\frac{1}{k!}\right)\cdot\left(\sum_{s_0=0}^{k}\left(\sum_{s_1=0}^{s_0}\left(\sum_{s_2=0}^{s_1}\left(\sum_{s_3=0}^{s_2}\frac{s_3^{\,s_1}\,(-1)^{(s_3+s_2)}}{\Gamma(s_2-s_3+1)\,\Gamma(s_3+1)}\right)\mathrm{C}(n,s_2)\, s_2!\right.\right.\right.$$
$$\left.\left.\left. r_{k-s_0+1}(k)\, \mathrm{C}(s_0, s_0-s_1)\,(-z)^{(s_0-s_1)}\left(\frac{\partial^{n-s_2}}{\partial x^{n-s_2}}\left(\left(\frac{\partial^{s_2}}{\partial x^{s_2}}u\right)v\right)\right)\right)\right)\right) \qquad (1.11)$$

**Proof**: in accordance with the definition of Leibniz operator (1.11), (1.12) can be written down as follows:

$$\sum_{i=0}^{n} C(i-z,k)\,C(n,i)\left(\frac{\partial^{n-i}}{\partial x^{n-i}}v\right)\left(\frac{\partial^{i}}{\partial x^{i}}u\right)=\left(\frac{1}{k!}\right)\cdot\left(\sum_{s_0=0}^{k}\left(\sum_{s_1=0}^{s_0}\left(\sum_{s_2=0}^{s_1}\left(\sum_{s_3=0}^{s_2}\frac{s_3^{s_1}(-1)^{(s_3+s_2)}}{\Gamma(s_2-s_3+1)\,\Gamma(s_3+1)}\right)C(n,s_2)\,s_2!\,r_{k-s_0+1}(k)\,C(s_0,s_0-s_1)\,(-z)^{(s_0-s_1)}\left(\frac{\partial^{n-s_2}}{\partial x^{n-s_2}}\left(\left(\frac{\partial^{s_2}}{\partial x^{s_2}}u\right)v\right)\right)\right)\right)\right)$$

where

$$r_1(k)=1$$

$$=-C(k,2)$$

$$r_3(k)=\sum_{i_2=z}^{z+k-1}(i_2-z)\left(\sum_{i_1=i_2+1}^{z+k-1}(i_1-z)\right)=\frac{k\,(k-1)\,(k-2)\,(3\,k-1)}{24}$$

$$r_4(k)=-\left(\sum_{i_3=z}^{z+k-1}(i_3-z)\left(\sum_{i_2=i_3+1}^{z+k-1}(i_2-z)\left(\sum_{i_1=i_2+1}^{z+k-1}(i_1-z)\right)\right)\right)=-\frac{k^2\,(k-2)\,(k-3)\,(k-1)^2}{48}$$

.............................................................................................................................................

.

$$=\Gamma(k)$$

**Proof**: As soon as

$$=\frac{(i-z)!}{(i-z-k)!}=(i-z)\,(i-z-1)\,(i-z-2)\,(i-z-3)\,..\,i-z-k+1$$

the expression on the right defines the algebraic equation for the $i$ variable of $k$ degree where the roots are the numbers .

In this algebraic equation, unknown values are the $r_s(k)$ coefficients. These coefficients are defined on the basis of Vieta theorem for algebraic equations of $k$ degree, so as to obtain the following equation:

$$\sum_{s=1}^{k+1} r_s(k)\,(i-z)^{(k+1-s)}=(i-z)\,(i-z-1)\,(i-z-2)\,(i-z-3)\,..\,i-z-k+1 \qquad (1.13)$$

Let's define the values of the $r_s(k)$ coefficients. Since the roots of (1.13) are

$$i_1 = z, i_2 = z + 1, i_3 = z + 2 \; ....$$

then, in accordance with Vieta theorem for algebraic equations of $k$ degree, coefficients are defined by (1.14):

$$r_1(k) = 1$$

$$= -\mathrm{C}(k, 2)$$

$$r_3(k) = \sum_{i_2 = z}^{z+k-1} (i_2 - z)\left(\sum_{i_1 = i_2 + 1}^{z+k-1} (i_1 - z)\right) = \frac{k\,(k-1)\,(k-2)\,(3\,k-1)}{24}$$

$$r_4(k) = -\left(\sum_{i_3 = z}^{z+k-1} (i_3 - z)\left(\sum_{i_2 = i_3 + 1}^{z+k-1} (i_2 - z)\left(\sum_{i_1 = i_2 + 1}^{z+k-1} (i_1 - z)\right)\right)\right) = -\frac{k^2\,(k-2)\,(k-3)\,(k-1)^2}{48}$$

..............................................................................................................................................
..........

$$= \Gamma(k)$$

Thus, it is proved that the following equality holds true:

$$\mathrm{C}(i - z, k) = \left(\frac{1}{k!}\right) \cdot \left(\sum_{s=1}^{k+1} r_s(k)\,(i - z)^{(k+1-s)}\right) \qquad (1.15)$$

where the $r_s(k)$ coefficients are defined by (1.14).

Substituting (1.15) to the following expression:

$$\sum_{i=0}^{n} \mathrm{C}(i - z, k)\,\mathrm{C}(n, i)\left(\frac{\partial^{\,n-i}}{\partial x^{n-i}} v\right)\left(\frac{\partial^{\,i}}{\partial x^{i}} u\right)$$

we obtain:

$$\sum_{i=0}^{n} \left(\left(\frac{1}{k!}\right) \cdot \left(\sum_{s=1}^{k+1} r_s(k)\,(i - z)^{(k+1-s)}\right)\right) \mathrm{C}(n, i)\left(\frac{\partial^{\,n-i}}{\partial x^{n-i}} v\right)\left(\frac{\partial^{\,i}}{\partial x^{i}} u\right)$$

Then, after re-grouping we obtain the following:

$$\left(\frac{1}{k!}\right) \cdot \left(\sum_{s=1}^{k+1} r_s(k)\left(\sum_{i=0}^{n} (i - z)^{(k+1-s)}\,\mathrm{C}(n, i)\left(\frac{\partial^{\,n-i}}{\partial x^{n-i}} v\right)\left(\frac{\partial^{\,i}}{\partial x^{i}} u\right)\right)\right)$$

In accordance with Newton binomial formula, we obtain:

$$\sum_{i=0}^{n} (i-z)^{(k+1-s)} \mathrm{C}(n,i) \left(\frac{\partial^{n-i}}{\partial x^{n-i}} v\right)\left(\frac{\partial^{i}}{\partial x^{i}} u\right) = \sum_{i=0}^{n} \left(\sum_{l=0}^{k+1-s} \mathrm{C}(k+1-s,l)\, i^{l} (-z)^{(k+1-s-l)}\right) \mathrm{C}(n,i) \left(\frac{\partial^{n-i}}{\partial x^{n-i}} v\right)\left(\frac{\partial^{i}}{\partial x^{i}} u\right)$$

Let's represent this expression as follows:

$$\sum_{l=0}^{k+1-s} \mathrm{C}(k+1-s,l)\,(-z)^{(k+1-s-l)} \left(\sum_{i=0}^{n} i^{l}\, \mathrm{C}(n,i) \left(\frac{\partial^{n-i}}{\partial x^{n-i}} v\right)\left(\frac{\partial^{i}}{\partial x^{i}} u\right)\right) = \sum_{l=0}^{k+1-s} \mathrm{C}(k+1-s,l)\,(-z)^{(k+1-s-l)} L[i^{l}]$$

Since, in accordance with (1.7):

$$L[i^{l}] = \sum_{i=0}^{l} \left(\sum_{s_0=0}^{i} \frac{s_0^{\,l}\,(-1)^{(s_0+i)}}{\Gamma(i-s_0+1)\,\Gamma(s_0+1)}\right) \mathrm{C}(n,i)\, i! \left(\frac{\partial^{n-i}}{\partial x^{n-i}}\left(\left(\frac{\partial^{i}}{\partial x^{i}} u\right) v\right)\right)$$

by combining the results, we obtain:

$$\sum_{i=0}^{n} \mathrm{C}(i-z,k)\, \mathrm{C}(n,i) \left(\frac{\partial^{n-i}}{\partial x^{n-i}} v\right)\left(\frac{\partial^{i}}{\partial x^{i}} u\right) = \left(\frac{1}{k!}\right) \cdot \left(\sum_{s=1}^{k+1} r_s(k) \left(\sum_{l=0}^{k+1-s} \mathrm{C}(k+1-s,l)\right.\right.$$

$$\left.\left.(-z)^{(k+1-s-l)} \left(\sum_{i=0}^{l} \left(\sum_{s_0=0}^{i} \frac{s_0^{\,l}\,(-1)^{(s_0+i)}}{\Gamma(i-s_0+1)\,\Gamma(s_0+1)}\right) \mathrm{C}(n,i)\, i! \left(\frac{\partial^{n-i}}{\partial x^{n-i}}\left(\left(\frac{\partial^{i}}{\partial x^{i}} u\right) v\right)\right)\right)\right)\right)$$

After the change of summation index $s$ to $s+1$ and renumbering, we eventually obtain (1.11). **Proposition 1.4. is proved.**

Let's discuss properties of some Leibniz operators.

**Let's prove the Proposition 1.5.: Let's assume that** $b_0(x), b_1(x), b_2(x), b_3(x) \,..\, b_m(x)$ **is a sequence of functions, and** $L_{i_0,k_0}[\ ]$ is **Leibniz operator for** $i$ **index of** , $b_{k_0}(x)$ **functions and** $n$ **parameter. Then the following equation holds true:**

$$\sum_{i_0=0}^{m-1} L_{i_0,k_0}[\mathrm{C}(i-z,i_0)] = \sum_{i_1=0}^{m-1} \left(\sum_{i_0=0}^{m-1-i_1} (-1)^{i_0}\, \mathrm{C}(k_0-1+i_0,k_0-1)\, \mathrm{C}(n,i_1) \left(\frac{\partial^{n-i_1}}{\partial x^{n-i_1}}\left([b_{i_1+i_0}(x)]_{i_0+k_0}\ \ b_{k_0}(x)\right)\right)\right) \qquad (1.16)$$

**Proof**: Considering the Proposition 1.4 proved above, (1.16) transforms as follows:

$$\sum_{i_0=0}^{m-1}\left(\left(\frac{1}{i_0!}\right)\cdot\left(\sum_{s_0=0}^{i_0}\left(\sum_{s_1=0}^{s_0}\left(\sum_{s_2=0}^{s_1}\left(\sum_{s_3=0}^{s_2}\frac{s_3^{s_1}(-1)^{(s_3+s_2)}}{\Gamma(s_2-s_3+1)\Gamma(s_3+1)}\right)C(n,s_2)\,s_2!\,r_{1+i_0-s_0}(i_0)\right.\right.\right.\right.$$

$$\left.\left.\left.\left.C(s_0,s_0-s_1)(-k_0)^{(s_0-s_1)}\left(\frac{\partial^{n-s_2}}{\partial x^{n-s_2}}\left(\left(\frac{\partial^{s_2}}{\partial x^{s_2}}[b_{i_0}(x)]_{i_0+k_0}\right)b_{k_0}(x)\right)\right)\right)\right)\right)\right) =$$

$$\sum_{i_1=0}^{m-1}\left(\sum_{i_0=0}^{m-1-i_1}(-1)^{i_0}C(k_0-1+i_0,k_0-1)\,C(n,i_1)\left(\frac{\partial^{n-i_1}}{\partial x^{n-i_1}}\left([b_{i_1+i_0}(x)]_{i_0+k_0}\,b_{k_0}(x)\right)\right)\right)$$

**Proof**: With $m=0$ the equality is satisfied identically. Let's check the validity of this formula for $m=1$. Indeed, in this case:

$$\frac{\partial^n}{\partial x^n}\left([b_0(x)]_{k_0}\,b_{k_0}(x)\right)=\frac{\partial^n}{\partial x^n}\left([b_0(x)]_{k_0}\,b_{k_0}(x)\right)$$

Let's use mathematical induction. Assume that this formula holds true for $m=p$. Let's prove that (1.16) holds true for $m=p+1$ as well. Indeed:
(1.17)

- $$\sum_{i_1=0}^{p}\left(\sum_{i_0=0}^{p-i_1}(-1)^{i_0}C(k_0-1+i_0,k_0-1)\,C(n,i_1)\left(\frac{\partial^{n-i_1}}{\partial x^{n-i_1}}\left([b_{i_1+i_0}(x)]_{i_0+k_0}\,b_{k_0}(x)\right)\right)\right)=\Bigg($$

$$\sum_{i_1=0}^{p-1}\left(\sum_{i_0=0}^{p-i_1-1}(-1)^{i_0}C(k_0-1+i_0,k_0-1)\,C(n,i_1)\left(\frac{\partial^{n-i_1}}{\partial x^{n-i_1}}\left([b_{i_1+i_0}(x)]_{i_0+k_0}\,b_{k_0}(x)\right)\right)\right)$$

$$\Bigg)+C(n,p)\left(\frac{\partial^{n-p}}{\partial x^{n-p}}\left([b_p(x)]_{p+1}\,b_{p+1}(x)\right)\right)$$

Similarly (1.5.67)

$$\sum_{i_0=0}^{p}\left(\left(\frac{1}{i_0!}\right)\cdot\left(\sum_{s_0=0}^{i_0}\left(\sum_{s_1=0}^{s_0}\left(\sum_{s_2=0}^{s_1}\left(\sum_{s_3=0}^{s_2}\frac{s_3^{s_1}(-1)^{(s_3+s_2)}}{\Gamma(s_2-s_3+1)\Gamma(s_3+1)}\right)C(n,s_2)\,s_2!\,r_{1+i_0-s_0}(i_0)\right.\right.\right.\right.$$

$$\left.\left.\left.\left.C(s_0,s_0-s_1)(-k_0)^{(s_0-s_1)}\left(\frac{\partial^{n-s_2}}{\partial x^{n-s_2}}\left(\left(\frac{\partial^{s_2}}{\partial x^{s_2}}[b_{i_0}(x)]_{i_0+k_0}\right)b_{k_0}(x)\right)\right)\right)\right)\right)\right) = \left(\sum_{i_0=0}^{p-1}\left(\left(\frac{1}{i_0!}\right)\right.\right.$$

$$\cdot\left(\sum_{s_0=0}^{i_0}\left(\sum_{s_1=0}^{s_0}\left(\sum_{s_2=0}^{s_1}\left(\sum_{s_3=0}^{s_2}\frac{s_3^{s_1}(-1)^{(s_3+s_2)}}{\Gamma(s_2-s_3+1)\Gamma(s_3+1)}\right)C(n,s_2)\,s_2!\,r_{1+i_0-s_0}(i_0)\right.\right.\right.$$

$$\mathrm{C}(s_0, s_0 - s_1)\,(-k_0)^{(s_0 - s_1)}\left(\frac{\partial^{n-s_2}}{\partial x^{n-s_2}}\left(\left(\frac{\partial^{s_2}}{\partial x^{s_2}}[b_{i_0}(x)]_{i_0+k_0}\right) b_{k_0}(x)\right)\right)\Bigg)\Bigg)\Bigg)\Bigg)\Bigg)\Bigg)$$

$$+\ \mathrm{C}(n,p)\left(\frac{\partial^{n-p}}{\partial x^{n-p}}\left([b_p(x)]_{p+1}\ b_{p+1}(x)\right)\right)$$

Since, by definition

$$\sum_{i_1=0}^{p-1}\left(\sum_{i_0=0}^{p-i_1-1}(-1)^{i_0}\,\mathrm{C}(k_0-1+i_0, k_0-1)\,\mathrm{C}(n,i_1)\left(\frac{\partial^{n-i_1}}{\partial x^{n-i_1}}\left([b_{i_1+i_0}(x)]_{i_0+k_0}\ b_{k_0}(x)\right)\right)\right) =$$

$$\sum_{i_0=0}^{p-1}\left(\left(\frac{1}{i_0!}\right)\cdot\left(\sum_{s_0=0}^{i_0}\left(\sum_{s_1=0}^{s_0}\left(\sum_{s_2=0}^{s_1}\left(\sum_{s_3=0}^{s_2}\frac{s_3^{\,s_1}\,(-1)^{(s_3+s_2)}}{\Gamma(s_2-s_3+1)\,\Gamma(s_3+1)}\right)\mathrm{C}(n,s_2)\,s_2!\right.\right.\right.\right.$$

$$\left.\left.\left.\left. r_{1+i_0-s_0}(i_0)\,\mathrm{C}(s_0,s_0-s_1)\,(-k_0)^{(s_0-s_1)}\left(\frac{\partial^{n-s_2}}{\partial x^{n-s_2}}\left(\left(\frac{\partial^{s_2}}{\partial x^{s_2}}[b_{i_0}(x)]_{i_0+k_0}\right) b_{k_0}(x)\right)\right)\right)\right)\right)\right)$$

then, comparing the results of (1.16) and (1.17) we obtain the identical equality. **The proposition is proved.**

# 2. The general theory of integration of linear second-order ODEs.

Since the most common are the linear second-order ODEs, let's present the integration theory for them (for greatest clarity), which is adequate to the general approach.
In accordance with Definition 1.1, **the general algorithm for the integration of linear homogeneous second-order ODEs:**

$$g_0(x)\left(\frac{\partial^2}{\partial x^2}y\right) + g_1(x)\left(\frac{\partial}{\partial x}y\right) + g_2(x)\,y = 0 \qquad (2.1)$$

**where are the unknown functions, and $g_p(x)$ are the given functions differentiable and integrable to the required degree, shall be understood as the definition of formulas:**

(2.2)

**where $\alpha_k$ are in the form of functional series:**

$$\alpha_k = \sum_{j=0}^{\infty}\xi_{k,j}(g_0(x), g_1(x), g_2(x)),\ k = 1, 2 \qquad (2.3)$$

**of certain functions, which are determined through the initial coefficients:** .

**Thus, the task of obtaining a general algorithm to solve linear ODEs of (2.1) type consists of two main tasks:**

**- to define the presentation of $F_i(x, \alpha_i)$ functions for the desired partial solutions (2.2)**

**- to define the functions.**

## 2.1. Derivation of the presentation of partial solutions.

Later on we shall consider homogeneous linear ODEs of second order only in their reduced form:

$$\frac{\partial^2}{\partial x^2}y = a_1(x)\left(\frac{\partial}{\partial x}y\right) + a_2(x)\,y \qquad (2.4)$$

where $a_1(x)$, $a_2(x)$ are the specified functions differentiable and integrable to the required degree.

Let's differentiate (2.4). Then we obtain:

$$\frac{\partial^3}{\partial x^3}y = \left(\frac{d}{dx}a_1(x)\right)\left(\frac{\partial}{\partial x}y\right) + a_1(x)\left(\frac{\partial^2}{\partial x^2}y\right) + \left(\frac{d}{dx}a_2(x)\right)y + a_2(x)\left(\frac{\partial}{\partial x}y\right) \qquad (2.5)$$

Since the right side of this ODE contains the $\frac{\partial^2}{\partial x^2}y$ value, which is defined by (2.4), then substituting it to (2.5), after combining similar terms, we get:

$$\frac{\partial^3}{\partial x^3}y = \left(\left(\frac{d}{dx}a_1(x)\right) + a_1(x)^2 + a_2(x)\right)\left(\frac{\partial}{\partial x}y\right) + \left(a_1(x)\,a_2(x) + \left(\frac{d}{dx}a_2(x)\right)\right)y \qquad (2.6)$$

Let's differentiate (2.6). Then we obtain:

$$\begin{aligned}\frac{\partial^4}{\partial x^4}y = &\left(\left(\frac{d^2}{dx^2}a_1(x)\right) + 2\,a_1(x)\left(\frac{d}{dx}a_1(x)\right) + \left(\frac{d}{dx}a_2(x)\right)\right)\left(\frac{\partial}{\partial x}y\right)\\ &+ \left(\left(\frac{d}{dx}a_1(x)\right) + a_1(x)^2 + a_2(x)\right)\left(\frac{\partial^2}{\partial x^2}y\right)\\ &+ \left(\left(\frac{d}{dx}a_1(x)\right)a_2(x) + a_1(x)\left(\frac{d}{dx}a_2(x)\right) + \left(\frac{d^2}{dx^2}a_2(x)\right)\right)y\\ &+ \left(a_1(x)\,a_2(x) + \left(\frac{d}{dx}a_2(x)\right)\right)\left(\frac{\partial}{\partial x}y\right)\end{aligned}$$

Again, we establish that the right side of this equation contains $\frac{\partial^2}{\partial x^2}y$, which is defined by (2.4). Substituting it to this equation, after combining similar terms, we get: (2.7)

$$\begin{aligned}\frac{\partial^4}{\partial x^4}y = &\\ &\left(\left(\frac{d^2}{dx^2}a_1(x)\right) + 3\,a_1(x)\left(\frac{d}{dx}a_1(x)\right) + 2\left(\frac{d}{dx}a_2(x)\right) + a_1(x)^3 + 2\,a_1(x)\,a_2(x)\right)\left(\frac{\partial}{\partial x}y\right)\\ &+ \left(2\left(\frac{d}{dx}a_1(x)\right)a_2(x) + a_2(x)\,a_1(x)^2 + a_2(x)^2 + a_1(x)\left(\frac{d}{dx}a_2(x)\right) + \left(\frac{d^2}{dx^2}a_2(x)\right)\right)y\end{aligned}$$

Proceeding in the same way, we obtain: (2.8)

$$\frac{\partial^5}{\partial x^5}y = \left(\left(\frac{d^3}{dx^3}a_1(x)\right) + 3\left(\frac{d}{dx}a_1(x)\right)^2 + 4\,a_1(x)\left(\frac{d^2}{dx^2}a_1(x)\right) + 3\left(\frac{d^2}{dx^2}a_2(x)\right)\right.$$
$$+ 6\,a_1(x)^2\left(\frac{d}{dx}a_1(x)\right) + 4\left(\frac{d}{dx}a_1(x)\right)a_2(x) + 5\,a_1(x)\left(\frac{d}{dx}a_2(x)\right) + a_1(x)^4$$
$$\left.+ 3\,a_2(x)\,a_1(x)^2 + a_2(x)^2\right)\left(\frac{\partial}{\partial x}y\right) + \left(3\left(\frac{d^2}{dx^2}a_1(x)\right)a_2(x) + 5\,a_2(x)\,a_1(x)\left(\frac{d}{dx}a_1(x)\right)\right.$$
$$+ 4\,a_2(x)\left(\frac{d}{dx}a_2(x)\right) + a_2(x)\,a_1(x)^3 + 2\,a_1(x)\,a_2(x)^2 + 3\left(\frac{d}{dx}a_1(x)\right)\left(\frac{d}{dx}a_2(x)\right)$$
$$\left.+ \left(\frac{d}{dx}a_2(x)\right)a_1(x)^2 + a_1(x)\left(\frac{d^2}{dx^2}a_2(x)\right) + \left(\frac{d^3}{dx^3}a_2(x)\right)\right)y$$

and so on.
Analyzing (2.4) (2.6) (2.7) (2.8), we can establish a common pattern: for the ODEs of any $n$ order the general structure of the right-hand side remains unchanged and defined by the following formula:

$$\frac{\partial^{2+n}}{\partial x^{2+n}}y = \alpha(x,n)\left(\frac{\partial}{\partial x}y\right) + \beta(x,n)\,y \tag{2.9}$$

Here $\alpha(x,n)$, $\beta(x,n)$ are certain functions of an independent ***x*** variable and *n* parameter, which is a natural number.
Later on, for the sake of simplicity let's assume that:

$$\alpha(x,n) = \alpha(n), \beta(x,n) = \beta(n) \tag{2.10}$$

Then, particularly for $n = 1$ we get:

$$\frac{\partial^3}{\partial x^3}y = \alpha(1)\left(\frac{\partial}{\partial x}y\right) + \beta(1)\,y \tag{2.11}$$

Comparing this equality with (2.6), we establish that:

$$\alpha(1) = \left(\frac{d}{dx}a_1(x)\right) + a_1(x)^2 + a_2(x) \tag{2.12}$$

$$\beta(1) = a_1(x)\,a_2(x) + \left(\frac{d}{dx}a_2(x)\right) \tag{2.13}$$

Again, assuming that $n = 2$ in (2.9), we have:

$$\frac{\partial^4}{\partial x^4}y = \alpha(2)\left(\frac{\partial}{\partial x}y\right) + \beta(2)\,y \tag{2.14}$$

Comparing this ODE with (2.7), we establish that:

$$\alpha(2) = \left(\frac{d^2}{dx^2}a_1(x)\right) + 3\,a_1(x)\left(\frac{d}{dx}a_1(x)\right) + 2\left(\frac{d}{dx}a_2(x)\right) + a_1(x)^3 + 2\,a_1(x)\,a_2(x) \tag{2.15}$$

$$\beta(2) = 2\left(\frac{d}{dx}a_1(x)\right)a_2(x) + a_2(x)\,a_1(x)^2 + a_2(x)^2 + a_1(x)\left(\frac{d}{dx}a_2(x)\right) + \left(\frac{d^2}{dx^2}a_2(x)\right) \tag{2.16}$$

Continuing in the same way, we can define other required values for the $\alpha(n)$, $\beta(n)$ coefficients.

**Definition 2.1.: Equality (2.9) is called the *n*-image of (2.4), and $\alpha(n)$, $\beta(n)$ functions are called the coefficients of the *n*-image.**

## 2.2. Properties of *n*-image.

**Property 1: *n*-image (2.9) is a differential equation of $n+2$ order and contains a solution of (2.4) ODE.**

*This property is obvious, because the equation for* **n**-*image* (2.9) *is derived, provided that the initial equation* (2.4) *is a basis for a consistent series of differentiations for both sides of this equation, taking into account its own value* (2.4), *without use of other relations.*

**Property 2:** The $\alpha(n)$, $\beta(n)$ coefficients for the ***n***-image are defined by the initial values of:

, (2.17)

*Indeed, assuming that* $n = 0$ *in (2.9) we obtain:*

$$\frac{\partial^2}{\partial x^2} y = \alpha(0)\left(\frac{\partial}{\partial x} y\right) + \beta(0)\, y$$

*Comparing this equation with the initial one* - (2.4) - *we obtain the condition* (2.17).

**Property 3: The $\alpha(n)$, $\beta(n)$ coefficients for the *n*-image are related by the following equations:**

$$\alpha(n+1) = \left(\frac{\partial}{\partial x}\alpha(n)\right) + \alpha(n)\, a_1(x) + \beta(n) \qquad (2.18)$$

$$\beta(n+1) = \alpha(n)\, a_2(x) + \left(\frac{\partial}{\partial x}\beta(n)\right) \qquad (2.19)$$

**Proof.** ***Indeed, substituting*** $n+1$ ***for n parameter*** *in (2.9) we obtain:*

$$\frac{\partial^{3+n}}{\partial x^{3+n}} y = \alpha(n+1)\left(\frac{\partial}{\partial x} y\right) + \beta(n+1)\, y$$

*On the other hand, upon differentiation of (2.9) with consideration of (2.4) we obtain:*

$$\frac{\partial^{3+n}}{\partial x^{3+n}} y = \left(\left(\frac{\partial}{\partial x}\alpha(n)\right) + \beta(n)\right)\left(\frac{\partial}{\partial x} y\right) + \left(\alpha(n)\, a + \left(\frac{\partial}{\partial x}\beta(n)\right)\right) y$$

*Since the left sides of these equations are equal, we obtain (2.18), (2.19) from the equality of their right sides.*

**Property 4: If partial solutions of $y_1(x)$, $y_2(x)$ are known , then the $\alpha(n)$, $\beta(n)$ coefficients of *n*- image are uniquely defined by the following equations:**

$$\alpha(n) = \left(\frac{1}{W}\right)\cdot\left(\left(\frac{\partial^{n+2}}{\partial x^{n+2}} y_2\right) y_1 - y_2\left(\frac{\partial^{n+2}}{\partial x^{n+2}} y_1\right)\right) \qquad (2.20)$$

$$\beta(n) = \left(\frac{1}{W}\right) \cdot \left( \left(\frac{\partial^{n+2}}{\partial x^{n+2}} y_1\right)\left(\frac{\partial}{\partial x} y_2\right) - \left(\frac{\partial}{\partial x} y_1\right)\left(\frac{\partial^{n+2}}{\partial x^{n+2}} y_2\right)\right) \qquad (2.21)$$

***Proof**: Since* $y_1(x)$ *,* $y_2(x)$ *are independent partial solutions, they must satisfy the following equations:*

$$\frac{\partial^{2+n}}{\partial x^{2+n}} y_1 = \alpha(n)\left(\frac{\partial}{\partial x} y_1\right) + \beta(n)\, y_1$$

$$\frac{\partial^{2+n}}{\partial x^{2+n}} y_2 = \alpha(n)\left(\frac{\partial}{\partial x} y_2\right) + \beta(n)\, y_2$$

*Then, solving this algebraic system of equations, we obtain* **(2.20), (2.21).**

## 2.3. The structure of the general solution for the linear ODE (2.4).

Analyzing (2.18), we can establish that:

$$\beta(n) = \alpha(n+1) - \left(\frac{\partial}{\partial x}\alpha(n)\right) - \alpha(n)\, a_1(x) \quad (2.22)$$

Consequently, the equation of ***n***-image (2.9) can be represented as follows:

$$\frac{\partial^{2+n}}{\partial x^{2+n}} y = \alpha(n)\left(\frac{\partial}{\partial x} y\right) + \left(\alpha(n+1) - \left(\frac{\partial}{\partial x}\alpha(n)\right) - \alpha(n)\, a_1(x)\right) y \qquad (2.23)$$

or, similarly, establishing from (2.19) that

$$\alpha(n) = -\frac{-\beta(n+1) + \left(\frac{\partial}{\partial x}\beta(n)\right)}{a_2(x)}$$

we obtain a second view:

$$\frac{\partial^{2+n}}{\partial x^{2+n}} y = \left(-\frac{-\beta(n+1) + \left(\frac{\partial}{\partial x}\beta(n)\right)}{a_2(x)}\right)\left(\frac{\partial}{\partial x} y\right) + \beta(n)\, y \qquad (2.24)$$

As follows from the above, the equation (2.9) of n-image is completely determined by the value of $\alpha(n)$ coefficient, or the $\beta(n)$ coefficient. It is obvious that the representation (2.24) may have certain specificity with $a_2(x) = 0$ , and, therefore, in general, the representation (2.23) is more preferable as it does not depend on the value of the $a_2(x)$ coefficient and excludes special cases. Thus, later on we shall use the representation (2.23) as a determinant.

Let's define a recurrence equation to obtain the $a_2(x)$ coefficient. For this purpose, let's determine the $\beta(n)$ value defined by (2.22) and substitute the equality

$$\beta(n+1) = \alpha(n+2) - \left(\frac{\partial}{\partial x}\alpha(n+1)\right) - \alpha(n+1)\, a_1(x)$$

which follows from it, to (2.19).
Then we obtain the desired recurrence equation:

$$\alpha(2+n) = \alpha(n)\left(a_2(x) - \left(\frac{d}{dx}a_1(x)\right)\right) + 2\left(\frac{\partial}{\partial x}\alpha(n+1)\right) + \alpha(n+1)\,a_1(x) - \left(\frac{\partial}{\partial x}\alpha(n)\right)a_1(x) - \left(\frac{\partial^2}{\partial x^2}\alpha(n)\right) \tag{2.25}$$

With the use of this recurrence equation, let's define all values of $\alpha(n)$ function for arbitrary natural values of $n$ paramenter.

However, in this case consideration of the $\alpha(n)$ function does not give any new information useful for the obtainment of new mathematical results.
A new mathematical result requires a formal deviation from the established restrictions, because only in this case we can obtain non-trivial patterns. Therefore, based on this hypothesis, let's formally **extend the domain of $n$ parameter to negative integer numbers.** It turns out that such non-trivial, formal approach makes it possible to define the general structure for partial solutions of linear ODEs (2.4).

Thus, let's assume for (2.25) that $n = -1$. Then we obtain:

$$\alpha(1) = \alpha(-1)\left(a_2(x) - \left(\frac{d}{dx}a_1(x)\right)\right) + 2\left(\frac{d}{dx}\alpha(0)\right) + \alpha(0)\,a_1(x) - \left(\frac{d}{dx}\alpha(-1)\right)a_1(x) - \left(\frac{d^2}{dx^2}\alpha(-1)\right)$$

With consideration of (2.12) and initial conditions (2.17) we obtain a non-homogeneous ODE for the unknown $\alpha(-1)$ function:

$$\frac{d^2}{dx^2}\alpha(-1) = -a_1(x)\left(\frac{d}{dx}\alpha(-1)\right) + \left(a_2(x) - \left(\frac{d}{dx}a_1(x)\right)\right)\alpha(-1) + 2\left(\frac{d}{dx}a_1(x)\right) - \left(\frac{d}{dx}a_1(x)\right) - a_2(x)$$

Let's make a substitution in this equation:

(2.26)

where $V(x)$ is a certain function.

Then, upon re-expression, we obtain a homogeneous ODE to solve the $V(x)$ function:

$$\frac{d^2}{dx^2}V(x) = -\left(\frac{d}{dx}V(x)\right)a_1(x) + \left(a_2(x) - \left(\frac{d}{dx}a_1(x)\right)\right)V(x) \tag{2.27}$$

Let's prove **Proposition 2.1: A homogeneous ODE (2.27) is adjoined with the ODE (2.4).**

**Proof**: By definition, [Mathematical Encyclopedia. Publishing House "Soviet Encyclopedia", Moscow, 1979] (2.4) has the following adjoined equation

$$\frac{\partial^2}{\partial x^2}Y = -\left(\frac{\partial}{\partial x}(a_1(x)\,Y)\right) + a_2(x)\,Y$$

Upon re-expression of this ODE, we obtain:

$$\frac{\partial^2}{\partial x^2} Y = -a_1(x)\left(\frac{\partial}{\partial x} Y\right) + \left(a_2(x) - \left(\frac{d}{dx} a_1(x)\right)\right) Y \tag{2.28}$$

Comparing (2.28) and (2.27) we confirm identical structure of these equations, which proves this Proposition. Thus, if we introduce new functions:

(2.29)

$$b_2(x) = a_2(x) - \left(\frac{d}{dx} a_1(x)\right) \tag{2.30}$$

then (2.27) implies that **V(x)** is the solution of the following equation:

$$\frac{d^2}{dx^2} \mathrm{V}(x) = \left(\frac{d}{dx} \mathrm{V}(x)\right) b_1(x) + b_2(x)\,\mathrm{V}(x) \tag{2.31}$$

However, the structure of this ODE is identical to (2.4), and, since in accordance with the Proposition 2.1 proved above, (2.4) is adjoined with (2.27), then the following holds true:

**Theorem 2.1.: Assuming that**

$$\frac{d^{2+n}}{dx^{2+n}} \mathrm{V}(x) = \varepsilon(n)\left(\frac{d}{dx} \mathrm{V}(x)\right) + \rho(n)\,\mathrm{V}(x) \tag{2.32}$$

**is an *n*-image equation of the ODE (2.27), and $\varepsilon(n)$, $\rho(n)$ are the *n*-image coefficients where $b_1(x)$, $b_2(x)$ are defined by (2.29), (2.30).**
**Then the solution for (2.4) shall be the following function:**

(2.33)

**Proof**: Since the structural form of equations (2.4) and (2.31) is the same, by definition $\varepsilon(-1)$ satisfies the following equation:

$$\frac{d^2}{dx^2}\varepsilon(-1) = -\left(\frac{d}{dx} b_1(x)\right) - b_2(x) + \varepsilon(-1)\left(b_2(x) - \left(\frac{d}{dx} b_1(x)\right)\right) + 2\left(\frac{d}{dx} b_1(x)\right) - \left(\frac{d}{dx}\varepsilon(-1)\right) b_1(x) \tag{2.34}$$

Then, by expressing the $\varepsilon(-1)$ function from (2.33) and substituting its value to this equation, we obtain:

$$\frac{d^2}{dx^2} \mathrm{y}(x) = -b_1(x)\left(\frac{d}{dx} \mathrm{y}(x)\right) + \left(b_2(x) - \left(\frac{d}{dx} b_1(x)\right)\right) \mathrm{y}(x)$$

If we substitute here the values of $b_1(x)$ and $b_2(x)$ coefficients, in accordance with (2.29), (2.30) upon appropriate transformations we obtain:

$$\frac{d^2}{dx^2} \mathrm{y}(x) = a_1(x)\left(\frac{d}{dx} \mathrm{y}(x)\right) + a_2(x)\,\mathrm{y}(x)$$

Since this equation coincides with (2.4), **the theorem** is proved. Thus, we obtain the structural form to present the first partial solution of ODE (2.4).

Formally, the structural form for the second partial solution can be obtained with the use of a well-known formula, which presents the second partial solution through the first one, in the following form:

$$y_2(x) = (-1 + \varepsilon(-1, x)) \int \frac{\mathbf{e}^{\left(-\int a_1(x)\, dx\right)}}{(-1 + \varepsilon(-1, x))^2}\, dx \qquad (2.35)$$

However, let's obtain a different presentation.
Let's prove the **Theorem 2.2**:

**The second partial solution** $y_2(x)$ **of (2.4) is determined by the following structural form:**

$$y_2(x) = \varepsilon(-2) - x\,\varepsilon(-1) + x \qquad (2.36)$$

**where** $\varepsilon(-2)$ is a certain new function, and $\varepsilon(-1)$ is a function, that defines the structural form of the first partial solution (2.33).

**Proof**: Since (2.32) in its structural form is fully identical to (2.4), the following recurrence relation coinciding with (2.25) holds true for the definition of $\varepsilon(n)$:

$$\begin{aligned}\varepsilon(2 + n) = {} & \varepsilon(n)\left(b_2(x) - \left(\frac{d}{dx} b_1(x)\right)\right) + 2\left(\frac{\partial}{\partial x}\varepsilon(n + 1)\right) + \varepsilon(n + 1)\, b_1(x) \\ & - \left(\frac{\partial}{\partial x}\varepsilon(n)\right) b_1(x) - \left(\frac{\partial^2}{\partial x^2}\varepsilon(n)\right)\end{aligned} \qquad (2.37)$$

Let's assume for this recurrence equation that $n = -2$. Then we obtain:

$$\begin{aligned}\varepsilon(0) = {} & \varepsilon(-2)\left(b_2(x) - \left(\frac{d}{dx} b_1(x)\right)\right) + 2\left(\frac{d}{dx}\varepsilon(-1)\right) + \varepsilon(-1)\, b_1(x) - \left(\frac{d}{dx}\varepsilon(-2)\right) b_1(x) \\ & - \left(\frac{d^2}{dx^2}\varepsilon(-2)\right)\end{aligned} \qquad (2.38)$$

Since , and $\varepsilon(-1)$ is determined from (2.33) as

(2.39)

where the $y_1(x)$ function is the solution of ODE (2.4), substituting these equalities to (2.38), upon transformations we have:

$$\begin{aligned}b_1(x) = {} & \varepsilon(-2)\left(b_2(x) - \left(\frac{d}{dx} b_1(x)\right)\right) + 2\left(\frac{d}{dx}(1 + y_1(x))\right) + (1 + y_1(x))\, b_1(x) \\ & - \left(\frac{d}{dx}\varepsilon(-2)\right) b_1(x) - \left(\frac{d^2}{dx^2}\varepsilon(-2)\right)\end{aligned} \qquad (2.40)$$

Let's make a substitution in this equation:

$$\varepsilon(-2) = y_2(x) + x\,\varepsilon(-1) - x \qquad (2.41)$$

where $y_2(x)$ is a certain new function. Then, upon transformation, (2.40) reduces to the following form:

$$x\left(\left(\frac{d^2}{dx^2}y_1(x)\right)+b_1(x)\left(\frac{d}{dx}y_1(x)\right)-\left(b_2(x)-\left(\frac{d}{dx}b_1(x)\right)\right)y_1(x)\right)+\left(\frac{d^2}{dx^2}y_2(x)\right)$$
$$+b_1(x)\left(\frac{d}{dx}y_2(x)\right)+\left(\left(\frac{d}{dx}b_1(x)\right)-b_2(x)\right)y_2(x)=0$$

Substituting (2.29), (2.30) we get:

$$\left(-a_2(x)\,y_1(x)+\left(\frac{d^2}{dx^2}y_1(x)\right)-a_1(x)\left(\frac{d}{dx}y_1(x)\right)\right)x-a_1(x)\left(\frac{d}{dx}y_2(x)\right)-a_2(x)\,y_2(x)$$
$$+\left(\frac{d^2}{dx^2}y_2(x)\right)=0$$

Since, by definition, $y_1(x)$ is a solution of (2.4), we get the following equality:

$$\frac{d^2}{dx^2}y_2(x)=a_1(x)\left(\frac{d}{dx}y_2(x)\right)+a_2(x)\,y_2(x) \qquad (2.42)$$

However, this equality is also identical to (2.4), so the $y_2(x)$ function is a solution of ODE (2.4) as well. Thus, **the Theorem 2.2. is proved**.
Obviously, the obtained partial solutions (2.4)

(2.43)

$$y_2(x)=\varepsilon(-2)-x\,\varepsilon(-1)+x \qquad (2.44)$$

are linearly independent functions. Indeed, the Wronskian of these functions is:

$$\begin{bmatrix}-1+\varepsilon(-1) & \varepsilon(-2)-x\,\varepsilon(-1)+x\\ \frac{d}{dx}\varepsilon(-1) & \frac{d}{dx}(\varepsilon(-2)-x\,\varepsilon(-1)+x)\end{bmatrix}=$$
$$2\,\varepsilon(-1)-1+(\varepsilon(-1)-1)\left(\frac{d}{dx}\varepsilon(-2)\right)-\varepsilon(-1)^2-\left(\frac{d}{dx}\varepsilon(-1)\right)\varepsilon(-2)$$

This implies the proof of the independence of the obtained partial solutions. Thus, we have solved the first task: to construct an algorithm for the ODE (2.4). Namely, we have defined the structural form of partial solutions represented by (2.43) (2.44). Then the general solution (2.4) is defined by the following formula:

$$\mathrm{y}(x)=C_1\,(-1+\varepsilon(-1))+C_2\,(\varepsilon(-2)-x\,\varepsilon(-1)+x)$$

Or, equivalently:

$$\mathrm{y}(x)=(C_1-C_2\,x)\,(-1+\varepsilon(-1))+C_2\,\varepsilon(-2) \qquad (2.45)$$

Thus, in order to obtain a general solution (2.45) for the ODE (2.4), it is required to construct an adjoined ODE with coefficients determined by (2.29), (2.30), and then define the values of $\varepsilon(-1)$ and $\varepsilon(-2)$.
Essentially, it can appear that the resulting formula for the general solution (2.45) is not representative, because in principle it is possible to assume arbitrary form for the desired solution.
However, in this case a significant factor is the possibility of a strict and consistent definition of only one $\varepsilon(n)$ function, because the desired $\varepsilon(-1)$ and $\varepsilon(-2)$ values are obtained automatically, then $n$ assumes values .

Moreover, a recurrence equation (2.37) is available to determine the $\varepsilon(n)$ function. In principle, thanks to (2.35), the $\varepsilon(-1)$ and $\varepsilon(-2)$ functions are interrelated by the following formula:

$$\varepsilon(-2,x) = (-1+\varepsilon(-1,x))\left(\int \frac{\mathbf{e}^{\left(-\int a_1(x)\,dx\right)}}{(-1+\varepsilon(-1,x))^2}\,dx + x\right) \tag{2.46}$$

**2.4. Non-homogeneous ODEs.** In order to fully consider the calculation of a general solution for linear second-order ODEs, let's discuss the finding of partial solution - $Y(x)$ - for the non-homogeneous equation:

$$\left(\frac{d^2}{dx^2}Y(x)\right) + a_1(x)\left(\frac{d}{dx}Y(x)\right) + a_2(x)\,Y(x) = F(x) \tag{2.47}$$

where $F(x)$ is a specified function.
This task can be solved in different ways, but all of them require either a pre-calculation of a special function, such as the Green function, or availability of two partial solutions, i.e. they require preparatory operations. A number of other methods are employed as well. However, it turns out that the solution for this problem is much simpler, i.e. no preliminary operations are required to calculate a partial solution of a non-homogenous ODE, and all that is needed is one partial solution - $y(x)$ - of a homogeneous ODE:

$$\left(\frac{d^2}{dx^2}y(x)\right) + a_1(x)\left(\frac{d}{dx}y(x)\right) + a_2(x)\,y(x) = 0 \tag{2.48}$$

Let's prove the following
**Theorem: The linear non-homogenous differential equation (2.47) has a partial solution**

$$Y(x) = y(x)\int \frac{\mathbf{e}^{\left(-\int a_1(x)\,dx\right)}\int \mathbf{e}^{\left(\int a_1(x)\,dx\right)}\,y(x)\,F(x)\,dx}{y(x)^2}\,dx \tag{2.49}$$

**where $y(x)$ is a non-trivial partial solution of the homogeneous linear ODE (2.48).**
**Proof**: Since, by definition, the $Y(x)$ function is the solution of (2.47), then (2.47) must be true.
By the substitution of (2.49), upon transformations we obtain:

$$\left(\left(\frac{d^2}{dx^2}y(x)\right) + a_1(x)\left(\frac{d}{dx}y(x)\right) + a_2(x)\,y(x)\right)\int \frac{1\int \mathbf{e}^{\left(\int a_1(x)\,dx\right)}\,y(x)\,F(x)\,dx}{\mathbf{e}^{\left(\int a_1(x)\,dx\right)}\,y(x)^2}\,dx + F(x) = F(x)$$

Since, by definition, $y(x)$ satisfies (2.48), we get the following identity:

**The theorem is proved.**

## 3. Calculation of the coefficients of *n*-image.

In accordance with the Definition 2.1, to construct an algorithm for the finding of general solution of the ODE (2.4) we need to define formulas to calculate the functions. Eventually, they must be determined by (2.3). For the sake of simplicity, we shall accomplish this task, without loss of generality, for the equation of the following form:

$$\frac{\partial^2}{\partial x^2} y = a\, y \tag{3.1}$$

which is derived from (2.4) upon the following transformation

where

$$a = a_2(x) + \frac{a_1(x)^2}{4} - \frac{\frac{d}{dx} a_1(x)}{2}$$

Since in this case $a_1(x) = 0$, from (2.29) and (2.30) we get:

$b_1(x) = 0$ ;

Consequently, the adjoined equation is equal to the initial one (3.1) and, thus, in order to find a general solution for this equation it is sufficient to obtain the value of $\alpha(n)$ coefficient for the *n*-image of this equation:

$$\frac{\partial^{2+n}}{\partial x^{2+n}} y = \alpha(n) \left( \frac{\partial}{\partial x} y \right) + \beta(n)\, y \tag{3.2}$$

**Note**. The case where the $a_1(x)$ coefficient of (1.1) is not zero shall be automatically considered as a special case in the theory of finding a general solution of $m$ -order ODEs.
Let's discuss the approach to the finding of $\alpha(n)$ coefficient.
Since, in accordance with Leibniz theorem:

$$\frac{\partial^n}{\partial x^n}(a\, y) = \sum_{i=0}^{n} \left( \left[ \begin{matrix} i \\ n \end{matrix} \right] \cdot \left( \frac{\partial^{n-i}}{\partial x^{n-i}} a \right) \right) \left( \frac{\partial^i}{\partial x^i} y \right) \tag{3.3}$$

then the equation of ***n***-image assumes the following form:

$$\sum_{i=0}^{n} \left( \left[ \begin{matrix} i \\ n \end{matrix} \right] \cdot \left( \frac{\partial^{n-i}}{\partial x^{n-i}} a \right) \right) \left( \frac{\partial^i}{\partial x^i} y \right) = \tag{3.4}$$

Since the right side of this equation defines the $\alpha(n), \beta(n)$ coefficients as certain expressions at $\frac{\partial}{\partial x} y$ and $y$ , then we can represent the left side of this equation in the following form:

$$\left( \sum_{i=2}^{n} \left( \left[ \begin{matrix} i \\ n \end{matrix} \right] \cdot \left( \frac{\partial^{n-i}}{\partial x^{n-i}} a \right) \right) \left( \frac{\partial^i}{\partial x^i} y \right) \right) + \left( \left[ \begin{matrix} 1 \\ n \end{matrix} \right] \cdot \left( \frac{\partial^{n-1}}{\partial x^{n-1}} a \right) \right) \left( \frac{\partial}{\partial x} y \right) + \left( \left[ \begin{matrix} 0 \\ n \end{matrix} \right] \cdot \left( \frac{\partial^n}{\partial x^n} a \right) \right) y \tag{3.5}$$

Let's substitute $i+2$ for the summation index $i$ in the following expression:

$$\sum_{i=2}^{n}\left(\begin{bmatrix} i \\ n \end{bmatrix}\cdot\left(\frac{\partial^{n-i}}{\partial x^{n-i}}a\right)\right)\left(\frac{\partial^{i}}{\partial x^{i}}y\right) \tag{3.6}$$

Then it shall take the following form:

$$\sum_{i=0}^{n-2}\left(\begin{bmatrix} i+2 \\ n \end{bmatrix}\cdot\left(\frac{\partial^{n-i-2}}{\partial x^{n-i-2}}a\right)\right)\left(\frac{\partial^{i+2}}{\partial x^{i+2}}y\right) \tag{3.7}$$

Since

$$= \sum_{i_1=0}^{i}\left(\begin{bmatrix} i_1 \\ i \end{bmatrix}\cdot\left(\frac{\partial^{i-i_1}}{\partial x^{i-i_1}}a\right)\right)\left(\frac{\partial^{i_1}}{\partial x^{i_1}}y\right) \tag{3.8}$$

the expression (3.7) equals:

$$\sum_{i=0}^{n-2}\left(\begin{bmatrix} i+2 \\ n \end{bmatrix}\cdot\left(\frac{\partial^{n-i-2}}{\partial x^{n-i-2}}a\right)\right)\left(\frac{\partial^{i+2}}{\partial x^{i+2}}y\right) = \sum_{i=0}^{n-2}\left(\sum_{i_1=0}^{i}\left(\begin{bmatrix} i+2 \\ n \end{bmatrix}\cdot\begin{bmatrix} i_1 \\ i \end{bmatrix}\cdot\left(\frac{\partial^{i-i_1}}{\partial x^{i-i_1}}a\right)\right)\left(\frac{\partial^{n-i-2}}{\partial x^{n-i-2}}a\right)\left(\frac{\partial^{i_1}}{\partial x^{i_1}}y\right)\right)$$

The right-hand side of this equation is represented as follows:

$$\begin{aligned}&\left(\sum_{i=2}^{n-2}\left(\sum_{i_1=2}^{i}\left(\left(\begin{bmatrix} i+2 \\ n \end{bmatrix}\begin{bmatrix} i_1 \\ i \end{bmatrix}\right)\cdot\left(\frac{\partial^{i-i_1}}{\partial x^{i-i_1}}a\right)\right)\left(\frac{\partial^{n-i-2}}{\partial x^{n-i-2}}a\right)\left(\frac{\partial^{i_1}}{\partial x^{i_1}}y\right)\right)\right)\\ &\quad+\left(\sum_{i=0}^{n-2}\left(\left(\begin{bmatrix} i+2 \\ n \end{bmatrix}\begin{bmatrix} 0 \\ i \end{bmatrix}\right)\cdot\left(\frac{\partial^{i}}{\partial x^{i}}a\right)\right)\left(\frac{\partial^{n-i-2}}{\partial x^{n-i-2}}a\right)\right)y\\ &\quad+\left(\sum_{i=0}^{n-2}\left(\left(\begin{bmatrix} i+2 \\ n \end{bmatrix}\begin{bmatrix} 1 \\ i \end{bmatrix}\right)\cdot\left(\frac{\partial^{i-1}}{\partial x^{i-1}}a\right)\right)\left(\frac{\partial^{n-i-2}}{\partial x^{n-i-2}}a\right)\right)\left(\frac{\partial}{\partial x}y\right)\end{aligned} \tag{3.9}$$

Thus, with consideration of (3.6), (3.5) becomes equal to (3.9):

$$\begin{aligned}&\left(\sum_{i=2}^{n-2}\left(\sum_{i_1=2}^{i}\left(\left(\begin{bmatrix} i+2 \\ n \end{bmatrix}\begin{bmatrix} i_1 \\ i \end{bmatrix}\right)\cdot\left(\frac{\partial^{i-i_1}}{\partial x^{i-i_1}}a\right)\right)\left(\frac{\partial^{n-i-2}}{\partial x^{n-i-2}}a\right)\left(\frac{\partial^{i_1}}{\partial x^{i_1}}y\right)\right)\right)\\ &\quad+\left(\left(\sum_{i=0}^{n-2}\left(\left(\begin{bmatrix} i+2 \\ n \end{bmatrix}\begin{bmatrix} 0 \\ i \end{bmatrix}\right)\cdot\left(\frac{\partial^{i}}{\partial x^{i}}a\right)\right)\left(\frac{\partial^{n-i-2}}{\partial x^{n-i-2}}a\right)\right)+\left(\begin{bmatrix} 0 \\ n \end{bmatrix}\cdot\left(\frac{\partial^{n}}{\partial x^{n}}a\right)\right)\right)y\\ &\quad+\left(\left(\sum_{i=0}^{n-2}\left(\left(\begin{bmatrix} i+2 \\ n \end{bmatrix}\begin{bmatrix} 1 \\ i \end{bmatrix}\right)\cdot\left(\frac{\partial^{i-1}}{\partial x^{i-1}}a\right)\right)\left(\frac{\partial^{n-i-2}}{\partial x^{n-i-2}}a\right)\right)+\left(\begin{bmatrix} 1 \\ n \end{bmatrix}\cdot\left(\frac{\partial^{n-1}}{\partial x^{n-1}}a\right)\right)\right)\left(\frac{\partial}{\partial x}y\right)\end{aligned}$$

Again, let's substitute $i_1+2$ for the summation index $i_1$ in of the following expression:

$$\sum_{i=2}^{n-2}\left(\sum_{i_1=2}^{i}\left(\left(\begin{bmatrix} i+2 \\ n \end{bmatrix}\begin{bmatrix} i_1 \\ i \end{bmatrix}\right)\cdot\left(\frac{\partial^{i-i_1}}{\partial x^{i-i_1}}a\right)\right)\left(\frac{\partial^{n-i-2}}{\partial x^{n-i-2}}a\right)\left(\frac{\partial^{i_1}}{\partial x^{i_1}}y\right)\right) \tag{3.10}$$

Then we obtain: (3.11)

$$\sum_{i=2}^{n-2}\left(\sum_{i_1=0}^{i-2}\left(\left(\begin{bmatrix}i+2\\n\end{bmatrix}\begin{bmatrix}i_1+2\\i\end{bmatrix}\right)\cdot\left(\frac{\partial^{i-i_1-2}}{\partial x^{i-i_1-2}}a\right)\right)\left(\frac{\partial^{n-i-2}}{\partial x^{n-i-2}}a\right)\left(\frac{\partial^{i_1+2}}{\partial x^{i_1+2}}y\right)\right)$$

This expression, taking into account the ***n***-image, is represented as follows:

$$\sum_{i=2}^{n-2}\left(\sum_{i_1=0}^{i-2}\left(\left(\begin{bmatrix}i+2\\n\end{bmatrix}\begin{bmatrix}i_1+2\\i\end{bmatrix}\right)\cdot\left(\frac{\partial^{i-i_1-2}}{\partial x^{i-i_1-2}}a\right)\right)\left(\frac{\partial^{n-i-2}}{\partial x^{n-i-2}}a\right)\left(\frac{\partial^{i_1}}{\partial x^{i_1}}(a\,y)\right)\right)$$

Again using Leibniz formula for the $\frac{\partial^{i_1}}{\partial x^{i_1}}(a\,y)$ component, we have:

$$\sum_{i=2}^{n-2}\left(\sum_{i_1=0}^{i-2}\left(\sum_{i_2=0}^{i_1}\left(\left(\left(\begin{bmatrix}i+2\\n\end{bmatrix}\cdot\begin{bmatrix}i_1+2\\i\end{bmatrix}\cdot\begin{bmatrix}i_2\\i_1\end{bmatrix}\cdot\left(\frac{\partial^{i-i_1-2}}{\partial x^{i-i_1-2}}a\right)\right)\left(\frac{\partial^{n-i-2}}{\partial x^{n-i-2}}a\right)\cdot\left(\frac{\partial^{i_1-i_2}}{\partial x^{i_1-i_2}}a\right)\right)\left(\frac{\partial^{i_2}}{\partial x^{i_2}}y\right)\right)\right)\right)$$

This expression can be converted to an equivalent form:

$$\sum_{i_2=2}^{i_1}\Bigg(\Bigg(\Bigg(\begin{bmatrix}i+2\\n\end{bmatrix}\cdot\begin{bmatrix}i_1+2\\i\end{bmatrix}\cdot\begin{bmatrix}i_2\\i_1\end{bmatrix}\cdot\left(\frac{\partial^{i-i_1-2}}{\partial x^{i-i_1-2}}a\right)\Bigg)\left(\frac{\partial^{n-i-2}}{\partial x^{n-i-2}}a\right)\cdot\left(\frac{\partial^{i_1-i_2}}{\partial x^{i_1-i_2}}a\right)\Bigg)\left(\frac{\partial^{i_2}}{\partial x^{i_2}}y\right)\Bigg)$$

$$+\left(\sum_{i=2}^{n-2}\left(\sum_{i_1=0}^{i-2}\left(\left(\begin{bmatrix}i+2\\n\end{bmatrix}\cdot\begin{bmatrix}i_1+2\\i\end{bmatrix}\cdot\begin{bmatrix}0\\i_1\end{bmatrix}\cdot\left(\frac{\partial^{i-i_1-2}}{\partial x^{i-i_1-2}}a\right)\right)\left(\frac{\partial^{n-i-2}}{\partial x^{n-i-2}}a\right)\cdot\left(\frac{\partial^{i_1}}{\partial x^{i_1}}a\right)\right)\right)\right)y$$

$$+\left(\sum_{i=2}^{n-2}\left(\sum_{i_1=0}^{i-2}\left(\left(\begin{bmatrix}i+2\\n\end{bmatrix}\cdot\begin{bmatrix}i_1+2\\i\end{bmatrix}\cdot\begin{bmatrix}1\\i_1\end{bmatrix}\cdot\left(\frac{\partial^{i-i_1-2}}{\partial x^{i-i_1-2}}a\right)\right)\left(\frac{\partial^{n-i-2}}{\partial x^{n-i-2}}a\right)\cdot\left(\frac{\partial^{i_1-1}}{\partial x^{i_1-1}}a\right)\right)\right)\right)\left(\frac{\partial}{\partial x}y\right)$$

Thus, (3.9) can be reduced as follows:

$$\sum_{i_2=2}^{i_1}\Bigg(\Bigg(\Bigg(\begin{bmatrix}i+2\\n\end{bmatrix}\cdot\begin{bmatrix}i_1+2\\i\end{bmatrix}\cdot\begin{bmatrix}i_2\\i_1\end{bmatrix}\cdot\left(\frac{\partial^{i-i_1-2}}{\partial x^{i-i_1-2}}a\right)\Bigg)\left(\frac{\partial^{n-i-2}}{\partial x^{n-i-2}}a\right)\cdot\left(\frac{\partial^{i_1-i_2}}{\partial x^{i_1-i_2}}a\right)\Bigg)\left(\frac{\partial^{i_2}}{\partial x^{i_2}}y\right)\Bigg)+\Bigg($$

$$\left(\sum_{i=2}^{n-2}\left(\sum_{i_1=0}^{i-2}\left(\left(\begin{bmatrix}i+2\\n\end{bmatrix}\cdot\begin{bmatrix}i_1+2\\i\end{bmatrix}\cdot\begin{bmatrix}0\\i_1\end{bmatrix}\cdot\left(\frac{\partial^{i-i_1-2}}{\partial x^{i-i_1-2}}a\right)\right)\left(\frac{\partial^{n-i-2}}{\partial x^{n-i-2}}a\right)\cdot\left(\frac{\partial^{i_1}}{\partial x^{i_1}}a\right)\right)\right)\right)$$

$$+\left(\sum_{i=0}^{n-2}\left(\left(\begin{bmatrix}i+2\\n\end{bmatrix}\begin{bmatrix}0\\i\end{bmatrix}\right)\cdot\left(\frac{\partial^{i}}{\partial x^{i}}a\right)\right)\left(\frac{\partial^{n-i-2}}{\partial x^{n-i-2}}a\right)+\left(\begin{bmatrix}0\\n\end{bmatrix}\cdot\left(\frac{\partial^{n}}{\partial x^{n}}a\right)\right)\right)y+\Bigg($$

$$\left(\sum_{i=2}^{n-2}\left(\sum_{i_1=0}^{i-2}\left(\left(\left(\begin{bmatrix} i+2\\ n\end{bmatrix}\cdot\begin{bmatrix} i_1+2\\ i\end{bmatrix}\cdot\begin{bmatrix}1\\ i_1\end{bmatrix}\cdot\left(\frac{\partial^{i-i_1-2}}{\partial x^{i-i_1-2}}a\right)\right)\left(\frac{\partial^{n-i-2}}{\partial x^{n-i-2}}a\right)\cdot\left(\frac{\partial^{i_1-1}}{\partial x^{i_1-1}}a\right)\right)\right)\right)$$

$$+\left(\sum_{i=0}^{n-2}\left(\left(\begin{bmatrix} i+2\\ n\end{bmatrix}\begin{bmatrix} i_1\\ i\end{bmatrix}\right)\cdot\left(\frac{\partial^{i-1}}{\partial x^{i-1}}a\right)\right)\left(\frac{\partial^{n-i-2}}{\partial x^{n-i-2}}a\right)+\left(\begin{bmatrix}1\\ n\end{bmatrix}\cdot\left(\frac{\partial^{n-1}}{\partial x^{n-1}}a\right)\right)\right)\left(\frac{\partial}{\partial x}y\right)$$

Continuing in the same manner, at the ***k***-th step we obtain: (3.12)

$$\sum_{i=2k}^{n-2}\left(\sum_{i_1=2(k-1)}^{i-2}\left(\sum_{i_2=2(k-1)}^{i_1-2}\left(\sum_{i_{k-1}=2}^{i_{k-2}-2}\left(\sum_{i_k=2}^{i_{k-1}}\left(\right.\right.\right.\right.\right.$$

$$\begin{bmatrix} i+2\\ n\end{bmatrix}\cdot\begin{bmatrix} i_1+2\\ i\end{bmatrix}\cdot\begin{bmatrix} i_{k-1}+2\\ i_{k-2}\end{bmatrix}\cdot\begin{bmatrix} i_k\\ i_{k-1}\end{bmatrix}\cdot\left(\frac{\partial^{n-i-2}}{\partial x^{n-i-2}}a\right)\cdot\left(\frac{\partial^{i-i_1-2}}{\partial x^{i-i_1-2}}a\right)\,..$$

$$\left(\left(\frac{\partial^{i_{k-2}-i_{k-1}-2}}{\partial x^{i_{k-2}-i_{k-1}-2}}a\right)\cdot\left(\frac{\partial^{i_{k-1}-i_k}}{\partial x^{i_{k-1}-i_k}}a\right)\right)\left(\frac{\partial^{i_k}}{\partial x^{i_k}}y\right)\left.\left.\left.\left.\left.\right)\right)\right)\right)\right) + \left(\left(\sum_{i=2k}^{n-2}\left(\sum_{i_1=2(k-1)}^{i-2}\left(\sum_{i_{k-1}=2}^{i_{k-2}-2}\left(\sum_{i_k=0}^{i_{k-1}-2}\left(\right.\right.\right.\right.\right.\right.$$

$$\left(\begin{bmatrix} i+2\\ n\end{bmatrix}\cdot\begin{bmatrix} i_1+2\\ i\end{bmatrix}\cdot\begin{bmatrix} i_k+2\\ i_{k-1}\end{bmatrix}\cdot\begin{bmatrix} 0\\ i_k\end{bmatrix}\cdot\left(\frac{\partial^{n-i-2}}{\partial x^{n-i-2}}a\right)\right)\left(\frac{\partial^{i-i_1-2}}{\partial x^{i-i_1-2}}a\right)\,..$$

$$\left(\frac{\partial^{i_{k-1}-i_k-2}}{\partial x^{i_{k-1}-i_k-2}}a\right)\cdot\left(\frac{\partial^{i_k}}{\partial x^{i_k}}a\right)\left.\left.\left.\left.\left.\left.\right)\right)\right)\right)\right)\right) + ?\,..$$

$$\left(\sum_{i=2}^{n-2}\left(\sum_{i_1=0}^{i-2}\left(\left(\left(\begin{bmatrix} i+2\\ n\end{bmatrix}\cdot\begin{bmatrix} i_1+2\\ i\end{bmatrix}\cdot\begin{bmatrix}0\\ i_1\end{bmatrix}\cdot\left(\frac{\partial^{i-i_1-2}}{\partial x^{i-i_1-2}}a\right)\right)\left(\frac{\partial^{n-i-2}}{\partial x^{n-i-2}}a\right)\cdot\left(\frac{\partial^{i_1}}{\partial x^{i_1}}a\right)\right)\right)\right)\right)$$

$$+\left(\sum_{i=0}^{n-2}\left(\left(\begin{bmatrix} i+2\\ n\end{bmatrix}\begin{bmatrix} 0\\ i\end{bmatrix}\right)\cdot\left(\frac{\partial^{i}}{\partial x^{i}}a\right)\right)\left(\frac{\partial^{n-i-2}}{\partial x^{n-i-2}}a\right)+\left(\begin{bmatrix}0\\ n\end{bmatrix}\cdot\left(\frac{\partial^{n}}{\partial x^{n}}a\right)\right)\right)y\;+$$

$$\left(\left(\sum_{i=2k}^{n-2}\left(\sum_{i_1=2(k-1)}^{i-2}\left(\sum_{i_{k-1}=2}^{i_{k-2}-2}\left(\sum_{i_k=0}^{i_{k-1}-2}\left(\left(\begin{bmatrix} i+2\\ n\end{bmatrix}\cdot\begin{bmatrix} i_1+2\\ i\end{bmatrix}\cdot\begin{bmatrix} i_k+2\\ i_{k-1}\end{bmatrix}\cdot\begin{bmatrix} 1\\ i_k\end{bmatrix}\cdot\left(\frac{\partial^{n-i-2}}{\partial x^{n-i-2}}a\right)\right)\left(\frac{\partial^{i-i_1-2}}{\partial x^{i-i_1-2}}a\right)\,..\right.\right.\right.\right.\right.\right.$$

$$\left(\frac{\partial^{i_{k-1}-i_k-2}}{\partial x^{i_{k-1}-i_k-2}}a\right)\cdot\left(\frac{\partial^{i_k-1}}{\partial x^{i_k-1}}a\right)\left.\left.\left.\left.\left.\left.\right)\right)\right)\right)\right)\right) + ?\,..$$

$$\left(\sum_{i=2}^{n-2}\left(\sum_{i_1=0}^{i-2}\left(\left(\left(\begin{bmatrix} i+2\\ n\end{bmatrix}\cdot\begin{bmatrix} i_1+2\\ i\end{bmatrix}\cdot\begin{bmatrix}1\\ i_1\end{bmatrix}\cdot\left(\frac{\partial^{i-i_1-2}}{\partial x^{i-i_1-2}}a\right)\right)\left(\frac{\partial^{n-i-2}}{\partial x^{n-i-2}}a\right)\cdot\left(\frac{\partial^{i_1-1}}{\partial x^{i_1-1}}a\right)\right)\right)\right)\right)$$

$$+\left(\sum_{i=0}^{n-2}\left(\left(\begin{bmatrix} i+2\\ n\end{bmatrix}\begin{bmatrix} 1\\ i\end{bmatrix}\right)\cdot\left(\frac{\partial^{i-1}}{\partial x^{i-1}}a\right)\right)\left(\frac{\partial^{n-i-2}}{\partial x^{n-i-2}}a\right)+\left(\begin{bmatrix}1\\ n\end{bmatrix}\cdot\left(\frac{\partial^{n-1}}{\partial x^{n-1}}a\right)\right)\right)\left(\frac{\partial}{\partial x}y\right)$$

Let's select such value of $k$ that the summand in (3.12): (3.13)

$$\sum_{i=2k}^{n-2}\Bigg(\sum_{i_1=2(k-1)}^{i-2}\Bigg(\sum_{i_2=2(k-1)}^{i_1-2}\Bigg(\sum_{i_{k-1}=2}^{i_{k-2}-2}\Bigg(\sum_{i_k=2}^{i_{k-1}}\Bigg( \begin{bmatrix} i+2\\ n\end{bmatrix}\cdot\begin{bmatrix} i_1+2\\ i\end{bmatrix}\cdot\begin{bmatrix} i_{k-1}+2\\ i_{k-2}\end{bmatrix}\cdot\begin{bmatrix} i_k\\ i_{k-1}\end{bmatrix}\cdot\Big(\frac{\partial^{n-i-2}}{\partial x^{n-i-2}}a\Big)\cdot\Big(\frac{\partial^{i-i_1-2}}{\partial x^{i-i_1-2}}a\Big)\cdot\cdot \Big(\Big(\frac{\partial^{i_{k-2}-i_{k-1}-2}}{\partial x^{i_{k-2}-i_{k-1}-2}}a\Big)\cdot\Big(\frac{\partial^{i_{k-1}-i_k}}{\partial x^{i_{k-1}-i_k}}a\Big)\Big)\Big(\frac{\partial^{i_k}}{\partial x^{i_k}}y\Big)\Bigg)\Bigg)\Bigg)\Bigg)\Bigg)\Bigg)$$

would assume a certain univalent finite value. For this purpose it is required to have the $n = 2p$ parameter where $p$ is a natural number.
In this case, the value of $k$ is determined from the equation . From here it follows that $k = p - 1$ . For the given value of $k$ the expression (3.13) is:

$$\frac{\partial^2}{\partial x^2}y = a\,y$$

Therefore, the remaining two summands from (3.12) are defined by the equations: (3.14)

$$\Bigg(\Bigg(\sum_{i=2(p-1)}^{2p-2}\Bigg(\sum_{i_1=2(p-2)}^{i-2}\Bigg(\sum_{i_{p-2}=2}^{i_{p-3}-2}\Bigg(\sum_{i_{p-1}=0}^{i_{p-2}-2}\Bigg( \Bigg(\begin{bmatrix} i+2\\ 2p\end{bmatrix}\cdot\begin{bmatrix} i_1+2\\ i\end{bmatrix}\cdot\begin{bmatrix} i_{p-1}+2\\ i_{p-2}\end{bmatrix}\cdot\begin{bmatrix} 0\\ i_{p-1}\end{bmatrix}\cdot\Big(\frac{\partial^{2p-i-2}}{\partial x^{2p-i-2}}a\Big)\Bigg)\Bigg)\Big(\frac{\partial^{i-i_1-2}}{\partial x^{i-i_1-2}}a\Big)\cdot\cdot \Big(\frac{\partial^{i_{p-2}-i_{p-1}-2}}{\partial x^{i_{p-2}-i_{p-1}-2}}a\Big)\cdot\Big(\frac{\partial^{i_{p-1}}}{\partial x^{i_{p-1}}}a\Big)\Bigg)\Bigg)\Bigg)\Bigg)\Bigg)\Bigg) + ?\cdot\cdot$$

$$\Bigg(\sum_{i=2}^{2p-2}\Bigg(\sum_{i_1=0}^{i-2}\Bigg(\Bigg(\Bigg(\Bigg(\begin{bmatrix} i+2\\ 2p\end{bmatrix}\cdot\begin{bmatrix} i_1+2\\ i\end{bmatrix}\cdot\begin{bmatrix} 0\\ i_1\end{bmatrix}\cdot\Big(\frac{\partial^{i-i_1-2}}{\partial x^{i-i_1-2}}a\Big)\Bigg)\Bigg)\Big(\frac{\partial^{2p-i-2}}{\partial x^{2p-i-2}}a\Big)\Bigg)\cdot\Big(\frac{\partial^{i_1}}{\partial x^{i_1}}a\Big)\Bigg)\Bigg)\Bigg)\Bigg)$$

$$+\Bigg(\sum_{i=0}^{2p-2}\Bigg(\Bigg(\Bigg(\begin{bmatrix} i+2\\ 2p\end{bmatrix}\begin{bmatrix} 0\\ i\end{bmatrix}\Bigg)\cdot\Big(\frac{\partial^{i}}{\partial x^{i}}a\Big)\Bigg)\Big(\frac{\partial^{2p-i-2}}{\partial x^{2p-i-2}}a\Big)\Bigg) + \Bigg(\begin{bmatrix} 0\\ 2p\end{bmatrix}\cdot\Big(\frac{\partial^{2p}}{\partial x^{2p}}a\Big)\Bigg) + a\Bigg)y\ +$$

$$\Bigg(\Bigg(\sum_{i=2(p-1)}^{2p-2}\Bigg(\sum_{i_1=2(p-2)}^{i-2}\Bigg(\sum_{i_{p-2}=2}^{i_{p-3}-2}\Bigg(\sum_{i_{p-1}=0}^{i_{p-2}-2}\Bigg( \Bigg(\begin{bmatrix} i+2\\ 2p\end{bmatrix}\cdot\begin{bmatrix} i_1+2\\ i\end{bmatrix}\cdot\begin{bmatrix} i_{p-1}+2\\ i_{p-2}\end{bmatrix}\cdot\begin{bmatrix} 1\\ i_{p-1}\end{bmatrix}\cdot\Big(\frac{\partial^{2p-i-2}}{\partial x^{2p-i-2}}a\Big)\Bigg)\Bigg)\Big(\frac{\partial^{i-i_1-2}}{\partial x^{i-i_1-2}}a\Big)\cdot\cdot \Big(\frac{\partial^{i_{p-2}-i_{p-1}-2}}{\partial x^{i_{p-2}-i_{p-1}-2}}a\Big)\cdot\Big(\frac{\partial^{i_{p-1}-1}}{\partial x^{i_{p-1}-1}}a\Big)\Bigg)\Bigg)\Bigg)\Bigg)\Bigg)\Bigg) + ?\cdot\cdot\Bigg($$

$$\sum_{i=2}^{2p-2}\left(\sum_{i_1=0}^{i-2}\left(\left(\left(\begin{bmatrix}i+2\\2p\end{bmatrix}\cdot\begin{bmatrix}i_1+2\\i\end{bmatrix}\cdot\begin{bmatrix}1\\i_1\end{bmatrix}\cdot\left(\frac{\partial^{i-i_1-2}}{\partial x^{i-i_1-2}}a\right)\right)\left(\frac{\partial^{2p-i-2}}{\partial x^{2p-i-2}}a\right)\right)\cdot\left(\frac{\partial^{i_1-1}}{\partial x^{i_1-1}}a\right)\right)\right)$$
$$+\left(\sum_{i=0}^{2p-2}\left(\left(\begin{bmatrix}i+2\\2p\end{bmatrix}\begin{bmatrix}0\\i\end{bmatrix}\right)\cdot\left(\frac{\partial^{i}}{\partial x^{i}}a\right)\right)\left(\frac{\partial^{2p-i-2}}{\partial x^{2p-i-2}}a\right)+\left(\begin{bmatrix}1\\2p\end{bmatrix}\cdot\left(\frac{\partial^{2p-1}}{\partial x^{2p-1}}a\right)\right)\right)\left(\frac{\partial}{\partial x}y\right)$$

In accordance with the formula of $n$-image, the expression (3.13) is taken into account as a summand in (3.14), so (3.14) must be equal to the following expression:

$$\alpha(2p)\left(\frac{\partial}{\partial x}y\right)+\beta(2p)\,y$$

Thus, equating the coefficients with $\frac{\partial}{\partial x}y$ and $y$, we obtain the desired values of the coefficients for the ***n***-image: (3.15)

$$\alpha(2p)=\left(\sum_{i=2(p-1)}^{2p-2}\left(\sum_{i_1=2(p-2)}^{i-2}\left(\sum_{i_{p-2}=2}^{i_{p-3}-2}\left(\sum_{i_{p-1}=0}^{i_{p-2}-2}\left(\right.\right.\right.\right.\right.$$
$$\left(\begin{bmatrix}i+2\\2p\end{bmatrix}\cdot\begin{bmatrix}i_1+2\\i\end{bmatrix}\cdot\begin{bmatrix}i_{p-1}+2\\i_{p-2}\end{bmatrix}\cdot\begin{bmatrix}1\\i_{p-1}\end{bmatrix}\cdot\left(\frac{\partial^{2p-i-2}}{\partial x^{2p-i-2}}a\right)\right)\left(\frac{\partial^{i-i_1-2}}{\partial x^{i-i_1-2}}a\right)..$$
$$\left.\left.\left.\left.\left.\left.\left(\frac{\partial^{i_{p-2}-i_{p-1}-2}}{\partial x^{i_{p-2}-i_{p-1}-2}}a\right)\cdot\left(\frac{\partial^{i_{p-1}-1}}{\partial x^{i_{p-1}-1}}a\right)\right)\right)\right)\right)\right)\right)+\,?..\left(\right.$$
$$\sum_{i=2}^{2p-2}\left(\sum_{i_1=0}^{i-2}\left(\left(\left(\begin{bmatrix}i+2\\2p\end{bmatrix}\cdot\begin{bmatrix}i_1+2\\i\end{bmatrix}\cdot\begin{bmatrix}1\\i_1\end{bmatrix}\cdot\left(\frac{\partial^{i-i_1-2}}{\partial x^{i-i_1-2}}a\right)\right)\left(\frac{\partial^{2p-i-2}}{\partial x^{2p-i-2}}a\right)\right)\cdot\left(\frac{\partial^{i_1-1}}{\partial x^{i_1-1}}a\right)\right)\right)$$
$$+\left(\sum_{i=1}^{2p-2}\left(\left(\begin{bmatrix}i+2\\2p\end{bmatrix}\begin{bmatrix}1\\i\end{bmatrix}\right)\cdot\left(\frac{\partial^{i-1}}{\partial x^{i-1}}a\right)\right)\left(\frac{\partial^{2p-i-2}}{\partial x^{2p-i-2}}a\right)\right)+\left(\begin{bmatrix}1\\2p\end{bmatrix}\cdot\left(\frac{\partial^{2p-1}}{\partial x^{2p-1}}a\right)\right)$$

$$\beta(2p)=\left(\sum_{i=2(p-1)}^{2p-2}\left(\sum_{i_1=2(p-2)}^{i-2}\left(\sum_{i_{p-2}=2}^{i_{p-3}-2}\left(\sum_{i_{p-1}=0}^{i_{p-2}-2}\left(\right.\right.\right.\right.\right.$$
$$\left(\begin{bmatrix}i+2\\2p\end{bmatrix}\cdot\begin{bmatrix}i_1+2\\i\end{bmatrix}\cdot\begin{bmatrix}i_{p-1}+2\\i_{p-2}\end{bmatrix}\cdot\begin{bmatrix}0\\i_{p-1}\end{bmatrix}\cdot\left(\frac{\partial^{2p-i-2}}{\partial x^{2p-i-2}}a\right)\right)\left(\frac{\partial^{i-i_1-2}}{\partial x^{i-i_1-2}}a\right)..$$
$$\left.\left.\left.\left.\left.\left.\left(\frac{\partial^{i_{p-2}-i_{p-1}-2}}{\partial x^{i_{p-2}-i_{p-1}-2}}a\right)\cdot\left(\frac{\partial^{i_{p-1}}}{\partial x^{i_{p-1}}}a\right)\right)\right)\right)\right)\right)\right)+\,?..$$
$$\left(\sum_{i=2}^{2p-2}\left(\sum_{i_1=0}^{i-2}\left(\left(\left(\begin{bmatrix}i+2\\2p\end{bmatrix}\cdot\begin{bmatrix}i_1+2\\i\end{bmatrix}\cdot\begin{bmatrix}0\\i_1\end{bmatrix}\cdot\left(\frac{\partial^{i-i_1-2}}{\partial x^{i-i_1-2}}a\right)\right)\left(\frac{\partial^{2p-i-2}}{\partial x^{2p-i-2}}a\right)\right)\cdot\left(\frac{\partial^{i_1}}{\partial x^{i_1}}a\right)\right)\right)\right.$$
$$\left.+\left(\sum_{i=0}^{2p-2}\left(\left(\begin{bmatrix}i+2\\2p\end{bmatrix}\begin{bmatrix}0\\i\end{bmatrix}\right)\cdot\left(\frac{\partial^{i}}{\partial x^{i}}a\right)\right)\left(\frac{\partial^{2p-i-2}}{\partial x^{2p-i-2}}a\right)\right)+\left(\begin{bmatrix}0\\2p\end{bmatrix}\cdot\left(\frac{\partial^{2p}}{\partial x^{2p}}a\right)\right)\right)+a$$

By way of summation index substitution, the obtained coefficients can be represented in equivalent form:

$$\alpha(2p+1) = 2p\left(\frac{\partial^{2p-1}}{\partial x^{2p-1}}a\right) + \left(\sum_{i_1=0}^{2p-2} i_1\left(\frac{\partial^{i_1-1}}{\partial x^{i_1-1}}a\right)\mathrm{C}(2p, i_1+2)\left(\frac{\partial^{2p-2-i_1}}{\partial x^{2p-2-i_1}}a\right)\right) +$$

$$\left(\sum_{i_2=0}^{2p-4}\left(\sum_{i_1=0}^{i_2} i_1\left(\frac{\partial^{i_1-1}}{\partial x^{i_1-1}}a\right)\mathrm{C}(i_2+2, i_1+2)\left(\frac{\partial^{i_2-i_1}}{\partial x^{i_2-i_1}}a\right)\right)\mathrm{C}(2p, i_2+4)\left(\frac{\partial^{2p-4-i_2}}{\partial x^{2p-4-i_2}}a\right)\right)$$

$$\left(\sum_{i_2=0}^{i_3}\left(\sum_{i_1=0}^{i_2} i_1\left(\frac{\partial^{i_1-1}}{\partial x^{i_1-1}}a\right)\mathrm{C}(i_2+2, i_1+2)\left(\frac{\partial^{i_2-i_1}}{\partial x^{i_2-i_1}}a\right)\right)\mathrm{C}(i_3+4, i_2+4)\left(\frac{\partial^{i_3-i_2}}{\partial x^{i_3-i_2}}a\right)\right)$$

$$\mathrm{C}(2p, i_3+6)\left(\frac{\partial^{2p-6-i_3}}{\partial x^{2p-6-i_3}}a\right)\Bigg) + \left(\sum_{i_4=0}^{2p-8}\left(\sum_{i_3=0}^{i_4}\right.\right.$$

$$\left(\sum_{i_2=0}^{i_3}\left(\sum_{i_1=0}^{i_2} i_1\left(\frac{\partial^{i_1-1}}{\partial x^{i_1-1}}a\right)\mathrm{C}(i_2+2, i_1+2)\left(\frac{\partial^{i_2-i_1}}{\partial x^{i_2-i_1}}a\right)\right)\mathrm{C}(i_3+4, i_2+4)\left(\frac{\partial^{i_3-i_2}}{\partial x^{i_3-i_2}}a\right)\right)$$

$$\mathrm{C}(i_4+6, i_3+6)\left(\frac{\partial^{i_4-i_3}}{\partial x^{i_4-i_3}}a\right)\Bigg)\mathrm{C}(2p, i_4+8)\left(\frac{\partial^{2p-8-i_4}}{\partial x^{2p-8-i_4}}a\right)\Bigg) + \ldots.$$

$$\beta(2p+1) = \left(\frac{\partial^{2p}}{\partial x^{2p}}a\right) + \left(\sum_{i_1=0}^{2p-2}\left(\frac{\partial^{i_1}}{\partial x^{i_1}}a\right)\mathrm{C}(2p, i_1+2)\left(\frac{\partial^{2p-2-i_1}}{\partial x^{2p-2-i_1}}a\right)\right)$$

$$+\left(\sum_{i_2=0}^{2p-4}\left(\sum_{i_1=0}^{i_2}\left(\frac{\partial^{i_1}}{\partial x^{i_1}}a\right)\mathrm{C}(i_2+2, i_1+2)\left(\frac{\partial^{i_2-i_1}}{\partial x^{i_2-i_1}}a\right)\right)\mathrm{C}(2p, i_2+4)\left(\frac{\partial^{2p-4-i_2}}{\partial x^{2p-4-i_2}}a\right)\right) +$$

$$\left(\sum_{i_3=0}^{2p-6}\left(\sum_{i_2=0}^{i_3}\left(\sum_{i_1=0}^{i_2}\left(\frac{\partial^{i_1}}{\partial x^{i_1}}a\right)\mathrm{C}(i_2+2, i_1+2)\left(\frac{\partial^{i_2-i_1}}{\partial x^{i_2-i_1}}a\right)\right)\mathrm{C}(i_3+4, i_2+4)\left(\frac{\partial^{i_3-i_2}}{\partial x^{i_3-i_2}}a\right)\right)\right.$$

$$\mathrm{C}(2p, i_3+6)\left(\frac{\partial^{2p-6-i_3}}{\partial x^{2p-6-i_3}}a\right)\Bigg) + \left(\sum_{i_4=0}^{2p-8}\left(\sum_{i_3=0}^{i_4}\right.\right.$$

$$\left(\sum_{i_2=0}^{i_3}\left(\sum_{i_1=0}^{i_2}\left(\frac{\partial^{i_1}}{\partial x^{i_1}}a\right)\mathrm{C}(i_2+2, i_1+2)\left(\frac{\partial^{i_2-i_1}}{\partial x^{i_2-i_1}}a\right)\right)\mathrm{C}(i_3+4, i_2+4)\left(\frac{\partial^{i_3-i_2}}{\partial x^{i_3-i_2}}a\right)\right)$$

$$\mathrm{C}(i_4+6, i_3+6)\left(\frac{\partial^{i_4-i_3}}{\partial x^{i_4-i_3}}a\right)\Bigg)\mathrm{C}(2p, i_4+8)\left(\frac{\partial^{2p-8-i_4}}{\partial x^{2p-8-i_4}}a\right)\Bigg) + \ldots.$$

These coefficients fully determine the equation of ***n***-image in the following form:

$$\frac{\partial^{2p+2}}{\partial x^{2p+2}}y = \alpha(2p)\left(\frac{\partial}{\partial x}y\right) + \beta(2p)\,y \qquad (3.16)$$

## 4. Investigation of the structure of $\alpha(2p)$ coefficients

Previously, it was proved that in order to define a general solution of a linear ODE, it is necessary and sufficient to know only the $\alpha(2p)$ or $\alpha(2p+1)$ function. Therefore, further on we shall consider only these functions; in particular, it shall be sufficient to consider the $\alpha(2p)$ function, which is defined as a function of a given coefficient of the initial ODE and its derivatives, which solves the problem of their complete definition. However, it is extremely difficult to use these coefficients in this form, because they are defined as a sum of summands, each (except for the first) being a partial summable sequence. The task is to determine the sum of these partial summable sequences as functions of ***p***, because they can be used in the scope of this theory only in this way. Therefore, let's represent the $\alpha(2p)$ coefficients as follows:

(4.1)

Here

$$\xi_0(2p) = 2p\left(\frac{\partial^{2p-1}}{\partial x^{2p-1}}a\right) \tag{4.2}$$

$$\xi_1(2p) = \sum_{i_1=0}^{2p-2} i_1\left(\frac{\partial^{i_1-1}}{\partial x^{i_1-1}}a\right)C(2p, i_1+2)\left(\frac{\partial^{2p-2-i_1}}{\partial x^{2p-2-i_1}}a\right) \tag{4.3}$$

$$\xi_2(2p) = \sum_{i_2=0}^{2p-4}\left(\sum_{i_1=0}^{i_2} i_1\left(\frac{\partial^{i_1-1}}{\partial x^{i_1-1}}a\right)C(i_2+2, i_1+2)\left(\frac{\partial^{i_2-i_1}}{\partial x^{i_2-i_1}}a\right)\right)C(2p, i_2+4)\left(\frac{\partial^{2p-4-i_2}}{\partial x^{2p-4-i_2}}a\right) \tag{4.4}$$

$$\left(\sum_{i_2=0}^{i_3}\left(\sum_{i_1=0}^{i_2} i_1\left(\frac{\partial^{i_1-1}}{\partial x^{i_1-1}}a\right)C(i_2+2, i_1+2)\left(\frac{\partial^{i_2-i_1}}{\partial x^{i_2-i_1}}a\right)\right)C(i_3+4, i_2+4)\left(\frac{\partial^{i_3-i_2}}{\partial x^{i_3-i_2}}a\right)\right) \tag{4.5}$$

$$C(2p, i_3+6)\left(\frac{\partial^{2p-6-i_3}}{\partial x^{2p-6-i_3}}a\right)$$

$$\xi_4(2p) = \sum_{i_4=0}^{2p-8}\left(\sum_{i_3=0}^{i_4}\left(\sum_{i_2=0}^{i_3}\left(\sum_{i_1=0}^{i_2} i_1\left(\frac{\partial^{i_1-1}}{\partial x^{i_1-1}}a\right)C(i_2+2, i_1+2)\left(\frac{\partial^{i_2-i_1}}{\partial x^{i_2-i_1}}a\right)\right)\right.\right.$$

$$\left.\left.C(i_3+4, i_2+4)\left(\frac{\partial^{i_3-i_2}}{\partial x^{i_3-i_2}}a\right)C(i_4+6, i_3+6)\left(\frac{\partial^{i_4-i_3}}{\partial x^{i_4-i_3}}a\right)\right)C(2p, i_4+8)\left(\frac{\partial^{2p-8-i_4}}{\partial x^{2p-8-i_4}}a\right)\right) \tag{4.6}$$

and so on.

$$\xi_s(2p) = \sum_{i_1=[i_s, i_1]}^{[2p-2s, i_s, i_2]} i_1\left(\frac{\partial^{i_1-1}}{\partial x^{i_1-1}}a\right)\left(\prod_{k=2}^{s-1} C(i_k+2(k-1), i_1+2(k-1))\left(\frac{\partial^{i_k-i_{k-1}}}{\partial x^{i_k-i_{k-1}}}a\right)\right) \tag{4.7}$$

Analysis of the obtained equations makes it possible to establish their main property, which is defined as **the Proposition 4.1.: Each $\xi_k(2p)$ summand, starting from $\xi_1(2p)$, is determined by the previous value of $\xi_{k-1}(2p)$ according to the following formula:**

$$\xi_k(2p) = \sum_{i_k=0}^{2p-2k} \xi_{k-1}(i_k+2k-2)\, C(2p, i_k+2k) \left(\frac{\partial^{2p-2k-i_k}}{\partial x^{2p-2k-i_k}} a\right) \tag{4.8}$$

**Proof**: Upon the replacement of index in (4.2) according to $p = \frac{i_1}{2}$ we obtain:

(4.9)

Then $\xi_1(2p)$ can be represented as follows:

$$\xi_1(2p) = \sum_{i_1=0}^{2p-2} \xi_0(i_1)\, C(2p, i_1+2) \left(\frac{\partial^{2p-3-i_1}}{\partial x^{2p-3-i_1}} a\right) \tag{4.10}$$

Again, making the parameter substitution according to $p = \frac{i_2}{2}$ in (7.3) we obtain:

$$\xi_1(i_2) = \sum_{i_1=0}^{i_2} i_1 \left(\frac{\partial^{i_1-1}}{\partial x^{i_1-1}} a\right) C(i_2+2, i_1+2) \left(\frac{\partial^{i_2-i_1}}{\partial x^{i_2-i_1}} a\right)$$

Then $\xi_2(2p)$ can be represented as follows:

$$\xi_2(2p) = \sum_{i_2=0}^{2p-4} \xi_1(i_2+2)\, C(2p, i_2+4) \left(\frac{\partial^{2p-5-i_2}}{\partial x^{2p-5-i_2}} a\right) \tag{4.11}$$

Quite similarly, for the representations (7.5), (7.6) we have:

$$\xi_3(2p) = \sum_{i_3=0}^{2p-6} \xi_2(i_3+4)\, C(2p, i_3+6) \left(\frac{\partial^{2p-6-i_3}}{\partial x^{2p-6-i_3}} a\right) \tag{4.12}$$

$$\xi_4(2p) = \sum_{i_4=0}^{2p-8} \xi_3(i_4+6)\, C(2p, i_4+8) \left(\frac{\partial^{2p-8-i_4}}{\partial x^{2p-8-i_4}} a\right) \tag{4.13}$$

Then, for the representation (7.8) we have:

$$\xi_s(2p) = \sum_{i_s=0}^{2p-2s} \xi_{s-1}(i_s+2s-2)\, C(2p, i_s+2s) \left(\frac{\partial^{2p-2s-i_s}}{\partial x^{2p-2s-i_s}} a\right) \tag{4.14}$$

It is obvious that

$$\xi_{s-1}(2p) = \sum_{i_{s-1}=0}^{2p-2s+2} \xi_{s-2}(i_{s-1}+2s-4)\, C(2p, i_{s-1}+2s-2)\left(\frac{\partial^{2p-2s+2-i_{s-1}}}{\partial x^{2p-2s+2-i_{s-1}}} a\right) \quad (4.15)$$

$$\xi_{s-2}(2p) = \sum_{i_{s-2}=0}^{2p-2s+4} \xi_{s-3}(i_{s-2}+2s-6)\, C(2p, i_{s-2}+2s-4)\left(\frac{\partial^{2p-2s+4-i_{s-2}}}{\partial x^{2p-2s+4-i_{s-2}}} a\right) \quad (4.16)$$

---

$$\xi_{2}(2p) = \sum_{i_{2}=0}^{2p-4} \xi_{1}(i_{2}+2)\, C(2p, i_{2}+4)\left(\frac{\partial^{2p-5-i_{2}}}{\partial x^{2p-5-i_{2}}} a\right)$$

$$\xi_{1}(i_2) = \sum_{i_{1}=0}^{i_2} i_1 \left(\frac{\partial^{i_1-1}}{\partial x^{i_1-1}} a\right) C(i_2+2, i_1+2)\left(\frac{\partial^{i_2-i_1}}{\partial x^{i_2-i_1}} a\right)$$

By way of successive substitution of (4.15), (4.16) and so on up to (4.13), (4.12), (4.11), (4.10) and (4.9) to (4.14), we finally obtain (4.7).
**Proposition 4.1 is proved.**

If we change $p$ parameter to $p+\frac{1}{2}$, (4.8) takes the following form:

$$\xi_{k}(2p+1) = \sum_{i_{k}=0}^{2p-2k+1} \xi_{k-1}(i_{k}+2k-2)\, C(2p+1, i_{k}+2k)\left(\frac{\partial^{2p-2k-i_{k}+1}}{\partial x^{2p-2k-i_{k}+1}} a\right) \quad (4.17)$$

That is, we obtained a formula to find the odd values of $\xi_k(2p+1)$.
Further on, we shall call the form of $\xi_k(2p)$ summands definition in the representation (4.8) **integral**. Obviously, there is an opposite representation as well. After the replacement of the summation index: and, further, after the substitution of in (4.8) and (4.17) we obtain the following equations:

$$\xi_{p-s}(2p) = \sum_{i=0}^{2s} \xi_{p-s-1}(2p-i-2)\, C(2p, i)\left(\frac{\partial^{i}}{\partial x^{i}} a\right) \quad (4.18)$$

$$\xi_{p-s}(-1+2p) = \sum_{i=0}^{-1+2s} \xi_{p-s-1}(2p-3-i)\, C(-1+2p, 2p-1-i)\left(\frac{\partial^{i}}{\partial x^{i}} a\right) \quad (4.19)$$

......

This representation shall be further called **differential**. While the integral representation (4.8) allows to determine $\xi_0(2p)$, $\xi_1(2p)$, $\xi_2(2p)$ ... $\xi_p(2p)$ in sequence, the differential representation (4.18), (4.19) allows sequential definition of $\xi_p(2p)$, $\xi_{p-1}(2p)$, $\xi_{p-2}(2p)$ ... $\xi_0(2p)$, $\xi_p(-1+2p)$, $\xi_{p-1}(-1+2p)$ .... . **Thus, the coefficients of *n*-image can be calculated in two principally equal ways.**

**Let's analyze the obtained results for the $\xi_k(2p)$ summands:**

1) Number of the $\xi_k(2p)$ summands equals $p-1$,
2) Only the value of the zero, that is, first summand is known:

$$\xi_0(2p) = 2p\left(\frac{\partial^{2p-1}}{\partial x^{2p-1}} a\right),$$

3) For the calculation of the $\xi_k(2p)$, $k$ **-1,2,3....** summands it is required to find partial sums, in accordance with (4.8), and for the calculation of $\xi_{p-1}(2p)$, $\xi_{p-2}(2p)$ ... $\xi_1(2p)$ in accordance with (4.17) it is required to solve recurrence equations successively, based on the following initial condition: $\xi_p(2p) = 0$.
**The integral method** is as follows:
- calculation of the summands from the known zero $\xi_0(2p)$ to $\xi_1(2p)$ according to the following formula:

$$\xi_1(2p) = \sum_{i_1=0}^{2p-2} \xi_0(i_1)\, \mathrm{C}(2p, i_1+2)\left(\frac{\partial^{2p-2-i_1}}{\partial x^{2p-2-i_1}} a\right)$$

from $\xi_1(2p)$ to $\xi_2(2p)$ according to the following formula:

$$\xi_2(2p) = \sum_{i_2=0}^{2p-4} \xi_1(i_2+2)\, \mathrm{C}(2p, i_2+4)\left(\frac{\partial^{2p-4-i_2}}{\partial x^{2p-4-i_2}} a\right)$$

and so on, in accordance with (4.8).

**The differential method** is inverse, as follows:

- with the use of (4.18) we demonstrate that

$$\xi_p(2p) = 0 \qquad (4.20)$$

**Indeed**: Let's assume in this formula $s = 0$. Then it takes the following form:

$$\xi_p(2p) = \xi_{p-1}(-2+2p)\, a$$

Assuming that

$$\xi_p(2p) = \mathrm{b}(p), \quad = \mathrm{b}(p-1)$$

we obtain the following recurrence equation:

Solution for this equation is the following function:

Since

$$\xi_0(2p) = 2p\left(\frac{\partial^{-1+2p}}{\partial x^{-1+2p}} a\right)$$

then

$$\xi_0(0) = 0$$

and, thus,

$$\xi_p(2p) = 0$$

which we needed to prove.

Let's assume $\xi_p(2p) = 0$ in (4.18). Then, upon transformation, this formula shall take the following form:

$$\xi_{p-1}(2p) = \xi_{p-2}(-2+2p)\,a + 2.\,\xi_{p-2}(-3+2p)\,p\left(\frac{\partial}{\partial x} a\right) + \xi_{p-2}(-4+2p)\,\mathrm{C}(2.p, 2.)\left(\frac{\partial^2}{\partial x^2} a\right) \tag{4.21}$$

Let's replace $p$ parameter by $p-2$ in (4.20). Then we obtain the following equation:

$= 0$

Since

$=$

then (4.21) takes the following form:

$$\xi_{p-1}(2p) = [\xi_{p-1}(2p)]_{p=p-1}\, a + 2\,\xi_{p-2}(-3+2p)\,p\left(\frac{\partial}{\partial x} a\right) \tag{4.22}$$

For the calculation of we shall use a formula, which follows from (4.19). Assuming that $s = 1$, we obtain:

$$\xi_p(2p+1) = \xi_{p-1}(-2+2p)\,(2p+1)\left(\frac{\partial}{\partial x} a\right) + \xi_{p-1}(-1+2p)\,a$$

Since, according to (4.20), , we have the following recurrence equation:

$$\xi_p(2p+1) = \xi_{p-1}(-1+2p)\,a$$

General solution for this equation is presented in the following form:

(4.23)

As follows from the formula

$$\xi_0(2p) = 2p\left(\frac{\partial^{-1+2p}}{\partial x^{-1+2p}}a\right) \tag{4.24}$$

with $p = \frac{1}{2}$

$$\xi_0(1) = a$$

and then, from (4.23), we finally obtain:

(4.25)

here, replacing $p$ parameter by $p-2$ we have:

and, thus, from (4.22) we obtain the following recurrence equation:

$$\xi_{p-1}(2p) = \xi_{p-2}(-2+2p)\,a + 2.\,a^{(p-1)}\,p\left(\frac{\partial}{\partial x}a\right)$$

The solution of this equation is the following function:

$$\xi_{p-1}(2p) = a^{(p-1)}\,p\,(p+1)\left(\frac{\partial}{\partial x}a\right) + \xi_{-1}(0)\,a^p \tag{4.26}$$

Assuming that $p = 1$ in this equation we obtain:

$$\xi_0(2) = 2\left(\frac{\partial}{\partial x}a\right) + \xi_{-1}(0)\,a \tag{4.27}$$

Since, according to (4.24)

from (4.27) we find:

$$\xi_{-1}(0) = -\frac{-2\left(\frac{\partial}{\partial x}a\right) + 2\left(\frac{\partial}{\partial x}a\right)}{a}$$

Thus, the general formula (4.26) has the following form:

$$\xi_{p-1}(2p) = a^{(p-1)}\,p\,(p+1)\left(\frac{\partial}{\partial x}a\right) \tag{4.28}$$

In the same way we get: (4.29)

$$\xi_{p-2}(2p) = \frac{a^{(p-2)} p^2 (p-1)(p+1)\left(\frac{\partial^3}{\partial x^3} a\right)}{3} + \frac{2\, a^{(p-3)} p^2 (p-1)(p-2)(p+1)\left(\frac{\partial}{\partial x} a\right)\left(\frac{\partial^2}{\partial x^2} a\right)}{3}$$
$$+ \frac{1\, a^{(p-4)} p^2 (p-1)(p-2)(p-3)(p+1)\left(\frac{\partial}{\partial x} a\right)^3}{6}$$

and so on.

**Conclusions: The obtained formulas, opposite to the integral approach, enable us to establish the value of the coefficient in only one given point, for each function of the sequence:**

$$\xi_1(2p),\ \xi_2(2p),\ \xi_3(2p)\ \ldots.\ \ \xi_{p-1}(2p)$$

**so here $p$ parameter can not have negative values, since in such case we receive non-existing coefficients (for example: assuming $p = -1$ in $\xi_{p-1}(2p)$ $\xi_{p-1}(2p)$ we obtain $\xi_{-2}(-2)$, but such values of the coefficient are not required for further exposition of the theory, i.e. they do not make sense). Thus, this approach can be used to develop a theory of differential equations, only provided that the general term is defined as follows: ...... In this case, by way of $p = k + s$ replacement, we arrive to the standard expression: , which is also meaningful for negative values of the $2k - 2s$ argument. This basically defines a task of the same complexity as with the integral approach, which shall be preferred.**

**Note**: As soon as (4.8) and (4.18) are equivalent, the following formula holds true:

$$\sum_{i=0}^{2p-2k} \left( \xi_{k-1}(i + 2k - 2)\, \mathrm{C}(2p, i + 2k)\left(\frac{\partial^{2p-2k-i}}{\partial x^{2p-2k-i}} a\right) - \xi_{k-1}(2p - i - 2)\, \mathrm{C}(2p, 2p - i)\left(\frac{\partial^{i}}{\partial x^{i}} a\right)\right) = 0$$

## 5. Integral method for finding $\xi_k(2p)$.

Let's give a practical view of the algorithm to determine the $\xi_k(2p)$, $k = 1, 2, 3$ ..... summands.

**5.1. Finding the $\xi_1(2p)$ summand**. Since $\xi_0(2p)$ is known, in accordance with (4.8) for $k = 1$ we have:

$$\xi_1(2p) = \sum_{i=0}^{-2+2p} \xi_0(i)\, \mathrm{C}(2p, i+2)\left(\frac{\partial^{2p-2-i}}{\partial x^{2p-2-i}} a\right) \qquad (5.1)$$

Since

$$\xi_0(2p) = 2p\left(\frac{\partial^{2p-1}}{\partial x^{2p-1}} a\right)$$

then

Then

$$\xi_1(2p) = \sum_{i=0}^{2p-2} i \left( \frac{\partial^{i-1}}{\partial x^{i-1}} a \right) \mathrm{C}(2p, i+2) \left( \frac{\partial^{2p-2-i}}{\partial x^{2p-2-i}} a \right) \tag{5.2}$$

Let's calculate the value of this sum. For this we shall change the variables for the summation index $i$ to $i-2$. Then (5.2) takes the following form:

$$\xi_1(2p) = \sum_{i=2}^{2p} (i-2) \left( \frac{\partial^{i-3}}{\partial x^{i-3}} a \right) \mathrm{C}(2p, i) \left( \frac{\partial^{2p-i}}{\partial x^{2p-i}} a \right)$$

This equation can be represented as follows: (5.3)

$$\xi_1(2p) = \left( \sum_{i=0}^{2p} (i-2) \left( \frac{\partial^{i-3}}{\partial x^{i-3}} a \right) \mathrm{C}(2p, i) \left( \frac{\partial^{2p-i}}{\partial x^{2p-i}} a \right) \right) - \left( \sum_{i=0}^{1} (i-2) \left( \frac{\partial^{i-3}}{\partial x^{i-3}} a \right) \mathrm{C}(2p, i) \left( \frac{\partial^{2p-i}}{\partial x^{2p-i}} a \right) \right)$$

Let's transform the first summand:

$$\sum_{i=0}^{2p} (i-2) \left( \frac{\partial^{i-3}}{\partial x^{i-3}} a \right) \mathrm{C}(2p, i) \left( \frac{\partial^{2p-i}}{\partial x^{2p-i}} a \right) = \left( \sum_{i=0}^{2p} i \left( \frac{\partial^{i-3}}{\partial x^{i-3}} a \right) \mathrm{C}(2p, i) \left( \frac{\partial^{2p-i}}{\partial x^{2p-i}} a \right) \right) - 2 \left( \sum_{i=0}^{2p} \left( \frac{\partial^{i-3}}{\partial x^{i-3}} a \right) \mathrm{C}(2p, i) \left( \frac{\partial^{2p-i}}{\partial x^{2p-i}} a \right) \right)$$

$$\left( \sum_{i=0}^{2p} i \left( \frac{\partial^{i}}{\partial x^{i}} \left( \int\int\int a\, dx\, dx\, dx \right) \right) \mathrm{C}(2p, i) \left( \frac{\partial^{2p-i}}{\partial x^{2p-i}} a \right) \right) - 2 \left( \sum_{i=0}^{2p} \left( \frac{\partial^{i}}{\partial x^{i}} \left( \int\int\int a\, dx\, dx\, dx \right) \right) \mathrm{C}(2p, i) \left( \frac{\partial^{2p-i}}{\partial x^{2p-i}} a \right) \right)$$

The second summand from the resulting equation is obtained with the use of Leibniz formula:

$$\sum_{i=0}^{2p} \left( \frac{\partial^{i}}{\partial x^{i}} \left( \int\int\int a\, dx\, dx\, dx \right) \right) \mathrm{C}(2p, i) \left( \frac{\partial^{2p-i}}{\partial x^{2p-i}} a \right) = \frac{\partial^{2p}}{\partial x^{2p}} \left( a \int\int\int a\, dx\, dx\, dx \right) \tag{5.4}$$

Let's establish the value of the following expression:

$$\sum_{i=0}^{2p} i \left( \frac{\partial^{i}}{\partial x^{i}} \left( \int\int\int a\, dx\, dx\, dx \right) \right) \mathrm{C}(2p, i) \left( \frac{\partial^{2p-i}}{\partial x^{2p-i}} a \right) \tag{5.5}$$

Since the expression (5.5) represents the impact of Leibniz operator on functions and $2p$ on $i$ index, that is:

$$\sum_{i=0}^{2p} i \left( \frac{\partial^{i}}{\partial x^{i}} \left( \int\int\int a\, dx\, dx\, dx \right) \right) \mathrm{C}(2p, i) \left( \frac{\partial^{2p-i}}{\partial x^{2p-i}} a \right) = L[i]$$

in accordance with (1.7) we can establish that:

$$\sum_{i=0}^{2p} i\left(\frac{\partial^i}{\partial x^i}(\int\int\int a\,dx\,dx\,dx)\right) C(2p,i)\left(\frac{\partial^{2p-i}}{\partial x^{2p-i}}a\right) = 2p\left(\frac{\partial^{2p-1}}{\partial x^{2p-1}}(a\int\int a\,dx\,dx)\right) \tag{5.10}$$

Thus, the final formula for the second coefficient takes the following form:

$$\begin{aligned}\xi_1(2p) = {} & 2p\left(\frac{\partial^{2p-1}}{\partial x^{2p-1}}(a\int\int a\,dx\,dx)\right) - 2\left(\frac{\partial^{2p}}{\partial x^{2p}}(a\int\int\int a\,dx\,dx\,dx)\right) \\ & + 2\int\int\int a\,dx\,dx\,dx\left(\frac{\partial^{2p}}{\partial x^{2p}}a\right) + 2\int\int a\,dx\,dx\,p\left(\frac{\partial^{2p-1}}{\partial x^{2p-1}}a\right)\end{aligned} \tag{5.11}$$

**Note**: Formally, it also gives the value of the partial sum:

$$\begin{aligned}\sum_{i_1=0}^{2p-2} i_1\left(\frac{\partial^{i_1-1}}{\partial x^{i_1-1}}a\right) C(2p,i_1+2)\left(\frac{\partial^{2p-2-i_1}}{\partial x^{2p-2-i_1}}a\right) = {} & 2p\left(\frac{\partial^{2p-1}}{\partial x^{2p-1}}(a\int\int a\,dx\,dx)\right) \\ & - 2\left(\frac{\partial^{2p}}{\partial x^{2p}}(a\int\int\int a\,dx\,dx\,dx)\right) + 2\int\int\int a\,dx\,dx\,dx\left(\frac{\partial^{2p}}{\partial x^{2p}}a\right) \\ & + 2\int\int a\,dx\,dx\,p\left(\frac{\partial^{2p-1}}{\partial x^{2p-1}}a\right)\end{aligned}$$

**5.2. Calculation of $\xi_2(2p)$.**

With the use of the recurrence formula (4.8) we have:

$$\xi_2(2p) = \sum_{i_2=0}^{2p-4} \xi_1(i_2+2)\,C(2p,i_2+4)\left(\frac{\partial^{2p-4-i_2}}{\partial x^{2p-4-i_2}}a\right) \tag{5.12}$$

Let's change of the summation index $i_2$ to $i-4$. This allows to transform this formula as follows:

$$\begin{aligned}\xi_2(2p) = {} & \sum_{i=4}^{2p} \xi_1(i-2)\,C(2p,i)\left(\frac{\partial^{2p-i}}{\partial x^{2p-i}}a\right) = \\ & \left(\sum_{i=0}^{2p} \xi_1(i-2)\,C(2p,i)\left(\frac{\partial^{2p-i}}{\partial x^{2p-i}}a\right)\right) - \left(\sum_{i=0}^{3} \xi_1(i-2)\,C(2p,i)\left(\frac{\partial^{2p-i}}{\partial x^{2p-i}}a\right)\right)\end{aligned} \tag{5.13}$$

With the replacement of

$$p = \frac{i-2}{2}$$

parameter, (5.11) takes the following form:

$$\xi_1(i-2) = \left(\frac{\partial^{i-3}}{\partial x^{i-3}}\left(a\int\int a\,dx\,dx\right)\right) i - 2\left(\frac{\partial^{i-3}}{\partial x^{i-3}}\left(a\int\int a\,dx\,dx\right)\right)$$
$$- 2\left(\frac{\partial^{i-2}}{\partial x^{i-2}}\left(a\int\int\int a\,dx\,dx\,dx\right)\right) + 2\int\int\int a\,dx\,dx\,dx\left(\frac{\partial^{i-2}}{\partial x^{i-2}}a\right)$$
$$+ \int\int a\,dx\,dx\left(\frac{\partial^{i-3}}{\partial x^{i-3}}a\right) i - 2\int\int a\,dx\,dx\left(\frac{\partial^{i-3}}{\partial x^{i-3}}a\right)$$

or

$$\xi_1(i-2) = \left(\frac{\partial^{i}}{\partial x^{i}}\left(\int\int\int a\int\int a\,dx\,dx\,dx\,dx\,dx\right)\right) i - 2\left(\frac{\partial^{i}}{\partial x^{i}}\left(\int\int\int a\int\int a\,dx\,dx\,dx\,dx\,dx\right)\right)$$
$$- 2\left(\frac{\partial^{i}}{\partial x^{i}}\left(\int\int a\int\int\int a\,dx\,dx\,dx\,dx\,dx\right)\right) + 2\int\int\int a\,dx\,dx\,dx\left(\frac{\partial^{i}}{\partial x^{i}}\left(\int\int a\,dx\,dx\right)\right)$$
$$+ \int\int a\,dx\,dx\left(\frac{\partial^{i}}{\partial x^{i}}\left(\int\int\int a\,dx\,dx\,dx\right)\right) i - 2\int\int a\,dx\,dx\left(\frac{\partial^{i}}{\partial x^{i}}\left(\int\int\int a\,dx\,dx\,dx\right)\right)$$

Thus, (5.13) becomes: (5.14)

$$\xi_2(2p) = \left(\sum_{i=0}^{2p}\left(\left(\frac{\partial^{i}}{\partial x^{i}}\left(\int\int\int a\int\int a\,dx\,dx\,dx\,dx\,dx\right)\right) i\right.\right.$$
$$- 2\left(\frac{\partial^{i}}{\partial x^{i}}\left(\int\int\int a\int\int a\,dx\,dx\,dx\,dx\,dx\right)\right) - 2\left(\frac{\partial^{i}}{\partial x^{i}}\left(\int\int a\int\int\int a\,dx\,dx\,dx\,dx\,dx\right)\right)$$
$$+ 2\int\int\int a\,dx\,dx\,dx\left(\frac{\partial^{i}}{\partial x^{i}}\left(\int\int a\,dx\,dx\right)\right) + \int\int a\,dx\,dx\left(\frac{\partial^{i}}{\partial x^{i}}\left(\int\int\int a\,dx\,dx\,dx\right)\right) i$$
$$\left.\left.- 2\int\int a\,dx\,dx\left(\frac{\partial^{i}}{\partial x^{i}}\left(\int\int\int a\,dx\,dx\,dx\right)\right)\right) \mathrm{C}(2p,i)\left(\frac{\partial^{2p-i}}{\partial x^{2p-i}}a\right)\right) - \left(\sum_{i=0}^{3}\left(\right.\right.$$
$$\left(\frac{\partial^{i}}{\partial x^{i}}\left(\int\int\int a\int\int a\,dx\,dx\,dx\,dx\,dx\right)\right) i - 2\left(\frac{\partial^{i}}{\partial x^{i}}\left(\int\int\int a\int\int a\,dx\,dx\,dx\,dx\,dx\right)\right)$$
$$- 2\left(\frac{\partial^{i}}{\partial x^{i}}\left(\int\int a\int\int\int a\,dx\,dx\,dx\,dx\,dx\right)\right) + 2\int\int\int a\,dx\,dx\,dx\left(\frac{\partial^{i}}{\partial x^{i}}\left(\int\int a\,dx\,dx\right)\right)$$
$$+ \int\int a\,dx\,dx\left(\frac{\partial^{i}}{\partial x^{i}}\left(\int\int\int a\,dx\,dx\,dx\right)\right) i - 2\int\int a\,dx\,dx\left(\frac{\partial^{i}}{\partial x^{i}}\left(\int\int\int a\,dx\,dx\,dx\right)\right)\Bigg)$$
$$\left.\left.\mathrm{C}(2p,i)\left(\frac{\partial^{2p-i}}{\partial x^{2p-i}}a\right)\right)\right)$$

The first summand on the right side of (5.14) is represented as follows:

$$\sum_{i=0}^{2p}\left(\left(\frac{\partial^i}{\partial x^i}\left(\int\int\int a\int\int a\,dx\,dx\,dx\,dx\,dx\right)\right)i-2\left(\frac{\partial^i}{\partial x^i}\left(\int\int\int a\int\int a\,dx\,dx\,dx\,dx\,dx\right)\right)\right.$$

$$-2\left(\frac{\partial^i}{\partial x^i}\left(\int\int a\int\int\int a\,dx\,dx\,dx\,dx\,dx\right)\right)+2\int\int\int a\,dx\,dx\,dx\left(\frac{\partial^i}{\partial x^i}(\int\int a\,dx\,dx)\right)$$

$$\left.+\int\int a\,dx\,dx\left(\frac{\partial^i}{\partial x^i}(\int\int\int a\,dx\,dx\,dx)\right)i-2\int\int a\,dx\,dx\left(\frac{\partial^i}{\partial x^i}(\int\int\int a\,dx\,dx\,dx)\right)\right)$$

$$=$$

$$C(2p,i)\left(\frac{\partial^{2p-i}}{\partial x^{2p-i}}a\right)$$

$$\sum_{i=0}^{2p}\left(\frac{\partial^i}{\partial x^i}\left(\int\int\int a\int\int a\,dx\,dx\,dx\,dx\,dx\right)\right)i\,C(2p,i)\left(\frac{\partial^{2p-i}}{\partial x^{2p-i}}a\right)$$

$$-2\left(\sum_{i=0}^{2p}\left(\frac{\partial^i}{\partial x^i}\left(\int\int\int a\int\int a\,dx\,dx\,dx\,dx\,dx\right)\right)C(2p,i)\left(\frac{\partial^{2p-i}}{\partial x^{2p-i}}a\right)\right)\ \mathbf{-2}$$

$$\sum_{i=0}^{2p}\left(\frac{\partial^i}{\partial x^i}\left(\int\int a\int\int\int a\,dx\,dx\,dx\,dx\,dx\right)\right)C(2p,i)\left(\frac{\partial^{2p-i}}{\partial x^{2p-i}}a\right)\ \mathbf{+2}$$

$$\int\int\int a\,dx\,dx\,dx\left(\sum_{i=0}^{2p}\left(\frac{\partial^i}{\partial x^i}(\int\int a\,dx\,dx)\right)C(2p,i)\left(\frac{\partial^{2p-i}}{\partial x^{2p-i}}a\right)\right)\ \mathbf{+}$$

$$\int\int a\,dx\,dx\left(\sum_{i=0}^{2p}\left(\frac{\partial^i}{\partial x^i}(\int\int\int a\,dx\,dx\,dx)\right)i\,C(2p,i)\left(\frac{\partial^{2p-i}}{\partial x^{2p-i}}a\right)\right)\ \mathbf{-2}$$

$$\int\int a\,dx\,dx\left(\sum_{i=0}^{2p}\left(\frac{\partial^i}{\partial x^i}(\int\int\int a\,dx\,dx\,dx)\right)C(2p,i)\left(\frac{\partial^{2p-i}}{\partial x^{2p-i}}a\right)\right)$$

Using (5.10) and Leibniz formula we can easily establish the validity of the following equations:

$$\sum_{i=0}^{2p}\left(\frac{\partial^i}{\partial x^i}\left(\int\int\int a\int\int a\,dx\,dx\,dx\,dx\,dx\right)\right)i\,C(2p,i)\left(\frac{\partial^{2p-i}}{\partial x^{2p-i}}a\right)=2p\left(\frac{\partial^{2p-1}}{\partial x^{2p-1}}\left(a\int\int a\int\int a\,dx\,dx\,dx\,dx\right)\right)$$

$$\sum_{i=0}^{2p}\left(\frac{\partial^i}{\partial x^i}\left(\int\int\int a\int\int a\,dx\,dx\,dx\,dx\,dx\right)\right)C(2p,i)\left(\frac{\partial^{2p-i}}{\partial x^{2p-i}}a\right)=\frac{\partial^{2p}}{\partial x^{2p}}\left(a\int\int\int a\int\int a\,dx\,dx\,dx\,dx\,dx\right)$$

$$\sum_{i=0}^{2p}\left(\frac{\partial^i}{\partial x^i}\left(\int\int a\int\int\int a\,dx\,dx\,dx\,dx\,dx\right)\right)C(2p,i)\left(\frac{\partial^{2p-i}}{\partial x^{2p-i}}a\right)=\frac{\partial^{2p}}{\partial x^{2p}}\left(a\int\int a\int\int\int a\,dx\,dx\,dx\,dx\,dx\right)$$

$$\sum_{i=0}^{2p}\left(\frac{\partial^i}{\partial x^i}(\int\int a\,dx\,dx)\right)C(2p,i)\left(\frac{\partial^{2p-i}}{\partial x^{2p-i}}a\right)=$$

$$\sum_{i=0}^{2p}\left(\frac{\partial^i}{\partial x^i}(\int\int\int a\,dx\,dx\,dx)\right)i\,C(2p,i)\left(\frac{\partial^{2p-i}}{\partial x^{2p-i}}a\right)=2p\left(\frac{\partial^{2p-1}}{\partial x^{2p-1}}(a\int\int a\,dx\,dx)\right)$$

$$\sum_{i=0}^{2p}\left(\frac{\partial^i}{\partial x^i}(\int\int\int a\,dx\,dx\,dx)\right)C(2p,i)\left(\frac{\partial^{2p-i}}{\partial x^{2p-i}}a\right)=$$

Substituting the obtained values to (5. 14) we eventually establish the required value of $\xi_2(2p)$. (5.15)

$$\xi_2(2p) = 2\left(\frac{\partial^{2p}}{\partial x^{2p}}\left(p\int a\iint a\iint a\,dx\,dx\,dx\,dx\,dx - a\iiint a\iint a\,dx\,dx\,dx\,dx\,dx\right.\right.$$

$$\left.\left. - a\iint a\iiint a\,dx\,dx\,dx\,dx\,dx\right)\right) + 2p\left((\iint a\,dx\,dx)^2 + \iint a\iint a\,dx\,dx\,dx\,dx\right.$$

$$\left. + 2\int a\iiint a\,dx\,dx\,dx\,dx - 2\iiint a\,dx\,dx\,dx\int a\,dx\right)\left(\frac{\partial^{2p}}{\partial x^{2p}}(\int a\,dx)\right)$$

$$+ 2\iint a\,dx\,dx\left(\frac{\partial^{2p}}{\partial x^{2p}}\left(p\int a\iint a\,dx\,dx\,dx - a\iiint a\,dx\,dx\,dx\right)\right)$$

$$+ 2\left(\iint a\iiint a\,dx\,dx\,dx\,dx\,dx + \iiint a\iint a\,dx\,dx\,dx\,dx\,dx\right)\left(\frac{\partial^{2p}}{\partial x^{2p}}a\right)$$

$$+ 2\iiint a\,dx\,dx\,dx\left(\frac{\partial^{2p}}{\partial x^{2p}}(a\iint a\,dx\,dx)\right)$$

## 5. Finding the $\xi_3(2p)$ summand

Taking in account the importance of $\xi_2(2p)$, let's introduce the following notations:

$$k_1 = \int a\iint a\iint a\,dx\,dx\,dx\,dx\,dx\ ,\ k_2 = -a\iiint a\iint a\,dx\,dx\,dx\,dx\,dx - a\iint a\iiint a\,dx\,dx\,dx\,dx\,dx$$

$$k_3 = (\iint a\,dx\,dx)^2 + \iint a\iint a\,dx\,dx\,dx\,dx + 2\int a\iiint a\,dx\,dx\,dx\,dx$$

$$- 2\iiint a\,dx\,dx\,dx\int a\,dx$$

$$k_4 = \int a\,dx\ ,\ ,\ ,\ ,\ k_8 = \iint a\iiint a\,dx\,dx\,dx\,dx\,dx + \iiint a\iint a\,dx\,dx\,dx\,dx\,dx$$

,

In this case $\xi_2(2p)$ takes the following form: (5.16)

$$\xi_2(2p) = 2\left(\frac{\partial^{2p}}{\partial x^{2p}}(p\,k_1 + k_2)\right) + 2p\,k_3\left(\frac{\partial^{2p}}{\partial x^{2p}}k_4\right) + 2\,k_5\left(\frac{\partial^{2p}}{\partial x^{2p}}(p\,k_6 + k_7)\right) + 2\,k_8\left(\frac{\partial^{2p}}{\partial x^{2p}}a\right)$$

$$+ 2\,k_9\left(\frac{\partial^{2p}}{\partial x^{2p}}k_{10}\right)$$

As soon as, in accordance with the recurrence formula (4.8)

$$\xi_3(2p) = \sum_{i_3=0}^{2p-6}\xi_2(i_3+4)\,\mathrm{C}(2p, i_3+6)\left(\frac{\partial^{2p-6-i_3}}{\partial x^{2p-6-i_3}}a\right) \tag{5.17}$$

by the replacement of summation index $i_3$ to $i-6$ it can be transformed as follows : (5.18)

$$\xi_3(2p) = \left(\sum_{i=0}^{2p} \xi_2(i-2)\, C(2p,i) \left(\frac{\partial^{2p-i}}{\partial x^{2p-i}} a\right)\right) - \left(\sum_{i=0}^{5} \xi_2(i-2)\, C(2p,i) \left(\frac{\partial^{2p-i}}{\partial x^{2p-i}} a\right)\right)$$

In this case, the equality (5.16), by the substitution of $p = \frac{i-2}{2}$ , can be transformed as follows:

$$\xi_2(i-2) = 2\left(\frac{\partial^{i-2}}{\partial x^{i-2}}\left(\frac{1 k_1 i}{2} - k_1 + k_2\right)\right) + k_3\left(\frac{\partial^{i-2}}{\partial x^{i-2}} k_4\right) i - 2 k_3\left(\frac{\partial^{i-2}}{\partial x^{i-2}} k_4\right)$$
$$+ 2 k_5\left(\frac{\partial^{i-2}}{\partial x^{i-2}}\left(\frac{1 k_6 i}{2} - k_6 + k_7\right)\right) + 2 k_8\left(\frac{\partial^{i-2}}{\partial x^{i-2}} a\right) + 2 k_9\left(\frac{\partial^{i-2}}{\partial x^{i-2}} k_{10}\right)$$

or the following form for calculation:

$$\xi_2(i-2) = i\left(\frac{\partial^i}{\partial x^i}\left(\int\!\!\int k_1\, dx\, dx\right)\right) + 2\left(\frac{\partial^i}{\partial x^i}\left(\int\!\!\int -k_1 + k_2\, dx\, dx\right)\right) + k_3\left(\frac{\partial^i}{\partial x^i}\left(\int\!\!\int k_4\, dx\, dx\right)\right) i$$
$$- 2 k_3\left(\frac{\partial^i}{\partial x^i}\left(\int\!\!\int k_4\, dx\, dx\right)\right) + k_5\, i\left(\frac{\partial^i}{\partial x^i}\left(\int\!\!\int k_6\, dx\, dx\right)\right) + 2 k_5\left(\frac{\partial^i}{\partial x^i}\left(\int\!\!\int -k_6 + k_7\, dx\, dx\right)\right)$$
$$+ 2 k_8\left(\frac{\partial^i}{\partial x^i}\left(\int\!\!\int a\, dx\, dx\right)\right) + 2 k_9\left(\frac{\partial^i}{\partial x^i}\left(\int\!\!\int k_{10}\, dx\, dx\right)\right)$$

Substituting this equality to (5.18) we get:

$$\xi_3(2p) = \left(\sum_{i=0}^{2p}\left(i\left(\frac{\partial^i}{\partial x^i}\left(\int\!\!\int k_1\, dx\, dx\right)\right) + 2\left(\frac{\partial^i}{\partial x^i}\left(\int\!\!\int -k_1 + k_2\, dx\, dx\right)\right)\right.\right.$$
$$+ k_3\left(\frac{\partial^i}{\partial x^i}\left(\int\!\!\int k_4\, dx\, dx\right)\right) i - 2 k_3\left(\frac{\partial^i}{\partial x^i}\left(\int\!\!\int k_4\, dx\, dx\right)\right) + k_5\, i\left(\frac{\partial^i}{\partial x^i}\left(\int\!\!\int k_6\, dx\, dx\right)\right)$$
$$+ 2 k_5\left(\frac{\partial^i}{\partial x^i}\left(\int\!\!\int -k_6 + k_7\, dx\, dx\right)\right) + 2 k_8\left(\frac{\partial^i}{\partial x^i}\left(\int\!\!\int a\, dx\, dx\right)\right)$$
$$\left.\left.+ 2 k_9\left(\frac{\partial^i}{\partial x^i}\left(\int\!\!\int k_{10}\, dx\, dx\right)\right)\right) C(2p,i)\left(\frac{\partial^{2p-i}}{\partial x^{2p-i}} a\right)\right) - \left(\sum_{i=0}^{5}\left(i\left(\frac{\partial^i}{\partial x^i}\left(\int\!\!\int k_1\, dx\, dx\right)\right)\right.\right.$$
$$+ 2\left(\frac{\partial^i}{\partial x^i}\left(\int\!\!\int -k_1 + k_2\, dx\, dx\right)\right) + k_3\left(\frac{\partial^i}{\partial x^i}\left(\int\!\!\int k_4\, dx\, dx\right)\right) i - 2 k_3\left(\frac{\partial^i}{\partial x^i}\left(\int\!\!\int k_4\, dx\, dx\right)\right)$$
$$+ k_5\, i\left(\frac{\partial^i}{\partial x^i}\left(\int\!\!\int k_6\, dx\, dx\right)\right) + 2 k_5\left(\frac{\partial^i}{\partial x^i}\left(\int\!\!\int -k_6 + k_7\, dx\, dx\right)\right) + 2 k_8\left(\frac{\partial^i}{\partial x^i}\left(\int\!\!\int a\, dx\, dx\right)\right)$$
$$\left.\left.+ 2 k_9\left(\frac{\partial^i}{\partial x^i}\left(\int\!\!\int k_{10}\, dx\, dx\right)\right)\right) C(2p,i)\left(\frac{\partial^{2p-i}}{\partial x^{2p-i}} a\right)\right)$$

As soon as

$$\sum_{i=0}^{2p} i\left(\frac{\partial^i}{\partial x^i}\left(\int\!\!\int k_1\, dx\, dx\right)\right) C(2p,i)\left(\frac{\partial^{2p-i}}{\partial x^{2p-i}} a\right) = 2p\left(\frac{\partial^{2p-1}}{\partial x^{2p-1}}\left(a \int k_1\, dx\right)\right)$$

$$2\left(\sum_{i=0}^{2p}\left(\frac{\partial^i}{\partial x^i}\left(\int -\int k_1\,dx+\int k_2\,dx\,dx\right)\right)C(2p,i)\left(\frac{\partial^{2p-i}}{\partial x^{2p-i}}a\right)\right)= \qquad 2\left(\frac{\partial^{2p}}{\partial x^{2p}}\left(a\int -\int k_1\,dx+\int k_2\,dx\,dx\right)\right)$$

$$k_3\left(\sum_{i=0}^{2p}\left(\frac{\partial^i}{\partial x^i}\left(\int\int k_4\,dx\,dx\right)\right)i\,C(2p,i)\left(\frac{\partial^{2p-i}}{\partial x^{2p-i}}a\right)\right)=2\,k_3\,p\left(\frac{\partial^{2p-1}}{\partial x^{2p-1}}\left(a\int k_4\,dx\right)\right)$$

$$-2\,k_3\left(\sum_{i=0}^{2p}\left(\frac{\partial^i}{\partial x^i}\left(\int\int k_4\,dx\,dx\right)\right)C(2p,i)\left(\frac{\partial^{2p-i}}{\partial x^{2p-i}}a\right)\right)=-2\,k_3\left(\frac{\partial^{2p}}{\partial x^{2p}}\left(a\int\int k_4\,dx\,dx\right)\right)$$

$$k_5\left(\sum_{i=0}^{2p}i\left(\frac{\partial^i}{\partial x^i}\left(\int\int k_6\,dx\,dx\right)\right)C(2p,i)\left(\frac{\partial^{2p-i}}{\partial x^{2p-i}}a\right)\right)=2\,k_5\,p\left(\frac{\partial^{2p-1}}{\partial x^{2p-1}}\left(a\int k_6\,dx\right)\right)$$

$$2\,k_5\left(\sum_{i=0}^{2p}\left(\frac{\partial^i}{\partial x^i}\left(\int -\int k_6\,dx+\int k_7\,dx\,dx\right)\right)C(2p,i)\left(\frac{\partial^{2p-i}}{\partial x^{2p-i}}a\right)\right)= \qquad 2\,k_5\left(\frac{\partial^{2p}}{\partial x^{2p}}\left(a\int -\int k_6\,dx+\int k_7\,dx\,dx\right)\right)$$

$$2\,k_8\left(\sum_{i=0}^{2p}\left(\frac{\partial^i}{\partial x^i}\left(\int\int a\,dx\,dx\right)\right)C(2p,i)\left(\frac{\partial^{2p-i}}{\partial x^{2p-i}}a\right)\right)=2\,k_8\left(\frac{\partial^{2p}}{\partial x^{2p}}\left(a\int\int a\,dx\,dx\right)\right)$$

$$2\,k_9\left(\sum_{i=0}^{2p}\left(\frac{\partial^i}{\partial x^i}\left(\int\int k_{10}\,dx\,dx\right)\right)C(2p,i)\left(\frac{\partial^{2p-i}}{\partial x^{2p-i}}a\right)\right)=2\,k_9\left(\frac{\partial^{2p}}{\partial x^{2p}}\left(a\int\int k_{10}\,dx\,dx\right)\right)$$

upon transformations we obtain: (5.19)

$$\begin{aligned}\xi_3(2p) = &-2\left(\frac{\partial^{2p}}{\partial x^{2p}}(k_1-p\,k_2)\right)+2\,k_3\left(\frac{\partial^{2p}}{\partial x^{2p}}\left(\int a\int\int a\,dx\,dx\,dx\,p-a\int\int\int a\,dx\,dx\,dx\right)\right)\\ &+2\,k_4\left(\frac{\partial^{2p}}{\partial x^{2p}}\left(a\int\int a\,dx\,dx\right)\right)+2\,p\,k_5\left(\frac{\partial^{2p}}{\partial x^{2p}}\left(\int a\,dx\right)\right)+2\,k_6\left(\frac{\partial^{2p}}{\partial x^{2p}}a\right)\\ &+2\int\int a\,dx\,dx\left(\frac{\partial^{2p}}{\partial x^{2p}}(k_7\,p+k_8)\right)\\ &+2\int\int\int a\,dx\,dx\,dx\left(\frac{\partial^{2p}}{\partial x^{2p}}\left(a\int\int a\int\int a\,dx\,dx\,dx\,dx\right)\right)\end{aligned}$$

$$\begin{aligned}k_1 = &\ a\int\int a\int\int\int a\int\int a\,dx\,dx\,dx\,dx\,dx\,dx\,dx+a\int\int a\int\int a\int\int\int a\,dx\,dx\,dx\,dx\,dx\,dx\,dx\\ &+a\int\int\int a\int\int a\int\int a\,dx\,dx\,dx\,dx\,dx\,dx\,dx\end{aligned}$$

$$k_2=\int a\int\int a\int\int a\int\int a\,dx\,dx\,dx\,dx\,dx\,dx\,dx$$

$$k_3 = \left(\int\int a\,dx\,dx\right)^2 + \int\int a \int\int a\,dx\,dx\,dx\,dx + 2\int a \int\int\int a\,dx\,dx\,dx\,dx$$
$$- 2\int\int\int a\,dx\,dx\,dx \int a\,dx$$

$$k_4 = \int\int a \int\int\int a\,dx\,dx\,dx\,dx\,dx + \int\int\int a \int\int a\,dx\,dx\,dx\,dx\,dx$$

$$k_5 = 2\int a \int\int a \int\int\int a\,dx\,dx\,dx\,dx\,dx\,dx + 2\int a \int\int\int a \int\int a\,dx\,dx\,dx\,dx\,dx\,dx$$
$$+ \left(\int\int a\,dx\,dx\right)^3 + 2\int\int a \int\int a\,dx\,dx\,dx\,dx \int\int a\,dx\,dx$$
$$+ 4\int a \int\int\int a\,dx\,dx\,dx\,dx \int\int a\,dx\,dx - 2\int a\,dx \int\int\int a \int\int a\,dx\,dx\,dx\,dx\,dx$$
$$- 2\int a \int\int a\,dx\,dx\,dx \int\int\int a\,dx\,dx\,dx - 2\int a\,dx \int\int a \int\int\int a\,dx\,dx\,dx\,dx\,dx$$
$$- 2\int\int\int a\,dx\,dx\,dx \int a\,dx \int\int a\,dx\,dx + \int\int a \int\int a \int\int a\,dx\,dx\,dx\,dx\,dx\,dx$$

$$k_6 = \left(\int\int a\,dx\,dx\right)^2 \int\int\int a\,dx\,dx\,dx - 2\int a\,dx \left(\int\int\int a\,dx\,dx\,dx\right)^2$$
$$+ 2\int a \int\int\int a\,dx\,dx\,dx\,dx \int\int\int a\,dx\,dx\,dx + \int\int\int a \int\int a \int\int a\,dx\,dx\,dx\,dx\,dx\,dx\,dx$$
$$+ \int\int a \int\int\int a \int\int a\,dx\,dx\,dx\,dx\,dx\,dx\,dx + \int\int a \int\int a \int\int\int a\,dx\,dx\,dx\,dx\,dx\,dx\,dx$$

$$k_7 = \int a \int\int a \int\int a\,dx\,dx\,dx\,dx\,dx$$

$$k_8 = -a \int\int a \int\int\int a\,dx\,dx\,dx\,dx\,dx - a \int\int\int a \int\int a\,dx\,dx\,dx\,dx\,dx$$

Proceeding in a similar way, we can obtain all following summands:

As we can see, with with the increase of $k$ number $\xi_k(2p)$ is becoming more cumbersome, which evidences a substantial information content.

## 5.4. General formula for the definition of coefficients.

Let's prove the **Theorem 5.4. The** $\xi_k(p)$ **,** $k = 1, 2, 3$ **... functions are determined by the following formula**:

$$\xi_k(p) = \left(\frac{\partial^p}{\partial x^p}\left(\eth \int G_k(\grave{a})\,dx - 2\,P_k(a)\right)\right) - \left(\sum_{s=0}^{k-1} \xi_{k-s-1}(-1)\left(\frac{\partial^p}{\partial x^p}\left(\eth \int G_s(\grave{a})\,dx - 2\,P_s(a)\right)\right)\right) - \left(\sum_{s=0}^{k-1} \xi_{k-s-1}(-2)\left(\frac{\partial^p}{\partial x^p} G_s(\grave{a})\right)\right)$$

(5.20)

**where ;**

$G_i(f)$ **is an integral operator, which operates on the specified** $f = f(x)$ **function in accordance with the following rules**:

$$G_0(f) = f \tag{5.21}$$

**=** (5.22)

**=** (5.23)

$$G_3(f) = G_1(G_1(G_1(G_0(f)))) = a\int\int a\int\int a\int\int f\,dx\,dx\,dx\,dx\,dx\,dx$$

**and so on.**

**The** $P_k(a)$ **operator with** $k = 0$ **is determined by the following formula:**

$$P_0(a) = 0 \tag{5.24}$$

**and operates on the specified *a* function according to the following rule:**

$$P_k(a) = \sum_{i=0}^{k-1} G_{k-i}\left(\int G_i(a)\,dx\right), k = 1, 2, 3 \tag{5.25}$$

**Proof**: we shall proceed in two phases.

Phase 1. Let's demonstrate that (5.20) gives the $\xi_k(p)$ values**,** which are the same as those defined with the use of integral method. In particular, we can obtain (5.11), (5.15), (5.19) from (5.20).

Indeed, assuming that $k = 0$ in (5.20) we get:

$$\xi_0(p) = \frac{\partial^p}{\partial x^p}\left(\eth \int G_0(\grave{a})\,dx - 2\,P_0(a)\right) \tag{5.26}$$

As soon as, in accordance with (5.21) $G_0(\grave{a}) = a$ , and in accordance with (5.24) $P_0(a) = 0$ , (5.26) takes the following form:

(5.27)

Because the operation of integration is written down through the operation of differentiation in accordance with the following formula:

(5.27) takes the following form:

(5.28)

This is fully identical with the known representation of (4.2):

In (5.20) let's assume $k = 1$ **:**

$$\xi_1(p) = \left(\frac{\partial^p}{\partial x^p}\left(\eth\int G_1(\grave{a})\,dx - 2\,P_1(a)\right)\right) - \xi_0(-1)\left(\frac{\partial^p}{\partial x^p}\left(\eth\int G_0(\grave{a})\,dx - 2\,P_0(a)\right)\right) - \xi_0(-2)\left(\frac{\partial^p}{\partial x^p}G_0(\grave{a})\right) \qquad (5.29)$$

As soon as, in accordance with (5.22)

$$G_0(\grave{a}) = a$$

and in accordance with (5.24), (5.25):

$$P_0(a) = 0 \qquad =$$

successively assuming in (5.28), we obtain the following equalities:

= =

substituting these to (5.29) we obtain:

$$\xi_1(p) = \left(\frac{\partial^p}{\partial x^p}\left(\eth\int\left(a\int\int \grave{a}\,dx\,dx\,dx - 2\,a\int\int\int \grave{a}\,dx\,dx\,dx\right)\right)\right) + \int\int \grave{a}\,dx\,dx\left(\frac{\partial^p}{\partial x^p}(\eth\int \grave{a}\,dx)\right) + 2\int\int\int \grave{a}\,dx\,dx\,dx\left(\frac{\partial^p}{\partial x^p}\grave{a}\right)$$

With consideration of the following equalities:

$$\int a\int\int \grave{a}\,dx\,dx\,dx = \frac{\partial^{-1}}{\partial x^{-1}}(a\int\int \grave{a}\,dx\,dx)\ ;$$

we obtain the following formula:

$$\xi_1(p) = p\left(\frac{\partial^{p-1}}{\partial x^{p-1}}(a\int\int a\,dx\,dx)\right) - 2\left(\frac{\partial^p}{\partial x^p}(a\int\int\int a\,dx\,dx\,dx)\right)$$
$$+ 2\int\int\int a\,dx\,dx\,dx\left(\frac{\partial^p}{\partial x^p}a\right) + \int\int a\,dx\,dx\,p\left(\frac{\partial^{p-1}}{\partial x^{p-1}}a\right) \qquad (5.30)$$

which is fully identical with (5.11), if we substitute $2p$ for $p$ parameter.

Again, assuming that $k = 2$ in (5.20):

$$\xi_2(p) = \left(\frac{\partial^p}{\partial x^p}\left(\eth\int G_2(\grave{a})\,dx - 2\,P_2(a)\right)\right)$$
$$- \left(\sum_{s=0}^{1}\xi_{k-s-1}(-1)\left(\frac{\partial^p}{\partial x^p}\left(\eth\int G_s(\grave{a})\,dx - 2\,P_s(a)\right)\right)\right) - \left(\sum_{s=0}^{1}\xi_{k-s-1}(-2)\left(\frac{\partial^p}{\partial x^p}G_s(\grave{a})\right)\right) \qquad (5.31)$$

As soon as, in accordance with (5.22)

$$G_0(\grave{a}) = a \; : \quad : \; G_2(\grave{a}) = a\int\int a\int\int \grave{a}\,dx\,dx\,dx\,dx$$

and, in accordance with (5.24), (5.25):

$$P_0(a) = 0 \; : \; = \; :$$

$$P_2(a) = G_1\left(\int G_1(\grave{a})\,dx\right) + G_2\left(\int G_0(\grave{a})\,dx\right) = a\int\int\int a\int\int \grave{a}\,dx\,dx\,dx\,dx\,dx + a\int\int a\int\int\int \grave{a}\,dx\,dx\,dx\,dx\,dx$$

subseqnetly assuming in (5.30) we obtain the following equalities:

$$=$$ $$=$$

$$\xi_1(-1) = -\int\int a\int\int a\,dx\,dx\,dx\,dx - 2\int a\int\int\int a\,dx\,dx\,dx\,dx + 2\int\int\int a\,dx\,dx\,dx\int a\,dx$$
$$- \int\int a\,dx\,dx\int\int a\,dx\,dx$$

$$\xi_1(-2) = -2\int\int\int a\int\int a\,dx\,dx\,dx\,dx\,dx - 2\int\int a\int\int\int a\,dx\,dx\,dx\,dx\,dx$$
$$+ 2\int\int\int a\,dx\,dx\,dx\int\int a\,dx\,dx - 2\int\int a\,dx\,dx\int\int\int a\,dx\,dx\,dx$$

which we substitute to (5.31) and obtain:

$$\xi_2(p) = \left(\frac{\partial^p}{\partial x^p}\left(p\int a\int\int a\int\int a\,dx\,dx\,dx\,dx\,dx - 2\,a\int\int\int a\int\int a\,dx\,dx\,dx\,dx\,dx\right.\right.$$

$$-2a\int\int a\int\int\int a\,dx\,dx\,dx\,dx\,dx\Bigg)\Bigg) - \Bigg(-\int\int a\int\int a\,dx\,dx\,dx\,dx$$

$$-2\int a\int\int\int a\,dx\,dx\,dx\,dx - \left(\int\int a\,dx\,dx\right)^2 + 2\int\int\int a\,dx\,dx\,dx\int a\,dx\Bigg)$$

$$\left(\frac{\partial^p}{\partial x^p}\left(p\int a\,dx\right)\right) + \int\int a\,dx\,dx\left(\frac{\partial^p}{\partial x^p}\left(p\int a\int\int a\,dx\,dx\,dx - 2a\int\int\int a\,dx\,dx\,dx\right)\right)$$

$$-\left(-2\int\int\int a\int\int a\,dx\,dx\,dx\,dx\,dx - 2\int\int a\int\int\int a\,dx\,dx\,dx\,dx\,dx\right)\left(\frac{\partial^p}{\partial x^p}a\right)$$

$$+2\int\int\int a\,dx\,dx\,dx\left(\frac{\partial^p}{\partial x^p}\left(a\int\int a\,dx\,dx\right)\right)$$

which is fully identical with (5.15), if we replace $p$ parameter by $2p$.

Quite similarly, it can be proved that (5.19) and other formulas follow from (5.20). For this purpose we can use:

**THE PROGRAM FOR THE CALCULATION OF SUMMANDS according to the following formula:**

$$\xi_k(p) = \left(\frac{\partial^p}{\partial x^p}\left(\delta\int G_k(\grave{a})\,dx - 2P_k(a)\right)\right) - \left(\sum_{s=0}^{k-1}\xi_{k-s-1}(-1)\left(\frac{\partial^p}{\partial x^p}\left(\delta\int G_s(\grave{a})\,dx - 2P_s(a)\right)\right)\right) - \left(\sum_{s=0}^{k-1}\xi_{k-s-1}(-2)\left(\frac{\partial^p}{\partial x^p}G_s(\grave{a})\right)\right)$$

```
> restart:alias(y=y(x),v=v(x),a=a(x)):k:=1:# here the number of the
coefficient is specified

> dn:= proc (n::integer, w,c) local k, Ds; global resd, x: option remember;
if n = 0 then resd := w elif 0 < n then resd := diff(w,`$`(x,n)) else
Ds  := w:for k to -n do Ds:=int(Ds,x) end do:Ds:=Ds+sum(c[i]*x^(-n-i)/(-n-
i)!,i=1..-n): resd := Ds end if: RETURN(resd) end proc:

> for k from 0 to k do

> R:=z->a*int(int(z,x),x):h[0]:=a:for l from 0 to k do h[l+1]:=R(h[l]) od:
l:='l':
```

> **for** $l$ **from** 0 **to** $k+1$ **do** $f_l := \int h_l\,dx$ **end do**; $l := 'l'$

**for** $m$ **from** 0 **to** $k$ **do for** $j$ **from** 0 **to** $k-m$ **do** $q_0 := f_m$; $q_{j+1} := \mathrm{R}(q_j)$ **end do end do**;

$G_{k+1} := \mathrm{add}(q_i, i = 1..k+1)$;

> $G_0 := 0$;

$G_{-1} := 0$;

$$\xi_0(-2) := -2\int\int\int a\,dx\,dx\,dx$$

> $\xi_k(-1) := \mathrm{expand}\Big(-\int f_k\,dx - 2\int G_k\,dx + \mathrm{add}\Big(\xi_{k-s-1}(-1)\int f_s\,dx, s = 0..k-1\Big) + 2\,\mathrm{add}\Big(\xi_{k-s-1}(-1)\int G_s\,dx, s = 0..k-1\Big) - \mathrm{add}\Big(\xi_{k-s-1}(-2)\int h_s\,dx, s = 0..k-1\Big)\Big)$

> $\xi_k(-2) := \mathrm{expand}\Big(-2\int\int f_k\,dx\,dx - 2\int\int G_k\,dx\,dx + 2\,\mathrm{add}\Big(\xi_{k-s-1}(-1)\int\int f_s\,dx\,dx, s = 0..k-1\Big) + 2\,\mathrm{add}\Big(\xi_{k-s-1}(-1)\int\int G_s\,dx\,dx, s = 0..k-1\Big) - \mathrm{add}\Big(\xi_{k-s-1}(-2)\int\int h_s\,dx\,dx, s = 0..k-1\Big)\Big);$
$s := \text{'s'};$

$s := \text{'s'}$

> $\xi_k(p) = \left(\frac{\partial^p}{\partial x^p}(p f_k - 2\,G_k)\right) - \mathrm{add}\left(\xi_{k-s-1}(-1)\left(\frac{\partial^p}{\partial x^p}(p f_s - 2\,G_s)\right), s = 0..k-1\right) - \mathrm{add}\left(\xi_{k-s-1}(-2)\left(\frac{\partial^p}{\partial x^p}h_s\right), s = 0..k-1\right)$

```
> od:
> k := k-1:# prints the value of the coefficient
```

> $\xi_k(p) = \left(\frac{\partial^p}{\partial x^p}(p f_k - 2\,G_k)\right) - \mathrm{add}\left(\xi_{k-s-1}(-1)\left(\frac{\partial^p}{\partial x^p}(p f_s - 2\,G_s)\right), s = 0..k-1\right) - \mathrm{add}\left(\xi_{k-s-1}(-2)\left(\frac{\partial^p}{\partial x^p}h_s\right), s = 0..k-1\right)$

**Phase 2. Finding the recurrence formula.**
Analyzing the obtained formulas, we can come to the conclusion that formulas found for $\xi_k(2\,p)$ with the use of integral method can be generally presented as follows:

$$\xi_k(2\,p) = \sum_{s=0}^{k} r_{k,s}(x)\left(\frac{d^{2\,p}}{dx^{2\,p}}f_{s,k}(x)\right)$$

where $r_{k,s}(x)$, $f_{s,k}(x)$ are the required functions.

As follows from the above, the task should be divided in two parts. The first one is to identify the $f_{s,k}(x)$ functions, and the second one is to identify the $r_{k,s}(x)$ functions.

The first task, the task of finding the $f_{s,k}(x)$ functions, is solved by analyzing the functions under the $\frac{d^{2\,p}}{dx^{2\,p}}[\ ]$ sign:

Let's present the functional coefficients $f_{s,k}(x)$ from the expression $\frac{d^{2\,p}}{dx^{2\,p}}f_{s,k}(x)$ for the first three consecutive values of $k$.

**for** $k = 0$ $\quad p\int a\,dx$

**for** $k = 1$ $\quad p\int a\,dx \qquad p\int a\int\int a\,dx\,dx\,dx - 2\,a\int\int\int a\,dx\,dx\,dx$

**for** $k = 2$ $\quad p\int a\,dx \quad p\int a\iint a\,dx\,dx\,dx - 2a\iiint a\,dx\,dx\,dx$

$$p\int a\iint a\iint a\,dx\,dx\,dx\,dx\,dx - 2a\iiint a\iint a\,dx\,dx\,dx\,dx\,dx \quad - 2a\iint a\iiint a\,dx\,dx\,dx\,dx\,dx$$

**for** $k = 3$ $\quad p\int a\,dx \quad p\int a\iint a\,dx\,dx\,dx - 2a\iiint a\,dx\,dx\,dx$

$$p\int a\iint a\iint a\,dx\,dx\,dx\,dx\,dx - 2a\iiint a\iint a\,dx\,dx\,dx\,dx\,dx \quad - 2a\iint a\iiint a\,dx\,dx\,dx\,dx\,dx$$

**As we can see, for each new value of $k$ all old values are repeated and new values are introduced in accordance with the rule of attachment of the integral operator $G_i(f)$, where $f$ is a certain function:**

$$G_0(f) = f \qquad (5.32)$$

(5.33)

$$G_2(f) = a\iint a\iint f\,dx\,dx\,dx\,dx = G_1(G_1(f))$$

$$G_3(f) = a\iint a\iint a\iint f\,dx\,dx\,dx\,dx\,dx\,dx =$$

**and so on.**

Moreover, let's introduce the derivative operator $P_k(a)$ determined by the following formula:

$$P_0(f) = 0 \qquad (5.34)$$

$$P_k(f) = \qquad (5.35)$$

In this case the following equalities hold true:

:

$$p\int a\iint a\,dx\,dx\,dx - 2a\iiint a\,dx\,dx\,dx = p\int G_1(a)\,dx - 2P_1(a)$$

$$p\int a\iint a\iint a\,dx\,dx\,dx\,dx\,dx - 2a\iiint a\iint a\,dx\,dx\,dx\,dx\,dx$$

$$- 2a\iint a\iiint a\,dx\,dx\,dx\,dx\,dx = p\int G_2(a)\,dx - P_2(a)$$

$$p\int a\int\int a\int\int a\int\int a\,dx\,dx\,dx\,dx\,dx\,dx\,dx - 2\,a\int\int\int a\int\int a\int\int a\,dx\,dx\,dx\,dx\,dx\,dx\,dx$$

$$-\,2\,a\int\int a\int\int\int a\int\int a\,dx\,dx\,dx\,dx\,dx\,dx\,dx$$

$$-\,2\,a\int\int a\int\int a\int\int\int a\,dx\,dx\,dx\,dx\,dx\,dx\,dx = p\int G_3(a)\,dx - P_3(a)$$

Thus, formally, the required $\xi_k(2p)$ value can be represented by the following formula:

$$\xi_k(2p) = \left(\frac{\partial^{2p}}{\partial x^{2p}}\left(\eth\int G_k(\grave{a})\,dx - P_k(a)\right)\right) + 2\left(\sum_{s=0}^{k-1} t_{k,s}\left(\frac{\partial^{2p}}{\partial x^{2p}}\left(\eth\int G_s(\grave{a})\,dx - P_s(a)\right)\right)\right)$$
$$+\,2\left(\sum_{s=0}^{k-1} e_{k,s}\left(\frac{\partial^{2p}}{\partial x^{2p}}G_s(\grave{a})\right)\right) \qquad (5.36)$$

where $t_{k,s}$, $e_{k,s}$ are the required functional coefficients.
In order to find the above coefficients, we need to compare the known obtained formulas with those obtained previously with the use of conventional method.
Indeed, as follows from (5.36):

$$\xi_1(2p) = 2\,t_{1,1}\left(\frac{\partial^{2p}}{\partial x^{2p}}\left(p\int a\int\int a\,dx\,dx\,dx - a\int\int\int a\,dx\,dx\,dx\right)\right) + 2\,e_{1,0}\left(\frac{\partial^{2p}}{\partial x^{2p}}a\right)$$
$$+\,2\,t_{1,0}\,p\left(\frac{\partial^{2p}}{\partial x^{2p}}\left(\int a\,dx\right)\right)$$

in fact, according to (5.11)

$$\xi_1(2p) = 2\left(\frac{\partial^{2p}}{\partial x^{2p}}\left(p\int a\int\int a\,dx\,dx\,dx - a\int\int\int a\,dx\,dx\,dx\right)\right) + 2\int\int\int a\,dx\,dx\,dx\left(\frac{\partial^{2p}}{\partial x^{2p}}a\right)$$
$$+\,2\int\int a\,dx\,dx\,p\left(\frac{\partial^{2p}}{\partial x^{2p}}\left(\int a\,dx\right)\right)$$

Therefore,

$$t_{1,1} = 1,\quad t_{1,0} = \int\int a\,dx\,dx\,,\quad e_{1,0} =$$

Again, from (5.36) we receive for $k = 2$ :

$$\xi_2(2p) = 2\,t_{2,2}\left(\frac{\partial^{2p}}{\partial x^{2p}}\left(p\int a\int\int a\int\int a\,dx\,dx\,dx\,dx\,dx - a\int\int\int a\int\int a\,dx\,dx\,dx\,dx\,dx\right.\right.$$

$$
- a \int\!\!\int a \int\!\!\int\!\!\int a\, dx\, dx\, dx\, dx\, dx \Bigg)\Bigg)
$$

$$
+ 2\, t_{2,1} \left( \frac{\partial^{2p}}{\partial x^{2p}} \left( p \int a \int\!\!\int a\, dx\, dx\, dx - a \int\!\!\int\!\!\int a\, dx\, dx\, dx \right) \right) + 2\, p\, t_{2,0} \left( \frac{\partial^{2p}}{\partial x^{2p}} \left( \int a\, dx \right) \right)
$$

$$
+ 2\, e_{2,0} \left( \frac{\partial^{2p}}{\partial x^{2p}} a \right) + 2\, e_{2,1} \left( \frac{\partial^{2p}}{\partial x^{2p}} \left( a \int\!\!\int a\, dx\, dx \right) \right)
$$

In fact, according to (5.15):

$$
\xi_2(2p) = 2 \left( \frac{\partial^{2p}}{\partial x^{2p}} \left( p \int a \int\!\!\int a \int\!\!\int a\, dx\, dx\, dx\, dx\, dx - a \int\!\!\int\!\!\int a \int\!\!\int a\, dx\, dx\, dx\, dx\, dx \right.\right.
$$

$$
\left.\left. - a \int\!\!\int a \int\!\!\int\!\!\int a\, dx\, dx\, dx\, dx\, dx \right)\right)
$$

$$
+ 2 \int\!\!\int a\, dx\, dx \left( \frac{\partial^{2p}}{\partial x^{2p}} \left( p \int a \int\!\!\int a\, dx\, dx\, dx - a \int\!\!\int\!\!\int a\, dx\, dx\, dx \right) \right) + 2\, p \left( \left( \int\!\!\int a\, dx\, dx \right)^2 \right.
$$

$$
\left. + \int\!\!\int a \int\!\!\int a\, dx\, dx\, dx\, dx + 2 \int a \int\!\!\int\!\!\int a\, dx\, dx\, dx\, dx - 2 \int\!\!\int\!\!\int a\, dx\, dx\, dx \int a\, dx \right)
$$

$$
\left( \frac{\partial^{2p}}{\partial x^{2p}} \left( \int a\, dx \right) \right)
$$

$$
+ 2 \left( \int\!\!\int\!\!\int a \int\!\!\int\!\!\int a\, dx\, dx\, dx\, dx\, dx + \int\!\!\int\!\!\int a \int\!\!\int a\, dx\, dx\, dx\, dx\, dx \right) \left( \frac{\partial^{2p}}{\partial x^{2p}} a \right)
$$

$$
+ 2 \int\!\!\int\!\!\int a\, dx\, dx\, dx \left( \frac{\partial^{2p}}{\partial x^{2p}} \left( a \int\!\!\int a\, dx\, dx \right) \right)
$$

Therefore,

$$
t_{2,2} = 1, \quad t_{2,1} = \int\!\!\int a\, dx\, dx, \quad t_{2,0} = \begin{array}{c} \left( \int\!\!\int a\, dx\, dx \right)^2 + \int\!\!\int a \int\!\!\int a\, dx\, dx\, dx\, dx + 2 \int a \int\!\!\int\!\!\int a\, dx\, dx\, dx\, dx \\ - 2 \int\!\!\int\!\!\int a\, dx\, dx\, dx \int a\, dx \end{array}
$$

$$
e_{2,0} = \int\!\!\int a \int\!\!\int\!\!\int a\, dx\, dx\, dx\, dx\, dx + \int\!\!\int\!\!\int a \int\!\!\int a\, dx\, dx\, dx\, dx\, dx, \qquad e_{2,1} =
$$

Quite similarly, let's establish the values of these coefficients for $k = 3$ :

$$
t_{3,3} = 1, \quad t_{3,2} = \int\!\!\int a\, dx\, dx, \quad t_{3,1} = \begin{array}{c} \left( \int\!\!\int a\, dx\, dx \right)^2 + \int\!\!\int a \int\!\!\int a\, dx\, dx\, dx\, dx + 2 \int a \int\!\!\int\!\!\int a\, dx\, dx\, dx\, dx \\ - 2 \int\!\!\int\!\!\int a\, dx\, dx\, dx \int a\, dx \end{array},
$$

$$t_{3,0} = 2\int a \iint a \iiint a \, dx\, dx\, dx\, dx\, dx\, dx + 2\int a \iiint a \iint a \, dx\, dx\, dx\, dx\, dx\, dx + \left(\iint a \, dx\, dx\right)^3$$
$$+ 2\iint a \iint a \, dx\, dx\, dx\, dx \iint a \, dx\, dx + 4\int a \iiint a \, dx\, dx\, dx\, dx \iint a \, dx\, dx$$
$$- 2\int a\, dx \iiint a \iint a \, dx\, dx\, dx\, dx\, dx - 2\int a \iint a \, dx\, dx\, dx \iiint a \, dx\, dx\, dx$$
$$- 2\int a\, dx \iint a \iiint a \, dx\, dx\, dx\, dx\, dx - 2\iiint a \, dx\, dx\, dx \int a\, dx \iint a \, dx\, dx$$
$$+ \iint a \iint a \iint a \, dx\, dx\, dx\, dx\, dx\, dx$$

$$e_{3,0} = \left(\iint a \, dx\, dx\right)^2 \iiint a \, dx\, dx\, dx - 2\int a\, dx \left(\iiint a \, dx\, dx\, dx\right)^2$$
$$+ 2\int a \iiint a \, dx\, dx\, dx\, dx \iiint a \, dx\, dx\, dx + \iiint a \iint a \iint a \, dx\, dx\, dx\, dx\, dx\, dx\, dx$$
$$+ \iint a \iiint a \iint a \, dx\, dx\, dx\, dx\, dx\, dx\, dx + \iint a \iint a \iiint a \, dx\, dx\, dx\, dx\, dx\, dx\, dx$$

$$e_{3,1} = \iint a \iiint a \, dx\, dx\, dx\, dx\, dx + \iiint a \iint a \, dx\, dx\, dx\, dx\, dx, \qquad e_{3,2} =$$

Analyzing the obtained results, we come to the following conclusions:

1). $t_{k,k} = 1$ coefficient for all $k = 1, 2, 3$ ...

2). The following equalities hold true:

(5.37)

for all $k = 2, 3$ .....

This means that for each new value of $\xi_k(2p)$ it is suffieicent to find only the values of $t_{k,0}$ and $e_{k,0}$, because other coefficients are known from the formulas for $\xi_{k-1}(2p)$ defined previously.

As soon as $t_{k,0}$ is the coefficient for , and $e_{k,0}$ applies for $\frac{\partial^{2p}}{\partial x^{2p}} a$, in order to define the relevant formulas we shall use the following known equality:

$$\xi_k(2p) = \sum_{i_k=0}^{2p-2k} \xi_{k-1}(i_k + 2k - 2)\, \mathrm{C}(2p, i_k + 2k) \left( \frac{\partial^{2p-2k-i_k}}{\partial x^{2p-2k-i_k}} a \right)$$

By the substitution of the summation index

$i_k = i - 2k$

we can represent it in equivalent form:

$$\xi_k(2p) = \sum_{i=2k}^{2p} \xi_{k-1}(i-2)\, C(2p,i) \left( \frac{\partial^{2p-i}}{\partial x^{2p-i}} a \right)$$

or

$$\xi_k(2p) = \left( \sum_{i=0}^{2p} \xi_{k-1}(i-2)\, C(2p,i) \left( \frac{\partial^{2p-i}}{\partial x^{2p-i}} a \right) \right) - \left( \sum_{i=0}^{2k-1} \xi_{k-1}(i-2)\, C(2p,i) \left( \frac{\partial^{2p-i}}{\partial x^{2p-i}} a \right) \right)$$

Because the required values of $t_{k,0}$ and $e_{k,0}$ coefficients are defined only for and $\frac{\partial^{2p}}{\partial x^{2p}} a$, the first summand in the obtained formula is discarded, and the following is obtained from the second summand:

$$\xi_{k-1}(i-2)\, C(2p,i) \left( \frac{\partial^{2p-i}}{\partial x^{2p-i}} a \right)$$

With $i = 0$ this expression assumes the following form:

By definition, we obtain from here the required formula:

Similarly, assuming that $i = 1$ we receive:

$$\xi_{k-s-1}(-1)\, 2p \left( \frac{\partial^{2p-1}}{\partial x^{2p-1}} a \right)$$

and, thus,

Thus, the required formula for $\xi_k(2p)$ takes the following form:

$$\xi_k(2p) = 2\left(\left(\frac{\partial^{2p}}{\partial x^{2p}}\left(\eth\int G_k(\grave{a})\,dx - P_k(a)\right)\right) - \left(\sum_{s=0}^{k-1}\xi_{k-s-1}(-1)\left(\frac{\partial^{2p}}{\partial x^{2p}}\left(\eth\int G_s(\grave{a})\,dx - P_s(a)\right)\right)\right)\right) - \left(\sum_{s=0}^{k-1}\xi_{k-s-1}(-2)\left(\frac{\partial^{2p}}{\partial x^{2p}}G_s(\grave{a})\right)\right)$$

By the substitution of $p$ for the $2p$ expression, we can transform this formula as follows:

$$\xi_k(p) = \left(\frac{\partial^{p}}{\partial x^{p}}\left(\eth\int G_k(\grave{a})\,dx - 2P_k(a)\right)\right) - \left(\sum_{s=0}^{k-1}\xi_{k-s-1}(-1)\left(\frac{\partial^{p}}{\partial x^{p}}\left(\eth\int G_s(\grave{a})\,dx - 2P_s(a)\right)\right)\right) - \left(\sum_{s=0}^{k-1}\xi_{k-s-1}(-2)\left(\frac{\partial^{p}}{\partial x^{p}}G_s(\grave{a})\right)\right)$$

**The theorem is proved.**

Thus, we have obtained the formula for the successive calculaton of the required $\xi_k(p)$ coefficients.

## 6. Formula for the common solution of a linear second-order ODE.

As demonstrated above, for a linear ODE

$$\frac{\partial^2}{\partial x^2}y = a\,y \qquad (6.1)$$

a general solution is determined by the following formula:

$$y(x) = (C_1 - C_2 x)(-1 + \alpha(-1)) + C_2\,\alpha(-2) \qquad (6.2)$$

As soon as

(6.3)

where

and we arrive to the task of determining the sum of a finite series

$$\sum_{k=0}^{p-1}\xi_k(p)$$

from there we can obtain the values for (6.3) equalities.

This task can be accomplished similarly to classical method employed by Newton for the expansion of, say, $(1+x)^m$ functions, where $m$ is a certain real number, to power series like:

It turnes out that in the cases where $x$ belongs to the convergence interval of

the following equality holds true:

.

According to this hypothesis, the formula for the *n*-image coefficient of $\alpha(p)$ can be written down as follows:

(6.4)

Therefore, the desired formulas (6.3) take the following form:

(6.5)

(6.6)

Let’s specify the values for functions.

Earlier, in the summary of the theory of indefinite integral calculation, it was shown that:

Thus, the desired formulas for the successive calculation of coeffieicents are as follows:

$$\xi_0(-2) = -2\int\int\int \grave{a}\,dx\,dx\,dx \qquad (6.7)$$

$$\xi_k(-1) = -\int\int G_k(\grave{a})\,dx\,dx - 2\int P_k(a)\,dx$$
$$+\left(\sum_{s=0}^{k-1}\xi_{k-s-1}(-1)\left(\int\int G_s(\grave{a})\,dx\,dx + 2\int P_s(a)\,dx\right)\right)$$
$$-\left(\sum_{s=0}^{k-1}\xi_{k-s-1}(-2)\int G_s(\grave{a})\,dx\right), k = 1, 2, 3 .. N \qquad (6.8)$$

$$\xi_k(-2) = -2\left(\int\int\int G_k(\grave{a})\,dx\,dx\,dx + \int\int P_k(a)\,dx\,dx\right)$$
$$+2\left(\sum_{s=0}^{k-1}\xi_{k-s-1}(-1)\left(\int\int\int G_s(\grave{a})\,dx\,dx\,dx + \int\int P_s(a)\,dx\,dx\right)\right)$$
$$-\left(\sum_{s=0}^{k-1}\xi_{k-s-1}(-2)\int\int G_s(\grave{a})\,dx\,dx\right) \qquad (6.9)$$

**Proposition 6.1. Series (6.5), (6.6) have a finite sum, that is, they converge:**

**Proof:** by definition, the $\xi_k(2p)$ summands, where are determined by the following formulas:

$$\xi_0(2p) = 2p\left(\frac{\partial^{2p-1}}{\partial x^{2p-1}}a\right)$$

$$\xi_1(2p) = \sum_{i_1=0}^{2p-2} i_1\left(\frac{\partial^{i_1-1}}{\partial x^{i_1-1}}a\right)C(2p, i_1+2)\left(\frac{\partial^{2p-2-i_1}}{\partial x^{2p-2-i_1}}a\right)$$

$$\xi_2(2p) = \sum_{i_2=0}^{2p-4}\left(\sum_{i_1=0}^{i_2} i_1\left(\frac{\partial^{i_1-1}}{\partial x^{i_1-1}}a\right)C(i_2+2, i_1+2)\left(\frac{\partial^{i_2-i_1}}{\partial x^{i_2-i_1}}a\right)\right)C(2p, i_2+4)$$

$$\left(\sum_{i_2=0}^{i_3}\left(\sum_{i_1=0}^{i_2} i_1\left(\frac{\partial^{i_1-1}}{\partial x^{i_1-1}}a\right)C(i_2+2, i_1+2)\left(\frac{\partial^{i_2-i_1}}{\partial x^{i_2-i_1}}a\right)\right)C(i_3+4, i_2+4)\left(\frac{\partial^{i_3-i_2}}{\partial x^{i_3-i_2}}a\right)\right)$$
$$C(2p, i_3+6)\left(\frac{\partial^{2p-6-i_3}}{\partial x^{2p-6-i_3}}a\right)$$

$$\xi_4(2p) = \sum_{i_4=0}^{2p-8}\left(\sum_{i_3=0}^{i_4}\left(\sum_{i_2=0}^{i_3}\left(\sum_{i_1=0}^{i_2} i_1\left(\frac{\partial^{i_1-1}}{\partial x^{i_1-1}}a\right)C(i_2+2, i_1+2)\left(\frac{\partial^{i_2-i_1}}{\partial x^{i_2-i_1}}a\right)\right)\right.\right.$$
$$\left.\left.C(i_3+4, i_2+4)\left(\frac{\partial^{i_3-i_2}}{\partial x^{i_3-i_2}}a\right)C(i_4+6, i_3+6)\left(\frac{\partial^{i_4-i_3}}{\partial x^{i_4-i_3}}a\right)\right)C(2p, i_4+8)\left(\frac{\partial^{2p-8-i_4}}{\partial x^{2p-8-i_4}}a\right)\right)$$

$$\xi_s(2p) = \sum_{i_1=[i_s, i_1]}^{[2p-2s, i_s, i_2]} i_1\left(\frac{\partial^{i_1-1}}{\partial x^{i_1-1}}a\right)\left(\prod_{k=2}^{s-1} C(i_k+2(k-1), i_1+2(k-1))\left(\frac{\partial^{i_k-i_{k-1}}}{\partial x^{i_k-i_{k-1}}}a\right)\right)$$
$$C(2p, i_s+2s)\left(\frac{\partial^{2p-2s-i_s}}{\partial x^{2p-2s-i_s}}a\right)$$

In this example, the following series:

=

can be represented as follows:

$$\alpha(-1) = \lim_{p\to\left(-\frac{1}{2}\right)}\left(\sum_{k=0}^{p-1}\xi_k(2p)\right) + \left(\sum_{k=p}^{\infty}\xi_k(2p)\right)$$

As soon as, by definition, all $\xi_k(2p)$ starting from $k = p$ are zero, then and, thus,

$$\sum_{k=0}^{\infty} \xi_k(-1) = \lim_{p \to \left(-\frac{1}{2}\right)} \sum_{k=0}^{p-1} \xi_k(2p)$$

Since, in accordance with (5.38):

$$\xi_{p-1}(2p) = a^{(p-1)} p\,(p+1)\left(\frac{\partial}{\partial x} a\right)$$

the value of $\xi_{p-2}(2p)$ is defined by (5.39), and all values of $\xi_{p-s}(2p)$, $s = 1 .. p$ are defined as well, the following series:

$$\sum_{k=0}^{\infty} \xi_k(-1)$$

is defined by a finite sum of strictly specified summands.

Obviously, in the $a < 1$ interval

$$\lim_{p \to \infty} a^{(p-1)} p\,(p+1)\left(\frac{\partial}{\partial x} a\right) = 0$$

Consequently,

Thus, the necessary condition of convergence is satisfied. And, as soon as this series is a finite expression, this means that it has a certain sum, that is, converges. **The proposition is proved.**
The general solution in the form of (6 2) is analytical. It can be used to obtain with the desired accuracy a general solution of the ODE (6.1) in a specified variation interval of the independent variable $x$ .

**Conclusions: A linear ODE**

$$\frac{\partial^2}{\partial x^2} y = a\,y \tag{6.8}$$

**where *a* is an arbitrary differentiable and integrable function, has the following general solution:**

$$\text{y}(x) = (C_1 - C_2\,x)\left(-1 + \left(\sum_{k=0}^{\infty} \xi_k(-1)\right)\right) + C_2\left(\sum_{k=0}^{\infty} \xi_k(-2)\right) \tag{6.9}$$

**where the desired values of** $\xi_k(-1), \xi_k(-2), k = 1, 2, 3 .. \infty$ **functions are defined by the following formulas**: **(6.10)**

$$\xi_0(-2) = -2 \int\!\!\int\!\!\int à\; dx\, dx\, dx$$

$$\xi_k(-1) = -\int\int G_k(\grave{a})\,dx\,dx - 2\int P_k(a)\,dx$$
$$+\left(\sum_{s=0}^{k-1} \xi_{k-s-1}(-1)\left(\int\int G_s(\grave{a})\,dx\,dx + 2\int P_s(a)\,dx\right)\right)$$
$$-\left(\sum_{s=0}^{k-1} \xi_{k-s-1}(-2)\int G_s(\grave{a})\,dx\right), k = 1, 2, 3 .. \infty$$

$$\xi_k(-2) = -2\left(\int\int\int G_k(\grave{a})\,dx\,dx\,dx + \int\int P_k(a)\,dx\,dx\right)$$
$$+ 2\left(\sum_{s=0}^{k-1} \xi_{k-s-1}(-1)\left(\int\int\int G_s(\grave{a})\,dx\,dx\,dx + \int\int P_s(a)\,dx\,dx\right)\right)$$
$$-\left(\sum_{s=0}^{k-1} \xi_{k-s-1}(-2)\int\int G_s(\grave{a})\,dx\,dx\right)$$

**where $G_i(a)$ is an integral operator, which operates on the specified $a$ function according to the following rules:**

$$G_0(a) = a$$

=

=

**and so on. With $k = 0$ the $P_k(a)$ operator is determined by thefollowing formula:**

$$P_0(a) = 0$$

**and operates on the specified *a* function according to the following rule**:

$$P_k(a) = \sum_{i=0}^{k-1} G_{k-i}\left(\int G_i(a)\,dx\right), k = 1, 2, 3$$

**Thus, all conditions required to define a general integration algorithm for a linear ODE (6.8), which are identified in Definition 1.1., are met.**

## 7. Examples

Let's discuss some specific examples of the construction of analytical solution. First of all, in order to demonstrate efficiency and generality of the above-described theory, let's consider the well-known "classroom" examples.

**Example No. 1. Find a general solution for the following linear ODE:**

$$\frac{\partial^2}{\partial x^2} y = x\,y \qquad (7.1)$$

**Solution:** in this example $\grave{a} = x$ .

Let's calculate the $\xi_k(-1)$ , $\xi_k(-2)$ $k = 0, 1, 2$ ..... summands with the use of (6.10).

Then we obtain:

$$= -\frac{x^3}{6} \qquad \xi_0(-2) = -2\int\int\int x\,dx\,dx\,dx = -\frac{x^4}{12}$$

$$\xi_1(-1) = -\int\int G_1(\dot a)\,dx\,dx + 2\int P_1(a)\,dx + \xi_0(-1)\int\int G_0(\dot a)\,dx\,dx - \xi_0(-2)\int G_0(\dot a)\,dx$$

$$\xi_1(-2) = -2\int\int\int G_1(\dot a)\,dx\,dx\,dx - 2\int\int P_1(a)\,dx\,dx$$
$$+ 2\xi_0(-1)\left(\int\int\int G_0(\dot a)\,dx\,dx\,dx + \int\int P_0(a)\,dx\,dx\right) - \xi_0(-2)\int\int G_0(\dot a)\,dx\,dx$$

or, in expanded presentation:

$$\xi_1(-1) = -\int\int x\int\int x\,dx\,dx\,dx\,dx - 2\int x\int\int\int x\,dx\,dx\,dx\,dx + \left(-\frac{x^3}{6}\right)\int\int x\,dx\,dx - \left(-\frac{x^4}{12}\right)\int x\,dx = -\frac{x^6}{180}$$

$$\xi_1(-2) = -2\int\int\int x\int\int x\,dx\,dx\,dx\,dx\,dx - 2\int\int x\int\int\int x\,dx\,dx\,dx\,dx\,dx + 2\left(-\frac{x^3}{6}\right)\int\int\int x\,dx\,dx\,dx - \left(-\frac{x^4}{12}\right)\int\int x\,dx\,dx = -\frac{x^7}{280}$$

Let's calculate $\xi_2(-1)$, $\xi_2(-2)$:

$$\xi_2(-1) = -\int\int G_2(\dot a)\,dx\,dx - 2\int P_2(a)\,dx + \left(\sum_{s=0}^{1}\xi_{1-s}(-1)\left(\int\int G_s(\dot a)\,dx\,dx + 2\int P_s(a)\,dx\right)\right) - \left(\sum_{s=0}^{1}\xi_{1-s}(-2)\int G_s(\dot a)\,dx\right) = -\frac{x^9}{12960}$$

$$\xi_2(-2) = -2\left(\int\int\int G_2(\dot a)\,dx\,dx\,dx + \int\int P_2(a)\,dx\,dx\right) + 2\left(\sum_{s=0}^{1}\xi_{1-s}(-1)\left(\int\int\int G_s(\dot a)\,dx\,dx\,dx + \int\int P_s(a)\,dx\,dx\right)\right) - \left(\sum_{s=0}^{1}\xi_{1-s}(-2)\int\int G_s(\dot a)\,dx\,dx\right) = -\frac{x^{10}}{18144}$$

In the similar manner we can obtain:

$$;\ ;\ \xi_5(-1) = -\frac{x^{18}}{109930867200}\ ;\ \xi_6(-1) = -\frac{x^{21}}{46170964224000}\ ;\ \xi_7(-1) = -\frac{x^{24}}{25486372251648000}$$

an so on.

$$;\ \xi_4(-2) = -\frac{41\,x^{16}}{18681062400}\ ;\ \xi_5(-2) = -\frac{89\,x^{19}}{12067966310400}\ ;\ \xi_6(-2) = -\frac{409\,x^{22}}{22808456326656000}\ ;$$

$$\xi_7(-2) = -\frac{x^{25}}{30279624560640000}$$ an so on.

In this case partial solutions are defined by the following formulals:

$$y_1(x) = -1 + \left( \sum_{k=0}^{N} \xi_k(-1) \right) \qquad (7.2)$$

$$y_2(x) = \left( \sum_{k=0}^{N} \xi_k(-2) \right) - x \left( \sum_{k=0}^{N} \xi_k(-1) \right) + x \qquad (7.3)$$

where **N** is the number of summands required to obtain the desired solution with the specified accuracy, within a defined range of independent argument ***x*** .

Assuming, in particular, $N = 7$ , we get:

$$y_1(x) = -\frac{x^{24}}{25486372251648000} - \frac{x^{21}}{46170964224000} - \frac{x^{18}}{109930867200} - \frac{x^{15}}{359251200} - \frac{x^{12}}{1710720} - \frac{x^{9}}{12960} - \frac{x^{6}}{180} - \frac{x^{3}}{6} - 1$$

$$y_2(x) = \frac{x^{25}}{161000868188160000} + \frac{x^{22}}{268334780313600} + \frac{x^{19}}{580811212800} + \frac{x^{16}}{1698278400} + \frac{x^{13}}{7076160} + \frac{x^{10}}{45360} + \frac{x^{7}}{504} + \frac{x^{4}}{12} + x$$

As follows from the above, the desired analytical solutions can be represented by the following formulas:

(7.4)

(7.5)

where $\mathrm{B}(k)$ is a numeric function of $k$ parameter. In particular:
$\mathrm{B}(0) = 1, \mathrm{B}(1) = 6, \mathrm{B}(2) = 180, \mathrm{B}(3) = 12960$ ...... and so on.

$\mathrm{H}(k)$ is a numeric function of $k$ parameter. In particular: $\mathrm{H}(0) = 1, \mathrm{H}(1) = 12, \mathrm{H}(2) = 504, \mathrm{H}(3) = 45360$ ....... and so on.

Substituting the obtained representations to the initial ODE (7.1), we can define recurrence formulas for $\mathrm{B}(k)$ and $\mathrm{H}(k)$ . Indeed, in particular upon the substitution of $y_1(x)$ to the ODE (7.1), we get:

$$\frac{\partial^2}{\partial x^2} \left( -\left( \sum_{k=0}^{\infty} \frac{x^{(3k)}}{\mathrm{B}(k)} \right) \right) = x \left( -\left( \sum_{k=0}^{\infty} \frac{x^{(3k)}}{\mathrm{B}(k)} \right) \right)$$

or

$$-\left( \sum_{k=0}^{\infty} \left( \frac{9 x^{(3k-2)} k^2}{\mathrm{B}(k)} - \frac{3 x^{(3k-2)} k}{\mathrm{B}(k)} \right) \right) + \left( \sum_{k=0}^{\infty} \frac{x^{(3k+1)}}{\mathrm{B}(k)} \right) = 0 \qquad (7.6)$$

Upon the replacement of $k$ index by $k+1$ in the summand:

we obtain the following equality:

$$\sum_{k=0}^{\infty} \frac{x^{(3k+1)}}{B(k)} = \sum_{k=1}^{\infty} \frac{x^{(-2+3k)}}{B(-1+k)}$$

Substituting this to (7.6), upon transformations we get:

$$-\left(\sum_{k=1}^{\infty} \left(\frac{9k^2}{B(k)} - \frac{3k}{B(k)} - \frac{1}{B(-1+k)}\right) x^{(3k-2)}\right) = 0$$

From here we obtain reccurence equation to find $B(k)$:

$$B(k+1) = (9(k+1)^2 - 3k - 3)\,B(k)$$

The solution, taking into account the initial condition of $B(0) = 1$, is as follows:

$$B(k) = \frac{9^k \Gamma\left(k + \frac{2}{3}\right) \Gamma(k+1)}{\Gamma\left(\frac{2}{3}\right)}$$

Thus, the first partial solution (7.4) takes the following form:

$$y_1(x) = -\left(\sum_{k=0}^{\infty} \frac{\Gamma\left(\frac{2}{3}\right) x^{(3k)}}{9^k \Gamma\left(k + \frac{2}{3}\right) \Gamma(k+1)}\right) = -\frac{\sqrt{x}\,\mathrm{BesselI}\left(-\frac{1}{3}, \frac{2\sqrt{x^3}}{3}\right) \Gamma\left(\frac{2}{3}\right) 3^{\left(\frac{2}{3}\right)}}{3}$$

We can check that this solution satisifies the initial ODE (7.1).
Proceeding in a similar way, by the substitution of the second partial solution (7.4) to the initial ODE (7.1) we obtain:

$$\frac{\partial^2}{\partial x^2}\left(\sum_{k=0}^{\infty} \frac{x^{(3k+1)}}{H(k)}\right) = x\left(\sum_{k=0}^{\infty} \frac{x^{(3k+1)}}{H(k)}\right)$$

Upon transformations:

$$\left(\sum_{k=0}^{\infty} \left(-\frac{x^{(3k-1)}(3k+1)}{H(k)} + \frac{x^{(3k-1)}(3k+1)^2}{H(k)}\right)\right) - \left(\sum_{k=0}^{\infty} \frac{x^{(2+3k)}}{H(k)}\right) = 0$$

Substituting $k+1$ for the summation index $k$ in the second summand, we get:

$$\sum_{k=0}^{\infty} \frac{x^{(2+3k)}}{H(k)} = \sum_{k=1}^{\infty} \frac{x^{(3k-1)}}{H(k-1)}$$

Thus,

$$\sum_{k=0}^{\infty} \frac{x^{(3k-1)}(3\,H(k-1)\,k + 9\,H(k-1)\,k^2 - H(k))}{H(k)\,H(k-1)} = 0$$

And we obtain the following recurrence equation:

$$\mathrm{H}(k+1) = (3k + 3 + 9(k+1)^2)\,\mathrm{H}(k)$$

With consideration of the initial condition, $\mathrm{H}(0) = 1$, the solution is:

$$\mathrm{H}(k) = \frac{3\,9^k\,\Gamma\left(k+\frac{4}{3}\right) 3^{\left(\frac{1}{2}\right)} \Gamma\left(\frac{2}{3}\right) \Gamma(k+1)}{2\pi}$$

Thus, the second partial solution of the initial ODE (7.1) is equal to:

$$y_2(x) = 2\pi \left( \sum_{k=0}^{\infty} \frac{x^{(3k+1)}}{3^{(2k+1)}\,\Gamma\left(k+\frac{4}{3}\right)\sqrt{3}\,\Gamma\left(\frac{2}{3}\right)\Gamma(k+1)} \right) = \frac{2\pi\sqrt{x}\,\mathrm{BesselI}\left(\frac{1}{3}, \frac{2\sqrt{x^3}}{3}\right) 3^{\left(\frac{5}{6}\right)}}{9\,\Gamma\left(\frac{2}{3}\right)}$$

As proved by a check, this solution satisfies the ODE (7.1).

Thus, the general solution of ODE (7.1) is equal to:

$$\mathrm{y}(x) = \frac{C_1\,\mathrm{BesselI}\left(-\frac{1}{3}, \frac{2\sqrt{x^3}}{3}\right)\Gamma\left(\frac{2}{3}\right) 3^{\left(\frac{2}{3}\right)} x^{\left(\frac{1}{2}\right)}}{3} + \frac{C_2\,2\pi\sqrt{x}\,\mathrm{BesselI}\left(\frac{1}{3}, \frac{2\sqrt{x^3}}{3}\right) 3^{\left(\frac{5}{6}\right)}}{9\,\Gamma\left(\frac{2}{3}\right)}$$

**The task is solved.**

**Example No. 2. Find a general solution for the following linear ODE:**

$$\frac{\partial^2}{\partial x^2} y = \mathbf{e}^x\, y \tag{7.7}$$

**Solution:** in this example $\dot{a} = \mathbf{e}^x$. Let's calculate the $\xi_k(-1)$, $\xi_k(-2)$ $k = 0, 1, 2$ ..... summands. We obtain the following:

**a) Zero summands $\xi_0(-1)$, $\xi_0(-2)$:**

$$= -\mathbf{e}^x \qquad \xi_0(-2) = -2\int\!\!\int\!\!\int \mathbf{e}^x\,dx\,dx\,dx = -2\,\mathbf{e}^x$$

**b) First summands $\xi_1(-1)$, $\xi_1(-2)$:**

$$\xi_1(-1) = -\int\!\!\int G_1(\dot{a})\,dx\,dx + 2\int P_1(a)\,dx + \xi_0(-1)\int\!\!\int G_0(\dot{a})\,dx\,dx - \xi_0(-2)\int G_0(\dot{a})\,dx$$

$$\xi_1(-2) = -2\int\!\!\int\!\!\int G_1(\dot{a})\,dx\,dx\,dx - 2\int\!\!\int P_1(a)\,dx\,dx$$
$$+ 2\,\xi_0(-1)\left(\int\!\!\int\!\!\int G_0(\dot{a})\,dx\,dx\,dx + \int\!\!\int P_0(a)\,dx\,dx\right) - \xi_0(-2)\int\!\!\int G_0(\dot{a})\,dx\,dx$$

or

$$\xi_1(-1) = -\int\int \mathbf{e}^x \int\int \mathbf{e}^x\, dx\, dx\, dx\, dx - 2\int \mathbf{e}^x \int\int\int \mathbf{e}^x\, dx\, dx\, dx\, dx + (-\mathbf{e}^x)\int\int \mathbf{e}^x\, dx\, dx - (-2\,\mathbf{e}^x)\int \mathbf{e}^x\, dx = -\frac{\mathbf{e}^{(2x)}}{4}$$

$$\xi_1(-2) = -2\int\int\int \mathbf{e}^x \int\int \mathbf{e}^x\, dx\, dx\, dx\, dx\, dx - 2\int\int \mathbf{e}^x \int\int\int \mathbf{e}^x\, dx\, dx\, dx\, dx\, dx + 2(-\mathbf{e}^x)\int\int\int \mathbf{e}^x\, dx\, dx\, dx - (-2\,\mathbf{e}^x)\int\int \mathbf{e}^x\, dx\, dx = -\frac{3\,\mathbf{e}^{(2x)}}{4}$$

**c) Second summands** $\xi_2(-1)$, $\xi_2(-2)$:

$$\xi_2(-1) = -\int\int G_2(\grave{a})\, dx\, dx - 2\int P_2(a)\, dx + \left(\sum_{s=0}^{1} \xi_{1-s}(-1)\left(\int\int G_s(\grave{a})\, dx\, dx + 2\int P_s(a)\, dx\right)\right) - \left(\sum_{s=0}^{1} \xi_{1-s}(-2)\int G_s(\grave{a})\, dx\right) = -\frac{\mathbf{e}^{(3x)}}{36}$$

$$\xi_2(-2) = -2\left(\int\int\int G_2(\grave{a})\, dx\, dx\, dx + \int\int P_2(a)\, dx\, dx\right) + 2\left(\sum_{s=0}^{1} \xi_{1-s}(-1)\left(\int\int\int G_s(\grave{a})\, dx\, dx\, dx + \int\int P_s(a)\, dx\, dx\right)\right) - \left(\sum_{s=0}^{1} \xi_{1-s}(-2)\int\int G_s(\grave{a})\, dx\, dx\right) = -\frac{11\,\mathbf{e}^{(3x)}}{108}$$

Similary we can find:

$$;\ ;\ ;\ ;\ \xi_7(-1) = -\frac{\mathbf{e}^{(8x)}}{1625702400} \quad \text{an so on.}$$

$$;\ ;\ ;\ ;\ \xi_7(-2) = -\frac{761\,\mathbf{e}^{(8x)}}{227598336000} \quad \text{an so on.}$$

In this case partial solutions are determined by the following formulas:

$$y_1(x) = -1 + \left(\sum_{k=0}^{N} \xi_k(-1)\right) \tag{7.8}$$

$$y_2(x) = \left(\sum_{k=0}^{N} \xi_k(-2)\right) - x\left(\sum_{k=0}^{N} \xi_k(-1)\right) + x \tag{7.9}$$

where **N** is the number of summands required to obtain the desired solution with specified accuracy in a given interval of independent argument ***x***.

Assuming that $N = 7$ in (7.8), (7.9), we obtain:

$$y_1(x) = -1 - \mathbf{e}^x - \frac{\mathbf{e}^{(2x)}}{4} - \frac{\mathbf{e}^{(3x)}}{36} - \frac{\mathbf{e}^{(4x)}}{576} - \frac{\mathbf{e}^{(5x)}}{14400} - \frac{\mathbf{e}^{(6x)}}{518400} - \frac{\mathbf{e}^{(7x)}}{25401600} - \frac{\mathbf{e}^{(8x)}}{1625702400}$$

$$y_2(x) = (-2+x)\,\mathbf{e}^x + \left(\frac{x}{4} - \frac{3}{4}\right)\mathbf{e}^{(2x)} + \left(\frac{x}{36} - \frac{11}{108}\right)\mathbf{e}^{(3x)} + \left(\frac{x}{576} - \frac{25}{3456}\right)\mathbf{e}^{(4x)}$$
$$+ \left(-\frac{137}{432000} + \frac{x}{14400}\right)\mathbf{e}^{(5x)} + \left(\frac{x}{518400} - \frac{49}{5184000}\right)\mathbf{e}^{(6x)}$$
$$+ \left(\frac{x}{25401600} - \frac{121}{592704000}\right)\mathbf{e}^{(7x)} + \left(-\frac{761}{227598336000} + \frac{x}{1625702400}\right)\mathbf{e}^{(8x)} + x$$

This means that the desired analytical solutions can be represented as follows:

(7.10)

$$y_2(x) = \sum_{k=0}^{\infty}\left(\frac{x}{\mathrm{H}(k)} - \mathrm{F}(k)\right)\mathbf{e}^{(kx)} \qquad (7.11)$$

where $\mathrm{B}(k)$ is a certain numeric function of $k$ parameter. In particular:
B(0) = 1, B(1) = 1, B(2) = 4, B(3) = 36; ...... and so on.

$\mathrm{H}(k)$ is a certain numeric function of $k$ parameter. In particular: $\mathrm{H}(0) = 1, \mathrm{H}(1) = 1, \mathrm{H}(2) = 4, \mathrm{H}(3) = 36$ ...... and so on.

F*(k)* is a certain numeric function of $k$ parameter. In particular:
F(0) = 1, F(1) = 2, F(2) = 3/4, F(3) = 25/3456...... and so on.

Substituting the obtained representations to the initial ODE (7.7), we can obtain recurring formulas for $\mathrm{B}(k)$ and $\mathrm{H}(k)$. Indeed, particularly for $y_1(x)$, we get:

$$\frac{\partial^2}{\partial x^2}\left(-\left(\sum_{k=0}^{\infty}\frac{\mathbf{e}^{(kx)}}{\mathrm{B}(k)}\right)\right) = \mathbf{e}^x\left(-\left(\sum_{k=0}^{\infty}\frac{\mathbf{e}^{(kx)}}{\mathrm{B}(k)}\right)\right)$$

or

$$\sum_{k=0}^{\infty}\left(-\frac{k^2\,\mathbf{e}^{(kx)}}{\mathrm{B}(k)} + \frac{\mathbf{e}^{(x+kx)}}{\mathrm{B}(k)}\right) = 0 \qquad (7.12)$$

Let’s substitute $k+1$ for $k$ index in the following summand:

We receive the following equality:

$$\sum_{k=0}^{\infty}\frac{\mathbf{e}^{(x+kx)}}{\mathrm{B}(k)} = \sum_{k=1}^{\infty}\frac{\mathbf{e}^{(kx)}}{\mathrm{B}(-1+k)}$$

Substituting it to (7.12), upon transformation we obtain:

$$\sum_{k=1}^{\infty}\left(-\frac{k^2}{\mathrm{B}(k)} + \frac{1}{\mathrm{B}(-1+k)}\right)\mathbf{e}^{(kx)} = 0$$

From here follows the reccurence equation to obtain $B(k)$.

The solution, with consideration of initial condition $B(0) = 1$, is as follows:

Thus, partial solution of (7.10) takes the following form:

$$y_1(x) = -\left(\sum_{k=0}^{\infty} \frac{\mathbf{e}^{(k\,x)}}{\Gamma(k+1)^2}\right) =$$

We can verify that this solution satisfies the initial ODE (7.7).

Quite similarly, substituting the second partial solution (7.11) to the initial ODE (7.7) we obtain:

$$\frac{\partial^2}{\partial x^2}\left(\sum_{k=0}^{\infty}\left(\frac{x}{H(k)} - F(k)\right)\mathbf{e}^{(k\,x)}\right) = \mathbf{e}^x\left(\sum_{k=0}^{\infty}\left(\frac{x}{H(k)} - F(k)\right)\mathbf{e}^{(k\,x)}\right)$$

Upon transformations:

$$\sum_{k=0}^{\infty}\left(\frac{2\,k\,\mathbf{e}^{(k\,x)}}{H(k)} + \left(\frac{x}{H(k)} - F(k)\right)k^2\,\mathbf{e}^{(k\,x)}\right) = \mathbf{e}^x\left(\sum_{k=0}^{\infty}\left(\frac{x}{H(k)} - F(k)\right)\mathbf{e}^{(k\,x)}\right)$$

Replacing the second summand of the summation index $k$ by $k+1$ we get:

$$\mathbf{e}^x\left(\sum_{k=0}^{\infty}\left(\frac{x}{H(k)} - F(k)\right)\mathbf{e}^{(k\,x)}\right) = \sum_{k=1}^{\infty}\frac{\mathbf{e}^{(k\,x)}\,(x - F(k-1)\,H(k-1))}{H(k-1)}$$

Thus

$$\sum_{k=1}^{\infty}(-\,\mathbf{e}^{(k\,x)}\,(-2\,k\,H(k-1) - k^2\,H(k-1)\,x + k^2\,H(k-1)\,F(k)\,H(k) + H(k)\,x$$
$$-\,H(k)\,F(k-1)\,H(k-1))/(H(k)\,H(k-1))) = 0$$

From here we obtain the following recurrence equations:

$$-2\,k\,H(k-1) + k^2\,H(k-1)\,F(k)\,H(k) - H(k)\,F(k-1)\,H(k-1) = 0$$

If the initial conditions: $H(0) = 1$, **F(0) = 0** are satisfied, they take the following form:

$$H(k+1) = (k+1)^2\,H(k)$$

$$-2\,(k+1)\,H(k) + (k+1)^2\,H(k)\,F(k+1)\,H(k+1) - H(k+1)\,F(k)\,H(k) = 0$$

Solution of the first recurrence equation is the following function:

In this case, the second equation assumes the following form:

$$F(k+1) = \frac{2k+2+F(k)\,\Gamma(k+2)^2}{\Gamma(k+2)^2\,(k+1)^2}$$

The solution is presented by the following function:

Thus, the second partial solution of the initial ODE (7.7) is as follows:

$$y_2(x) = \sum_{k=0}^{\infty} \frac{(x-2\,\Psi(k+1)-2\gamma)\,\mathbf{e}^{(kx)}}{\Gamma(k+1)^2} \qquad (7.13)$$

Let’s prove that this series converges with the use of D'Alembert criterion. Indeed,

$$\lim_{k\to\infty} \frac{(x-2\,\Psi(k+2)-2\gamma)\,\mathbf{e}^{((k+1)x)}\,\Gamma(k+1)^2}{\Gamma(k+2)^2\,(x-2\,\Psi(k+1)-2\gamma)\,\mathbf{e}^{(kx)}} = 0$$

Now let’s prove that $y_2(x)$ is a solution of the ODE (7.7).
For this purpose, we shall write it down as follows:

$$y_2(x) = \lim_{N\to\infty} \sum_{k=0}^{N} \frac{(x-2\,\Psi(k+1)-2\gamma)\,\mathbf{e}^{(kx)}}{\Gamma(k+1)^2}$$

Substituting this function to ODE (7.7), we get:

$$\lim_{N\to\infty} \left(\frac{\partial^2}{\partial x^2}\left(\sum_{k=0}^{N} \frac{(x-2\,\Psi(k+1)-2\gamma)\,\mathbf{e}^{(kx)}}{\Gamma(k+1)^2}\right)\right) - \mathbf{e}^x\left(\sum_{k=0}^{N} \frac{(x-2\,\Psi(k+1)-2\gamma)\,\mathbf{e}^{(kx)}}{\Gamma(k+1)^2}\right) = 0$$

Upon differentiation we obtain:

$$\lim_{N\to\infty} \left(\sum_{k=0}^{N}\left(\frac{2\,\mathbf{e}^{(kx)}\,k}{\Gamma(k+1)^2} + \frac{(x-2\,\Psi(k+1)-2\gamma)\,\mathbf{e}^{(kx)}\,k^2}{\Gamma(k+1)^2}\right)\right) - \mathbf{e}^x\left(\sum_{k=0}^{N} \frac{(x-2\,\Psi(k+1)-2\gamma)\,\mathbf{e}^{(kx)}}{\Gamma(k+1)^2}\right) = 0 \qquad (7.14)$$

Replacing the summation index $k$ by $k-1$ in the second summand, we get

$$\mathbf{e}^x\left(\sum_{k=0}^{N} \frac{(x-2\,\Psi(k+1)-2\gamma)\,\mathbf{e}^{(kx)}}{\Gamma(k+1)^2}\right)$$

and define that:

$$\sum_{k=0}^{N} \frac{(-x+2\,\Psi(k+1)+2\gamma)\,\mathbf{e}^{((k+1)x)}}{\Gamma(k+1)^2} = \sum_{k=1}^{N+1}\left(-\frac{\mathbf{e}^{(kx)}\,(x-2\,\Psi(k)-2\gamma)}{\Gamma(k)^2}\right)$$

Consequently, equality (7.14) takes the following form:

$$\lim_{N\to\infty}\left(\sum_{k=0}^{N}\left(-\frac{\mathbf{e}^{(k\,x)}(x-2\,\Psi(k)-2\,\gamma)}{\Gamma(k)^2}+\frac{\mathbf{e}^{(k\,x)}\left(x-2\,\Psi(k+1)+\frac{2}{k}-2\,\gamma\right)}{\Gamma(k)^2}\right)\right)$$
$$-\frac{\mathbf{e}^{((N+1)\,x)}(x-2\,\Psi(N+1)-2\,\gamma)}{\Gamma(N+1)^2}-\left(\lim_{k\to 0}-\frac{\mathbf{e}^{(k\,x)}(x-2\,\Psi(k)-2\,\gamma)}{\Gamma(k)^2}\right)=0$$

Upon transformations, with consideration of the following formula:

$$\Psi(k)\,k-k\,\Psi(k+1)+1=0\,,\quad k=1,2,3\ \ldots.$$

we get:

$$-\left(\lim_{N\to\infty}\frac{\mathbf{e}^{((N+1)\,x)}(x-2\,\Psi(N+1)-2\,\gamma)}{\Gamma(N+1)^2}\right)=0$$

Since the left-hand expression has zero limit, the resulting identity proves that (10.13) is the solution of ODE (10.7).
Thus, the general solution of ODE (7.7) is as follows:

$$\mathrm{y}(x)=C_1\left(-\mathrm{BesselI}\left(0,2\,(\mathbf{e}^x)^{\left(\frac{1}{2}\right)}\right)\right)+C_2\left(\sum_{k=0}^{\infty}\frac{(x-2\,\Psi(k+1)-2\,\gamma)\,\mathbf{e}^{(k\,x)}}{\Gamma(k+1)^2}\right)$$

**The task is solved.**

**Example No. 3. Find a general solution for the following linear ODE**

(7.15)

**Solution:** in this example $\grave{a}=\sin(x)$. Let's calculate the $\xi_k(-1)$, $\xi_k(-2)$ $k=0,1,2\ \ldots\ldots$ summands.
**a) Zero summands** $\xi_0(-1)$, $\xi_0(-2)$**:**

$$\xi_0(-1)=-\int\!\!\int\sin(x)\,dx\,dx=\sin(x)\qquad\qquad \xi_0(-2)=-2\int\!\!\int\!\!\int\sin(x)\,dx\,dx\,dx=-2\cos(x)$$

**b) first summands** $\xi_1(-1)$, $\xi_1(-2)$**:**

$$\xi_1(-1)=-\int\!\!\int G_1(\grave{a})\,dx\,dx+2\int P_1(a)\,dx+\xi_0(-1)\int\!\!\int G_0(\grave{a})\,dx\,dx-\xi_0(-2)\int G_0(\grave{a})\,dx=$$
$$-\frac{9\sin(x)^2}{4}+\frac{x^2}{4}-2\cos(x)^2$$

$$\xi_1(-2)=-2\int\!\!\int\!\!\int G_1(\grave{a})\,dx\,dx\,dx-2\int\!\!\int P_1(a)\,dx\,dx$$
$$+2\,\xi_0(-1)\left(\int\!\!\int\!\!\int G_0(\grave{a})\,dx\,dx\,dx+\int\!\!\int P_0(a)\,dx\,dx\right)-\xi_0(-2)\int\!\!\int G_0(\grave{a})\,dx\,dx=$$
$$\frac{3\sin(x)\cos(x)}{4}-\frac{3\,x}{4}+\frac{x^3}{6}$$

**c) second summands** $\xi_2(-1)$, $\xi_2(-2)$:

$$\xi_2(-1) = -\int\int G_2(\grave{a})\,dx\,dx - 2\int P_2(a)\,dx + \left(\sum_{s=0}^{1} \xi_{1-s}(-1)\left(\int\int G_s(\grave{a})\,dx\,dx + 2\int P_s(a)\,dx\right)\right) - \left(\sum_{s=0}^{1} \xi_{1-s}(-2)\int G_s(\grave{a})\,dx\right) =$$

$$\frac{34\sin(x)^3}{9} - \frac{\sin(x)}{12} - \frac{x^2\sin(x)}{4} - \cos(x)\,x + \frac{15\cos(x)^2\sin(x)}{4}$$

$$\xi_2(-2) = -2\left(\int\int\int G_2(\grave{a})\,dx\,dx\,dx + \int\int P_2(a)\,dx\,dx\right) + 2\left(\sum_{s=0}^{1} \xi_{1-s}(-1)\left(\int\int\int G_s(\grave{a})\,dx\,dx\,dx + \int\int P_s(a)\,dx\,dx\right)\right) - \left(\sum_{s=0}^{1} \xi_{1-s}(-2)\int\int G_s(\grave{a})\,dx\,dx\right) =$$

$$-\frac{443\cos(x)\sin(x)^2}{108} - \frac{28\cos(x)}{27} - \frac{1\cos(x)\,x^2}{2} + \frac{3\sin(x)\,x}{4} - \frac{1\sin(x)\,x^3}{6} - 4\cos(x)^3$$

Quite similarly, we can obtain the following summands:

$$\xi_3(-1) = -\frac{43\sin(x)^2}{192} - \frac{50\cos(x)^2}{27} + \frac{5x^2}{64} + \frac{7\sin(x)^2x^2}{8} + \frac{3\cos(x)\,x\sin(x)}{8} - \frac{11915\sin(x)^4}{1728} - \frac{2353\cos(x)^2\sin(x)^2}{216} + \frac{13x^2\cos(x)^2}{16} - \frac{x^4}{96} - 4\cos(x)^4;$$

$$\xi_4(-1) = \frac{99871\sin(x)}{103680} + \frac{80425\cos(x)^2\sin(x)^3}{3456} + \frac{x^3\cos(x)}{9} + \frac{151\sin(x)\,x^2}{576} + \frac{36733\cos(x)^2\sin(x)}{20736} - \frac{737\cos(x)\,x}{432} + \frac{7759\sin(x)^3}{11520} - \frac{113x^2\cos(x)^2\sin(x)}{72} - \frac{25x\cos(x)^3}{8} - \frac{343\sin(x)^2x\cos(x)}{108} + \frac{1060231\sin(x)^5}{86400} + \frac{x^4\sin(x)}{96} - \frac{x^3\cos(x)\sin(x)^2}{36} - \frac{x^3\cos(x)^3}{36} - \frac{227\sin(x)^3x^2}{144} + 11\sin(x)\cos(x)^4;$$ and so on.

$$\xi_3(-2) = -\frac{2927x}{6912} + \frac{3\cos(x)\sin(x)\,x^2}{16} - \frac{18977\sin(x)\cos(x)}{6912} + \frac{7x^3\cos(x)^2}{24} + \frac{x^3\sin(x)^2}{3} + \frac{35x^3}{96} + \frac{17305\cos(x)\sin(x)^3}{3456} - \frac{21\cos(x)^2x}{16} - \frac{1x^5}{120} - \frac{3\sin(x)^2x}{2} + 5\cos(x)^3\sin(x);$$

$$\xi_4(-2) = -\frac{13\cos(x)\sin(x)^2x^2}{48} - \frac{67x^3\cos(x)^2\sin(x)}{108} - \frac{1049893\cos(x)}{243000} + \frac{67\cos(x)^2x\sin(x)}{24} - \frac{25918963\cos(x)\sin(x)^2}{7776000} - \frac{733\cos(x)\,x^2}{288} + \frac{1343\sin(x)\,x}{6912} - \frac{271\sin(x)\,x^3}{864} - \frac{11255\cos(x)^3}{1944} - \frac{4667137\cos(x)\sin(x)^4}{432000} - \frac{53x^2\cos(x)^3}{216} + \frac{x^4\cos(x)}{16} + \frac{45\sin(x)^3x}{16} + \frac{1\sin(x)\,x^5}{120} - \frac{8123\cos(x)^3\sin(x)^2}{432} - \frac{5\sin(x)^3x^3}{8} - 8\cos(x)^5;$$ and so on.

In this case partial solutions are determined by the following formulas:

$$y_1(x) = -1 + \left( \sum_{k=0}^{N} \xi_k(-1) \right) \tag{7.16}$$

$$y_2(x) = \left( \sum_{k=0}^{N} \xi_k(-2) \right) - x \left( \sum_{k=0}^{N} \xi_k(-1) \right) + x \tag{7.17}$$

where **N** is the number of summands required to obtain the desired solution with the specified accuracy, within a defined range of independent argument ***x*** .

Assuming, in particular $N = 4$ , we get:

$$y_1(x) = \left( -\frac{1}{96} + \frac{\sin(x)}{96} \right) x^4 + \frac{x^3 \cos(x)}{12} + \left( -\frac{25 \sin(x)}{16} - \frac{\cos(2x)}{32} + \frac{\sin(3x)}{576} + \frac{75}{64} \right) x^2 + \left( \frac{11 \cos(3x)}{864} - \frac{187 \cos(x)}{32} + \frac{3 \sin(2x)}{16} \right) x - \frac{132853}{13824} + \frac{256465 \sin(x)}{13824} + \frac{437 \cos(2x)}{576} - \frac{4169 \sin(3x)}{82944} + \frac{\sin(5x)}{230400} - \frac{\cos(4x)}{4608}$$

$$y_2(x) = \left( -\frac{\sin(x)}{480} + \frac{1}{480} \right) x^5 - \frac{x^4 \cos(x)}{48} + \left( \frac{\cos(2x)}{96} + \frac{11 \sin(x)}{24} - \frac{\sin(3x)}{1728} - \frac{21}{64} \right) x^3 + \left( -\frac{3 \sin(2x)}{32} - \frac{11 \cos(3x)}{1728} + \frac{163 \cos(x)}{64} \right) x^2 + \left( \frac{\cos(4x)}{4608} + \frac{10799}{1536} - \frac{383 \cos(2x)}{576} - \frac{\sin(5x)}{230400} - \frac{68201 \sin(x)}{4608} + \frac{3737 \sin(3x)}{82944} \right) x - \frac{25 \sin(4x)}{27648} - \frac{174593 \cos(x)}{6912} + \frac{2599 \sin(2x)}{1728} - \frac{137 \cos(5x)}{6912000} + \frac{55853 \cos(3x)}{497664}$$

Analysis of the obtained solutions shows that they can be represented by the following formulas:

$$y_1(x) = \sum_{i=0}^{k} \left( \sum_{j=0}^{k-i} x^{(2i)} \left( x \sin((2j+1)x)\, b_{i,j} + c_{i,j}\, x \cos(2xj) + e_{i,j} \cos((2j+1)x) + g_{i,j} \sin(2xj) \right) \right) \tag{7.18}$$

$$y_2(x) = \sum_{i=0}^{k} \left( \sum_{j=0}^{k-i} x^{(2i)} \left( x \sin((2j+1)x)\, \rho_{i,j} + x\, \sigma_{i,j} \cos(2xj) + \tau_{i,j} \cos((2j+1)x) + \upsilon_{i,j} \sin(2xj) \right) \right) \tag{7.19}$$

where $k$ is a natural number.

$b_{i,j}, c_{i,s}, e_{i,j}, g_{i,s}, \rho_{i,j}, \sigma_{i,s}, \tau_{i,j}, \upsilon_{i,s}$ are numeric functions, for which we need to define recurrence formulas. For the definition of these functions, let's substitute (7.18), (7.19) to the initial ODE (7.15).
Indeed, in particular with $y_1(x)$ in ODE (7.15) we get:

$$\frac{\partial^2}{\partial x^2}\left(\sum_{i=0}^{k}\left(\sum_{j=0}^{k-i} x^{(2i)}\left(x\sin((2j+1)x)\,b_{i,j}+c_{i,j}\,x\cos(2xj)+e_{i,j}\cos((2j+1)x)+g_{i,j}\sin(2xj)\right)\right)\right)=\sin(x)\left(\sum_{i=0}^{k}\left(\sum_{j=0}^{k-i} x^{(2i)}\left(x\sin((2j+1)x)\,b_{i,j}+c_{i,j}\,x\cos(2xj)+e_{i,j}\cos((2j+1)x)+g_{i,j}\sin(2xj)\right)\right)\right)$$

Expanding this equality and making the necessary substitutions of the relevant summation indexes, we eventually obtain the following formulas for the required numeric coefficients :

$$-4(2i+3)(i+1)\,b_{i+1,j}+2(2j+1)^2\,b_{i,j}+8(2j+1)(i+1)\,e_{i+1,j}+c_{i,j}-c_{i,j+1}=0 \quad (7.20)$$

$$4(2j+1)(2i+1)\,b_{i,j}+4(i+1)(2i+1)\,e_{i+1,j}-2(2j+1)^2\,e_{i,j}+g_{i,j}-g_{i,j+1}=0$$
(7.21)

$$-4(i+1)\,g_{i+1,j}+8(i-j)(i+j)\,g_{i,j}-8j(2i+1)\,c_{i,j}-e_{i,j-1}+e_{i,j}=0$$
(7.22)

$$4(2i+3)(i+1)\,c_{i+1,j}+16j(i+1)\,g_{i+1,j}-8\,c_{i,j}\,j^2-b_{i,j}+b_{i,j-1}=0$$
(7.23)

Substituting the second solution, (7.19), to the initial ODE (7.15) we have:

$$\frac{\partial^2}{\partial x^2}\left(\sum_{i=0}^{k}\left(\sum_{j=0}^{k-i} x^{(2i)}\left(x\sin((2j+1)x)\,\rho_{i,j}+x\,\sigma_{i,j}\cos(2xj)+\tau_{i,j}\cos((2j+1)x)+\upsilon_{i,j}\sin(2xj)\right)\right)\right)$$
$$=\sin(x)\left(\sum_{i=0}^{k}\left(\sum_{j=0}^{k-i} x^{(2i)}\left(x\sin((2j+1)x)\,\rho_{i,j}+x\,\sigma_{i,j}\cos(2xj)+\tau_{i,j}\cos((2j+1)x)+\upsilon_{i,j}\sin(2xj)\right)\right)\right)$$

From here, with the use of similar transformations, we obtain the following recurrence equations for the numerical coefficients .

$$4(i+1)(2i+1)\,\upsilon_{i+1,j}-8j(2i+1)\,\sigma_{i,j}+\tau_{i,j}-8j^2\,\upsilon_{i,j}-\tau_{i,j-1}=0$$
(7.24)

$$4(i+1)(2i+1)\,\tau_{i+1,j}-2(2j+1)^2\,\tau_{i,j}+4(2j+1)(2i+1)\,\rho_{i,j}+\upsilon_{i,j}-\upsilon_{i,j+1}=0$$
(7.25)

$$-8(2j+1)(i+1)\,\tau_{i+1,j}+4(2i+3)(i+1)\,\rho_{i+1,j}-2(2j+1)^2\,\rho_{i,j}-\sigma_{i,j}+\sigma_{i,j+1}=0$$
(7.26)

$$16j(i+1)\,\upsilon_{i+1,j}+4(2i+3)(i+1)\,\sigma_{i+1,j}-8\,\sigma_{i,j}\,j^2+\rho_{i,j-1}-\rho_{i,j}=0$$
(7.27)

Thus, the general solution of ODE (7.15) is as follows:

$$y(x) = C_1 \left( \sum_{i=0}^{k} \left( \sum_{j=0}^{k-i} x^{(2i)} \left( x \sin((2j+1)x)\, b_{i,j} + c_{i,j}\, x \cos(2xj) + e_{i,j} \cos((2j+1)x) + g_{i,j} \sin(2xj) \right) \right) \right) + C_2 \left( \sum_{i=0}^{k} \left( \sum_{j=0}^{k-i} x^{(2i)} \left( x \sin((2j+1)x)\, \rho_{i,j} + x\, \sigma_{i,j} \cos(2xj) + \tau_{i,j} \cos((2j+1)x) + \upsilon_{i,j} \sin(2xj) \right) \right) \right)$$

where numerical coefficients are determined by (7.20) - (7.27).

**The task is solved.**

**Example No. 4. Find a general solution for the following linear ODE:**

(7.28)

**Solution:** In this example $\grave{a} = \ln(x)$. Let's calculate the $\xi_k(-1)$, $\xi_k(-2)$ $k = 0, 1, 2$ ......summands. Then we have:

$$\xi_0(-1) = -\int\int \ln(x)\, dx\, dx =$$

$$\xi_0(-2) = -2\int\int\int \ln(x)\, dx\, dx\, dx =$$

$$\xi_1(-1) = -\int\int G_1(\grave{a})\, dx\, dx + 2\int P_1(a)\, dx + \xi_0(-1)\int\int G_0(\grave{a})\, dx\, dx - \xi_0(-2)\int G_0(\grave{a})\, dx = -\frac{x^4 \ln(x)^2}{24} + \frac{x^4 \ln(x)}{9} - \frac{25\, x^4}{432}$$

$$\xi_1(-2) = -2\int\int\int G_1(\grave{a})\, dx\, dx\, dx - 2\int\int P_1(a)\, dx\, dx + 2\,\xi_0(-1)\left(\int\int\int G_0(\grave{a})\, dx\, dx\, dx + \int\int P_0(a)\, dx\, dx\right) - \xi_0(-2)\int\int G_0(\grave{a})\, dx\, dx = -\frac{x^5 \ln(x)^2}{30} + \frac{29\, x^5 \ln(x)}{300} - \frac{2819\, x^5}{54000}$$

$$\xi_2(-1) = -\int\int G_2(\grave{a})\, dx\, dx - 2\int P_2(a)\, dx + \left(\sum_{s=0}^{1} \xi_{1-s}(-1)\left(\int\int G_s(\grave{a})\, dx\, dx + 2\int P_s(a)\, dx\right)\right) - \left(\sum_{s=0}^{1} \xi_{1-s}(-2)\int G_s(\grave{a})\, dx\right) = -\frac{x^6 \ln(x)^3}{720} + \frac{113\, x^6 \ln(x)^2}{21600} - \frac{889\, x^6 \ln(x)}{162000} + \frac{2021\, x^6}{1215000}$$

$$\xi_2(-2) = -2\left(\iiint G_2(\grave{a})\,dx\,dx\,dx + \iint P_2(a)\,dx\,dx\right) + 2\left(\sum_{s=0}^{1} \xi_{1-s}(-1)\left(\iiint G_s(\grave{a})\,dx\,dx\,dx + \iint P_s(a)\,dx\,dx\right)\right) - \left(\sum_{s=0}^{1} \xi_{1-s}(-2)\iint G_s(\grave{a})\,dx\,dx\right) =$$

$$-\frac{x^7\ln(x)^3}{840} + \frac{2489\,x^7\ln(x)^2}{529200} - \frac{7489\,x^7\ln(x)}{1481760} + \frac{12091561\,x^7}{7779240000}$$

Quite similarly, we obtain:

$$\xi_3(-1) = -\frac{x^8\,1\ln(x)^4}{40320} + \frac{127\,x^8\ln(x)^3}{1058400} - \frac{84061\,x^8\ln(x)^2}{444528000} + \frac{22059049\,x^8\ln(x)}{186701760000} - \frac{17351633\,x^8}{697019904000};$$

$$\xi_4(-1) = -\frac{x^{10}\ln(x)^5}{3628800} + \frac{3713\,x^{10}\ln(x)^4}{2286144000} - \frac{153541\,x^{10}\ln(x)^3}{45008460000} + \frac{492518029\,x^{10}\ln(x)^2}{151228425600000} - \frac{542779529651\,x^{10}\ln(x)}{381095632512000000} + \frac{7830520197209\,x^{10}}{34298606926080000000};$$ and so on.

$$\xi_3(-2) = -\frac{x^9\ln(x)^4}{45360} + \frac{629\,x^9\ln(x)^3}{5715360} - \frac{2118401\,x^9\ln(x)^2}{12002256000} + \frac{120458131\,x^9\ln(x)}{1080203040000} - \frac{18046415923\,x^9}{762191265024000};$$

$$\xi_4(-2) = \frac{41668782325106407\,x^{11}}{189881482451791284000000} - \frac{x^{11}\ln(x)^5}{3991680} + \frac{69611\,x^{11}\ln(x)^4}{46103904000} - \frac{220052737\,x^{11}\ln(x)^3}{68464297440000} + \frac{137096461039\,x^{11}\ln(x)^2}{44282707584192000} - \frac{1393786594886281\,x^{11}\ln(x)}{1022930545194835200000};$$ an so on.

In this case (in practice), partial solutions are determined by the following formulas:

$$y_1(x) = -1 + \left(\sum_{k=0}^{N} \xi_k(-1)\right) \tag{7.29}$$

$$y_2(x) = \left(\sum_{k=0}^{N} \xi_k(-2)\right) - x\left(\sum_{k=0}^{N} \xi_k(-1)\right) + x \tag{7.30}$$

where **N** is the number of summands required to obtain the desired solution with the specified accuracy, within a defined range of independent argument ***x*** .

Assuming, in particular, $N = 4$ in (7.29), (7.30) we get:

$$y_1(x) = -\frac{x^{10}\ln(x)^5}{3628800} + \left(-\frac{x^8}{40320} + \frac{3713\,x^{10}}{2286144000}\right)\ln(x)^4$$
$$+\left(-\frac{x^6}{720} - \frac{153541\,x^{10}}{45008460000} + \frac{127\,x^8}{1058400}\right)\ln(x)^3$$
$$+\left(\frac{492518029\,x^{10}}{151228425600000} - \frac{x^4}{24} - \frac{84061\,x^8}{444528000} + \frac{113\,x^6}{21600}\right)\ln(x)^2$$
$$+\left(-\frac{542779529651\,x^{10}}{381095632512000000} + \frac{x^4}{9} - \frac{889\,x^6}{162000} - \frac{x^2}{2} + \frac{22059049\,x^8}{186701760000}\right)\ln(x) - 1$$
$$+\frac{7830520197209\,x^{10}}{342986069260800000000} + \frac{3\,x^2}{4} - \frac{25\,x^4}{432} - \frac{17351633\,x^8}{697019904000} + \frac{2021\,x^6}{1215000}$$

$$y_2(x) = \frac{x^{11}\ln(x)^5}{39916800} + \left(\frac{x^9}{362880} - \frac{31607\,x^{11}}{276623424000}\right)\ln(x)^4$$
$$+\left(\frac{x^7}{5040} + \frac{31511803\,x^{11}}{159750027360000} - \frac{71\,x^9}{7144200}\right)\ln(x)^3$$
$$+\left(-\frac{559\,x^7}{1058400} + \frac{75623\,x^9}{6001128000} - \frac{356133410639\,x^{11}}{221413537920960000} + \frac{x^5}{120}\right)\ln(x)^2 +$$
$$\left(\frac{x^3}{6} - \frac{13\,x^5}{900} + \frac{48179\,x^7}{111132000} + \frac{541141476663763\,x^{11}}{876797610167001600000} - \frac{20073827\,x^9}{3024568512000}\right)\ln(x)$$
$$+\frac{17\,x^5}{3000} - \frac{508937\,x^7}{4667544000} - \frac{2847798132559381\,x^{11}}{3214924570612339200000} + x + \frac{10601083\,x^9}{8710757314560}$$
$$-\frac{5\,x^3}{36}$$

Analysis of the obtained formulas demonstrates that the desired analytical solutions can be represented as follows:

$$y_1(x) = \sum_{k=0}^{\infty} \ln(x)^k \left(\sum_{i=k}^{\infty} b_{i,k}\, x^{(2i)}\right) \qquad (7.31)$$

$$y_2(x) = \sum_{k=0}^{\infty} \ln(x)^k \left(\sum_{i=k}^{\infty} c_{i,k}\, x^{(2i+1)}\right) \qquad (7.32)$$

where $b_{i,k}$, $c_{i,k}$ are the desired numeric functions.

Substituting (7.31), (7.32) to the initial ODE (7.28), we obtain recurrence formulas for $b_{i,k}$ and $c_{i,k}$.

Indeed, with $y_1(x)$ in ODE (7.28) we get:

$$\frac{\partial^2}{\partial x^2}\left(\sum_{k=0}^{\infty} \ln(x)^k \left(\sum_{i=k}^{\infty} b_{i,k}\, x^{(2i)}\right)\right) = \mathbf{e}^x \left(\sum_{k=0}^{\infty} \ln(x)^k \left(\sum_{i=k}^{\infty} b_{i,k}\, x^{(2i)}\right)\right)$$

or (10.33)

$$\left(\sum_{k=0}^{\infty}\left(\sum_{i=k}^{\infty}\left(\frac{b_{i,k}\ln(x)^k k^2 x^{(2i)}}{x^2\ln(x)^2}-\frac{b_{i,k}\ln(x)^k k x^{(2i)}}{x^2\ln(x)}-\frac{b_{i,k}\ln(x)^k k x^{(2i)}}{x^2\ln(x)^2}\right.\right.\right.$$
$$\left.\left.\left.+\frac{4 b_{i,k}\ln(x)^k k x^{(2i)} i}{x^2\ln(x)}+\frac{4 b_{i,k}\ln(x)^k x^{(2i)} i^2}{x^2}-\frac{2 b_{i,k}\ln(x)^k x^{(2i)} i}{x^2}\right)\right)\right)$$
$$-\ln(x)\left(\sum_{k=0}^{\infty}\left(\sum_{i=k}^{\infty} b_{i,k}\ln(x)^k x^{(2i)}\right)\right)=0$$

Let's replace $k$ index by $k+1$ and $i \rightarrow i+1$ in the following summand:

$$\ln(x)\left(\sum_{k=0}^{\infty}\left(\sum_{i=k}^{\infty} b_{i,k}\ln(x)^k x^{(2i)}\right)\right)$$

We obtain the following equality:

$$\ln(x)\left(\sum_{k=0}^{\infty}\left(\sum_{i=k}^{\infty} b_{i,k}\ln(x)^k x^{(2i)}\right)\right)=\sum_{k=1}^{\infty}\left(\sum_{i=k}^{\infty}\ln(x)^k x^{(2(i-1))} b_{i-1,-1+k}\right)$$

Similarly, for other summands we have:

$$\sum_{k=0}^{\infty}\left(\sum_{i=k}^{\infty}\frac{b_{i,k}\ln(x)^k k^2 x^{(2i)}}{x^2\ln(x)^2}\right)=\sum_{s=-2}^{\infty}\left(\sum_{i=2+s}^{\infty}\ln(x)^s(2+s)^2 b_{i,2+s} x^{(-2+2i)}\right)$$

$$\sum_{k=0}^{\infty}\left(\sum_{i=k}^{\infty}\left(-\frac{b_{i,k}\ln(x)^k k x^{(2i)}}{x^2\ln(x)}\right)\right)=\sum_{s=-1}^{\infty}\left(\sum_{i=1+s}^{\infty}(-1-s) b_{i,1+s}\ln(x)^s x^{(-2+2i)}\right)$$

$$\sum_{k=0}^{\infty}\left(\sum_{i=k}^{\infty}\left(-\frac{b_{i,k}\ln(x)^k k x^{(2i)}}{x^2\ln(x)^2}\right)\right)=\sum_{s=-2}^{\infty}\left(\sum_{i=2+s}^{\infty}(-\ln(x)^s(2+s) b_{i,2+s} x^{(-2+2i)})\right)$$

$$\sum_{k=0}^{\infty}\left(\sum_{i=k}^{\infty}\frac{4 b_{i,k}\ln(x)^k k x^{(2i)} i}{x^2\ln(x)}\right)=\sum_{s=-1}^{\infty}\left(\sum_{i=1+s}^{\infty}(4+4s) b_{i,1+s}\ln(x)^s x^{(-2+2i)} i\right)$$

Substituting these to (7.33), upon transformations we have:

$$\sum_{k=0}^{\infty}\left(\sum_{i=k}^{\infty}\frac{\ln(x)^k x^{(2i)}((k+3)(k+2) b_{i,k+2}+(4i-1)(k+1) b_{i,k+1}+2i(2i-1) b_{i,k}-b_{i}}{x^2}\right.$$
$$\left.\vphantom{\sum}\right)=0$$

Thus we obtain recurrence equation to find $\mathrm{B}(k)$ .

$$(k+3)(k+2) b_{i,k+2}+(4i-1)(k+1) b_{i,k+1}+2i(2i-1) b_{i,k}-b_{i-1,k-1}=0 \qquad (7.34)$$

Thus, partial solution (7.31) takes the following form:

$$y_1(x) = \sum_{k=0}^{\infty} \ln(x)^k \left( \sum_{i=k}^{\infty} b7_{i,k} x^{(2i)} \right) \tag{7.35}$$

We can check that this solution satisifies the initial ODE (7.28).
Proceeding in a similar way, we can obtain recurrence equation to find $c_{i,k}$ coefficients.

Indeed,

$$\frac{\partial^2}{\partial x^2}\left( \sum_{k=0}^{\infty} \left( \sum_{i=k}^{\infty} c_{i,k} \ln(x)^k x^{(2i+1)} \right) \right) = \ln(x) \left( \sum_{k=0}^{\infty} \left( \sum_{i=k}^{\infty} c_{i,k} \ln(x)^k x^{(2i+1)} \right) \right)$$

From there follows (10.36)

$$\left( \sum_{k=0}^{\infty} \left( \sum_{i=k}^{\infty} \left( \frac{c_{i,k} \ln(x)^k k^2 x^{(2i+1)}}{x^2 \ln(x)^2} - \frac{c_{i,k} \ln(x)^k k x^{(2i+1)}}{x^2 \ln(x)} - \frac{c_{i,k} \ln(x)^k k x^{(2i+1)}}{x^2 \ln(x)^2} \right.\right.\right.$$
$$+ \frac{2 c_{i,k} \ln(x)^k k x^{(2i+1)} (2i+1)}{x^2 \ln(x)} + \frac{c_{i,k} \ln(x)^k x^{(2i+1)} (2i+1)^2}{x^2}$$
$$\left.\left.\left. - \frac{c_{i,k} \ln(x)^k x^{(2i+1)} (2i+1)}{x^2} \right) \right) \right) - \ln(x) \left( \sum_{k=0}^{\infty} \left( \sum_{i=k}^{\infty} c_{i,k} \ln(x)^k x^{(2i+1)} \right) \right) = 0$$

When, upon the replacement of summation index, we obtain the following equalities:

$$\sum_{k=0}^{\infty} \left( \sum_{i=k}^{\infty} c_{i,k} \ln(x) \ln(x)^k x^{(2i+1)} \right) = \sum_{k=1}^{\infty} \left( \sum_{i=k}^{\infty} \ln(x)^k x^{(-1+2s)} c_{-1+i,k-1} \right)$$

$$\sum_{k=0}^{\infty} \left( \sum_{i=k}^{\infty} \frac{c_{i,k} \ln(x)^k k^2 x^{(2i+1)}}{x^2 \ln(x)^2} \right) = \sum_{k=-2}^{\infty} \left( \sum_{i=k+2}^{\infty} \ln(x)^k x^{(2i-1)} c_{i,k+2} (k+2)^2 \right)$$

$$\sum_{k=0}^{\infty} \left( \sum_{i=k}^{\infty} \left( - \frac{c_{i,k} \ln(x)^k k x^{(2i+1)}}{x^2 \ln(x)} \right) \right) = \sum_{k=-1}^{\infty} \left( \sum_{i=k+1}^{\infty} (-\ln(x)^k x^{(2i-1)} (k+1) c_{i,k+1}) \right)$$

$$\sum_{k=0}^{\infty} \left( \sum_{i=k}^{\infty} \left( - \frac{c_{i,k} \ln(x)^k k x^{(2i+1)}}{x^2 \ln(x)^2} \right) \right) = \sum_{k=-2}^{\infty} \left( \sum_{i=k+2}^{\infty} (-\ln(x)^k x^{(2i-1)} c_{i,k+2} (k+2)) \right)$$

$$\sum_{k=0}^{\infty} \left( \sum_{i=k}^{\infty} \frac{2 c_{i,k} \ln(x)^k k x^{(2i+1)} (2i+1)}{x^2 \ln(x)} \right) =$$
$$\sum_{k=-1}^{\infty} \left( \sum_{i=k+1}^{\infty} 2 \ln(x)^k x^{(2i-1)} (2i+1) (k+1) c_{i,k+1} \right)$$

Thus, substituting the obtained values to (7.36), we get:

$$\sum_{k=0}^{\infty} \left( \sum_{i=k}^{\infty} \ln(x)^k x^{(2i-1)} \right.$$
$$\left. ((k+2)(k+1) c_{i,k+2} + (k+1)(1+4i) c_{i,k+1} + 2i(2i+1) c_{i,k} - c_{-1+i,k-1}) \right) = 0$$

From here we receive the desired recurrence equation:

$$(k+2)(k+1)c_{i,k+2} + (k+1)(1+4i)c_{i,k+1} + 2i(2i+1)c_{i,k} - c_{-1+i,k-1} = 0 \qquad (7.36)$$

Thus, the second partial solution of the initial ODE (7.28) is as follows:

$$y_2(x) = \sum_{k=0}^{\infty} \ln(x)^k \left( \sum_{i=k}^{\infty} c_{i,k} x^{(2i+1)} \right)$$

Consequently, the general solution of ODE (7.28) is defined by the following formula:

$$y(x) = C_1 \left( \sum_{k=0}^{\infty} \ln(x)^k \left( \sum_{i=k}^{\infty} b_{i,k} x^{(2i)} \right) \right) + C_2 \left( \sum_{k=0}^{\infty} \ln(x)^k \left( \sum_{i=k}^{\infty} c_{i,k} x^{(2i+1)} \right) \right)$$

where $b_{i,k}$ is defined by the recurrence expression (7.34), and $c_{i,k}$ is defined by the recurrence expression (7.36). **The task is solved.**

Further, for the calculation of $\xi_k(-1)$, $\xi_k(-2)$ summands we shall use a **PROGRAM for the calculation of $\xi_k(-1)$, $\xi_k(-2)$ and partial solutions of $y_1(x)$ and $y_2(x)$, up to $k$, of the linear ODE $\frac{\partial^2}{\partial x^2} y = a\, y$**

```
#############################################################################

> restart:alias(y=y(x),v=v(x),a=a(x)):k:=7:a:=(x^3-1):

> a:=(x^3-1):# type here the functional coefficient

k:=7: # type here the number of coefficients, which should be calculated

> dn:= proc (n::integer, w,c) local k, Ds; global resd, x: option remember;
if n = 0 then resd := w elif 0 < n then resd := diff(w,`$`(x,n)) else
Ds := w:for k to -n do Ds:=int(Ds,x) end do:Ds:=Ds: resd := Ds end if:
RETURN(resd) end proc:

> for k from 0 to k do

> R:=z->a*int(int(z,x),x):h[0]:=a:for l from 0 to k do h[l+1]:=R(h[l]) od:
l:='l':
```

> **for** $l$ **from** 0 **to** $k+1$ **do** $f_l := \int h_l\, dx$ **end do**; $l := 'l'$

**for** $m$ **from** 0 **to** $k$ **do for** $j$ **from** 0 **to** $k - m$ **do** $q_0 := f_m$; $q_{j+1} := \mathrm{R}(q_j)$ **end do end do**;

> $G_{k+1} := \mathrm{add}(q_i, i = 1..k+1)$;
> $G_0 := 0$;
> $G_{-1} := 0$;

$\xi_0(-2) := -2 \int\int\int a\, dx\, dx\, dx$

$$\xi_k(-1) := \operatorname{expand}\left(-\int f_k\,dx - 2\int G_k\,dx + \operatorname{add}\left(\xi_{k-s-1}(-1)\int f_s\,dx, s = 0\,..\,k-1\right) + 2\operatorname{add}\left(\xi_{k-s-1}(-1)\int G_s\,dx, s = 0\,..\,k-1\right) - \operatorname{add}\left(\xi_{k-s-1}(-2)\int h_s\,dx, s = 0\,..\,k-1\right)\right)$$

$$\xi_k(-2) := \operatorname{expand}\left(-2\int\int f_k\,dx\,dx - 2\int\int G_k\,dx\,dx + 2\operatorname{add}\left(\xi_{k-s-1}(-1)\int\int f_s\,dx\,dx, s = 0\,..\,k-1\right) + 2\operatorname{add}\left(\xi_{k-s-1}(-1)\int\int G_s\,dx\,dx, s = 0\,..\,k-1\right) - \operatorname{add}\left(\xi_{k-s-1}(-2)\int\int h_s\,dx\,dx, s = 0\,..\,k-1\right)\right);$$

$s := 's';$

$s := 's'$

```
> a1:=add(xi[s](-1),s=0..k):
> a2:=add(xi[s](-2),s=0..k):
> od:
> Diff(y,`$`(x,2)) = a*y;
> k:=k-1:
> for l from 0 to k do Xi[l](-1):=xi[l](-1) od:for l from 0 to k do Xi[l](-2):=xi[l](-2) od:
> y[1](x)=combine(sort(-1+a1,x));
> y[2](x)=sort(expand(a2-x*a1+x),x);
> #################################################################
```

Let's define the method to solve linear ODEs with the use of this program.

**Example No. 5. Find a general solution for the following linear ODE in analytical form:**

(7.37)

**Solution**: If we enter the following equalities to the program:

, **for all** $k$,

we obtain the desired values of ; from $k = 0$ **to k=4** :

$$\Xi_1(-1) := -\frac{1\,x^{18}}{22032} + \frac{37\,x^{11}}{7920} - \frac{1\,x^4}{24}$$

$$\Xi_2(-1) := -\frac{1\,x^{27}}{15466464} + \frac{1429\,x^{20}}{115117200} - \frac{367\,x^{13}}{1235520} + \frac{1\,x^6}{720}$$

$$\Xi_3(-1) := -\frac{1\,x^{36}}{19487744640} + \frac{1008383\,x^{29}}{65619566812800} - \frac{1852423\,x^{22}}{2765575612800} + \frac{2083\,x^{15}}{259459200} - \frac{1\,x^{8}}{40320}$$

$$\Xi_4(-1) := \frac{1138219\,x^{38}}{103793749806146400} + \frac{1\,x^{10}}{3628800} - \frac{4259\,x^{17}}{35286451200} - \frac{1\,x^{45}}{38585734387200} - \frac{61821827\,x^{31}}{83911021061868000} + \frac{300689\,x^{24}}{19082471728320}$$

$$\Xi_0(-2) := -\frac{1\,x^{10}}{360} + \frac{1\,x^{3}}{3}$$

$$\Xi_1(-2) := -\frac{1\,x^{19}}{77520} + \frac{79\,x^{12}}{23760} - \frac{1\,x^{5}}{30}$$

$$\Xi_2(-2) := -\frac{223\,x^{28}}{10285198560} + \frac{2207\,x^{21}}{241746120} - \frac{703\,x^{14}}{2882880} + \frac{1\,x^{7}}{840}$$

$$\Xi_3(-2) := \frac{1811\,x^{16}}{259459200} - \frac{29\,x^{37}}{1522209386880} + \frac{11361673\,x^{30}}{984293502192000} - \frac{169889\,x^{23}}{304345641600} - \frac{1\,x^{9}}{45360}$$

$$\Xi_4(-2) := \frac{845507\,x^{25}}{61161768360000} + \frac{1\,x^{11}}{3991680} - \frac{1289\,x^{46}}{124778547861327360} + \frac{2741467\,x^{39}}{327769736229936000} - \frac{6385785223\,x^{32}}{102930852502558080000} - \frac{11461\,x^{18}}{105859353600}$$

---

$$y_1(x) = -\frac{1\,x^{45}}{38585734387200} + \frac{1138219\,x^{38}}{103793749806146400} - \frac{1\,x^{36}}{19487744640} - \frac{61821827\,x^{31}}{83911021061868000} + \frac{1008383\,x^{29}}{65619566812800} - \frac{1\,x^{27}}{15466464} + \frac{300689\,x^{24}}{19082471728320} - \frac{1852423\,x^{22}}{2765575612800} + \frac{1429\,x^{20}}{115117200} - \frac{1\,x^{18}}{22032} - \frac{4259\,x^{17}}{35286451200} + \frac{2083\,x^{15}}{259459200} - \frac{367\,x^{13}}{1235520} + \frac{37\,x^{11}}{7920} + \frac{1\,x^{10}}{3628800} - \frac{1\,x^{9}}{72} - \frac{1\,x^{8}}{40320} + \frac{1\,x^{6}}{720} - \frac{1\,x^{4}}{24} + \frac{1\,x^{2}}{2} - 1$$

$$y_2(x) = \frac{1\,x^{46}}{64160092483200} - \frac{73327\,x^{39}}{28179298589904000} + \frac{1\,x^{37}}{30995213760} + \frac{211357451\,x^{32}}{1816426808868672000} - \frac{2129\,x^{30}}{556727094000} + \frac{1\,x^{28}}{23269680} - \frac{33907\,x^{25}}{17539036515000} + \frac{26099\,x^{23}}{233853820200} - \frac{467\,x^{21}}{142203600} + \frac{1\,x^{19}}{30780} + \frac{47\,x^{18}}{3780691200} - \frac{17\,x^{16}}{16216200} + \frac{23\,x^{14}}{432432} - \frac{2\,x^{12}}{1485} - \frac{1\,x^{11}}{39916800} + \frac{1\,x^{10}}{90} + \frac{1\,x^{9}}{362880} - \frac{1\,x^{7}}{5040} + \frac{1\,x^{5}}{120} - \frac{1\,x^{3}}{6} + x$$

Reviewing the obtained solutions, we can see that it is very difficult to define their construction patterns. However, it is easy to define the $\Xi_k(-1)$ **and** $\Xi_k(-2)$ pattern. Indeed:

$$\Xi_k(-1) = \sum_{i=0}^{k+1} b_{i,k} x^{(2(k+1)+7i)} \qquad \Xi_k(-2) = \sum_{i=0}^{k+1} c_{i,k} x^{(2k+3+7i)}$$

In this case partial solutions are defined by the following formulas:

$$y_1(x) = -1 + \left( \sum_{k=0}^{N} \left( \sum_{i=0}^{k+1} b_{i,k} x^{(2(k+1)+7i)} \right) \right)$$

$$y_2(x) = \left( \sum_{k=0}^{N} \left( \sum_{i=0}^{k+1} c_{i,k} x^{(2k+3+7i)} \right) \right) - x \left( \sum_{k=0}^{N} \left( \sum_{i=0}^{k+1} b_{i,k} x^{(2(k+1)+7i)} \right) \right) + x$$

where $N$ is a parameter, which defines the number of summands and tends to infinity.

Analyzing (7.38) and (7.39), we eventually obtain the following views:

$$y_1(x) = -1 + \left( \sum_{k=0}^{N} \left( \sum_{i=0}^{k+1} b_{i,k} x^{(2(k+1)+7i)} \right) \right)$$

$$y_2(x) = \left( \sum_{k=0}^{N} \left( \sum_{i=0}^{k+1} r_{i,k} x^{(2k+3+7i)} \right) \right) + x$$

here we assume that

Substituting the obtained values of $y_1(x)$ and $y_2(x)$ functions to the initial equation (7.37), we receive the required recurrence equations to determine the numerical coefficients $r_{i,k}, b_{i,k}$:

$$b_{i-1,k-1} = (2k+2+7i)(2k+1+7i) b_{i,k} + b_{i,k-1} \tag{7.38}$$

$$r_{i-1,k-1} = (7i+2k+3)(2k+2+7i) r_{i,k} + r_{i,k-1} \tag{7.39}$$

Thus, the desired general solution of the initial equation assumes the following form:

$$y = C_1 \left( -1 + \left( \sum_{k=0}^{\infty} \left( \sum_{i=0}^{k+1} b_{i,k} x^{(2(k+1)+7i)} \right) \right) \right) + C_2 \left( \left( \sum_{k=0}^{\infty} \left( \sum_{i=0}^{k+1} r_{i,k} x^{(2k+3+7i)} \right) \right) + x \right)$$

where numerical coefficients $r_{i,k}, b_{i,k}$ are determined by formulas (7.38), (7.39). **The task is solved.**

**Example No. 6. Find a general solution for the following linear ODE:**

(7.40)

**Solution**: If we enter the following equalities to the program:

**, for all** $k$ **,**

we obtain the desired values of ; from $k = 0$ to $k = 6$ :

$$\Xi_0(-1) := -\frac{4\,x^{\left(\frac{5}{2}\right)}\,(-14+3\,x)}{105}$$

$$\Xi_1(-1) := -\frac{2\,x^{7}}{735} + \frac{8\,x^{6}}{315} - \frac{4\,x^{5}}{75}$$

$$\Xi_2(-1) := \frac{16\,x^{\left(\frac{19}{2}\right)}}{41895} - \frac{8\,x^{\left(\frac{21}{2}\right)}}{293265} - \frac{656\,x^{\left(\frac{17}{2}\right)}}{401625} + \frac{32\,x^{\left(\frac{15}{2}\right)}}{14625}$$

$$\Xi_3(-1) := \frac{2848\,x^{11}}{57432375} - \frac{32\,x^{10}}{658125} - \frac{4\,x^{14}}{26687115} + \frac{32\,x^{13}}{11437335} - \frac{32012\,x^{12}}{1762732125}$$

$$\Xi_4(-1) := -\frac{509312\,x^{\left(\frac{27}{2}\right)}}{581502796875} + \frac{7874656\,x^{\left(\frac{29}{2}\right)}}{17942850300375} + \frac{256\,x^{\left(\frac{25}{2}\right)}}{378421875} + \frac{32\,x^{\left(\frac{33}{2}\right)}}{2642024385} - \frac{16\,x^{\left(\frac{35}{2}\right)}}{30823617825} - \frac{15242768\,x^{\left(\frac{31}{2}\right)}}{144207352414125}$$

$$\Xi_5(-1) := -\frac{4\,x^{21}}{3236479871625} - \frac{245817524\,x^{17}}{38128556888296875} + \frac{405952\,x^{16}}{40123692984375} - \frac{256\,x^{15}}{39734296875} + \frac{16\,x^{20}}{462354267375} - \frac{2428744\,x^{19}}{6393192623692875} + \frac{74467936\,x^{18}}{35042386636632375}$$

$$\Xi_6(-1) := \frac{2048\,x^{\left(\frac{35}{2}\right)}}{45893112890625} - \frac{21408256\,x^{\left(\frac{37}{2}\right)}}{259800912073828125} + \frac{387460384\,x^{\left(\frac{39}{2}\right)}}{6074153637915287925} - \frac{167699021552\,x^{\left(\frac{41}{2}\right)}}{6266823657523538484375} + \frac{370906048\,x^{\left(\frac{43}{2}\right)}}{56666922594282225225} - \frac{241182752\,x^{\left(\frac{45}{2}\right)}}{259787382263759975625} + \frac{32\,x^{\left(\frac{47}{2}\right)}}{456343661899125} - \frac{16\,x^{\left(\frac{49}{2}\right)}}{7453613144352375}$$

$$\Xi_0(-2) := -\frac{16\,x^{\left(\frac{7}{2}\right)}\,(-6+x)}{315}$$

$$\Xi_1(-2) := -\frac{1\,x^{8}}{630} + \frac{16\,x^{7}}{945} - \frac{4\,x^{6}}{105}$$

$$\Xi_2(-2) := -\frac{496\,x^{\left(\frac{19}{2}\right)}}{401625} + \frac{724\,x^{\left(\frac{21}{2}\right)}}{2639385} + \frac{992\,x^{\left(\frac{17}{2}\right)}}{580125} - \frac{362\,x^{\left(\frac{23}{2}\right)}}{20235285}$$

$$\Xi_3(-2) := -\frac{1376\,x^{11}}{34459425} - \frac{83\,x^{15}}{789176115} + \frac{332\,x^{14}}{157835223} - \frac{6846506\,x^{13}}{481225870125} + \frac{1376\,x^{12}}{34459425}$$

$$\Xi_4(-2) := \frac{852224\, x^{\left(\frac{27}{2}\right)}}{1486062703125} - \frac{426112\, x^{\left(\frac{29}{2}\right)}}{581502796875} - \frac{149063048\, x^{\left(\frac{33}{2}\right)}}{1755936585277875} + \frac{11992\, x^{\left(\frac{35}{2}\right)}}{1276097777955} + \frac{281041808\, x^{\left(\frac{31}{2}\right)}}{778719703036275} - \frac{5996\, x^{\left(\frac{37}{2}\right)}}{15738539261445}$$

$$\Xi_5(-2) := \frac{1264\, x^{21}}{45989238101625} + \frac{49792\, x^{17}}{5731956140625} - \frac{24896\, x^{16}}{4458188109375} - \frac{158\, x^{22}}{168627206372625} - \frac{454673894\, x^{20}}{1462695175536469875} + \frac{1430928656\, x^{19}}{805974892642544625} - \frac{23959694116\, x^{18}}{4380298329579046875}$$

$$\Xi_6(-2) := \frac{378368\, x^{\left(\frac{37}{2}\right)}}{9622256002734375} - \frac{1618888\, x^{\left(\frac{51}{2}\right)}}{970482792234112282125} - \frac{189184\, x^{\left(\frac{39}{2}\right)}}{2624251637109375} + \frac{509590676527264\, x^{\left(\frac{41}{2}\right)}}{9211661065318008523078125} - \frac{2703289090900432\, x^{\left(\frac{43}{2}\right)}}{117759883348524811659890625} + \frac{99432991184\, x^{\left(\frac{45}{2}\right)}}{1795773411269482550375} - \frac{1601961931352\, x^{\left(\frac{47}{2}\right)}}{2078143185680721549014625} + \frac{3237776\, x^{\left(\frac{49}{2}\right)}}{57087223072594840125}$$

---

$$y_1(x) = -\frac{16\, x^{\left(\frac{49}{2}\right)}}{7453613144352375} + \frac{32\, x^{\left(\frac{47}{2}\right)}}{456343661899125} - \frac{241182752\, x^{\left(\frac{45}{2}\right)}}{2597873822637599756 25} + \frac{370906048\, x^{\left(\frac{43}{2}\right)}}{56666922594282225225} - \frac{4\, x^{21}}{3236479871625} - \frac{167699021552\, x^{\left(\frac{41}{2}\right)}}{626682365752353848437 5} + \frac{16\, x^{20}}{462354267375} + \frac{387460384\, x^{\left(\frac{39}{2}\right)}}{6074153637915287925} - \frac{2428744\, x^{19}}{6393192623692875} - \frac{21408256\, x^{\left(\frac{37}{2}\right)}}{259800912073828125} + \frac{74467936\, x^{18}}{35042386636632375} - \frac{141903184\, x^{\left(\frac{35}{2}\right)}}{299085416708203125} - \frac{245817524\, x^{17}}{38128556888296875} + \frac{32\, x^{\left(\frac{33}{2}\right)}}{2642024385} + \frac{405952\, x^{16}}{40123692984375} - \frac{15242768\, x^{\left(\frac{31}{2}\right)}}{144207352414125} - \frac{256\, x^{15}}{39734296875} + \frac{7874656\, x^{\left(\frac{29}{2}\right)}}{17942850300375} - \frac{4\, x^{14}}{26687115} - \frac{509312\, x^{\left(\frac{27}{2}\right)}}{581502796875} + \frac{32\, x^{13}}{11437335} + \frac{256\, x^{\left(\frac{25}{2}\right)}}{378421875} - \frac{32012\, x^{12}}{1762732125} + \frac{2848\, x^{11}}{57432375} - \frac{8\, x^{\left(\frac{21}{2}\right)}}{293265} - \frac{32\, x^{10}}{658125} + \frac{16\, x^{\left(\frac{19}{2}\right)}}{41895} - \frac{656\, x^{\left(\frac{17}{2}\right)}}{401625} + \frac{32\, x^{\left(\frac{15}{2}\right)}}{14625} - \frac{2\, x^{7}}{735} + \frac{8\, x^{6}}{315} - \frac{4\, x^{5}}{75} - \frac{4\, x^{\left(\frac{7}{2}\right)}}{35} + \frac{8\, x^{\left(\frac{5}{2}\right)}}{15} - 1$$

$$y_2(x) = \frac{8\,x^{\left(\frac{51}{2}\right)}}{16719489917031825} - \frac{765328\,x^{\left(\frac{49}{2}\right)}}{57087223072594840125} + \frac{1312571672\,x^{\left(\frac{47}{2}\right)}}{8332570912913879506875}$$
$$- \frac{97602780976\,x^{\left(\frac{45}{2}\right)}}{96796778631476966917875} + \frac{2\,x^{22}}{6690472155675} + \frac{459428946752\,x^{\left(\frac{43}{2}\right)}}{120779367536948524779375}$$
$$- \frac{25216\,x^{21}}{3541171333825125} - \frac{1008305818496\,x^{\left(\frac{41}{2}\right)}}{119069649492947231203125} + \frac{66170558\,x^{20}}{958317528799756125}$$
$$+ \frac{256\,x^{\left(\frac{39}{2}\right)}}{24825696328125} - \frac{93944624\,x^{19}}{268658297547514875} + \frac{101148332\,x^{\left(\frac{37}{2}\right)}}{761638571293359375}$$
$$+ \frac{2910672976\,x^{18}}{2978602864113751875} - \frac{3464\,x^{\left(\frac{35}{2}\right)}}{1276097777955} - \frac{64\,x^{17}}{44730984375}$$
$$+ \frac{4284008\,x^{\left(\frac{33}{2}\right)}}{205868427239475} + \frac{64\,x^{16}}{74551640625} - \frac{33732368\,x^{\left(\frac{31}{2}\right)}}{432622057242375} + \frac{1\,x^{15}}{22365315}$$
$$+ \frac{256\,x^{\left(\frac{29}{2}\right)}}{1789239375} - \frac{548\,x^{14}}{789176115} - \frac{512\,x^{\left(\frac{27}{2}\right)}}{4970109375} + \frac{34414\,x^{13}}{8749561275} - \frac{128\,x^{12}}{13253625}$$
$$+ \frac{2\,x^{\left(\frac{23}{2}\right)}}{213003} + \frac{64\,x^{11}}{7363125} - \frac{284\,x^{\left(\frac{21}{2}\right)}}{2639385} + \frac{32\,x^{\left(\frac{19}{2}\right)}}{80325} - \frac{64\,x^{\left(\frac{17}{2}\right)}}{133875} + \frac{1\,x^{8}}{882} - \frac{8\,x^{7}}{945} + \frac{8\,x^{6}}{525} + \frac{4\,x^{\left(\frac{9}{2}\right)}}{63} - \frac{8\,x^{\left(\frac{7}{2}\right)}}{35} + x$$

Reviewing the obtained data, we can see that it is very difficult to define construction patterns for both solutions. Thus, let's apply this analysis to the obtained summands $\Xi_3(-1)$ and $\Xi_3(-2)$. This way we can define the following patterns:

$$\Xi_k(-1) = \sum_{i=0}^{k+1} b_{i,k}\, x^{\left(\frac{5k}{2} + \frac{5}{2} + i\right)} \qquad \Xi_k(-2) = \sum_{i=0}^{k+1} c_{i,k}\, x^{\left(\frac{5k}{2} + \frac{7}{2} + i\right)}$$

In this case partial solutions are defined by the following formulas:

$$y_1(x) = -1 + \left( \sum_{k=0}^{N} \left( \sum_{i=0}^{k+1} b_{i,k}\, x^{\left(\frac{5k}{2} + \frac{5}{2} + i\right)} \right) \right)$$

$$y_2(x) = \left( \sum_{k=0}^{N} \left( \sum_{i=0}^{k+1} r_{i,k}\, x^{\left(\frac{5k}{2} + \frac{7}{2} + i\right)} \right) \right) + x$$

Substituting the obtained solutions to the initial equation (7.40), we can define recurrence equations to find the numerical coefficients $r_{i,k}, b_{i,k}$:

$$(5k + 5 + 2i)(5k + 3 + 2i)\, b_{i,k} - 4\, b_{i-1,k-1} + 8\, b_{i,k-1} = 0 \qquad (7.41)$$

$$(5k + 7 + 2i)(5k + 5 + 2i)\, r_{i,k} - 4\, r_{i-1,k-1} + 8\, r_{i,k-1} = 0 \qquad (7.42)$$

Thus, the desired general solution of the ODE (7.40) is as follows:

$$y = C_1\left(-1 + \left(\sum_{k=0}^{\infty}\left(\sum_{i=0}^{k+1} b_{i,k}\, x^{\left(\frac{5k}{2}+\frac{5}{2}+i\right)}\right)\right)\right) + C_2\left(\left(\sum_{k=0}^{\infty}\left(\sum_{i=0}^{k+1} r_{i,k}\, x^{\left(\frac{5k}{2}+\frac{7}{2}+i\right)}\right)\right) + x\right)$$

where numerical coefficients $r_{i,k}, b_{i,k}$ are defined by (7.41), (7.42). **The task is solved.**

**Example No. 7. Find a general solution for the following linear ODE**:

(7.43)

**is defined by the following formula:**

$$y = C_1\,(x+1)^{\left(\frac{1}{2}+\frac{\sqrt{5}}{2}\right)} + C_2\,(x+1)^{\left(-\frac{\sqrt{5}}{2}+\frac{1}{2}\right)} \tag{7.44}$$

**Define partial solution and evaluate the effectiveness of convergence of this series in comparison with the exact solution (7.44).**

**Solution**. If we enter the following equalities to the program:

**, for all** $k$,

we obtain the desired values of from $k = 0$ to $k = 6$:

$$\Xi_1(-1) := -\frac{1\ln(x+1)^2}{2} - \ln(x+1) - 2$$

$$\Xi_2(-1) := \frac{1\ln(x+1)^3}{6} + \ln(x+1)^2 + 4\ln(x+1) + 6$$

$$\Xi_3(-1) := -\frac{1\ln(x+1)^3}{2} - \frac{7\ln(x+1)^2}{2} - 13\ln(x+1) - 22 - \frac{1\ln(x+1)^4}{24}$$

$$\Xi_4(-1) := 46\ln(x+1) + \frac{1\ln(x+1)^5}{120} + 12\ln(x+1)^2 + \frac{11\ln(x+1)^3}{6} + 80 + \frac{1\ln(x+1)^4}{6}$$

$$\Xi_5(-1) := -\frac{1\ln(x+1)^5}{24} - 166\ln(x+1) - 43\ln(x+1)^2 - \frac{20\ln(x+1)^3}{3} - \frac{2\ln(x+1)^4}{3} - 296 - \frac{1\ln(x+1)^6}{720}$$

$$\Xi_6(-1) := 1106 + \frac{1\ln(x+1)^6}{120} + 610\ln(x+1) + 157\ln(x+1)^2 + \frac{11\ln(x+1)^5}{60} + \frac{74\ln(x+1)^3}{3} + \frac{31\ln(x+1)^4}{12} + \frac{1\ln(x+1)^7}{5040}$$

$$\Xi_0(-2) := 2\ln(x+1)(x+1) - 2x - 2$$

$$\Xi_1(-2) := -4\ln(x+1)x - 4\ln(x+1) + 4x + 4$$

$$\Xi_2(-2) := -14 - 14x + 14\ln(x+1) + 14\ln(x+1)x - \ln(x+1)^2 - \ln(x+1)^2 x + \frac{1\ln(x+1)^3}{3} + \frac{1\ln(x+1)^3 x}{3}$$

$$\Xi_3(-2) := 48 + 48x - 48\ln(x+1) - 48\ln(x+1)x + 4\ln(x+1)^2 + 4\ln(x+1)^2 x - \frac{4\ln(x+1)^3}{3} - \frac{4\ln(x+1)^3 x}{3}$$

$$\Xi_4(-2) := -172 - 172x + 172\ln(x+1) + 172\ln(x+1)x + \frac{1\ln(x+1)^5}{60} - 16\ln(x+1)^2 - 16\ln(x+1)^2 x + \frac{16\ln(x+1)^3}{3} + \frac{1\ln(x+1)^5 x}{60} + \frac{16\ln(x+1)^3 x}{3} - \frac{1\ln(x+1)^4}{12}$$

$$\Xi_5(-2) := 628 + 628x - 628\ln(x+1) - 628\ln(x+1)x - \frac{1\ln(x+1)^5}{10} + 62\ln(x+1)^2 + 62\ln(x+1)^2 x - \frac{62\ln(x+1)^3}{3} - \frac{1\ln(x+1)^5 x}{10} - \frac{62\ln(x+1)^3 x}{3} + \frac{1\ln(x+1)^4}{2}$$

$$\Xi_6(-2) := -2326 - 2326x + 2326\ln(x+1) + 2326\ln(x+1)x + \frac{29\ln(x+1)^5}{60} - 239\ln(x+1)^2 - 239\ln(x+1)^2 x + \frac{239\ln(x+1)^3}{3} + \frac{29\ln(x+1)^5 x}{60} + \frac{239\ln(x+1)^3 x}{3} - \frac{29\ln(x+1)^4}{12} - \frac{29\ln(x+1)^4 x}{12} + \frac{1\ln(x+1)^7 x}{2520} - \frac{1\ln(x+1)^6 x}{360} - \frac{1\ln(x+1)^6}{360} + \frac{1\ln(x+1)^7}{2520}$$

---

$$y_1(x) = 871 + 481\ln(x+1) + \frac{3\ln(x+1)^5}{20} + 123\ln(x+1)^2 + \frac{39\ln(x+1)^3}{2} + \frac{49\ln(x+1)^4}{24} + \frac{1\ln(x+1)^6}{144} + \frac{1\ln(x+1)^7}{5040}$$

$$y_2(x) = \left(-2705 + 1353\ln(x+1) - 313\ln(x+1)^2 + \frac{1\ln(x+1)^5}{4} + \frac{263\ln(x+1)^3}{6} - \frac{97\ln(x+1)^4}{24} + \frac{1\ln(x+1)^7}{5040} - \frac{7\ln(x+1)^6}{720}\right)x - 1834 - \frac{1\ln(x+1)^6}{360} + 1834\ln(x+1) + \frac{190\ln(x+1)^3}{3} + \frac{2\ln(x+1)^5}{5} - 190\ln(x+1)^2 - 2\ln(x+1)^4 + \frac{1\ln(x+1)^7}{2520}$$

Analyzing the obtained solutions, we establish that the formula for the first solution can be represented as follows:

(7.45)

where $b_i$ are the desired numerical coefficients.
In order to determine $b_i$, let's substitute (7.45) to (7.44):

$$\left(\frac{\partial^2}{\partial x^2}\left(\sum_{i=0}^{\infty} b_i \ln(x+1)^i\right)\right) - \left(\left(\frac{1}{(x+1)^2}\right) \cdot \left(\sum_{i=0}^{\infty} b_i \ln(x+1)^i\right)\right) = 0$$

Expanding this equality, we get:

$$\sum_{i=0}^{\infty} \left(\left(\left(\frac{1}{(x+1)^2 \ln(x+1)^2}\right) \cdot (-i\, b_i \ln(x+1)^{(i+1)} - b_i \ln(x+1)^i\, i + b_i \ln(x+1)^i\, i^2)\right) - \left(\left(\frac{1}{(x+1)^2}\right) \cdot (b_i)\right) \ln(x+1)^i\right)$$

From here, with the use of previously described method, we obtain the following recurrence equation:

$$i\,(i-1)\, b_i + (-i+1)\, b_{i-1} - b_{i-2} = 0 \qquad (7.46)$$

In order to define effectiveness, that is, convergence rate of the partial solution (7.45) with consideration of (7.46), we can use a

**PROGRAM**
**for the calculation of relative accuracy**

$$\left|\frac{\left(\frac{\partial^2}{\partial x^2} y\right) - \left(\left(\frac{1}{(x+1)^2}\right) \cdot y\right)}{\left(\frac{1}{(x+1)^2}\right) \cdot y}\right| \le \delta$$

**of partial solution (7.45), where $b_i$ coefficients are defined by (7.46).**

```
> restart:N:=26:

> alias(y=y(x),a=a(x)):a:=1/(x+1)^2:yr:=
diff(y,x$2)=a*y;b[0]:=871:b[1]:=481:
```

> **for** $i$ **from** 0 **to** $N$ **do** $b_{i+2} := \frac{-(-i-1)\, b_{i+1} + b_i}{(i+2)\,(i+1)}$ **end do**

```
> Y:=y =add(b[i]*ln(x+1)^i,i = 0 .. N);

> K0:=0:K1:=10;

> dy:=expand(subs(Y,(lhs(yr)-rhs(yr))/(a*y))):for z from K0 to K1 do
he(z):=evalf(subs(x=z,dy)) od:with(plots):

s := 's':for s from K0 to K1 do if abs(he(s))<1 then print(he(s),x=s) else
he(s):=(he(s)) fi od:s := 's': x:='x':print(yr);print(pointplot({seq([n,
(he(n))],n=K0..K1)},numpoints=30,axes=BOXED,color=red,xtickmarks=20,ytickm
arks=20,scaling=UNCONSTRAINED,title="RELATIVE ACCURACY GRAPH""")));
```

Relative accuracy for the specified values of $x$ in the range from -10 to 10 is as follows.

$$Y := y = 871 + 481\ln(x+1) + 676\ln(x+1)^2 + \frac{611}{2}\ln(x+1)^3 + \frac{3185}{24}\ln(x+1)^4 + \frac{2509}{60}\ln(x+1)^5 + \frac{8203}{720}\ln(x+1)^6 + \frac{1469}{560}\ln(x+1)^7 + \frac{1339}{2520}\ln(x+1)^8 + \frac{6929}{72576}\ln(x+1)^9 + \frac{56069}{3628800}\ln(x+1)^{10} + \frac{15119}{6652800}\ln(x+1)^{11} + \frac{20969}{68428800}\ln(x+1)^{12} + \frac{18269}{479001600}\ln(x+1)^{13} + \frac{739}{167650560}\ln(x+1)^{14} + \frac{15943}{33530112000}\ln(x+1)^{15} + \frac{77389}{1609445376000}\ln(x+1)^{16} + \frac{62609}{13680285696000}\ln(x+1)^{17} + \frac{202607}{492490285056000}\ln(x+1)^{18} + \frac{1457}{41588068515840}\ln(x+1)^{19} + \frac{37}{13054290480000}\ln(x+1)^{20} + \frac{858257}{3930072474746880000}\ln(x+1)^{21} + \frac{1388689}{86461594444431360000}\ln(x+1)^{22} + \frac{374491}{331436112036986880000}\ln(x+1)^{23} + \frac{727127}{9545360026665222144000}\ln(x+1)^{24} + \frac{5882581}{119317000333315276800000}\ln(x+1)^{25} + \frac{1189777}{387780251083274649600000}\ln(x+1)^{26}$$

$$K1 := 10$$

```
Warning, the name changecoords has been redefined
```

$-0., x = 0$

$-0.8239920116\ 10^{-24}, x = 1$

$-0.4782780812\ 10^{-19}, x = 2$

$-0.1049226610\ 10^{-16}, x = 3$

$-0.3115563645\ 10^{-15}, x = 4$

$-0.3434897659\ 10^{-14}, x = 5$

$-0.2124522542\ 10^{-13}, x = 6$

$-0.9049928483\ 10^{-13}, x = 7$

$-0.2979190031\ 10^{-12}, x = 8$

$-0.8132629957\ 10^{-12}, x = 9$

$-0.1927480134\ 10^{-11}, x = 10$

As we can see, with **N=26** the calculation accuracy is already sufficient for any practical applications.To compare the effectiveness of the obtained solution and the exact partial solution, for instance:

let’s present this solution in a series up to **N=26**

$$y = 1.618033988\,x + .4999999992\,x^2 + 1. + .1452538823\text{e-}3\,x^{22} - .6366100191\text{e-}1\,x^3 + .2199433523\text{e-}1\,x^4 - .1047795179\text{e-}1\,x^5 + .5906012806\text{e-}2\,x^6 - .3697135339\text{e-}2\,x^7 + .2487232093\text{e-}2\,x^8 - .1763714520\text{e-}2\,x^9 + .1301968064\text{e-}2\,x^{10} - .9920956424\text{e-}3\,x^{11} + .7756506330\text{e-}3\,x^{12} - .6194445006\text{e-}3\,x^{13} + .5036068751\text{e-}3\,x^{14} - .4157095473\text{e-}3\,x^{15} + .3476881895\text{e-}3\,x^{16} - .2941435132\text{e-}3\,x^{17} + .2513614180\text{e-}3\,x^{18} - .2167260109\text{e-}3\,x^{19} + .1883562077\text{e-}3\,x^{20} - .1648741623\text{e-}3\,x^{21} - .1287199866\text{e-}3\,x^{23} + .1146785991\text{e-}3\,x^{24} - .1026693003\text{e-}3\,x^{25}$$

Let's present a program that calculates relative error for this private solution with the same interval of ***x***: from -10 to 10, for an exact solution represented by a series, up to 26[th] term of the series:

```
> restart:

> alias(y=y(x),a=a(x)):a:=1/(x+1)^2:yr:= diff(y,x$2)=1/(x+1)^2*y;

> Y:=y = 1.618033988*x+.4999999992*x^2+1.+.1452538823e-3*x^22-.6366100191e-1*x^3+.2199433523e-1*x^4-.1047795179e-1*x^5+.5906012806e-2*x^6-.3697135339e-2*x^7+.2487232093e-2*x^8-.1763714520e-2*x^9+.1301968064e-2*x^10-.9920956424e-3*x^11+.7756506330e-3*x^12-.6194445006e-3*x^13+.5036068751e-3*x^14-.4157095473e-3*x^15+.3476881895e-3*x^16-.2941435132e-3*x^17+.2513614180e-3*x^18-.2167260109e-3*x^19+.1883562077e-3*x^20-.1648741623e-3*x^21-.1287199866e-3*x^23+.1146785991e-3*x^24-.1026693003e-3*x^25:

> K0:=0:K1:=10:

> dy:=expand(subs(Y,(lhs(yr)-rhs(yr))/(a*y))):for z from K0 to K1 do he(z):=evalf(subs(x=z,dy)) od:with(plots):

s := 's':for s from K0 to K1 do if abs(he(s))<1 then print(he(s),x=s) else he(s):=(he(s)) fi od:s := 's': x:='x':print(yr);print(pointplot({seq([n,(he(n))],n=K0..K1)},numpoints=30,axes=BOXED,color=red,xtickmarks=20,ytickmarks=20,scaling=UNCONSTRAINED,title="RELATIVE ACCURACY GRAPH""")));
```

$$K1 := 10$$

```
Warning, the name changecoords has been redefined
```

**As we can see, the effectiveness of the obtained solution is five time higher than the convergence rate of the exact solution. The task is solved**.

**Example No. 8. Obtain with the relative accuracy** $\delta = 10^{(-4)}$ :

**both partial solutions of linear ODE in the range of [-3,3]:**

$$\frac{\partial^2}{\partial x^2} y = (x^4 - 2x^2 + x - 3)\, y$$

**Solution**: On the basis of the theory of second-order linear ODEs, there is a
**PROGRAM**
**for the calculation of partial solutions of linear ODE:**

$$\frac{\partial^2}{\partial x^2} y = a\, y$$

**with specified relative accuracy:**

**for the [K0,K1] interval of the independent variable** ***x*** **.**

```
> restart:alias(C=binomial,ex=expand,y=y(x)):
> N:=16;a[1](x):=0:a[2](x):=x^4-2*x^2+x-3:
> yr:=diff(y,x$2)=a[1](x)*diff(y,x)+a[2](x)*y;
> b[1](x):= -a[1](x):b[2](x):= a[2](x)-diff(a[1](x),x):
> for k from 0 to N do
> dn := proc (n::integer, w) local k, Ds; global resd, x; option remember;
if n = 0 then resd := w elif 0 < n then resd := diff(w,`$`(x,n)) else
Ds := w; for k to -n do Ds := int(Ds,x) end do; resd := Ds end if;
RETURN(resd) end proc:
> R:=H->add(dn(l-2,H)*b[2-l](x),l = 0 .. 1):
> G[0]:=b[2](x):P(-1):=0:G[-1]:=0:
> for s to k+1 do  G[s]:=R(G[s-1]) od;
```

$> \quad P(k) := \mathrm{add}(\xi_{k-s-1}(-1)\, P(s), s = 0..k-1) - \mathrm{add}(\xi_{k-s-1}(-2)\, G_s, s = 0..k-1)$

$$\xi_k(-1) := \mathrm{dn}\left(-1, -\int G_k\, dx - P(k)\right) - \mathrm{add}\left(\xi_{k-s-1}(-1)\, \mathrm{dn}\left(-1, -\int G_s\, dx - P(s)\right), s = 0..k-1\right) - \mathrm{add}(\xi_{k-s-1}(-2)\, \mathrm{dn}(-1, G_s), s = 0..k-1)$$

$$\xi_k(-2) := \mathrm{dn}\left(-2, -2\int G_k\, dx - P(k)\right) - \mathrm{add}\left(\xi_{k-s-1}(-1)\, \mathrm{dn}\left(-2, -2\int G_s\, dx - P(s)\right), s = 0..k-1\right) - \mathrm{add}(\xi_{k-s-1}(-2)\, \mathrm{dn}(-2, G_s), s = 0..k-1)$$

```
> od:
> z1:=alpha(-1)=ex(ex(ex(ex(ex(ex(value(add(xi[k](-1),k=0..N)))))))):
> z2:=alpha(-2)=ex(ex(ex(ex(ex(ex(value(add(xi[k](-2),k=0..N)))))))):
> Y1:=y=-1+rhs(z1);
> K0:=-10.1:K:=10.1:
> dy:=subs(Y1,(lhs(yr)-rhs(yr))/(a[2](x)*y)):for s from K0 to K do
he(s):=evalf(subs(x=s,dy)) od:with(plots):
s := 's':for s from K0 to K do if abs(he(s))<1 then print(he(s),x=s) else
he(s):=(he(s)) fi od:s := 's': x:='x':print(yr);print(pointplot({seq([n,
(he(n))],n=K0..K)},numpoints=30,axes=BOXED,color=red,xtickmarks=20,ytickma
rks=20,scaling=UNCONSTRAINED,title=""""));
> print("############################################");
> Y2:=y=expand(rhs(z2)-rhs(z1)*x+x);
> dy:=subs(Y2,(lhs(yr)-rhs(yr))/(a[2](x)*y)):for s from K0 to K do
he(s):=evalf(subs(x=s,dy)) od:with(plots):
s := 's':for s from K0 to K do if abs(he(s))<1 then print(he(s),x=s) else
he(s):=(he(s)) fi od:s := 's': x:='x':print(yr);print(pointplot({seq([n,
(he(n))],n=K0..K)},numpoints=30,axes=BOXED,color=red,xtickmarks=20,ytickma
rks=20,scaling=UNCONSTRAINED,title="RELATIVE ERROR"""));
> print("############################################");
```

With the use of this program, we shall present both partial solutions of the initial ODE, which satisfy the specified conditions in the defined interval. Also we present the graphs of relative accuracy for each of the partial solutions in the defined interval.

$$N := 10$$

$$yr := \frac{\partial^2}{\partial x^2} y = (x^4 - 2x^2 + x - 3)\, y$$

**The first partial solution, which satisfies the above conditions, is defined by the following formula:**

$$Y1 := y = -1 - \frac{2287703133976036138263224944 31}{116463234040694957416296348646972571217024000000000000} x^{54} - \frac{5}{24} x^4 + \frac{3}{2} x^2 + \frac{1}{10} x^5 - \frac{1}{6} x^3 - \frac{17}{144} x^6 - \frac{1}{240} x^7 + \frac{4121}{2661120} x^{11} - \frac{17}{2592} x^9 + \frac{569}{13440} x^8 - \frac{746953}{479001600} x^{12} - \frac{4169}{3628800} x^{10} + \frac{5249}{259459200} x^{13} + \frac{15168047}{2117187072000} x^{17} - \frac{10030771}{186810624000} x^{15} + \frac{7031}{25159680} x^{14} - \frac{6535078051}{914624815104000} x^{18} + \frac{2868893}{597793996800} x^{16} + \frac{1198386590660250847}{611818352245381162033152000000} x^{31} - \frac{3184406869236619152583}{84798023621209829057794867200000000} x^{33} + \frac{171459194770539300 07}{117775032807235873691381760000000000} x^{35} + \frac{31002456404085277}{279096022529662309171200000000} x^{32}$$

$$+\frac{94728908994900254153204829461}{1151260746236440599001507173543058052428800000000000}x^{49}$$

$$+\frac{6592351}{20324995891200}x^{19}+\frac{247866275147}{1598764176801792000}x^{23}-\frac{45221359}{249443131392000}x^{21}$$

$$+\frac{12998835937}{16550353797120000}x^{20}-\frac{2531240104812571}{15511210043330985984000}x^{24}$$

$$+\frac{5928436247029}{16057153253965824000}x^{22}+\frac{76800386120663}{6803162299706572800000}x^{25}$$

$$+\frac{446503311556171}{2349532842724198010880000}x^{29}-\frac{17301110229257359}{523503338962420776960000000}x^{27}$$

$$+\frac{527945555677997}{442205549480927232000000}x^{26}-\frac{6661750398524497}{300244562051976622080000000}x^{30}$$

$$+\frac{845722629648471677}{9527760769116058140672000000}x^{28}$$

$$-\frac{143797828077266615053}{72078320078028354699125637120000000}x^{36}$$

$$+\frac{169294224814513897258 1}{144156640156056709398251274240000000}x^{34}$$

$$+\frac{886609034374377982017 1}{4242792022774850878880350003200000000}x^{37}$$

$$+\frac{7753522875390471889327474 1}{508081875683128818106231925504455680000000000}x^{43}$$

$$+\frac{1248839243378477124081 7}{18473823599165496535124857305600000000}x^{38}$$

$$-\frac{4200451590074680715917}{14445696648971516087616429772800000000}x^{39}$$

$$+\frac{22927438795888554899555 47}{3150895353073667089030895662043136000000000}x^{41}$$

$$+\frac{2280832109700523581263}{22603266521331901642976295997440000000}x^{40}$$

$$-\frac{60125303468413386376290347}{472634302961050063354634349306470400000000000}x^{42}$$

$$+\frac{475822813692271963742531353 3}{7712682872869895458852600629157637222400000000000}x^{46}$$

$$-\frac{164180737019940243365678697 79}{27284652313701028021102210827880243507200000000000}x^{48}$$

$$+\frac{9161200379129625976551 37}{334378670321375376018631267212334080000000000}x^{44}$$

$$-\frac{31451954072071797657864 49}{1927681913160922345118495110315008000000000000}x^{45}$$

$$+\frac{1661899192088128640740976 71}{7661451283738513303026888910002604608000000000000}x^{47}$$

$$+\frac{458869201805224336109446818486 77}{163561258876020025100999840584081604734348800000000000000}x^{52}$$

$$+\frac{50656872573109731324043443 9449}{863445559677330449251130380157293539321600000000000000}x^{50}$$

$$-\frac{1736061197369360705884636 8311}{2648765325927449799206475151159216271001600000000000}x^{51}$$

$$
- \frac{1901984868711583428353257463}{137207133909924997314862164465912077412480000000000000} x^{53}
+ \frac{26511372463031787 47407}{86382789636124477542046644011834961778483200000000} x^{61}
$$

$$
- \frac{2308440694527812281099357 7}{1177606828165036274913467612050372571019863040000000 00} x^{59}
$$

$$
- \frac{840174304494224632599745 1461}{620557065073100779692459049801583080568248320000000000} x^{57}
$$

$$
- \frac{18064778961158297540528951561}{1794771577728584126367078996057210155470080000000000000} x^{56}
$$

$$
- \frac{110573937741042522494202 2539}{363734180060443157103449759461712711285425152000000000 0} x^{60}
$$

$$
+ \frac{16604738469241700140815256 3517}{2218491507636335287400541103040659513031487744000000000 0} x^{58}
$$

$$
+ \frac{405122048456637008813619442 7801}{13954567891422391409836522984466279350262246400000000000} x^{55}
$$

$$
- \frac{1104804203}{17784836309792449730438243315758 08000000} x^{63}
$$

$$
- \frac{18175128474381512328345493}{5399563616199199354513348693748359650682129536000000000 0} x^{62}
$$

$$
+ \frac{127091359}{262521136524867206160578991287500800000} x^{64}
- \frac{1}{12380777361784688317711810560000 0} x^{66}
$$

```
Warning, the name changecoords has been redefined
```

The relative error of the obtained first partial solution with the specified values of independent variables is defined by the following values:

$0.00001440608857, x = -3.1$

$-0.2477354033\ 10^{-7}, x = -2.1$

$-0.8060318022\ 10^{-9}, x = -1.1$

$0.3777416540\ 10^{-9}, x = -0.1$

$-0.5694933139\ 10^{-8}, x = 0.9$

$-0.5933001644\ 10^{-8}, x = 1.9$

$-.4697685561e-7\ x = 2.9$

"#############################################"

The second partial solution, which satisfies the specified conditions, is defined by the following formula:

$$
Y2 := y = x + \frac{10561158690671752703312547573 17}{33341333540111774347511505965216634811232640000000000 00} x^{54}
+ \frac{1}{12} x^4 - \frac{1}{40} x^5 - \frac{1}{2} x^3 - \frac{1}{40} x^6 + \frac{37}{720} x^7 - \frac{199}{211200} x^{11} - \frac{151}{17280} x^9 - \frac{1}{480} x^8 - \frac{23}{84480} x^{12}
+ \frac{11}{5184} x^{10} + \frac{420919}{889574400} x^{13} - \frac{112207499}{16937496576000} x^{17} - \frac{2534017}{62270208000} x^{15}
- \frac{10513}{296524800} x^{14} - \frac{3910763}{4234374144000} x^{18} + \frac{721879}{53374464000} x^{16}
$$

$$+\frac{35512933519360518677}{917727528368071743049728000000}x^{31}-\frac{55325023776277433203}{2826600787373660968593162240000000}x^{33}$$

$$-\frac{72299784529615455845 33}{2102284335609160345391164416000000000}x^{35}-\frac{125587817571135947}{22247940081650224073932800000}x^{32}$$

$$+\frac{1506910074031843735200263112583}{20722693432255930782027129123775044943718400000000000}x^{49}$$

$$+\frac{1572395017}{914624815104000}x^{19}-\frac{182893787941}{9256003128852480000}x^{23}-\frac{132740873}{1662954209280000}x^{21}$$

$$-\frac{7098527}{40649991782400}x^{20}-\frac{2343053983}{1682909659791360000}x^{24}+\frac{24615911}{645617516544000}x^{22}$$

$$+\frac{982486114265413}{2982925008332881920000000}x^{25}-\frac{25082595173751690 1}{767514061956571350220800000 0}x^{29}$$

$$-\frac{877212252908071}{1193954983598503526400000 0}x^{27}-\frac{264947011061989}{646300418472124416000000}x^{26}$$

$$-\frac{65073848166403}{685280412461224419840000 00}x^{30}+\frac{6292532295109351}{1047006677924841553920000 00}x^{28}$$

$$-\frac{233927598947478011}{1413300393686830484296581120000 0000}x^{36}$$

$$+\frac{25531043030780747547 7}{4239901181060491452889743360000000}x^{34}$$

$$+\frac{71341029435090694031539}{2333535612526167983384192501760000000 0}x^{37}$$

$$+\frac{4389480113323438841220506 9}{2540409378415644090531159627522278400000000000}x^{43}$$

$$-\frac{36174475025069093939}{7071320037958084798133916672000000 0}x^{38}$$

$$+\frac{10375151693674830630821 9}{3842555308626423279305970319564800000000}x^{39}$$

$$-\frac{34331831630364667213981 7}{1369954501336376995230824200888320000000 0}x^{41}$$

$$+\frac{70054793362825675453 3}{168533127571334354355525014016000000000}x^{40}$$

$$+\frac{10077471290471511963 23}{19096335473173739933520579769958400000000}x^{42}$$

$$+\frac{23962024397517680705993 27}{1147281654768355395723749509203609600000000000}x^{46}$$

$$+\frac{27233322075980379662987 3}{520302294311613806657174119524794880000000000 0}x^{48}$$

$$-\frac{30918044554228779478722 7}{940892362376164477974503565748992000000000 0}x^{44}$$

$$+\frac{61118169833102882031698 3}{1667262068084420157433618506875904000000000 0}x^{45}$$

$$-\frac{82915541250207194972513629 7}{618292559740301073577608578701964582400000000000 0}x^{47}$$

$$+\frac{1}{340055462002209314737427251200000}x^{67}$$

$$+\frac{65610313140832156512823738 07}{859233630117928837303563890497989668398080000000000 0}x^{52}$$

$$
\begin{aligned}
&- \frac{5732861630819074150444802917}{363556025127297031263633844276755174451200000000000} x^{50} \\
&+ \frac{18325746322012402262546979557 21}{733928725725730881863460823133699508423360000000000000} x^{51} \\
&- \frac{117293868237993984393189281957}{215548547325246623019575663742164945142240000000000000} x^{53} \\
&+ \frac{8767828280414528840706070855 8469}{24935742153916671568753432742433795665265628481587200000000000} x^{61} \\
&- \frac{256419077250678630368869165 67}{174521331934058375942175900105865215025143702528000000000} x^{59} \\
&+ \frac{41493166312786735453382314506 8041}{3623536129139347636087550468299743871284763315200000000000000} x^{57} \\
&- \frac{770483329986123084938627357 17}{144714037392528503509415793912983637706423296000000000000} x^{56} \\
&+ \frac{753425192165647912508489 63}{8136192630958432444856685319620755945228144640000000000} x^{60} \\
&+ \frac{295229456449440237992272331}{1837803615793413847550744109027765277067504640000000000} x^{58} \\
&+ \frac{21723956551866404934425330262468983}{100124024620955658365577052413545554338131617920000000000000000} x^{55} \\
&- \frac{13589632709148467}{851973645278120460128638531119975008747520000000} x^{65} \\
&+ \frac{32783166995066815187495057 3}{14459558675110677490298473906012863180382238737920000000000} x^{63} \\
&- \frac{22115276444026338193741}{27850086617967234897351927883918570563290611200000000} x^{62} \\
&+ \frac{2486353501}{146603964720394301987193602155905024000000} x^{64}
\end{aligned}
$$

The relative error of the obtained second partial solution with the specified values of independent variable is defined by the following values:

$0.2391000929\ 10^{-5}, x = -3.1$
$-0.3276463183\ 10^{-7}, x = -2.1$
$-0.1012319411\ 10^{-8}, x = -1.1$
$-0.9046644911\ 10^{-47}, x = -0.1$
$0.8631757216\ 10^{-14}, x = 0.9$
$0.2323828797\ 10^{-7}, x = 1.9$
$-.5897879005e-6\ x = 2.9$
"################################################"

As can be seen, in the given interval the obtained solutions already ensure the required accuracy with ten terms of the series. **The task is solved.**

**Example No. 9. Obtain with the relative accuracy** $\delta = 10^{(-5)}$ :

**both partial solutions of linear ODE in the range of [1,3]:**

**Solution**: with the use of the presented PROGRAM we define for $N := 6$ , that the first partial solution satisfying the established conditions is as follows:

$$Y1 := y = -1 - \frac{1}{6} x^3 \ln(x) + \frac{5 x^3}{36} - \frac{8530899170945477354381\ x^{18}}{1934198611137401252759145676800000}$$
$$- \frac{1}{180} x^6 \ln(x)^2 + \frac{47}{5400} x^6 \ln(x) - \frac{457 x^6}{162000} + \frac{158087}{145050701206118400 0} x^{21} \ln(x)^6$$
$$- \frac{20819858459}{949356839394044928000 00} x^{21} \ln(x)^5$$
$$+ \frac{6868106068749683}{29824994466403315458048000000} x^{21} \ln(x)^4$$
$$- \frac{128091145386350 35309}{9369820261565265584300359680000 0} x^{21} \ln(x)^3$$
$$+ \frac{2559356429799391719648 67}{5575043055631333022658714009600000000} x^{21} \ln(x)^2$$
$$- \frac{160418588058346007732585249}{19902903708603858890891609014272000000000} x^{21} \ln(x)$$
$$+ \frac{227699}{5755980206592000} x^{18} \ln(x)^5 - \frac{831400613}{12557630150714880000} x^{18} \ln(x)^4$$
$$+ \frac{434196781485847}{78902101762971734016 00000} x^{18} \ln(x)^3$$
$$- \frac{148954726814081900 09}{619697107246379998961664000000} x^{18} \ln(x)^2$$
$$+ \frac{115695755979527481162 17}{22123186728695765962931404800000 00} x^{18} \ln(x) - \frac{1}{12960} x^9 \ln(x)^3$$
$$+ \frac{91}{518400} x^9 \ln(x)^2 - \frac{1199}{10368000} x^9 \ln(x) + \frac{5581 x^9}{248832000} - \frac{1}{1710720} x^{12} \ln(x)^4$$
$$+ \frac{3923}{2258150400} x^{12} \ln(x)^3 - \frac{429991}{248396544000} x^{12} \ln(x)^2 + \frac{91046557}{131153375232000} x^{12} \ln(x)$$
$$- \frac{328000063 x^{12}}{3462449106124800} - \frac{1}{46170964224000} x^{21} \ln(x)^7 - \frac{1}{109930867200} x^{18} \ln(x)^6$$
$$- \frac{1}{359251200} x^{15} \ln(x)^5 + \frac{33841}{3319481088000} x^{15} \ln(x)^4$$
$$- \frac{104356993}{7668001313280000} x^{15} \ln(x)^3 + \frac{394947097471}{4723488808980480000 0} x^{15} \ln(x)^2$$
$$- \frac{57504306348079}{242472425527664640000 00} x^{15} \ln(x) + \frac{3786437592072193\ x^{15}}{152757628082428723200000 00}$$
$$+ \frac{47497816195154206292068448 29\ x^{21}}{835921955761362073417447578599424000000000 0}$$

Warning, the name changecoords has been redefined

Relative accuracy for the first solution

$$0.5689094752\ 10^{-10}, x = 1.1$$
$$0.7943952838\ 10^{-9}, x = 2.1$$

$$0.1673119868\ 10^{-5}, x = 3.1$$

The second partial solution satisfying the established conditions is as follows:

$$Y2 := y = x + \frac{76539562714198210013282873}{485357121663778101182425496970854400000} x^{22} \ln(x)^3 + \frac{1}{12} x^4 \ln(x)$$
$$- \frac{7 x^4}{144} + \frac{13904886805637}{465900913312568219566080000} x^{22} \ln(x)^5$$
$$+ \frac{409044628581278169705788962367 5571}{498644801826421953911868901692485117411328000000 0} x^{22} \ln(x)$$
$$- \frac{16653492899020112418045510047}{336131346269399934688704616363720704000000 00} x^{22} \ln(x)^2 + \frac{1}{504} x^7 \ln(x)^2$$
$$- \frac{101}{42336} x^7 \ln(x) + \frac{1145 x^7}{1778112} + \frac{1}{268334780313600} x^{22} \ln(x)^7$$
$$+ \frac{1}{580811212800} x^{19} \ln(x)^6 - \frac{609763200636373 49737}{2126632672870016392302546124800000} x^{22} \ln(x)^4$$
$$- \frac{1206249332648687}{164684663286758041190400000} x^{19} \ln(x)^3$$
$$+ \frac{4801229682268587073}{16270844732731694469611520 00000} x^{19} \ln(x)^2 + \frac{1}{45360} x^{10} \ln(x)^3$$
$$- \frac{257}{6350400} x^{10} \ln(x)^2 + \frac{10123}{444528000} x^{10} \ln(x) - \frac{5791 x^{10}}{1481760000}$$
$$- \frac{797868478177609 x^{16}}{262846432759460659200 00000} + \frac{1}{1698278400} x^{16} \ln(x)^5$$
$$- \frac{68197}{37090400256000} x^{16} \ln(x)^4 + \frac{438441901}{202513585397760000} x^{16} \ln(x)^3$$
$$- \frac{36118962017}{300877326876672000 00} x^{16} \ln(x)^2 + \frac{34219823279399}{109519346983108608000000} x^{16} \ln(x)$$
$$+ \frac{1}{7076160} x^{13} \ln(x)^4 - \frac{13523}{38635833600} x^{13} \ln(x)^3 + \frac{21332947}{70317217152000} x^{13} \ln(x)^2$$
$$- \frac{2984473133}{27423714689280000} x^{13} \ln(x) + \frac{11594425933 x^{13}}{855619898305536000}$$
$$- \frac{634618227756454789951}{10599293140179503825918361600000 00} x^{19} \ln(x)$$
$$+ \frac{120578906286632276491549 x^{19}}{253747077775897321592485576704000000 0} - \frac{521701}{80337806954496000} x^{19} \ln(x)^5$$
$$+ \frac{26925528179}{27780813644864716800 00} x^{19} \ln(x)^4$$
$$- \frac{2529577527686026083063905502939 2449 x^{22}}{46074779688761388541456686516385624848806707 2000000}$$
$$- \frac{20190013}{1224830204828246016000} x^{22} \ln(x)^6$$

Relative accuracy for the second solution

$$0., x = 1.1$$
$$-0.2730613555\ 10^{-9}, x = 2.1$$
$$0.1017074219\ 10^{-6}, x = 3.1$$

"##################################################"

**The task is solved.**

**Example No. 10. Find a general solution for the following linear ODE**:

(7.47)

**Solution**. Let's apply the following transformation:

(7.48)

(7.47) can be reduced as follows:

$$\frac{d^2}{dx^2} Z(x) = \left( \frac{1}{4} + x\, \mathbf{e}^x \right) Z(x) \tag{7.49}$$

With the use of the program presented above for partial solutions of ODEs, for (7.49) we have:

$$\begin{aligned} Z_1(x) = & - \frac{x^8}{10321920} - \frac{x^6}{46080} - \frac{x^4}{384} - \frac{x^2}{8} \\ & + \left( - \frac{x^7}{46080} + \frac{7\,x^6}{23040} - \frac{23\,x^5}{3840} + \frac{11\,x^4}{192} - \frac{13\,x^3}{24} + 3\,x^2 - \frac{23\,x}{2} + 20 \right) \mathbf{e}^x \\ & + \left( - \frac{x^6}{1536} + \frac{x^5}{96} - \frac{391\,x^4}{3072} + \frac{117\,x^3}{128} - \frac{135\,x^2}{32} + \frac{1325\,x}{128} - \frac{1085}{128} \right) \mathbf{e}^{(2\,x)} \\ & + \left( - \frac{x^5}{288} + \frac{5\,x^4}{108} - \frac{1603\,x^3}{5184} + \frac{923\,x^2}{864} - \frac{6389\,x}{3888} + \frac{559}{648} \right) \mathbf{e}^{(3\,x)} \\ & + \left( - \frac{x^4}{576} + \frac{x^3}{72} - \frac{265\,x^2}{6912} + \frac{899\,x}{20736} - \frac{2801}{165888} \right) \mathbf{e}^{(4\,x)} - 1 \end{aligned}$$

$$\begin{aligned} Z_2(x) = & \frac{x^9}{92897280} + \frac{x^7}{322560} + \frac{x^5}{1920} + \frac{x^3}{24} + x \\ & + \left( \frac{x^6}{864} - \frac{17\,x^5}{864} + \frac{337\,x^4}{1728} - \frac{973\,x^3}{864} + \frac{9287\,x^2}{2592} - \frac{124525\,x}{23328} + \frac{97475}{34992} \right) \mathbf{e}^{(3\,x)} \\ & + \left( \frac{x^7}{7680} - \frac{19\,x^6}{7680} + \frac{199\,x^5}{5120} - \frac{571\,x^4}{1536} + \frac{335\,x^3}{128} - \frac{2969\,x^2}{256} + \frac{3539\,x}{128} - \frac{361}{16} \right) \mathbf{e}^{(2\,x)} \\ & + \left( \frac{x^8}{322560} - \frac{x^7}{20160} + \frac{3\,x^6}{2560} - \frac{13\,x^5}{960} + \frac{x^4}{6} - \frac{5\,x^3}{4} + \frac{31\,x^2}{4} - 28\,x + 48 \right) \mathbf{e}^x \\ & + \left( \frac{x^5}{576} - \frac{73\,x^4}{3456} + \frac{715\,x^3}{6912} - \frac{4949\,x^2}{20736} + \frac{125423\,x}{497664} - \frac{95729}{995328} \right) \mathbf{e}^{(4\,x)} \end{aligned}$$

From here we obtain the following representation for the desired solutions:

$$Z_1(x) = \sum_{i=0}^{k} \mathbf{e}^{(i\,x)} \left( \sum_{j=0}^{k-i} b_j\, x^j \right) \tag{7.50}$$

$$Z_2(x) = \sum_{i=0}^{k} \mathbf{e}^{(i\,x)} \left( \sum_{j=0}^{k-i} c_j\, x^{(j+1)} \right) \tag{7.51}$$

where $c_j$, $b_j$ are the desired numerical coefficients.

To define the specified numerical coefficients, we shall successively substitute their values (7.50), (7.51) to the ODE (7.49).

In particular, substituting (7.50) to (7.49) we get, upon transformations:

$$\sum_{i=0}^{k}\left(\sum_{j=0}^{k-i}\left(\frac{e^{(ix)}(2ijb_j x^{(j+1)}+i^2 b_j x^{(j+2)}-b_j x^j j+b_j x^j j^2)}{x^2}-b_j x^{(j+1)}e^{(x+ix)}-\frac{b_j x^j e^{(ix)}}{4}\right)\right)=0$$

From here we can easily find the desired representation to define $b_j$.

$$4(j+2)(j+1)b_{j+2}+8i(j+1)b_{j+1}+(2i-1)(2i+1)b_j-4b_{j-1}=0 \tag{7.52}$$

In a similar way, substituting (7.51) to (7.49), we get, upon transformations:

$$\sum_{i=0}^{k}\left(\sum_{j=0}^{k-i}\left(\frac{e^{(ix)}(2ic_j j x^{(j+2)}+2ic_j x^{(j+2)}+c_j j^2 x^{(j+1)}+c_j j x^{(j+1)}+i^2 c_j x^{(j+3)})}{x^2}-c_j x^{(j+2)}e^{(x+ix)}-\frac{c_j x^{(j+1)}e^{(ix)}}{4}\right)\right)=0$$

Hence we obtain the desired recurrence equation for the definition of $c_j$.

$$4(j+2)(j+1)c_{j+1}+8i(j+1)c_j+(2i-1)(2i+1)c_{j-1}-4c_{j-2}=0 \tag{7.53}$$

Thus, the general solution with consideration of (7.50), (7.51) and (7.48) takes the following form:

$$y=\left(C_1\left(\sum_{i=0}^{N}e^{(ix)}\left(\sum_{j=0}^{N-i}b_j x^j\right)\right)+C_2\left(\sum_{i=0}^{N}e^{(ix)}\left(\sum_{j=0}^{N-i}c_j x^{(j+1)}\right)\right)\right)e^{\left(\frac{x}{2}\right)}$$

where **N** is a sufficiently great natural number tending to infinity, and numerical coefficients $c_j, b_j$ are defined by the recurrence equations (7.52), (7.53).

**The task is solved.**

**The above examples lead to the following conclusions:**

**1) This method for solving linear ODEs is common for any representations.**

**2) In this method, the replacement of variables is basically needed to remove the first derivative, since in its reduced form (6.1) the initial ODE can be solved with a much higher accuracy because of the reduction of the number of summands in the formulas for and, thus, during the calculation the summands depending on $a_1(x)$ are zero.**

## 8. The common algorithm for linear ODEs with variable coefficients of *m* order.

## 8.1 Introduction.

**In accordance with the Definition 1.1., general solution of a linear homogeneous ODE of *m* order**

$$\left(\frac{\partial^m}{\partial x^m} y\right) + \left(\sum_{p=1}^{m} a_p(x) \left(\frac{\partial^{m-p}}{\partial x^{m-p}} y\right)\right) = 0 \qquad (8.1)$$

**where** $y = \mathrm{y}(x)$ **is the desired function, and** $a_p(x)$ **are the specified functions, should be found in the category of analytical representations.**

$$y_i = \mathrm{F}(x, \alpha_k(x)), (i, k) = 1, 2, 3 .. m \qquad (8.2)$$

$\phi_k(x)$ **components of this common form are found in the form of a series converging throughout the range of independent variable *x*:**

$$\alpha_k(x) = \sum_{j=0}^{\infty} R_{k,j}(x), k = 1, 2, 3 .. m \qquad (8.3)$$

**following from certain** $R_{k,j}(x)$ **functions, which are defined by the initial coefficients .**

**Thus, the task of obtaining a common algorithm to solve linear ODEs (8.1) consists of two main tasks:**

**- to determine a linear representation of the desired partial solutions (8 .2)**
**- to determine** $R_{i,j}(x)$ **functions from (8.3).**

Let's consider each task individually. For this purpose, we shall introduce new concepts.

## 8.2 General classic definitions:

Let's assume a linear, generally non-homogenous LDE of ***m*** order.

$$\left(\frac{\partial^m}{\partial x^m} y\right) + \left(\sum_{p=1}^{m} a_p(x) \left(\frac{\partial^{m-p}}{\partial x^{m-p}} y\right)\right) = \mathrm{f}(x)$$

**where** $y = \mathrm{y}(x)$ **is the desired function, and** $a_p(x)$ **are the specified functions**.

The task of constructing a single algorithm to obtain $y$ function is solved, if we somehow define partial solutions for the following homogeneous ODE:

$$\left(\frac{\partial^m}{\partial x^m} y\right) + \left(\sum_{p=1}^{m} a_p(x) \left(\frac{\partial^{m-p}}{\partial x^{m-p}} y\right)\right) = 0$$

In this case, the general solution for the initial equation is defined by formula [1]:

$$y = \left(\sum_{i=1}^{m} C_i y_i\right) + \left(\sum_{i=1}^{m} y_i \int \frac{\mathrm{W}_i(x)}{\mathrm{W}(x)} dx\right)$$

where $\mathrm{W}(x)$ is Wronskian and $\mathrm{W}_i(x)$ is a determinant, resulting from Wronskian upon replacement of $i$ th column by starting at the top.

## 8.3. Equation of *n*-image

As soon as on the basis of second-order linear ODEs we have defined a close relation between the initial equation and its adjoined analogue, it makes sense to consider the general theory based on the condition that the equation of *n*-image can be defined for the adjoined equation.

Let's assume a linear homogeneous ODE of $m$ order in reduced form.

$$\frac{\partial^m}{\partial x^m} y = \sum_{p=1}^{m} b_p(x) \left( \frac{\partial^{m-p}}{\partial x^{m-p}} y \right) \tag{8.4}$$

where $b_p(x)$ are specified functions.

Let's differentiate this equation for variable ***x*.**

$$\frac{\partial^{m+1}}{\partial x^{m+1}} y = \left( \sum_{p=1}^{m} \left( \frac{d}{dx} b_p(x) \right) \left( \frac{\partial^{m-p}}{\partial x^{m-p}} y \right) \right) + \left( \sum_{p=1}^{m} b_p(x) \left( \frac{\partial^{m+1-p}}{\partial x^{m+1-p}} y \right) \right)$$

On the right side of this equality we get (with $p = 1$) the $\frac{\partial^m}{\partial x^m} y$ derivative, that is, we obtain the following equation:

$$\frac{\partial^{m+1}}{\partial x^{m+1}} y = \left( \sum_{p=1}^{m} \left( \frac{d}{dx} b_p(x) \right) \left( \frac{\partial^{m-p}}{\partial x^{m-p}} y \right) \right) + \left( \sum_{p=2}^{m} b_p(x) \left( \frac{\partial^{m+1-p}}{\partial x^{m+1-p}} y \right) \right) + b_1(x) \left( \frac{\partial^m}{\partial x^m} y \right)$$

Replacing here $\frac{\partial^m}{\partial x^m} y$ by an equivalent expression, in accordance with (11.4), upon transformation we get:

$$\begin{aligned} \frac{\partial^{m+1}}{\partial x^{m+1}} y = & \left( \sum_{p=1}^{m-1} \left( \left( \frac{d}{dx} b_p(x) \right) + b_{1+p}(x) + b_1(x)\, b_p(x) \right) \left( \frac{\partial^{m-p}}{\partial x^{m-p}} y \right) \right) \\ & + \left( \left( \frac{d}{dx} b_m(x) \right) + b_1(x)\, b_m(x) \right) y \end{aligned} \tag{8.5}$$

Again, let's differentiate the equality (8.5) for the variable *x*.

$$\begin{aligned} \frac{\partial^{2+m}}{\partial x^{2+m}} y = & \left( \sum_{p=1}^{m-1} \left( \right. \right. \\ & \left( \left( \frac{d^2}{dx^2} b_p(x) \right) + \left( \frac{d}{dx} b_{1+p}(x) \right) + \left( \frac{d}{dx} b_1(x) \right) b_p(x) + b_1(x) \left( \frac{d}{dx} b_p(x) \right) \right) \left( \frac{\partial^{m-p}}{\partial x^{m-p}} y \right) \\ & + \left. \left. \left( \left( \frac{d}{dx} b_p(x) \right) + b_{1+p}(x) + b_1(x)\, b_p(x) \right) \left( \frac{\partial^{1+m-p}}{\partial x^{1+m-p}} y \right) \right) \right) + \left( \left( \frac{d}{dx} b_m(x) \right) + b_1(x)\, b_m(x) \right) \left( \frac{\partial}{\partial x} y \right) \\ & + \left( \left( \frac{d^2}{dx^2} b_m(x) \right) + \left( \frac{d}{dx} b_1(x) \right) b_m(x) + b_1(x) \left( \frac{d}{dx} b_m(x) \right) \right) y \end{aligned}$$

As soon as in the first summand of this equality (with $p = 1$) we again have the $\frac{\partial^m}{\partial x^m} y$ derivate, replacing this by an equivalent expression, in accordance with (8.4), upon transformation we get (8.6)

$$\frac{\partial^{2+m}}{\partial x^{2+m}} y = \left( \sum_{p=1}^{m-1} \left( \right.\right.$$

$$\left( \left( \frac{d^2}{dx^2} b_p(x) \right) + \left( \frac{d}{dx} b_{1+p}(x) \right) + \left( \frac{d}{dx} b_1(x) \right) b_p(x) + b_1(x) \left( \frac{d}{dx} b_p(x) \right) \right) \left( \frac{\partial^{m-p}}{\partial x^{m-p}} y \right)$$

$$+ \left( \left( \frac{d}{dx} b_p(x) \right) + b_{1+p}(x) + b_1(x)\, b_p(x) \right) \left( \frac{\partial^{1+m-p}}{\partial x^{1+m-p}} y \right) \left.\left.\right)\right)$$

$$+ \left( \left( \frac{d^2}{dx^2} b_m(x) \right) + \left( \frac{d}{dx} b_1(x) \right) b_m(x) + b_1(x) \left( \frac{d}{dx} b_m(x) \right) \right) y$$

$$+ \left( \left( \frac{d}{dx} b_m(x) \right) + b_1(x)\, b_m(x) \right) \left( \frac{\partial}{\partial x} y \right)$$

$$+ \left( \left( \frac{d^2}{dx^2} b_1(x) \right) + \left( \frac{d}{dx} b_2(x) \right) + \left( \frac{d}{dx} b_1(x) \right) b_1(x) + b_1(x) \left( \frac{d}{dx} b_1(x) \right) \right) \left( \frac{\partial^{m-1}}{\partial x^{m-1}} y \right)$$

$$+ \left( \left( \frac{d}{dx} b_1(x) \right) + b_2(x) + b_1(x)^2 \right) \left( \sum_{p=1}^{m} b_p(x) \left( \frac{\partial^{m-p}}{\partial x^{m-p}} y \right) \right)$$

Thus, on the right side we get an expression, which is structurally identical to the equality (8.5).

where $V_p(k)$ are functional coefficients derived from the desired function. With $k = 1$ these are, respectively, the coefficients defined by (8.5), and with $k = 2$ these are the coefficients defined by (8.6).

Proceeding the same way, we eventually obtain the following linear ODE of $m + n$ order:

$$\frac{\partial^{m+n}}{\partial x^{m+n}} y = \sum_{p=1}^{m} \alpha_p(n) \left( \frac{\partial^{m-p}}{\partial x^{m-p}} y \right) \tag{8.7}$$

where $\alpha_p(n)$ are the functions defined as described above, and *n* is a natural number that specifies the type of $\alpha_p(n)$ and a new order of the ODE (8.7).

**Definition 8.2. The ODE (8.7) is called the equation of *n*-image for the ODE (8.4), and $\alpha_p(n)$ functions are called coefficients of *n*-image equation.**

Further on, we shall often call **the equation of *n*-image (8.7) for ODE (8.4) in a simpler way: *n*-image of ODE (8.4).**

## 8.4 Properties of *n*-image equation.

**Property 1: The equation of *n*-image (8.7) is a differential equation of $n + m$ order and contains all solutions for the initial ODE (8.4).**

*This property is obvious, because the equation of* ***n****-image* (8.7) *is derived, provided that the initial equation* (11.4) *is the basis for a successive series of differentiations of the right and left parts of this equation, taking* *into account its initial value* (8.4). *Thus, the new equation* (8.7) *is, essentially, the initial equation differentiated for **n** times with the replacement of all derivatives according to* (8.4).

**Property 2:** The coefficients of ***n***-image $\alpha_p(n)$ are defined by the initial values:

, (8.8)

*Indeed, assuming that* $n = 0$ *in* (8.7) *we get:*

$$\frac{\partial^m}{\partial x^m} y = \sum_{p=1}^{m} \alpha_p(0) \left( \frac{\partial^{m-p}}{\partial x^{m-p}} y \right)$$

*Comparing this equation with the initial equation* (8.4) *we obtain, from the condition of equality of the coefficients with identical* $\frac{\partial^{m-p}}{\partial x^{m-p}} y$ *derivatives,* the equality (8.8).

**Property 3: Coefficients of *n*-image (8.7) are related by recurrence equations**:

$$\alpha_p(n+1) = \left( \frac{\partial}{\partial x} \alpha_p(n) \right) + \alpha_{1+p}(n) + \alpha_1(n)\, b_p(x) \qquad (8.9)$$

$$\alpha_m(n+1) = \left( \frac{\partial}{\partial x} \alpha_m(n) \right) + b_m(x)\, \alpha_1(n) \qquad (8.10)$$

**Proof**. *Indeed, replacing **n** parameter by* $n+1$ *in* (8.7), *we have:*

$$\frac{\partial^{m+n+1}}{\partial x^{m+n+1}} y = \sum_{p=1}^{m} \alpha_p(n+1) \left( \frac{\partial^{m-p}}{\partial x^{m-p}} y \right)$$

*On the other hand, upon differentiation of* (8.7) *with consideration of* (8.4) *we have:*

$$\frac{\partial^{m+n+1}}{\partial x^{m+n+1}} y = \left( \sum_{p=1}^{m-1} \left( \left( \frac{\partial}{\partial x} \alpha_p(n) \right) + \alpha_{1+p}(n) + \alpha_1(n)\, b_p(x) \right) \left( \frac{\partial^{m-p}}{\partial x^{m-p}} y \right) \right)$$
$$+ \left( \left( \frac{\partial}{\partial x} \alpha_m(n) \right) + b_m(x)\, \alpha_1(n) \right) y$$

*As soon as the left parts of these equalities are equal, from the equality of their right parts we obtain* (8.9), (8.10).

## 9. Definition of the construction of general solution for homogeneous ODEs.

An equation of *n*-image is completely defined upon the definition of formulas for all $\alpha_p(n)$ coefficients. So let's analyze the recurrence relations (8.9) and (8.10). Obviously, with consideration of the initial conditions (8.8), recurrence relations (8.9) and (8.10) should allow to define, for any given natural values of *n*, all required values of $\alpha_p(n)$ coefficients. However, this does not make much sense, because it gives no new scientific results for the definition of the algorithm to solve differential equations. New results, which would be useful to answer the question about the structure of the general solution for a homogeneous ODE, can be obtained only with the use of unconventional approach. In this case this means that a system of recurrence equations is solved for the case of negative parameter values: .

**Definition 12.1. The  functions, in the case where $n$ parameter assumes the values , are called special functions**.

**9.1. The case of $n = -1$** . A system of recurrence equations (8.9), (8.10) takes the following form:

$$\alpha_p(0) = \left(\frac{d}{dx}\alpha_p(-1)\right) + \alpha_{1+p}(-1) + \alpha_1(-1)\, b_p(x) \qquad (9.1)$$

$$\alpha_m(0) = \left(\frac{d}{dx}\alpha_m(-1)\right) + b_m(x)\,\alpha_1(-1) \qquad (9.2)$$

It is obvious that $\alpha_p(-1)$ are functions of the $\boldsymbol{x}$ variable with specific, i.e. negative value of parameter $\boldsymbol{n}$.

As follows from (9.1):

$$\alpha_{1+p}(-1) = \alpha_p(0) - \left(\frac{d}{dx}\alpha_p(-1)\right) - \alpha_1(-1)\, b_p(x)$$

or, with consideration of initial conditions (8.8) :

$$\alpha_{1+p}(-1) = -\left(\frac{d}{dx}\alpha_p(-1)\right) - (-1 + \alpha_1(-1))\, b_p(x) \qquad (9.3)$$

Thus,

$$\alpha_m(-1) = -\left(\frac{d}{dx}\alpha_{m-1}(-1)\right) - (-1 + \alpha_1(-1))\, b_{m-1}(x) \qquad (9.4)$$

$$\alpha_{m-1}(-1) = -\left(\frac{d}{dx}\alpha_{m-2}(-1)\right) - (-1 + \alpha_1(-1))\, b_{m-2}(x) \qquad (9.5)$$

$$\alpha_{m-2}(-1) = -\left(\frac{d}{dx}\alpha_{m-3}(-1)\right) - (-1 + \alpha_1(-1))\, b_{m-3}(x) \qquad (9.6)$$

$$\alpha_{m-3}(-1) = -\left(\frac{d}{dx}\alpha_{m-4}(-1)\right) - (-1 + \alpha_1(-1))\, b_{m-4}(x) \qquad (9.7)$$

$$\alpha_{m-4}(-1) = -\left(\frac{d}{dx}\alpha_{m-5}(-1)\right) - (-1 + \alpha_1(-1))\, b_{m-5}(x) \qquad (9.8)$$

and so on.

Therefore, if we introduce a $Z_k[\ ]$ operator, operating on the $\mathrm{f}(x)$ function according to the following rule:

$$Z_k[\mathrm{f}(x)] = -\left(\frac{d}{dx}\mathrm{f}(x)\right) - (-1 + \alpha_1(-1))\, b_k(x) \qquad (9.9)$$

(9.3) can be written down as follows:

(9.10)

Let's note that $Z_k[\ ]$ has the following property:

$$Z_k[\mathrm{f}(x)] = Z_k[\mathrm{f}(x) + C] \tag{9.11}$$

where $C$ is an arbitrary constant.

Thus, with consideration of (9.9), (9.4)-(9.8) take the following form:

$$\alpha_{m-3}(-1) = Z_{m-4}[Z_{m-5}[\alpha_{m-5}(-1)]]$$

$$\alpha_{m-2}(-1) = Z_{m-3}[Z_{m-4}[Z_{m-5}[\alpha_{m-5}(-1)]]]$$

$$\alpha_{m-1}(-1) = Z_{m-2}[Z_{m-3}[Z_{m-4}[Z_{m-5}[\alpha_{m-5}(-1)]]]]$$

$$\alpha_{m}(-1) = Z_{m-1}[Z_{m-2}[Z_{m-3}[Z_{m-4}[Z_{m-5}[\alpha_{m-5}(-1)]]]]]$$

Therefore, in a general case, we can write down:

$$\alpha_{m}(-1) = Z_{m-1}[Z_{m-2}[Z_{m-3}[Z_{m-4}[Z_{m-p}[\alpha_{m-p}(-1)]]]]]$$

With $p = m - 1$ we get:

$$\alpha_{m}(-1) = Z_{m-1}[Z_{m-2}[Z_{m-3}[Z_{m-4}[\ ][\ ][\ ][Z_4[Z_3[Z_2[Z_1[\alpha_1(-1)]]]]]]]]$$

Substituting this value to (9.2), taking into account the initial condition (8.8), we have:

$$b_m(x) = \left(\frac{\partial}{\partial x}(Z_{m-1}[Z_{m-2}[Z_{m-3}[Z_{m-4}[\ ][\ ][\ ][Z_4[Z_3[Z_2[Z_1[\alpha_1(-1)]]]]]]]])\right) + b_m(x)\,\alpha_1(-1)$$

Following from that:

$$-\left(\frac{\partial}{\partial x}(Z_{m-1}[Z_{m-2}[Z_{m-3}[Z_{m-4}[\ ][\ ][\ ][Z_4[Z_3[Z_2[Z_1[\alpha_1(-1)]]]]]]]])\right) - (-1 + \alpha_1(-1))\,b_m(x) = 0$$

In operator form, with consideration of (9.9), the eventual equation is as follows:

$$Z_m[Z_{m-1}[Z_{m-2}[Z_{m-3}[Z_{m-4}[\ ][\ ][\ ][Z_4[Z_3[Z_2[Z_1[\alpha_1(-1)]]]]]]]]] = 0 \tag{9.12}$$

Let's use the property of $Z_k[\ ]$ operator (9.11). In such case the equation (9.12) can be represented as follows:

$$Z_m[Z_{m-1}[Z_{m-2}[Z_{m-3}[Z_{m-4}[\ ][\ ][\ ][Z_4[Z_3[Z_2[Z_1[-1 + \alpha_1(-1)]]]]]]]]] = 0 \tag{9.13}$$

Let's introduce a new variable:

(9.14)

In this case, the equality (9.13) takes the following form:

$$Z_m[Z_{m-1}[Z_{m-2}[Z_{m-3}[Z_{m-4}[\ ][\ ][\ ][Z_4[Z_3[Z_2[Z_1[\mathrm{Y}(x)]]]]]]]]]]=0 \qquad (9.15)$$

As soon as the operation on the $\mathrm{f}(x)$ function acts in accordance with (9.9), with consideration of (9.14) it is equal to:

$$Z_k[\mathrm{f}(x)]=-\left(\frac{d}{dx}\mathrm{f}(x)\right)-\mathrm{Y}(x)\,b_k(x) \qquad (9.16)$$

Therefore, the $\mathrm{Y}(x)$ function is the solution of the differential equation (9.15). Let's define this ODE. As soon as, with consideration of (9.16):

$$Z_1[\mathrm{Y}(x)]=-\left(\frac{d}{dx}\mathrm{Y}(x)\right)-\mathrm{Y}(x)\,b_1(x)$$

$$Z_2[Z_1[\mathrm{Y}(x)]]=-\left(\frac{d}{dx}\left(-\left(\frac{d}{dx}\mathrm{Y}(x)\right)-\mathrm{Y}(x)\,b_1(x)\right)\right)-\mathrm{Y}(x)\,b_2(x)=$$
$$\left(\frac{d^2}{dx^2}\mathrm{Y}(x)\right)+\left(\frac{d}{dx}(\mathrm{Y}(x)\,b_1(x))\right)-\mathrm{Y}(x)\,b_2(x)$$

$$Z_3[Z_2[Z_1[\mathrm{Y}(x)]]]=Z_3\left[\left(\frac{d^2}{dx^2}\mathrm{Y}(x)\right)+\left(\frac{d}{dx}(\mathrm{Y}(x)\,b_1(x))\right)-\mathrm{Y}(x)\,b_2(x)\right]=$$
$$-\left(\frac{d}{dx}\left(\left(\frac{d^2}{dx^2}\mathrm{Y}(x)\right)+\left(\frac{d}{dx}(\mathrm{Y}(x)\,b_1(x))\right)-\mathrm{Y}(x)\,b_2(x)\right)\right)-\mathrm{Y}(x)\,b_3(x)=$$
$$-\left(\frac{d^3}{dx^3}\mathrm{Y}(x)\right)-\left(\frac{d^2}{dx^2}(\mathrm{Y}(x)\,b_1(x))\right)+\left(\frac{d}{dx}(\mathrm{Y}(x)\,b_2(x))\right)-\mathrm{Y}(x)\,b_3(x)$$

$$Z_4[Z_3[Z_2[Z_1[\mathrm{Y}(x)]]]]=$$
$$Z_4\left[-\left(\frac{d^3}{dx^3}\mathrm{Y}(x)\right)-\left(\frac{d^2}{dx^2}(\mathrm{Y}(x)\,b_1(x))\right)+\left(\frac{d}{dx}(\mathrm{Y}(x)\,b_2(x))\right)-\mathrm{Y}(x)\,b_3(x)\right]=$$
$$\left(\frac{d^4}{dx^4}\mathrm{Y}(x)\right)+\left(\frac{d^3}{dx^3}(\mathrm{Y}(x)\,b_1(x))\right)-\left(\frac{d^2}{dx^2}(\mathrm{Y}(x)\,b_2(x))\right)+\left(\frac{d}{dx}(\mathrm{Y}(x)\,b_3(x))\right)$$
$$-\mathrm{Y}(x)\,b_4(x)$$

and so on, in the expanded view of the ODE (9.15) it equals: (9.17)

$$(-1)^m\left(\frac{d^m}{dx^m}\mathrm{Y}(x)\right)=(-1)^{(m-1)}\left(\frac{d^{m-1}}{dx^{m-1}}(\mathrm{Y}(x)\,b_1(x))\right)$$
$$+(-1)^{(m-2)}\left(\frac{d^{m-2}}{dx^{m-2}}(\mathrm{Y}(x)\,b_2(x))\right)+(-1)^{(m-3)}\left(\frac{d^{m-3}}{dx^{m-3}}(\mathrm{Y}(x)\,b_3(x))\right)..$$
$$-\left(\frac{d^3}{dx^3}(\mathrm{Y}(x)\,b_{m-3}(x))\right)+\left(\frac{d^2}{dx^2}(\mathrm{Y}(x)\,b_{m-2}(x))\right)-\left(\frac{d}{dx}(\mathrm{Y}(x)\,b_{m-1}(x))\right)$$
$$+\mathrm{Y}(x)\,b_m(x)$$

**However, this ODE is adjoined with (8.4).** Thus, their solutions are linked by Green formula [Mathematical Encyclopedia. Vol. 5. Publishing House "Soviet Encyclopedia", Moscow, 1979]:

$$\sum_{k=1}^{m}\left(\sum_{j=0}^{k-j}(-1)^j\left(\frac{d^j}{dx^j}(\mathrm{Y}(x)\,b_{m-k}(x))\right)\left(\frac{\partial^{m-k-j}}{\partial x^{m-k-j}}y\right)\right)=const$$

With consideration of (9.14), this formula assumes the following view:

$$\sum_{k=1}^{m}\left(\sum_{j=0}^{k-j}(-1)^j\left(\frac{d^j}{dx^j}((-1+\alpha_1(-1))\,b_{m-k}(x))\right)\left(\frac{\partial^{m-k-j}}{\partial x^{m-k-j}}y_1\right)\right)=const \qquad (9.18)$$

**Thus, we have established a mathematical functional relationship between the $\alpha_1(-1)$ function and a certain exact partial solution of the ODE (8.4) $y=y_1$. In order to establish this relationship through specific formulas, we need to reduce the adjoined ODE (9.17) to the corresponding ODE (8.4).**
As soon as, in accordance with Leibniz formula:

$$\frac{d^{m-k}}{dx^{m-k}}(\mathrm{Y}(x)\,b_k(x))=\sum_{i=0}^{m-k}\mathrm{C}(m-k,i)\left(\frac{d^i}{dx^i}\mathrm{Y}(x)\right)\left(\frac{d^{m-k-i}}{dx^{m-k-i}}b_k(x)\right)$$

the ODE (9.17) assumes the following form:

$$\begin{aligned}(-1)^m\left(\frac{d^m}{dx^m}\mathrm{Y}(x)\right)&=(-1)^{(m-1)}\left(\sum_{i=0}^{m-1}\mathrm{C}(m-1,i)\left(\frac{d^i}{dx^i}\mathrm{Y}(x)\right)\left(\frac{d^{m-1-i}}{dx^{m-1-i}}b_1(x)\right)\right)\\
&+(-1)^{(m-2)}\left(\sum_{i=0}^{m-2}\mathrm{C}(m-2,i)\left(\frac{d^i}{dx^i}\mathrm{Y}(x)\right)\left(\frac{d^{m-2-i}}{dx^{m-2-i}}b_2(x)\right)\right)\\
&+(-1)^{(m-3)}\left(\sum_{i=0}^{m-3}\mathrm{C}(m-3,i)\left(\frac{d^i}{dx^i}\mathrm{Y}(x)\right)\left(\frac{d^{m-3-i}}{dx^{m-3-i}}b_3(x)\right)\right)..-\Bigg(\\
&\left(\frac{d^3}{dx^3}\mathrm{Y}(x)\right)b_{m-3}(x)+3\left(\frac{d^2}{dx^2}\mathrm{Y}(x)\right)\left(\frac{d}{dx}b_{m-3}(x)\right)+3\left(\frac{d}{dx}\mathrm{Y}(x)\right)\left(\frac{d^2}{dx^2}b_{m-3}(x)\right)\\
&+\mathrm{Y}(x)\left(\frac{d^3}{dx^3}b_{m-3}(x)\right)\Bigg)+\left(\frac{d^2}{dx^2}\mathrm{Y}(x)\right)b_{m-2}(x)+2\left(\frac{d}{dx}\mathrm{Y}(x)\right)\left(\frac{d}{dx}b_{m-2}(x)\right)\\
&+\mathrm{Y}(x)\left(\frac{d^2}{dx^2}b_{m-2}(x)\right)-\left(\left(\frac{d}{dx}\mathrm{Y}(x)\right)b_{m-1}(x)+\mathrm{Y}(x)\left(\frac{d}{dx}b_{m-1}(x)\right)\right)+\mathrm{Y}(x)\,b_m(x)\end{aligned}$$

Upon successive transformations, we eventually obtain:

$$(-1)^m\left(\frac{d^m}{dx^m}\mathrm{Y}(x)\right)=\sum_{k=1}^{m}\left(\sum_{i=1}^{k}(-1)^{(-i)}\mathrm{C}(m-i,m-k)\left(\frac{d^{k-i}}{dx^{k-i}}b_i(x)\right)\right)\left(\frac{d^{m-k}}{dx^{m-k}}\mathrm{Y}(x)\right)$$

This ODE can be transformed as follows:

$$\frac{d^m}{dx^m}\mathrm{Y}(x)=\sum_{k=1}^{m}\left(\sum_{i=1}^{k}(-1)^{(-m-i)}\mathrm{C}(m-i,m-k)\left(\frac{d^{k-i}}{dx^{k-i}}b_i(x)\right)\right)\left(\frac{d^{m-k}}{dx^{m-k}}\mathrm{Y}(x)\right) \qquad (9.19)$$

Thus, introducing the $a_k(x)$ functions in accordance with the following formulas:

$$a_k(x)=\sum_{i=1}^{k}(-1)^{(-m-i)}\mathrm{C}(m-i,m-k)\left(\frac{d^{k-i}}{dx^{k-i}}b_i(x)\right) \qquad (9.20)$$

we obtain a new representation of the ODE (9.19):

$$\frac{d^m}{dx^m}\mathrm{Y}(x) = \sum_{k=1}^{m} a_k(x)\left(\frac{d^{m-k}}{dx^{m-k}}\mathrm{Y}(x)\right) \tag{9.21}$$

**This ODE is fully identical in its structure to the ODE (9.4) and, in that sense, these equations are equivalent. Hence the ODE (9.21) can be taken as an initial one, i.e. the task can be stipulated as follows: for the given functions we need to define partial solution of $\mathrm{Y}(x)$. In this case, such solution is defined by the formula (9.14).**

However, in order to determine $\mathrm{Y}(x)$ we need to know the $\alpha_1(-1)$ function determined by the values of functions, which depend upon the given functions. Therefore the task is to define the functions through the functions with the use of (9.20). This task is solved rather simple, because the ODEs (8.4) and (9.21) are adjoined. We simply need to interchange $b_i(x)$ and $a_k(x)$, and reverse the sign of $-m$ parameter. This way we can write down the following general formula:

$$b_k(x) = \sum_{i=1}^{k} (-1)^{(m-i)}\,\mathrm{C}(m-i, m-k)\left(\frac{d^{k-i}}{dx^{k-i}}a_i(x)\right) \tag{9.22}$$

By way of a direct check (substitution to (9.20)), we can verify the validity of the obtained formula.

**Thus, we have constructed the calculation of the first partial solution for the ODE (9.21) with the specified $a_k(x)$, $k = 1\,..\,m$ coefficients.**

**The calculation of this solution is as follows:**

**1) In accordance with the task, the solution of is obtained in the form determined by (9.14), that is**

**2) For the calculation of $\alpha_1(-1)$ function we need to define the coeffieicents for the adjoined ODE (9.4) with the use of (9.22). This task shall be sloved later.**

**9.2. The case of $n = -m$**. Let's prove the following **theorem.**

**If we define the special functions with the use of following equations**:

(9.23)

$$\alpha_1(-j) = \sum_{i=1}^{j} \frac{Y_{j-i+1}(x)\,x^{(i-1)}}{(i-1)!} \qquad j = 2\,..\,m \tag{9.24}$$

**where the functions are partial solutions of the ODE :**

$$\frac{d^m}{dx^m}\mathrm{Y}(x) = \sum_{k=1}^{m}\left(\sum_{i=1}^{k} (-1)^{(-m-i)}\,\mathrm{C}(m-i, m-k)\left(\frac{d^{k-i}}{dx^{k-i}}b_i(x)\right)\right)\left(\frac{d^{m-k}}{dx^{m-k}}\mathrm{Y}(x)\right) \tag{9.25}$$

**and special functions are defined from the coefficients of *n*-image of the ODE**:

$$\frac{d^m}{dx^m}\mathrm{R}(x) = \sum_{k=1}^{m} b_k(x)\left(\frac{d^{m-k}}{dx^{m-k}}\mathrm{R}(x)\right) \tag{9.26}$$

**then the system of recurrence equations**

$$\alpha_p(n+1) = \left(\frac{\partial}{\partial x}\alpha_p(n)\right) + \alpha_{1+p}(n) + \alpha_1(n)\, b_p(x) \tag{9.27}$$

$$\alpha_m(n+1) = \left(\frac{\partial}{\partial x}\alpha_m(n)\right) + b_m(x)\,\alpha_1(n) \tag{9.28}$$

**is solved identically.**

As a preliminary operation, let's reduce the specified system of recurrence equations (8.9), (8.10) to such form, where it is defined by a single recurrence equation of $\alpha_1(n)$ function. Indeed, if we assume $p = m - 1$ in (8.9), we get:

$$\alpha_{m-1}(n+1) = \left(\frac{\partial}{\partial x}\alpha_{m-1}(n)\right) + \alpha_m(n) + \alpha_1(n)\, b_{m-1}(x)$$

Let's express $\alpha_m(n)$ and $\alpha_m(n+1)$ from here:

$$\alpha_m(n) = \alpha_{m-1}(n+1) - \left(\frac{\partial}{\partial x}\alpha_{m-1}(n)\right) - \alpha_1(n)\, b_{m-1}(x)$$

$$\alpha_m(n+1) = \alpha_{m-1}(n+2) - \left(\frac{\partial}{\partial x}\alpha_{m-1}(n+1)\right) - \alpha_1(n+1)\, b_{m-1}(x)$$

Substituting these to (8.10), upon transformation we receive:

$$\alpha_{m-1}(n+2) - 2\left(\frac{\partial}{\partial x}\alpha_{m-1}(n+1)\right) - \alpha_1(n+1)\, b_{m-1}(x) + \left(\frac{\partial^2}{\partial x^2}\alpha_{m-1}(n)\right) - b_m(x)\,\alpha_1(n) + \left(\frac{\partial}{\partial x}(\alpha_1(n)\, b_{m-1}(x))\right) = 0 \tag{9.29}$$

As we can see, this recurrence equation does not contain the $\alpha_m(n)$ functions.

Again, assuming that $p = m - 2$ in (8.9), we get:

$$\alpha_{m-2}(n+1) = \left(\frac{\partial}{\partial x}\alpha_{m-2}(n)\right) + \alpha_{m-1}(n) + \alpha_1(n)\, b_{m-2}(x)$$

From here, let's express $\alpha_{m-1}(n)$, and

$$\alpha_{m-1}(n) = \alpha_{m-2}(n+1) - \left(\frac{\partial}{\partial x}\alpha_{m-2}(n)\right) - \alpha_1(n)\, b_{m-2}(x)$$

$$\alpha_{m-1}(n+1) = \alpha_{m-2}(n+2) - \left(\frac{\partial}{\partial x}\alpha_{m-2}(n+1)\right) - \alpha_1(n+1)\, b_{m-2}(x)$$

$$\alpha_{m-1}(n+2) = \alpha_{m-2}(n+3) - \left(\frac{\partial}{\partial x}\alpha_{m-2}(n+2)\right) - \alpha_1(n+2)\, b_{m-2}(x)$$

Substituting these equalities to (9.23), upon transformations we obtain:

$$\alpha_{m-2}(n+3) - 3\left(\frac{\partial}{\partial x}\alpha_{m-2}(n+2)\right) - \alpha_1(n+2)\,b_{m-2}(x) + 3\left(\frac{\partial^2}{\partial x^2}\alpha_{m-2}(n+1)\right)$$
$$- \alpha_1(n+1)\,b_{m-1}(x) - \left(\frac{\partial^3}{\partial x^3}\alpha_{m-2}(n)\right) - b_m(x)\,\alpha_1(n) + \left(\frac{\partial}{\partial x}(\alpha_1(n)\,b_{m-1}(x))\right)$$
$$+ 2\left(\frac{\partial}{\partial x}(\alpha_1(n+1)\,b_{m-2}(x))\right) - \left(\frac{\partial^2}{\partial x^2}(\alpha_1(n)\,b_{m-2}(x))\right) = 0$$

Procceding in a similar manner, we eventaually obtain at the $s$ step the following recurrence ODE:

$$\left(\sum_{i=0}^{s} (-1)^i \operatorname{C}(s+1, i)\left(\frac{\partial^{s-i+1}}{\partial x^{s-i+1}}\alpha_{m-s}(n+i)\right)\right)$$
$$+ \left(\sum_{j=0}^{s}\left(\sum_{i=0}^{s-j}(-1)^{(s+i-1)}\left(\frac{\partial^{s-i}}{\partial x^{s-i}}(\alpha_1(n+j)\,b_{m-s+i}(x))\right)\right)\right) = 0$$

Assuming here $s = m - 1$, we eventually get:

$$\left(\sum_{i=0}^{m} (-1)^{(i+m)} \operatorname{C}(m, i)\left(\frac{\partial^{m-i}}{\partial x^{m-i}}\alpha_1(n+i)\right)\right)$$
$$+ \left(\sum_{j=0}^{m}\left(\sum_{i=0}^{m-1-j}(-1)^{(j+i+m-2)}\operatorname{C}(m-1-i, j)\left(\frac{\partial^{m-1-i-j}}{\partial x^{m-1-i-j}}(\alpha_1(n+j)\,b_{1+i}(x))\right)\right)\right) = 0$$

(9.30)

Thus, we have received a certain linear ODE of the $\alpha_1(n)$ function. If we solve this, this means that we have solved the initial system of recurrence equations (9.27), (9.28).
Let’s consider individual cases, because they clearly demonstrate the correctness of the general theory.

**The variant of** $m = 2$ : Let's separate out the components from special summands, that is, the summands of the following type: $\alpha_1(-1)$, $\alpha_1(-2)$ (9.31)

$$(-1)^m\left(\frac{d^m}{dx^m}\alpha_1(-2)\right) + \left(\sum_{i=0}^{m-1}(-1)^{(-2+i+m)}\left(\frac{d^{m-1-i}}{dx^{m-1-i}}(\alpha_1(-2)\,b_{1+i}(x))\right)\right)$$
$$+ (-1)^{(1+m)}\,m\left(\frac{d^{m-1}}{dx^{m-1}}\alpha_1(-1)\right)$$
$$+ \left(\sum_{i=0}^{m-2}(-1)^{(-1+i+m)}(m-1-i)\left(\frac{d^{m-2-i}}{dx^{m-2-i}}(\alpha_1(-1)\,b_{1+i}(x))\right)\right)$$
$$+ \left(\sum_{j=2}^{m}\left(\sum_{i=0}^{m-1-j}(-1)^{(j+i+m-2)}\operatorname{C}(m-1-i, j)\left(\frac{\partial^{m-1-i-j}}{\partial x^{m-1-i-j}}(\alpha_1(-2+j)\,b_{1+i}(x))\right)\right)\right)$$
$$+ \left(\sum_{i=2}^{m}(-1)^{(i+m)}\operatorname{C}(m, i)\left(\frac{\partial^{m-i}}{\partial x^{m-i}}\alpha_1(-2+i)\right)\right) = 0$$

In this case we have received a certian linear ODE of $m$ order for the special $\alpha_1(-2)$ function and of $m - 1$ order for the special $\alpha_1(-1)$ function.

Let’s consider some specific cases because with their aid we can easily define the desired construction for the general solution.

Equation (9.24) takes the following form:

$$\left(\frac{d^2}{dx^2}\alpha_1(-2)\right) - 2\left(\frac{d}{dx}\alpha_1(-1)\right) + \left(\frac{d}{dx}(\alpha_1(-2)\,b_1(x))\right) - \alpha_1(-2)\,b_2(x) - \alpha_1(-1)\,b_1(x) + \alpha_1(0) = 0$$

As soon as, according to the property 1, , this equation assumes the following form:

$$\left(\frac{d^2}{dx^2}\alpha_1(-2)\right) - 2\left(\frac{d}{dx}\alpha_1(-1)\right) + \left(\frac{d}{dx}(\alpha_1(-2)\,b_1(x))\right) - \alpha_1(-2)\,b_2(x) - \alpha_1(-1)\,b_1(x) + b_1(x) = 0 \qquad (9.32)$$

Previously it was demonstrated that the special $\alpha_1(-1)$ function defines the construction of the first partial solution for the adjoined ODE (12.21), where the $a_k(x)$ coefficients are defined by (9.20).

In accordance with (9.14):

(9.33)

where (in this case with $m = 2$ ) the $Y_1(x)$ function is the solution of the following ODE:

$$\left(\frac{d^2}{dx^2}Y_1(x)\right) + \left(\frac{d}{dx}Y_1(x)\right)b_1(x) + \left(\left(\frac{d}{dx}b_1(x)\right) - b_2(x)\right)Y_1(x) = 0 \qquad (9.34)$$

Analysis of the ODE (9.25) demonstrates that, if we assume:

$$\alpha_1(-2) = Y_1(x)\,x + Y_2(x) \qquad (9.35)$$

where the $Y_2(x)$ function is the solution of the following ODE.

$$\left(\frac{d^2}{dx^2}Y_2(x)\right) + \left(\frac{d}{dx}Y_2(x)\right)b_1(x) + \left(\left(\frac{d}{dx}b_1(x)\right) - b_2(x)\right)Y_2(x) = 0 \qquad (9.36)$$

then the ODE (9.25) takes the following form:

$$\left(\left(\frac{d^2}{dx^2}Y_1(x)\right) + \left(\frac{d}{dx}Y_1(x)\right)b_1(x) + \left(\left(\frac{d}{dx}b_1(x)\right) - b_2(x)\right)Y_1(x)\right)x + \left(\frac{d^2}{dx^2}Y_2(x)\right) + \left(\frac{d}{dx}Y_2(x)\right)b_1(x) + \left(\left(\frac{d}{dx}b_1(x)\right) - b_2(x)\right)Y_2(x) = 0$$

By virtue of (9.27) and (9.29), this ODE is identically solved. Therefore, the above means that two partial solutions, , of the following ODE

$$\frac{d^2}{dx^2}\mathrm{Y}(x) = a_1(x)\left(\frac{d}{dx}\mathrm{Y}(x)\right) + a_2(x)\,\mathrm{Y}(x) \qquad (9.37)$$

where

(9.38)

$$a_2(x) = -\left(\frac{d}{dx} b_1(x)\right) + b_2(x) \tag{9.39}$$

are defined by the following formulas:

(9.40)

$$Y_2(x) = \alpha_1(-2) + x\,(1 - \alpha_1(-1)) \tag{9.41}$$

Here the special functions are derivatives from the coefficients of $n$- image ODE:

$$\frac{d^2}{dx^2} \mathrm{Y}(x) = b_1(x)\left(\frac{d}{dx} \mathrm{Y}(x)\right) + b_2(x)\, \mathrm{Y}(x) \tag{9.42}$$

**The variant of** $m = 3$. (9.54) takes the following form:

$$-\left(\frac{d^3}{dx^3} \alpha_1(-3)\right) + 3\left(\frac{d^2}{dx^2} \alpha_1(-2)\right) - 3\left(\frac{d}{dx} \alpha_1(-1)\right) - \left(\frac{d^2}{dx^2} (\alpha_1(-3)\, b_1(x))\right)$$
$$+ \left(\frac{d}{dx} (\alpha_1(-3)\, b_2(x))\right) - \alpha_1(-3)\, b_3(x) + 2\left(\frac{d}{dx} (\alpha_1(-2)\, b_1(x))\right) - \alpha_1(-2)\, b_2(x)$$
$$- \alpha_1(-1)\, b_1(x) + \alpha_1(0) = 0$$

With consideration of the initial condition: , this equation equals: (9.43)

$$-\left(\frac{d^3}{dx^3} \alpha_1(-3)\right) + 3\left(\frac{d^2}{dx^2} \alpha_1(-2)\right) - 3\left(\frac{d}{dx} \alpha_1(-1)\right) - \left(\frac{d^2}{dx^2} (\alpha_1(-3)\, b_1(x))\right)$$
$$+ \left(\frac{d}{dx} (\alpha_1(-3)\, b_2(x))\right) - \alpha_1(-3)\, b_3(x) + 2\left(\frac{d}{dx} (\alpha_1(-2)\, b_1(x))\right) - \alpha_1(-2)\, b_2(x)$$
$$- \alpha_1(-1)\, b_1(x) + b_1(x) = 0$$

As an operational hypothesis, let's assume:

(9.44)

$$\alpha_1(-2) = Y_1(x)\, x + Y_2(x) \tag{9.45}$$

$$\alpha_1(-3) = \frac{Y_1(x)\, x^2}{2} + Y_2(x)\, x + Y_3(x) \tag{9.46}$$

where the functions are the solution of the following ODE:

$$-\left(\frac{d^3}{dx^3} \mathrm{Y}(x)\right) - \left(\frac{d^2}{dx^2} \mathrm{Y}(x)\right) b_1(x) + \left(-2\left(\frac{d}{dx} b_1(x)\right) + b_2(x)\right)\left(\frac{d}{dx} \mathrm{Y}(x)\right)$$
$$+ \left(\left(\frac{d}{dx} b_2(x)\right) - \left(\frac{d^2}{dx^2} b_1(x)\right) - b_3(x)\right) \mathrm{Y}(x) = 0 \tag{9.47}$$

Substituting (9.56), (9.57) (9.58) to (9.55), we get:

$$\left(-\frac{1}{2}\left(\frac{d^3}{dx^3}Y_1(x)\right)-\frac{1}{2}\left(\frac{d^2}{dx^2}Y_1(x)\right)b_1(x)+\left(-\left(\frac{d}{dx}b_1(x)\right)+\frac{1}{2}b_2(x)\right)\left(\frac{d}{dx}Y_1(x)\right)+\left(\frac{1}{2}\left(\frac{d}{dx}b_2(x)\right)-\frac{1}{2}\left(\frac{d^2}{dx^2}b_1(x)\right)-\frac{1}{2}b_3(x)\right)Y_1(x)\right)x^2+\left(-\left(\frac{d^3}{dx^3}Y_2(x)\right)-\left(\frac{d^2}{dx^2}Y_2(x)\right)b_1(x)+\left(b_2(x)-2\left(\frac{d}{dx}b_1(x)\right)\right)\left(\frac{d}{dx}Y_2(x)\right)+\left(-\left(\frac{d^2}{dx^2}b_1(x)\right)-b_3(x)+\left(\frac{d}{dx}b_2(x)\right)\right)Y_2(x)\right)x-\left(\frac{d^3}{dx^3}Y_3(x)\right)-\left(\frac{d^2}{dx^2}Y_3(x)\right)b_1(x)+\left(b_2(x)-2\left(\frac{d}{dx}b_1(x)\right)\right)\left(\frac{d}{dx}Y_3(x)\right)+\left(-\left(\frac{d^2}{dx^2}b_1(x)\right)-b_3(x)+\left(\frac{d}{dx}b_2(x)\right)\right)Y_3(x)=0$$

As soon as the functions are the solution of (9.59), this equality is solved identically.

From the system of equations (9.56)-(9.58) we obtain:

(9.48)

$$Y_2(x)=(-\alpha_1(-1)+1)x+\alpha_1(-2) \quad (9.49)$$

$$Y_3(x)=\left(\frac{\alpha_1(-1)}{2}-\frac{1}{2}\right)x^2-\alpha_1(-2)x+\alpha_1(-3) \quad (9.50)$$

Thus, we have demonstrated that the functions defined by the equalities (9.60)-(9.62) are partial solutions of the following equation:

$$\frac{d^3}{dx^3}\mathrm{Y}(x)=a_1(x)\left(\frac{d^2}{dx^2}\mathrm{Y}(x)\right)+a_2(x)\left(\frac{d}{dx}\mathrm{Y}(x)\right)+a_3(x)\,\mathrm{Y}(x)$$

where

$$a_2(x)=-2\left(\frac{d}{dx}b_1(x)\right)+b_2(x)$$

$$a_2(x)=\left(\frac{d}{dx}b_2(x)\right)-\left(\frac{d^2}{dx^2}b_1(x)\right)-b_3(x)$$

and the special functions are defined by the coefficient of $n$- image $\alpha_1(n)$.

**9.4. The variant of** $n=-4$. In the equality (9.53) let's assume $n=-4$

$$\left(\sum_{i=0}^{m} (-1)^{(i+m)} C(m,i)\left(\frac{\partial^{m-i}}{\partial x^{m-i}} \alpha_1(-4+i)\right)\right) + \left(\sum_{j=0}^{m}\left(\sum_{i=0}^{m-1-j} (-1)^{(j+i+m-2)} C(m-1-i,j)\left(\frac{\partial^{m-1-i-j}}{\partial x^{m-1-i-j}}(\alpha_1(-4+j)\, b_{1+i}(x))\right)\right)\right) = 0$$

(9.51)

Let's separate out the components from special summands, that is, the summands of the following type:

$$(-1)^m C(m,0)\left(\frac{d^m}{dx^m}\alpha_1(-4)\right) + (-1)^{(1+m)} C(m,1)\left(\frac{d^{m-1}}{dx^{m-1}}\alpha_1(-3)\right) + (-1)^{(2+m)} C(m,2)\left(\frac{d^{m-2}}{dx^{m-2}}\alpha_1(-2)\right) + (-1)^{(3+m)} C(m,3)\left(\frac{d^{m-3}}{dx^{m-3}}\alpha_1(-1)\right)$$
$$+ \left(\sum_{i=0}^{m-1} (-1)^{(-2+i+m)} C(m-1-i,0)\left(\frac{d^{m-1-i}}{dx^{m-1-i}}(\alpha_1(-4)\, b_{1+i}(x))\right)\right)$$
$$+ \left(\sum_{i=0}^{m-2} (-1)^{(-1+i+m)} C(m-1-i,1)\left(\frac{d^{m-2-i}}{dx^{m-2-i}}(\alpha_1(-3)\, b_{1+i}(x))\right)\right)$$
$$+ \left(\sum_{i=0}^{m-3} (-1)^{(i+m)} C(m-1-i,2)\left(\frac{d^{m-3-i}}{dx^{m-3-i}}(\alpha_1(-2)\, b_{1+i}(x))\right)\right)$$
$$+ \left(\sum_{i=0}^{m-4} (-1)^{(1+i+m)} C(m-1-i,3)\left(\frac{d^{m-4-i}}{dx^{m-4-i}}(\alpha_1(-1)\, b_{1+i}(x))\right)\right)$$
$$+ \left(\sum_{i=0}^{m} (-1)^{(i+m)} C(m,i)\left(\frac{\partial^{m-i}}{\partial x^{m-i}}\alpha_1(-4+i)\right)\right)$$
$$+ \left(\sum_{j=4}^{m}\left(\sum_{i=4}^{m-1-j} (-1)^{(j+i+m-2)} C(m-1-i,j)\left(\frac{\partial^{m-1-i-j}}{\partial x^{m-1-i-j}}(\alpha_1(-4+j)\, b_{1+i}(x))\right)\right)\right) = 0$$

(9.64)

With consideration of the above representations for (9.56), (9.57), (9.58), for which, as it was demonstrated earlier, it is enough to add the representation for $\alpha_1(-4)$, let's define the desired system of equalities as follows:

(9.52)

$$\alpha_1(-2) = Y_1(x)\, x + Y_2(x) \qquad (9.53)$$

$$\alpha_1(-3) = \left(\left(\frac{1}{2}\right) \cdot (Y_1(x)\, x^2)\right) + Y_2(x)\, x + Y_3(x) \qquad (9.54)$$

$$\alpha_1(-4) = \left(\left(\frac{1}{6}\right) \cdot (Y_1(x)\, x^3)\right) + \left(\left(\frac{1}{2}\right) \cdot (Y_2(x)\, x^2)\right) + Y_3(x)\, x + Y_4(x) \qquad (9.55)$$

where the functions present the solution of the ODE: (9.56)

$$2\left(\frac{d^4}{dx^4}Y(x)\right)+\left(\frac{d^3}{dx^3}Y(x)\right)b_1(x)+\left(-b_2(x)+3\left(\frac{d}{dx}b_1(x)\right)\right)\left(\frac{d^2}{dx^2}Y(x)\right)$$
$$+\left(b_3(x)+3\left(\frac{d^2}{dx^2}b_1(x)\right)-2\left(\frac{d}{dx}b_2(x)\right)\right)\left(\frac{d}{dx}Y(x)\right)$$
$$+\left(-b_4(x)+\left(\frac{d^3}{dx^3}b_1(x)\right)-\left(\frac{d^2}{dx^2}b_2(x)\right)+\left(\frac{d}{dx}b_3(x)\right)\right)Y(x)=0$$

Let’s consider private cases that shall help us to define the presentation of the special $\alpha_1(-4)$ function.

**The variant of** $m=4$. Equation (9.54), with consideration of the initial condition , takes the following form: (9.57)

$$2\left(\frac{d^4}{dx^4}\alpha_1(-4)\right)-8\left(\frac{d^3}{dx^3}\alpha_1(-3)\right)+12\left(\frac{d^2}{dx^2}\alpha_1(-2)\right)-8\left(\frac{d}{dx}\alpha_1(-1)\right)$$
$$+\left(\frac{d^3}{dx^3}(\alpha_1(-4)\,b_1(x))\right)-\left(\frac{d^2}{dx^2}(\alpha_1(-4)\,b_2(x))\right)+\left(\frac{d}{dx}(\alpha_1(-4)\,b_3(x))\right)$$
$$-\alpha_1(-4)\,b_4(x)-3\left(\frac{d^2}{dx^2}(\alpha_1(-3)\,b_1(x))\right)+2\left(\frac{d}{dx}(\alpha_1(-3)\,b_2(x))\right)-\alpha_1(-3)\,b_3(x)$$
$$+3\left(\frac{d}{dx}(\alpha_1(-2)\,b_1(x))\right)-\alpha_1(-2)\,b_2(x)-\alpha_1(-1)\,b_1(x)+b_1(x)=0$$

As an operating hypothesis, let’s assume the values of the special functions defined by the equalities (9.54) (9.57) with the condition (9.59). In this case, upon substitution to (9.59), we obtain the equality, which is solved identically, since the functions are the solution of the ODE (9.58), and the summands of the left-hand expression are the equations of the functions.

From the system of equations (12.54) - (12.57) we obtain:

(9.58)

$$Y_2(x)=(-\alpha_1(-1)+1)\,x+\alpha_1(-2) \tag{9.59}$$

$$Y_3(x)=x^2\frac{1}{2}(\alpha_1(-1)-1)-\alpha_1(-2)\,x+\alpha_1(-3) \tag{9.60}$$

$$Y_4(x)=\frac{x^3\,(1-\alpha_1(-1))}{6}+\frac{x^2\,\alpha_1(-2)}{2}-\alpha_1(-3)\,x+\alpha_1(-4) \tag{9.61}$$

Thus, it is demonstrated that the functions defined by the equalities (9.60) - (9.63) are partial solutions of the following equation:

$$\frac{d^4}{dx^4}Y(x)=a_1(x)\left(\frac{d^3}{dx^3}Y(x)\right)+a_2(x)\left(\frac{d^2}{dx^2}Y(x)\right)+a_3(x)\left(\frac{d}{dx}Y(x)\right)+a_4(x)\,Y(x) \tag{9.62}$$

where

$$a_1(x)=-\left(\left(\frac{1}{2}\right)\cdot(b_1(x))\right) \tag{9.63}$$

$$a_2(x) = \left(-\frac{1}{2}\right)\cdot\left(-b_2(x) + 3\left(\frac{d}{dx}b_1(x)\right)\right)$$

$$a_3(x) = \left(-\frac{1}{2}\right)\cdot\left(b_3(x) + 3\left(\frac{d^2}{dx^2}b_1(x)\right) - 2\left(\frac{d}{dx}b_2(x)\right)\right)$$

$$a_4(x) = \left(-\frac{1}{2}\right)\cdot\left(-b_4(x) + \left(\frac{d^3}{dx^3}b_1(x)\right) - \left(\frac{d^2}{dx^2}b_2(x)\right) + \left(\frac{d}{dx}b_3(x)\right)\right)$$

and the special functions are defined by the coefficient of n-image $\alpha_1(n)$ for the following equation:

$$\frac{d^4}{dx^4}R(x) = b_1(x)\left(\frac{d^3}{dx^3}R(x)\right) + b_2(x)\left(\frac{d^2}{dx^2}R(x)\right) + b_3(x)\left(\frac{d}{dx}R(x)\right) + b_4(x)\,R(x) \qquad (9.64)$$

Similarly, we can prove that the functions defined by the equalities (9.70) - (9.73) are partial solutions of the equations and have greater order.

Substituting the values of special functions (9.64) - (9.67) to (9.64), upon transformations with consideration of (9.68), we obtain the following identity.

Thus, we arrive to the following hypothesis:

**If special functions are defined by the following equalities**:

(9.65)

$$\alpha_1(-j) = \sum_{i=1}^{j} \frac{Y_{j-i+1}(x)\,x^{(i-1)}}{(i-1)!} \qquad j = 2 .. m \qquad (9.66)$$

**where the functions are partial solutions of the following ODE:**

$$\frac{d^m}{dx^m}Y(x) = \sum_{k=1}^{m}\left(\sum_{i=1}^{k}(-1)^{(-m-i)}C(m-i, m-k)\left(\frac{d^{k-i}}{dx^{k-i}}b_i(x)\right)\right)\left(\frac{d^{m-k}}{dx^{m-k}}Y(x)\right) \qquad (9.67)$$

**and the special functions are derivatives from the *n*-image coefficients of ODE**:

$$\frac{d^m}{dx^m}R(x) = \sum_{k=1}^{m} b_k(x)\left(\frac{d^{m-k}}{dx^{m-k}}R(x)\right) \qquad (9.68)$$

**the equation (9.69)**

$$\left(\sum_{i=0}^{m}(-1)^{(i+m)}C(m,i)\left(\frac{\partial^{m-i}}{\partial x^{m-i}}\alpha_1(n+i)\right)\right)$$
$$+\left(\sum_{j=0}^{m}\left(\sum_{i=0}^{m-1-j}(-1)^{(j+i+m-2)}C(m-1-i,j)\left(\frac{\partial^{m-1-i-j}}{\partial x^{m-1-i-j}}(\alpha_1(n+j)\,b_{1+i}(x))\right)\right)\right) = 0$$

- **is solved identically.**
  **Proof.** Assuming that $n = -m$ in (12.81), we obtain (9.70)

$$\left(\sum_{i=0}^{m} (-1)^{(i+m)} C(m,i)\left(\frac{\partial^{m-i}}{\partial x^{m-i}}\alpha_1(-m+i)\right)\right) + \left(\sum_{j=0}^{m}\left(\sum_{i=0}^{m-1-j} (-1)^{(j+i+m-2)} C(m-1-i,j)\left(\frac{\partial^{m-1-i-j}}{\partial x^{m-1-i-j}}(\alpha_1(-m+j)\, b_{1+i}(x))\right)\right)\right) = 0$$

Substituting here (12.77), (12.78) and taking into account the initial condition: , we have:

$$\left(\sum_{i=0}^{m} (-1)^{(i+m)} C(m,i)\left(\frac{\partial^{m-i}}{\partial x^{m-i}}\alpha_1(-m+i)\right)\right) + b_1(x) + \left(\sum_{j=0}^{m}\left(\sum_{i=0}^{m-1-j} (-1)^{(j+i+m-2)} C(m-1-i,j)\left(\frac{\partial^{m-1-i-j}}{\partial x^{m-1-i-j}}(\alpha_1(-m+j)\, b_{1+i}(x))\right)\right)\right) = 0$$

As soon as all solutions of the given ODE are defined only by one special function, $\alpha_1(-m+j)$, let's introduce the following

**Definition:**

(9.71)

**Thus, the desired structural form for partial solutions of the initial ODE (8.4) shall take the following form:**

$$Y_i(x) = \frac{(-1)^i (1-\alpha(-1))\, x^{(i-1)}}{(i-1)!} + \left(\sum_{k=2}^{i} \frac{(-1)^{(k+i)}\alpha(-k)\, x^{(i-k)}}{(i-k)!}\right) \qquad i = 1, 2 .. m \qquad (9.72)$$

**where  are the desired functions**.

**The key property of this formula is that, in order to obtain the general solution it is enough to find formulas for the $\alpha(n)$ function, where $n$ is the parameter. If you know this function, the special functions are defined by a plain substitution of the appropriate negative value of $i$ for $n$ parameter.**

# 10. Calculation of the n- image coefficient $\alpha(n)$.

**10. General provisions**. Given the special status of a linear second-order ODE, let's define formulas for the general solution of the following ODE:

$$\frac{\partial^m}{\partial x^m} Y = \sum_{p=1}^{m} a_p(x)\left(\frac{\partial^{m-p}}{\partial x^{m-p}} Y\right) \qquad (10.1)$$

where $a_p(x)$ are the specified functions.

Therefore, the adjoined equation is as follows:

$$\frac{d^m}{dx^m} R(x) = \sum_{k=1}^{m} b_k(x)\left(\frac{d^{m-k}}{dx^{m-k}} R(x)\right)$$

where the $b_p(x)$ coefficients are defined through the specified $a_p(x)$ coefficients by the following formulas:

$$b_k(x) = \sum_{i=1}^{k} (-1)^i \, \mathrm{C}(m-i, m-k) \left( \frac{d^{k-i}}{dx^{k-i}} a_i(x) \right) \qquad (10.2)$$

In accordance with the **Theorem on the construction of ODE solution**, general solution for the ODE (10.1) is defined by the following formula:

(10.3)

where partial solutions are defined by the following formulas:

$$Y_i(x) = \frac{(-1)^i (1 - \alpha(-1)) \, x^{(i-1)}}{(i-1)!} + \left( \sum_{k=2}^{i} \frac{(-1)^{(k+i)} \alpha(-k) \, x^{(i-k)}}{(i-k)!} \right) \qquad (10.4)$$

As known, special $\alpha(-k)$ functions are defined through the value of the *n*-image coefficient $\alpha_1(n)$, which is, in turn, defined by the equation of *n*-image:

$$\frac{d^{m+n}}{dx^{m+n}} \mathrm{R}(x) = \sum_{p=1}^{m} \alpha_p(n) \left( \frac{d^{m-p}}{dx^{m-p}} \mathrm{R}(x) \right) \qquad (10.5)$$

for the adjoined equation (10.2).

## 10.2. The definition of representation form for the special $\alpha(n)$ function

As soon as the equation of *n*-image (10.5) can be represented as follows:

$$\mathrm{Diff}\left( \frac{d^m}{dx^m} \mathrm{R}(x), x \$ n \right) = \sum_{p=1}^{m} \alpha_p(n) \left( \frac{d^{m-p}}{dx^{m-p}} \mathrm{R}(x) \right) \qquad (10.6)$$

with consideration of the initial equality (10.2) we get:

$$\sum_{k=1}^{m} \left( \frac{\partial^n}{\partial x^n} \left( b_k(x) \left( \frac{d^{m-k}}{dx^{m-k}} \mathrm{R}(x) \right) \right) \right) = \sum_{p=1}^{m} \alpha_p(n) \left( \frac{d^{m-p}}{dx^{m-p}} \mathrm{R}(x) \right)$$

In accordance with Leibniz Theorem:

$$\frac{\partial^n}{\partial x^n} \left( b_k(x) \left( \frac{d^{m-k}}{dx^{m-k}} \mathrm{R}(x) \right) \right) = \sum_{i=0}^{n} \left( \left[ \begin{matrix} i \\ n \end{matrix} \right] \cdot \left( \frac{d^{n-i}}{dx^{n-i}} b_k(x) \right) \right) \left( \frac{\partial^i}{\partial x^i} \left( \frac{d^{m-k}}{dx^{m-k}} \mathrm{R}(x) \right) \right)$$

Then the equation of ***n***-image (10.6) takes the following form:

$$\sum_{i=0}^{n} \left( \sum_{k=1}^{m} \left( \left[ \begin{matrix} i \\ n \end{matrix} \right] \cdot \left( \frac{d^{n-i}}{dx^{n-i}} b_k(x) \right) \right) \left( \frac{d^{m+i-k}}{dx^{m+i-k}} \mathrm{R}(x) \right) \right) = \sum_{p=1}^{m} \alpha_p(n) \left( \frac{d^{m-p}}{dx^{m-p}} \mathrm{R}(x) \right) \qquad (10.7)$$

In the left part of the equality the derivatives above the order are replaced by the equality (10.2). Consequently,

since the left part shall be the sum of derivatives, where the accompanying coefficients present a sum of a certain sequence of functions, then the formula for the *n*-image coefficients of $\alpha_p(n)$ should be obtained in the following form :

(10.8)

where $\xi_{i,p}(n)$ are certain functions with ***n*** parameter. N is a natural number

Taking into account that

the first coefficient of *n*-image shall be obtained in the following form:

(10.9)

$\xi_i(n)$ are the desired functions.

Taking into account the equality (10.8), the equation of *n*-image (10.5) takes the following form:

$$\frac{d^{m+n}}{dx^{m+n}}\mathrm{R}(x)=\sum_{p=1}^{m}\left(\sum_{i=0}^{N}\xi_{i,p}(n)\right)\left(\frac{d^{m-p}}{dx^{m-p}}\mathrm{R}(x)\right)$$

Or, in equivalent form,

$$\frac{d^{m+n}}{dx^{m+n}}\mathrm{R}(x)=\sum_{i=0}^{N}\left(\sum_{p=1}^{m}\xi_{i,p}(n)\left(\frac{d^{m-p}}{dx^{m-p}}\mathrm{R}(x)\right)\right) \tag{10.10}$$

Let's represent the left-hand part of (10.7) in an equivalent form:

$$\sum_{k=1}^{m}\left(\sum_{i=0}^{n}\left(\begin{bmatrix} i \\ n \end{bmatrix}\cdot\left(\frac{d^{n-i}}{dx^{n-i}}b_k(x)\right)\right)\left(\frac{d^{m+i-k}}{dx^{m+i-k}}\mathrm{R}(x)\right)\right)=$$
$$\left(\sum_{i=1}^{n}\left(\begin{bmatrix} i \\ n \end{bmatrix}\cdot\left(\frac{d^{n-i}}{dx^{n-i}}b_1(x)\right)\right)\left(\frac{d^{m+i-1}}{dx^{m+i-1}}\mathrm{R}(x)\right)\right)+\left(\frac{d^{n}}{dx^{n}}b_1(x)\right)\left(\frac{d^{m-1}}{dx^{m-1}}\mathrm{R}(x)\right)$$
$$+\left(\sum_{i=2}^{n}\left(\begin{bmatrix} i \\ n \end{bmatrix}\cdot\left(\frac{d^{n-i}}{dx^{n-i}}b_2(x)\right)\right)\left(\frac{d^{m+i-2}}{dx^{m+i-2}}\mathrm{R}(x)\right)\right)$$
$$+\left(\begin{bmatrix} 0 \\ n \end{bmatrix}\cdot\left(\frac{d^{n}}{dx^{n}}b_2(x)\right)\right)\left(\frac{d^{m-2}}{dx^{m-2}}\mathrm{R}(x)\right)+\left(\begin{bmatrix} 1 \\ n \end{bmatrix}\cdot\left(\frac{d^{n-1}}{dx^{n-1}}b_2(x)\right)\right)\left(\frac{d^{m-1}}{dx^{m-1}}\mathrm{R}(x)\right)$$

$$+\left(\sum_{i=3}^{n}\left(\begin{bmatrix} i \\ n \end{bmatrix}\cdot\left(\frac{d^{n-i}}{dx^{n-i}}b_3(x)\right)\right)\left(\frac{d^{m+i-3}}{dx^{m+i-3}}\mathrm{R}(x)\right)\right)$$

$$+\left(\begin{bmatrix} 0 \\ n \end{bmatrix}\cdot\left(\frac{d^{n}}{dx^{n}}b_3(x)\right)\right)\left(\frac{d^{m-3}}{dx^{m-3}}\mathrm{R}(x)\right)+\left(\begin{bmatrix} 1 \\ n \end{bmatrix}\cdot\left(\frac{d^{n-1}}{dx^{n-1}}b_3(x)\right)\right)\left(\frac{d^{m-2}}{dx^{m-2}}\mathrm{R}(x)\right)$$

$$+\left(\begin{bmatrix} 2 \\ n \end{bmatrix}\cdot\left(\frac{d^{n-2}}{dx^{n-2}}b_3(x)\right)\right)\left(\frac{d^{m-1}}{dx^{m-1}}\mathrm{R}(x)\right) + \ldots\ldots +$$

$$\left(\sum_{i=m}^{n}\left(\begin{bmatrix} i \\ n \end{bmatrix}\cdot\left(\frac{d^{n-i}}{dx^{n-i}}b_m(x)\right)\right)\left(\frac{d^{i}}{dx^{i}}\mathrm{R}(x)\right)\right)+\sum_{i=0}^{m-1}\left(\begin{bmatrix} i \\ n \end{bmatrix}\cdot\left(\frac{d^{n-i}}{dx^{n-i}}b_m(x)\right)\right)\left(\frac{d^{i}}{dx^{i}}\mathrm{R}(x)\right)$$

Grouping in the sequence of orders, we obtain (10.11)

$$\left(\sum_{k=1}^{m}\left(\sum_{i=k}^{n}\left(\begin{bmatrix} i \\ n \end{bmatrix}\cdot\left(\frac{d^{n-i}}{dx^{n-i}}b_k(x)\right)\right)\left(\frac{d^{m+i-k}}{dx^{m+i-k}}\mathrm{R}(x)\right)\right)\right)$$

$$+\left(\sum_{k=0}^{m-1}\left(\sum_{i=0}^{k}\left(\begin{bmatrix} i \\ n \end{bmatrix}\cdot\left(\frac{d^{n-i}}{dx^{n-i}}b_{m+i-k}(x)\right)\right)\left(\frac{d^{k}}{dx^{k}}\mathrm{R}(x)\right)\right)\right)$$

Returning to (10.7), with consideration of the obtained result and the fact that the right-hand part of (13.7) equals the right-hand part of (13.10), we obtain the following equality:

$$\left(\sum_{k=1}^{m}\left(\sum_{i=k}^{n}\left(\begin{bmatrix} i \\ n \end{bmatrix}\cdot\left(\frac{d^{n-i}}{dx^{n-i}}b_k(x)\right)\right)\left(\frac{d^{m+i-k}}{dx^{m+i-k}}\mathrm{R}(x)\right)\right)\right)$$

$$+\left(\sum_{k=0}^{m-1}\left(\sum_{i=0}^{k}\left(\begin{bmatrix} i \\ n \end{bmatrix}\cdot\left(\frac{d^{n-i}}{dx^{n-i}}b_{m+i-k}(x)\right)\right)\left(\frac{d^{k}}{dx^{k}}\mathrm{R}(x)\right)\right)\right) = \sum_{i=0}^{N}\left(\sum_{p=1}^{m}\xi_{i,p}(n)\left(\frac{d^{m-p}}{dx^{m-p}}\mathrm{R}(x)\right)\right)$$

Let's present it in the following form: (10.11)

$$\left(\sum_{k=1}^{m}\left(\sum_{i=k}^{n}\left(\begin{bmatrix} i \\ n \end{bmatrix}\cdot\left(\frac{d^{n-i}}{dx^{n-i}}b_k(x)\right)\right)\left(\frac{d^{m+i-k}}{dx^{m+i-k}}\mathrm{R}(x)\right)\right)\right)$$

$$+\left(\sum_{k=0}^{m-1}\left(\sum_{i=0}^{k}\left(\begin{bmatrix} i \\ n \end{bmatrix}\cdot\left(\frac{d^{n-i}}{dx^{n-i}}b_{m+i-k}(x)\right)\right)\left(\frac{d^{k}}{dx^{k}}\mathrm{R}(x)\right)\right)\right) =$$

$$\left(\sum_{i=1}^{N}\left(\sum_{p=1}^{m}\xi_{i,p}(n)\left(\frac{d^{m-p}}{dx^{m-p}}\mathrm{R}(x)\right)\right)\right)+\left(\sum_{p=1}^{m}\xi_{0,p}(n)\left(\frac{d^{m-p}}{dx^{m-p}}\mathrm{R}(x)\right)\right)$$

As follows from here

$$\sum_{k=0}^{m-1}\left(\sum_{i=0}^{k}\left(\begin{bmatrix} i \\ n \end{bmatrix}\cdot\left(\frac{d^{n-i}}{dx^{n-i}}b_{m+i-k}(x)\right)\right)\left(\frac{d^{k}}{dx^{k}}\mathrm{R}(x)\right)\right) = \sum_{p=1}^{m}\xi_{0,p}(n)\left(\frac{d^{m-p}}{dx^{m-p}}\mathrm{R}(x)\right)$$

Thus, the desired zero, that is, initial value of the function is equal to:

$$\xi_0(n) = \sum_{i=0}^{m-1}\left(\begin{bmatrix} i \\ n \end{bmatrix}\cdot\left(\frac{d^{n-i}}{dx^{n-i}}b_{i+1}(x)\right)\right) \qquad (10.12)$$

From (10.11) we also obtain the second identity:

$$\sum_{k=1}^{m}\left(\sum_{i=k}^{n}\left(\left[\begin{matrix} i \\ n \end{matrix}\right]\cdot\left(\frac{d^{n-i}}{dx^{n-i}}b_k(x)\right)\right)\left(\frac{d^{m+i-k}}{dx^{m+i-k}}\mathrm{R}(x)\right)\right)=\sum_{i=1}^{N}\left(\sum_{p=1}^{m}\xi_{i,p}(n)\left(\frac{d^{m-p}}{dx^{m-p}}\mathrm{R}(x)\right)\right) \qquad (10.13)$$

If in the summand

$$\sum_{k=1}^{m}\left(\sum_{i=k}^{n}\left(\left[\begin{matrix} i \\ n \end{matrix}\right]\cdot\left(\frac{d^{n-i}}{dx^{n-i}}b_k(x)\right)\right)\left(\frac{d^{m+i-k}}{dx^{m+i-k}}\mathrm{R}(x)\right)\right)$$

we replace the summation index $i$ by $i-k$, we obtain

$$\sum_{k=1}^{m}\left(\sum_{i=k}^{n}\left(\left[\begin{matrix} i \\ n \end{matrix}\right]\cdot\left(\frac{d^{n-i}}{dx^{n-i}}b_k(x)\right)\right)\left(\frac{d^{m+i-k}}{dx^{m+i-k}}\mathrm{R}(x)\right)\right)=$$
$$\sum_{k=1}^{m}\left(\sum_{i=0}^{-k+n}\left(\left[\begin{matrix} k+i \\ n \end{matrix}\right]\cdot\left(\frac{d^{n-k-i}}{dx^{n-k-i}}b_k(x)\right)\right)\left(\frac{d^{m+i}}{dx^{m+i}}\mathrm{R}(x)\right)\right)$$

Substituting the value (10.2) to the right-hand part, we obtain:

$$\sum_{k=1}^{m}\left(\sum_{i=0}^{-k+n}\left(\left[\begin{matrix} k+i \\ n \end{matrix}\right]\cdot\left(\frac{d^{n-k-i}}{dx^{n-k-i}}b_k(x)\right)\right)\left(\frac{\partial^i}{\partial x^i}\left(\sum_{k_1=1}^{m}b_{k_1}(x)\left(\frac{d^{m-k_1}}{dx^{m-k_1}}\mathrm{R}(x)\right)\right)\right)\right)$$

Then: (10.14)

$$\sum_{k=1}^{m}\left(\sum_{i=0}^{-k+n}\left(\left[\begin{matrix} k+i \\ n \end{matrix}\right]\cdot\left(\frac{d^{n-k-i}}{dx^{n-k-i}}b_k(x)\right)\right)\left(\frac{\partial^i}{\partial x^i}\left(\sum_{k_1=1}^{m}b_{k_1}(x)\left(\frac{d^{m-k_1}}{dx^{m-k_1}}\mathrm{R}(x)\right)\right)\right)\right)=$$
$$\sum_{k=1}^{m}\left(\sum_{i=0}^{-k+n}\left(\sum_{k_1=1}^{m}\left(\left[\begin{matrix} k+i \\ n \end{matrix}\right]\cdot\left(\frac{d^{n-k-i}}{dx^{n-k-i}}b_k(x)\right)\right)\left(\frac{\partial^i}{\partial x^i}\left(b_{k_1}(x)\left(\frac{d^{m-k_1}}{dx^{m-k_1}}\mathrm{R}(x)\right)\right)\right)\right)\right)$$

Since, in accordance with Leibniz formula

$$\frac{\partial^i}{\partial x^i}\left(b_{k_1}(x)\left(\frac{d^{m-k_1}}{dx^{m-k_1}}\mathrm{R}(x)\right)\right)=\sum_{i_1=0}^{i}\left[\begin{matrix} i_1 \\ i \end{matrix}\right]\left(\frac{d^{m+i_1-k_1}}{dx^{m+i_1-k_1}}\mathrm{R}(x)\right)\left(\frac{d^{i-i_1}}{dx^{i-i_1}}b_{k_1}(x)\right)$$

we eventually obtain for the right-hand part of (10.14):

$$\sum_{k=1}^{m}\sum_{i=0}^{-k+n}\left(\left[\begin{matrix} k+i \\ n \end{matrix}\right]\cdot\left(\left(\frac{d^{n-k-i}}{dx^{n-k-i}}b_k(x)\right)\left(\sum_{k_1=1}^{m}\left(\sum_{i_1=0}^{i}\left[\begin{matrix} i_1 \\ i \end{matrix}\right]\left(\frac{d^{i-i_1}}{dx^{i-i_1}}b_{k_1}(x)\right)\left(\frac{d^{m+i_1-k_1}}{dx^{m+i_1-k_1}}\mathrm{R}(x)\right)\right)\right)\right)\right)$$

Grouping in the secuence of orders, we get an equivalent equality:

$$\sum_{k=1}^{m}\left(\sum_{i=0}^{-k+n}\left(\begin{bmatrix}k+i\\ n\end{bmatrix}\cdot\left(\left(\frac{d^{n-k-i}}{dx^{n-k-i}}b_k(x)\right)\left(\sum_{k_1=1}^{m}\left(\sum_{i_1=0}^{i}\begin{bmatrix}i_1\\ i\end{bmatrix}\left(\frac{d^{i-i_1}}{dx^{i-i_1}}b_{k_1}(x)\right)\left(\frac{d^{m+i_1-k_1}}{dx^{m+i_1-k_1}}\mathrm{R}(x)\right)\right)\right)\right)\right)\right)$$

$$=\sum_{k=1}^{m}\left(\sum_{i=0}^{-k+n}\left(\begin{bmatrix}k+i\\ n\end{bmatrix}\cdot\left(\left(\frac{d^{n-k-i}}{dx^{n-k-i}}b_k(x)\right)\left(\sum_{k_1=1}^{m}\left(\sum_{i_1=k_1}^{i}\begin{bmatrix}i_1\\ i\end{bmatrix}\left(\frac{d^{i-i_1}}{dx^{i-i_1}}b_{k_1}(x)\right)\left(\frac{d^{m+i_1-k_1}}{dx^{m+i_1-k_1}}\mathrm{R}(x)\right)\right)\right)\right.\right.\right.$$

$$\left.\left.\left.+\left(\sum_{k_1=0}^{m-1}\left(\sum_{i_1=0}^{k_1}\left(\begin{bmatrix}i_1\\ i\end{bmatrix}\cdot\left(\frac{d^{i-i_1}}{dx^{i-i_1}}b_{m+i_1-k_1}(x)\right)\right)\left(\frac{d^{k_1}}{dx^{k_1}}\mathrm{R}(x)\right)\right)\right)\right)\right)\right)$$

In this case, the equality (10.13) assumes the following form (10.15)

$$\left(\sum_{k=1}^{m}\left(\sum_{i=0}^{-k+n}\left(\begin{bmatrix}k+i\\ n\end{bmatrix}\cdot\left(\left(\frac{d^{n-k-i}}{dx^{n-k-i}}b_k(x)\right)\left(\sum_{k_1=1}^{m}\left(\sum_{i_1=k_1}^{i}\begin{bmatrix}i_1\\ i\end{bmatrix}\left(\frac{d^{i-i_1}}{dx^{i-i_1}}b_{k_1}(x)\right)\left(\frac{d^{m+i_1-k_1}}{dx^{m+i_1-k_1}}\mathrm{R}(x)\right)\right)\right)\right)\right)\right.\right.$$

$$\begin{bmatrix}k+i\\ n\end{bmatrix}\cdot\left(\left(\frac{d^{n-k-i}}{dx^{n-k-i}}b_k(x)\right)\left(\sum_{k_1=0}^{m-1}\left(\sum_{i_1=0}^{k_1}\left(\begin{bmatrix}i_1\\ i\end{bmatrix}\cdot\left(\frac{d^{i-i_1}}{dx^{i-i_1}}b_{m+i_1-k_1}(x)\right)\right)\left(\frac{d^{k_1}}{dx^{k_1}}\mathrm{R}(x)\right)\right)\right)\right)$$

$$\left.\left.\right)\right)=\left(\sum_{i=2}^{N}\left(\sum_{p=1}^{m}\xi_{i,p}(n)\left(\frac{d^{m-p}}{dx^{m-p}}\mathrm{R}(x)\right)\right)\right)+\left(\sum_{p=1}^{m}\xi_{1,p}(n)\left(\frac{d^{m-p}}{dx^{m-p}}\mathrm{R}(x)\right)\right)$$

From here we obtain the first desired value of $\xi_1(n)$ :

$$\xi_1(n)=\sum_{i_1=0}^{m-1}\left(\sum_{k=1}^{m}\left(\sum_{i=0}^{n-k}\left(\begin{bmatrix}k+i\\ n\end{bmatrix}\cdot\begin{bmatrix}i_1\\ i\end{bmatrix}\right)\left(\frac{d^{n-k-i}}{dx^{n-k-i}}b_k(x)\right)\left(\frac{d^{i-i_1}}{dx^{i-i_1}}b_{i_1+1}(x)\right)\right)\right)\qquad(10.16)$$

From (10.15), taking into account (10.16), we obtain the following equality: (10.17)

$$\sum_{k=1}^{m}\left(\sum_{i=0}^{-k+n}\left(\begin{bmatrix}k+i\\ n\end{bmatrix}\cdot\left(\left(\frac{d^{n-k-i}}{dx^{n-k-i}}b_k(x)\right)\left(\sum_{k_1=1}^{m}\left(\sum_{i_1=k_1}^{i}\begin{bmatrix}i_1\\ i\end{bmatrix}\left(\frac{d^{i-i_1}}{dx^{i-i_1}}b_{k_1}(x)\right)\left(\frac{d^{m+i_1-k_1}}{dx^{m+i_1-k_1}}\mathrm{R}(x)\right)\right)\right)\right)\right)\right.$$

$$\left.=\sum_{i=2}^{N}\left(\sum_{p=1}^{m}\xi_{i,p}(n)\left(\frac{d^{m-p}}{dx^{m-p}}\mathrm{R}(x)\right)\right)\right)$$

On the left side of this equality we shall again replace the summation index $i_1$ by $i_1-k_1$. Then we obtain the new equivalent representation:

$$\sum_{k=1}^{m}\left(\sum_{i=0}^{-k+n}\left(\begin{bmatrix}k+i\\ n\end{bmatrix}\cdot\left(\left(\frac{d^{n-k-i}}{dx^{n-k-i}}b_k(x)\right)\left(\sum_{k_1=1}^{m}\left(\sum_{i_1=0}^{i-k_1}\begin{bmatrix}i_1+k_1\\ i\end{bmatrix}\left(\frac{d^{i-i_1-k_1}}{dx^{i-i_1-k_1}}b_{k_1}(x)\right)\left(\frac{d^{m+i_1}}{dx^{m+i_1}}\mathrm{R}(x)\right)\right)\right)\right)\right)\right)$$

Substituting here the value (10 2) we obtain: (10.18)

$$\sum_{k=1}^{m}\Bigg(\sum_{i=0}^{-k+n}\Bigg(\begin{bmatrix}k+i\\n\end{bmatrix}\cdot\Bigg(\Bigg(\frac{d^{n-k-i}}{dx^{n-k-i}}b_k(x)\Bigg)$$
$$\Bigg(\sum_{k_1=1}^{m}\Bigg(\sum_{i_1=0}^{i-k_1}\begin{bmatrix}i_1+k_1\\i\end{bmatrix}\Bigg(\frac{d^{i-i_1-k_1}}{dx^{i-i_1-k_1}}b_{k_1}(x)\Bigg)\Bigg(\frac{\partial^{i_1}}{\partial x^{i_1}}\Bigg(\sum_{k_2=1}^{m}b_{k_2}(x)\Bigg(\frac{d^{m-k_2}}{dx^{m-k_2}}\mathrm{R}(x)\Bigg)\Bigg)\Bigg)\Bigg)\Bigg)\Bigg)\Bigg)\Bigg)$$

With consideration of Leibniz formula:

$$\frac{\partial^{i_1}}{\partial x^{i_1}}\Bigg(b_{k_2}(x)\Bigg(\frac{d^{m-k_2}}{dx^{m-k_2}}\mathrm{R}(x)\Bigg)\Bigg)=\sum_{i_2=0}^{i_1}\begin{bmatrix}i_2\\i_1\end{bmatrix}\Bigg(\frac{d^{m+i_2-k_2}}{dx^{m+i_2-k_2}}\mathrm{R}(x)\Bigg)\Bigg(\frac{d^{i_1-i_2}}{dx^{i_1-i_2}}b_{k_2}(x)\Bigg)$$

Thus, the expression (10.18) can be reduced as follows:

$$\sum_{k=1}^{m}\Bigg(\sum_{i=0}^{-k+n}\Bigg(\begin{bmatrix}k+i\\n\end{bmatrix}\cdot\Bigg(\Bigg(\frac{d^{n-k-i}}{dx^{n-k-i}}b_k(x)\Bigg)\Bigg(\sum_{k_1=1}^{m}\Bigg(\sum_{i_1=0}^{i-k_1}$$
$$\begin{bmatrix}i_1+k_1\\i\end{bmatrix}\Bigg(\frac{d^{i-i_1-k_1}}{dx^{i-i_1-k_1}}b_{k_1}(x)\Bigg)\Bigg(\sum_{k_2=1}^{m}\Bigg(\sum_{i_2=0}^{i_1}\begin{bmatrix}i_2\\i_1\end{bmatrix}\Bigg(\frac{d^{i_1-i_2}}{dx^{i_1-i_2}}b_{k_2}(x)\Bigg)\Bigg(\frac{d^{m+i_2-k_2}}{dx^{m+i_2-k_2}}\mathrm{R}(x)\Bigg)\Bigg)\Bigg)\Bigg)\Bigg)\Bigg)\Bigg)\Bigg)\Bigg)$$

This expression is equivalent to

$$\sum_{k=1}^{m}\Bigg(\sum_{i=0}^{-k+n}\Bigg(\begin{bmatrix}k+i\\n\end{bmatrix}\cdot\Bigg(\Bigg(\frac{d^{n-k-i}}{dx^{n-k-i}}b_k(x)\Bigg)\Bigg(\sum_{k_1=1}^{m}\Bigg(\sum_{i_1=0}^{i-k_1}\begin{bmatrix}i_1+k_1\\i\end{bmatrix}\Bigg(\frac{d^{i-i_1-k_1}}{dx^{i-i_1-k_1}}b_{k_1}(x)\Bigg)\Bigg($$
$$\Bigg(\sum_{k_2=1}^{m}\Bigg(\sum_{i_2=k_2}^{i_1}\begin{bmatrix}i_2\\i_1\end{bmatrix}\Bigg(\frac{d^{i_1-i_2}}{dx^{i_1-i_2}}b_{k_2}(x)\Bigg)\Bigg(\frac{d^{m+i_2-k_2}}{dx^{m+i_2-k_2}}\mathrm{R}(x)\Bigg)\Bigg)\Bigg)$$
$$+\Bigg(\sum_{k_2=0}^{m-1}\Bigg(\sum_{i_2=0}^{k_2}\begin{bmatrix}i_2\\i_1\end{bmatrix}\Bigg(\frac{d^{i_1-i_2}}{dx^{i_1-i_2}}b_{m+i_2-k_2}(x)\Bigg)\Bigg(\frac{d^{k_2}}{dx^{k_2}}\mathrm{R}(x)\Bigg)\Bigg)\Bigg)\Bigg)\Bigg)\Bigg)\Bigg)\Bigg)\Bigg)\Bigg)$$

In this case the equality (10.17) assumes the following form: (10.19)

$$\Bigg(\sum_{k=1}^{m}\Bigg(\sum_{i=0}^{-k+n}\Bigg(\begin{bmatrix}k+i\\n\end{bmatrix}\cdot\Bigg(\Bigg(\frac{d^{n-k-i}}{dx^{n-k-i}}b_k(x)\Bigg)\Bigg(\sum_{k_1=1}^{m}\Bigg(\sum_{i_1=0}^{i-k_1}$$
$$\begin{bmatrix}i_1+k_1\\i\end{bmatrix}\Bigg(\frac{d^{i-i_1-k_1}}{dx^{i-i_1-k_1}}b_{k_1}(x)\Bigg)\Bigg(\sum_{k_2=1}^{m}\Bigg(\sum_{i_2=k_2}^{i_1}\begin{bmatrix}i_2\\i_1\end{bmatrix}\Bigg(\frac{d^{i_1-i_2}}{dx^{i_1-i_2}}b_{k_2}(x)\Bigg)\Bigg(\frac{d^{m+i_2-k_2}}{dx^{m+i_2-k_2}}\mathrm{R}(x)\Bigg)\Bigg)\Bigg)$$
$$\Bigg)\Bigg)\Bigg)\Bigg)\Bigg)\Bigg)+\Bigg(\sum_{k=1}^{m}\Bigg(\sum_{i=0}^{-k+n}\Bigg(\begin{bmatrix}k+i\\n\end{bmatrix}\cdot\Bigg(\Bigg(\frac{d^{n-k-i}}{dx^{n-k-i}}b_k(x)\Bigg)\Bigg(\sum_{k_1=1}^{m}\Bigg(\sum_{i_1=0}^{i-k_1}$$

$$\begin{bmatrix} i_1+k_1 \\ i \end{bmatrix}\left(\frac{d^{i-i_1-k_1}}{dx^{i-i_1-k_1}} b_{k_1}(x)\right)\left(\sum_{k_2=0}^{m-1}\left(\sum_{i_2=0}^{k_2}\begin{bmatrix} i_2 \\ i_1\end{bmatrix}\left(\frac{d^{i_1-i_2}}{dx^{i_1-i_2}} b_{m+i_2-k_2}(x)\right)\left(\frac{d^{k_2}}{dx^{k_2}}\mathrm{R}(x)\right)\right)\right)$$

$$= \left(\sum_{i=3}^{N}\left(\sum_{p=1}^{m} \xi_{i,p}(n)\left(\frac{d^{m-p}}{dx^{m-p}}\mathrm{R}(x)\right)\right)\right)+\left(\sum_{p=1}^{m} \xi_{2,p}(n)\left(\frac{d^{m-p}}{dx^{m-p}}\mathrm{R}(x)\right)\right)$$

As follows from here, the desired second summand for the first coefficient of *n*-image is defined by the following formula: (10.20)

$$\xi_2(n) = \sum_{k=1}^{m}\left(\sum_{i=0}^{-k+n}\left(\begin{bmatrix} k+i \\ n \end{bmatrix}\cdot\left(\left(\frac{d^{n-k-i}}{dx^{n-k-i}} b_k(x)\right)\right.\right.\right.$$

$$\left(\sum_{k_1=1}^{m}\left(\sum_{i_1=0}^{i-k_1}\begin{bmatrix} i_1+k_1 \\ i \end{bmatrix}\left(\frac{d^{i-i_1-k_1}}{dx^{i-i_1-k_1}} b_{k_1}(x)\right)\left(\sum_{i_2=0}^{m-1}\begin{bmatrix} i_2 \\ i_1 \end{bmatrix}\left(\frac{d^{i_1-i_2}}{dx^{i_1-i_2}} b_{1+i_2}(x)\right)\right)\right)\right)$$

Again, similarly to the described algorithm, from (10.19) we obtain the following equality: (13.21)

$$\sum_{k=1}^{m}\left(\sum_{i=0}^{-k+n}\left(\begin{bmatrix} k+i \\ n \end{bmatrix}\cdot\left(\left(\frac{d^{n-k-i}}{dx^{n-k-i}} b_k(x)\right)\left(\sum_{k_1=1}^{m}\left(\sum_{i_1=0}^{i-k_1}\right.\right.\right.\right.\right.$$

$$\begin{bmatrix} i_1+k_1 \\ i \end{bmatrix}\left(\frac{d^{i-i_1-k_1}}{dx^{i-i_1-k_1}} b_{k_1}(x)\right)\left(\sum_{k_2=1}^{m}\left(\sum_{i_2=k_2}^{i_1}\begin{bmatrix} i_2 \\ i_1 \end{bmatrix}\left(\frac{d^{i_1-i_2}}{dx^{i_1-i_2}} b_{k_2}(x)\right)\left(\frac{d^{m+i_2-k_2}}{dx^{m+i_2-k_2}}\mathrm{R}(x)\right)\right)\right)$$

$$= \sum_{i=3}^{N}\left(\sum_{p=1}^{m} \xi_{i,p}(n)\left(\frac{d^{m-p}}{dx^{m-p}}\mathrm{R}(x)\right)\right)$$

Further we can similarly find: (10.21)

$$\xi_3(n) = \sum_{k=1}^{m}\left(\sum_{i=0}^{-k+n}\left(\begin{bmatrix} k+i \\ n \end{bmatrix}\cdot\left(\left(\frac{d^{n-k-i}}{dx^{n-k-i}} b_k(x)\right)\left(\sum_{k_1=1}^{m}\left(\sum_{i_1=0}^{i-k_1}\begin{bmatrix} i_1+k_1 \\ i \end{bmatrix}\left(\frac{d^{i-i_1-k_1}}{dx^{i-i_1-k_1}} b_{k_1}(x)\right)\right.\right.\right.\right.\right.$$

$$\left(\sum_{k_2=1}^{m}\left(\sum_{i_2=0}^{i_1-k_2}\begin{bmatrix} i_2+k_2 \\ i_1 \end{bmatrix}\left(\frac{d^{i_1-i_2-k_2}}{dx^{i_1-i_2-k_2}} b_{k_2}(x)\right)\left(\sum_{i_3=0}^{m-1}\begin{bmatrix} i_3 \\ i_2 \end{bmatrix}\left(\frac{d^{i_2-i_3}}{dx^{i_2-i_3}} b_{1+i_3}(x)\right)\right)\right)\right)$$

$$\xi_4(n) = \sum_{k=1}^{m}\left(\sum_{i=0}^{-k+n}\left(\begin{bmatrix} k+i \\ n \end{bmatrix}\cdot\left(\left(\frac{d^{n-k-i}}{dx^{n-k-i}} b_k(x)\right)\left(\sum_{k_1=1}^{m}\left(\sum_{i_1=0}^{i-k_1}\begin{bmatrix} i_1+k_1 \\ i \end{bmatrix}\left(\frac{d^{i-i_1-k_1}}{dx^{i-i_1-k_1}} b_{k_1}(x)\right)\right.\right.\right.\right.\right.$$

$$\sum_{k_2=1}^{m}\left(\sum_{i_2=0}^{i_1-k_2}\begin{bmatrix} i_2+k_2 \\ i_1 \end{bmatrix}\left(\frac{d^{i_1-i_2-k_2}}{dx^{i_1-i_2-k_2}} b_{k_2}(x)\right)\right.$$

$$\left(\sum_{k_3=1}^{m}\left(\sum_{i_3=0}^{i_2-k_3}\begin{bmatrix} i_3+k_3 \\ i_2 \end{bmatrix}\left(\frac{d^{i_2-i_3-k_3}}{dx^{i_2-i_3-k_3}} b_{k_3}(x)\right)\left(\sum_{i_4=0}^{m-1}\begin{bmatrix} i_4 \\ i_3 \end{bmatrix}\left(\frac{d^{i_3-i_4}}{dx^{i_3-i_4}} b_{1+i_4}(x)\right)\right)\right)\right)$$

$$\xi_5(n)=\sum_{k=1}^{m}\Bigg(\sum_{i=0}^{-k+n}\Bigg(\begin{bmatrix}k+i\\ n\end{bmatrix}\cdot\Bigg(\Bigg(\frac{d^{n-k-i}}{dx^{n-k-i}}b_k(x)\Bigg)\Bigg(\sum_{k_1=1}^{m}\Bigg(\sum_{i_1=0}^{i-k_1}\begin{bmatrix}i_1+k_1\\ i\end{bmatrix}\Bigg(\frac{d^{i-i_1-k_1}}{dx^{i-i_1-k_1}}b_{k_1}(x)\Bigg)\Bigg($$

$$\sum_{k_2=1}^{m}\Bigg(\sum_{i_2=0}^{i_1-k_2}\begin{bmatrix}i_2+k_2\\ i_1\end{bmatrix}\Bigg(\frac{d^{i_1-i_2-k_2}}{dx^{i_1-i_2-k_2}}b_{k_2}(x)\Bigg)\Bigg(\sum_{k_3=1}^{m}\Bigg(\sum_{i_3=0}^{i_2-k_3}\begin{bmatrix}i_3+k_3\\ i_2\end{bmatrix}\Bigg(\frac{d^{i_2-i_3-k_3}}{dx^{i_2-i_3-k_3}}b_{k_3}(x)\Bigg)$$

$$\Bigg(\sum_{k_4=1}^{m}\Bigg(\sum_{i_4=0}^{i_3-k_4}\begin{bmatrix}i_4+k_4\\ i_3\end{bmatrix}\Bigg(\frac{d^{i_3-i_4-k_4}}{dx^{i_3-i_4-k_4}}b_{k_4}(x)\Bigg)\Bigg(\sum_{i_5=0}^{m-1}\begin{bmatrix}i_5\\ i_4\end{bmatrix}\Bigg(\frac{d^{i_4-i_5}}{dx^{i_4-i_5}}b_{1+i_5}(x)\Bigg)\Bigg)\Bigg)\Bigg)\Bigg)\Bigg)\Bigg)\Bigg)\Bigg)\Bigg)\Bigg)\Bigg)\Bigg)\Bigg)\Bigg)$$

and so on.

As a general case, we can write down (10.22)

$$\xi_N(n)=\sum_{k=1}^{m}\Bigg(\sum_{i=0}^{-k+n}\Bigg(\begin{bmatrix}k+i\\ n\end{bmatrix}\cdot\Bigg(\Bigg(\frac{d^{n-k-i}}{dx^{n-k-i}}b_k(x)\Bigg)\Bigg(\sum_{k_1=1}^{m}\Bigg(\sum_{i_1=0}^{i-k_1}\begin{bmatrix}i_1+k_1\\ i\end{bmatrix}\Bigg(\frac{d^{i-i_1-k_1}}{dx^{i-i_1-k_1}}b_{k_1}(x)\Bigg)\Bigg($$

$$\sum_{k_2=1}^{m}\Bigg(\sum_{i_2=0}^{i_1-k_2}\begin{bmatrix}i_2+k_2\\ i_1\end{bmatrix}\Bigg(\frac{d^{i_1-i_2-k_2}}{dx^{i_1-i_2-k_2}}b_{k_2}(x)\Bigg)[[0,0,0]]\Bigg(\sum_{k_{N-2}=1}^{m}\Bigg(\sum_{i_{N-2}=0}^{i_{N-3}-k_{N-2}}$$

$$\begin{bmatrix}i_{N-2}+k_{N-2}\\ i_{N-3}\end{bmatrix}\Bigg(\frac{d^{i_{N-3}-i_{N-2}-k_{N-2}}}{dx^{i_{N-3}-i_{N-2}-k_{N-2}}}b_{k_{N-2}}(x)\Bigg)\Bigg(\sum_{k_{N-1}=1}^{m}\Bigg(\sum_{i_{N-1}=0}^{i_{N-2}-k_{N-1}}$$

$$\begin{bmatrix}i_{N-1}+k_{N-1}\\ i_{N-2}\end{bmatrix}\Bigg(\frac{d^{i_{N-2}-i_{N-1}-k_{N-1}}}{dx^{i_{N-2}-i_{N-1}-k_{N-1}}}b_{k_{N-1}}(x)\Bigg)\Bigg(\sum_{i_N=0}^{m-1}\begin{bmatrix}i_N\\ i_{N-1}\end{bmatrix}\Bigg(\frac{d^{i_{N-1}-i_N}}{dx^{i_{N-1}-i_N}}b_{1+i_N}(x)\Bigg)\Bigg)\Bigg)\Bigg)\Bigg)$$

## 10.3 The key property of serial summands for the special $\alpha(n)$ function

Let's write down the obtained values of

$$\xi_0(n)=\sum_{i=0}^{m-1}\Bigg(\begin{bmatrix}i\\ n\end{bmatrix}\cdot\Bigg(\frac{d^{n-i}}{dx^{n-i}}b_{i+1}(x)\Bigg)\Bigg)$$

$$\xi_1(n)=\sum_{i_1=0}^{m-1}\Bigg(\sum_{k=1}^{m}\Bigg(\sum_{i=0}^{n-k}\Bigg(\Bigg(\begin{bmatrix}k+i\\ n\end{bmatrix}\cdot\begin{bmatrix}i_1\\ i\end{bmatrix}\Bigg)\Bigg(\frac{d^{n-k-i}}{dx^{n-k-i}}b_k(x)\Bigg)\Bigg(\frac{d^{i-i_1}}{dx^{i-i_1}}b_{i_1+1}(x)\Bigg)\Bigg)\Bigg)\Bigg)$$

$$\xi_2(n)=\sum_{k=1}^{m}\Bigg(\sum_{i=0}^{-k+n}\Bigg(\begin{bmatrix}k+i\\ n\end{bmatrix}\cdot\Bigg(\Bigg(\frac{d^{n-k-i}}{dx^{n-k-i}}b_k(x)\Bigg)$$

$$\Bigg(\sum_{k_1=1}^{m}\Bigg(\sum_{i_1=0}^{i-k_1}\begin{bmatrix}i_1+k_1\\ i\end{bmatrix}\Bigg(\frac{d^{i-i_1-k_1}}{dx^{i-i_1-k_1}}b_{k_1}(x)\Bigg)\Bigg(\sum_{i_2=0}^{m-1}\begin{bmatrix}i_2\\ i_1\end{bmatrix}\Bigg(\frac{d^{i_1-i_2}}{dx^{i_1-i_2}}b_{1+i_2}(x)\Bigg)\Bigg)\Bigg)\Bigg)\Bigg)\Bigg)\Bigg)\Bigg)$$

$$\xi_3(n) = \sum_{k=1}^{m}\left(\sum_{i=0}^{-k+n}\left(\begin{bmatrix}k+i\\n\end{bmatrix}\cdot\left(\left(\frac{d^{n-k-i}}{dx^{n-k-i}}b_k(x)\right)\left(\sum_{k_1=1}^{m}\left(\sum_{i_1=0}^{i-k_1}\begin{bmatrix}i_1+k_1\\i\end{bmatrix}\left(\frac{d^{i-i_1-k_1}}{dx^{i-i_1-k_1}}b_{k_1}(x)\right)\right.\right.\right.\right.\right.$$
$$\left.\left.\left.\left.\left.\left(\sum_{k_2=1}^{m}\left(\sum_{i_2=0}^{i_1-k_2}\begin{bmatrix}i_2+k_2\\i_1\end{bmatrix}\left(\frac{d^{i_1-i_2-k_2}}{dx^{i_1-i_2-k_2}}b_{k_2}(x)\right)\left(\sum_{i_3=0}^{m-1}\begin{bmatrix}i_3\\i_2\end{bmatrix}\left(\frac{d^{i_2-i_3}}{dx^{i_2-i_3}}b_{1+i_3}(x)\right)\right)\right)\right)\right)\right)\right)\right)\right)$$

and so on.

Analyzing the obtained formulas, let's prove the following **Theorem: Any two successive summands are related by the following equality:**

$$\xi_{s+1}(n) = \sum_{k=1}^{m}\left(\sum_{i_0=0}^{n-k}\begin{bmatrix}i_0+k\\n\end{bmatrix}\left(\frac{d^{n-i_0-k}}{dx^{n-i_0-k}}b_k(x)\right)\xi_s(i_0)\right) \tag{10.23}$$

**Proof**: Let's check the validity of this formula for a number of first summands. Indeed, in accordance with (10.12):

$$\xi_0(n) = \sum_{i=0}^{m-1}\left(\begin{bmatrix}i\\n\end{bmatrix}\cdot\left(\frac{d^{n-i}}{dx^{n-i}}b_{i+1}(x)\right)\right) \tag{10.24}$$

and, in accordance with (10.16)

$$\xi_1(n) = \sum_{i_1=0}^{m-1}\left(\sum_{k=1}^{m}\left(\sum_{i=0}^{n-k}\left(\begin{bmatrix}k+i\\n\end{bmatrix}\cdot\begin{bmatrix}i_1\\i\end{bmatrix}\right)\left(\frac{d^{n-k-i}}{dx^{n-k-i}}b_k(x)\right)\left(\frac{d^{i-i_1}}{dx^{i-i_1}}b_{i_1+1}(x)\right)\right)\right) \tag{10.25}$$

Let’s assume $s=0$ in the equality (10.23). Then we get:

$$\xi_1(n) = \sum_{k=1}^{m}\left(\sum_{i_0=0}^{n-k}\begin{bmatrix}i_0+k\\n\end{bmatrix}\left(\frac{d^{n-i_0-k}}{dx^{n-i_0-k}}b_k(x)\right)\xi_0(i_0)\right) \tag{10.26}$$

Substituting $n=i_0$ to (10.24), we get:

$$\xi_0(i_0) = \sum_{i=0}^{m-1}\left(\begin{bmatrix}i\\i_0\end{bmatrix}\cdot\left(\frac{d^{i_0-i}}{dx^{i_0-i}}b_{i+1}(x)\right)\right)$$

Substituting this value to (10.26) we obtain:

$$\xi_1(n) = \sum_{k=1}^{m}\left(\sum_{i_0=0}^{n-k}\begin{bmatrix}i_0+k\\n\end{bmatrix}\left(\frac{d^{n-i_0-k}}{dx^{n-i_0-k}}b_k(x)\right)\left(\sum_{i=0}^{m-1}\left(\begin{bmatrix}i\\i_0\end{bmatrix}\cdot\left(\frac{d^{i_0-i}}{dx^{i_0-i}}b_{i+1}(x)\right)\right)\right)\right)$$

Comparing this formula with (10.25), we can see that they are identical. Again, assuming that $n=i_0$ in (10.25), we obtain:

$$\xi_1(i_0) = \sum_{i_1=0}^{m-1}\left(\sum_{k_0=1}^{m}\left(\sum_{i_2=0}^{i_0-k_0}\left(\begin{bmatrix} k_0+i_2 \\ i_0 \end{bmatrix}\cdot\begin{bmatrix} i_1 \\ i_2 \end{bmatrix}\right)\left(\frac{d^{i_0-k_0-i_2}}{dx^{i_0-k_0-i_2}}b_{k_0}(x)\right)\left(\frac{d^{i_2-i_1}}{dx^{i_2-i_1}}b_{i_1+1}(x)\right)\right)\right) \qquad (10.27)$$

Assuming that $n = 1$ in (10.23) we obtain:

$$\xi_2(n) = \sum_{k=1}^{m}\left(\sum_{i_0=0}^{n-k}\begin{bmatrix} i_0+k \\ n \end{bmatrix}\left(\frac{d^{n-i_0-k}}{dx^{n-i_0-k}}b_k(x)\right)\xi_1(i_0)\right)$$

Substituting here (10.27) we have: (10.28)

$$\xi_2(n) = \sum_{k=1}^{m}\left(\sum_{i_0=0}^{n-k}\begin{bmatrix} i_0+k \\ n \end{bmatrix}\left(\frac{d^{n-i_0-k}}{dx^{n-i_0-k}}b_k(x)\right)\right.$$
$$\left.\left(\sum_{i_1=0}^{m-1}\left(\sum_{k_0=1}^{m}\left(\sum_{i_2=0}^{i_0-k_0}\left(\begin{bmatrix} k_0+i_2 \\ i_0 \end{bmatrix}\cdot\begin{bmatrix} i_1 \\ i_2 \end{bmatrix}\right)\left(\frac{d^{i_0-k_0-i_2}}{dx^{i_0-k_0-i_2}}b_{k_0}(x)\right)\left(\frac{d^{i_2-i_1}}{dx^{i_2-i_1}}b_{i_1+1}(x)\right)\right)\right)\right)\right)$$

Comparing this value of $\xi_2(n)$ with the similar one obtained previously and defined by (10.20), we can confirm their identity.

Let's assume that (10.23) is solved up to $s = p$. Let’s prove that it holds true for $s = p$ as well. Indeed, in this case (10.23) takes the following form:

$$\xi_{p+2}(n) = \sum_{k=1}^{m}\left(\sum_{i_0=0}^{n-k}\begin{bmatrix} i_0+k \\ n \end{bmatrix}\left(\frac{d^{n-i_0-k}}{dx^{n-i_0-k}}b_k(x)\right)\xi_{p+1}(i_0)\right) \qquad (10.29)$$

The value of $\xi_{p+1}(i_0)$ is defined with $s = p$. Therefore, from (10.23) we have:

$$\xi_{p+1}(n) = \sum_{k=1}^{m}\left(\sum_{i_0=0}^{n-k}\begin{bmatrix} i_0+k \\ n \end{bmatrix}\left(\frac{d^{n-i_0-k}}{dx^{n-i_0-k}}b_k(x)\right)\xi_p(i_0)\right)$$

The value of $\xi_p(i_0)$ is also defined with $s = p-1$. Therefore, from (10.23) we have:

$$\xi_p(n) = \sum_{k=1}^{m}\left(\sum_{i_0=0}^{n-k}\begin{bmatrix} i_0+k \\ n \end{bmatrix}\left(\frac{d^{n-i_0-k}}{dx^{n-i_0-k}}b_k(x)\right)\xi_{p-1}(i_0)\right)$$

Since, by definition, (10.23) is defined up to $s = p$, we eventually obtain on the $p+1$ step:

$$\xi_1(n) = \sum_{k=1}^{m}\left(\sum_{i_0=0}^{n-k}\begin{bmatrix} i_0+k \\ n \end{bmatrix}\left(\frac{d^{n-i_0-k}}{dx^{n-i_0-k}}b_k(x)\right)\xi_0(i_0)\right)$$

Successively substituting all these values to (10.29), we get:

$$\xi_{p+2}(n)=\sum_{k=1}^{m}\left(\sum_{i=0}^{-k+n}\left(\begin{bmatrix}k+i\\ n\end{bmatrix}\cdot\left(\left(\frac{d^{n-k-i}}{dx^{n-k-i}}b_k(x)\right)\left(\sum_{k_1=1}^{m}\left(\sum_{i_1=0}^{i-k_1}\begin{bmatrix}i_1+k_1\\ i\end{bmatrix}\left(\frac{d^{i-i_1-k_1}}{dx^{i-i_1-k_1}}b_{k_1}(x)\right)\right)\right(\right.\right.\right.$$

$$\sum_{k_2=1}^{m}\left(\sum_{i_2=0}^{i_1-k_2}\begin{bmatrix}i_2+k_2\\ i_1\end{bmatrix}\left(\frac{d^{i_1-i_2-k_2}}{dx^{i_1-i_2-k_2}}b_{k_2}(x)\right)[[0,0,0]]\left(\sum_{k_p=1}^{m}\left(\sum_{i_p=0}^{i_{p-1}-k_p}\begin{bmatrix}i_{N-2}+k_{N-2}\\ i_{N-3}\end{bmatrix}\right.\right.\right.$$

$$\left(\frac{d^{i_{p-1}-i_p-k_p}}{dx^{i_{p-1}-i_p-k_p}}b_{k_p}(x)\right)\left(\sum_{k_{p+1}=1}^{m}\left(\sum_{i_{p+1}=0}^{i_p-k_{p+1}}\right.\right.$$

$$\left.\left.\left.\begin{bmatrix}i_{p+1}+k_{p+1}\\ i_p\end{bmatrix}\left(\frac{d^{i_p-i_{p+1}-k_{p+1}}}{dx^{i_p-i_{p+1}-k_{p+1}}}b_{k_{p+1}}(x)\right)\left(\sum_{i_{p+2}=0}^{m-1}\begin{bmatrix}i_{p+2}\\ i_{p+1}\end{bmatrix}\left(\frac{d^{i_{p+1}-i_{p+2}}}{dx^{i_{p+1}-i_{p+2}}}b_{1+i_{p+2}}(x)\right)\right)\right)\right)\right)$$

However, this formula is fully identical to (10.22) with $N=p+2$.
**The theorem is proved.**

Thus, if we introduce the differential operator $\Xi_n[\;]$, which operates on the $\mathrm{f}(x)$ function and depends upon the ***n*** parameter in accordance with the following formula:

$$\Xi_n[\mathrm{f}(x)]=\sum_{k=1}^{m}\left(\sum_{i_0=0}^{n-k}\begin{bmatrix}i_0+k\\ n\end{bmatrix}\left(\frac{d^{n-i_0-k}}{dx^{n-i_0-k}}b_k(x)\right)[\mathrm{f}(x)]\right)\tag{10.30}$$

then

(10.31)

(10.32)

$$\xi_3(n)=\Xi_n[\Xi_{i_1}[\Xi_{i_0}[\xi_0(n)]]]\tag{10.33}$$

$$\xi_4(n)=\Xi_n[\Xi_{i_2}[\Xi_{i_1}[\Xi_{i_0}[\xi_0(n)]]]]\tag{10.34}$$

---

$$\xi_k(n)=\Xi_n[\Xi_{i_{k-2}}[\Xi_{i_{k-3}}[[[0,0,0]]\Xi_{i_2}[\Xi_{i_1}[\Xi_{i_0}[\xi_0(n)]]]]]]\tag{10.35}$$

Consequently, the formula for the first coefficient of *n*-image takes the following form:

$$\alpha(n)=$$
$$\xi_0(n)+\Xi_n[\xi_0(n)]+\left(\sum_{k=2}^{N}\Xi_n[\Xi_{i_{k-2}}[\Xi_{i_{k-3}}[[[0,0,0]]\Xi_{i_2}[\Xi_{i_1}[\Xi_{i_0}[\xi_0(n)]]]]]]\right)\tag{10.36}$$

**Note: it is obvious that**

$$\Xi_{i_0}[\mathrm{f}(x)]=\sum_{k=1}^{m}\left(\sum_{i_1=0}^{i_0-k}\begin{bmatrix}i_1+k\\ i_0\end{bmatrix}\left(\frac{d^{i_0-i_1-k}}{dx^{i_0-i_1-k}}b_k(x)\right)[\mathrm{f}(x)]\right)\tag{10.37}$$

and so on.

## 10.4. The definition of serial components for the special $\alpha(n)$ function.

The summands are defined by the formulas (13.31)...(13.37), represented as partial summable sequences. However, in such form these coefficients cannot be used to determine the desired special functions: $\alpha(n), n = -1, -2, -3 .. -m$, because we cannot identify their values without calculation of the sum of $\alpha(n)$. Thus, **the calculation of partial solutions for linear ODEs with variable coefficients** is impossible without the calculation of the sums of (13.31)...(13.37) series. Let's consider this solution successively for the summands, up to the definition of the common member of this sequence.

**10.4.1. The $\xi_0(n)$ coefficient** is defined by (10.12):

$$\xi_0(n) = \sum_{i=0}^{m-1}\left(\left[\begin{matrix} i \\ n \end{matrix}\right]\cdot\left(\frac{d^{n-i}}{dx^{n-i}} b_{i+1}(x)\right)\right)$$

Let's represent it as follows:

$$\xi_0(n) = \frac{\partial^n}{\partial x^n}\left(\sum_{i=0}^{m-1}\left[\begin{matrix} i \\ n \end{matrix}\right][b_{i+1}(x)]_i\right) \qquad (10.38)$$

Because the ***n*** parameter does not fall within the limits of summation, this formula is definitive and not subject to alteration.

**10.4.2. The $\xi_1(n)$ coefficient** is defined by (10.16)

$$\xi_1(n) = \sum_{i_1=0}^{m-1}\left(\sum_{k=1}^{m}\left(\sum_{i=0}^{n-k}\left(\left[\begin{matrix} k+i \\ n \end{matrix}\right]\cdot\left[\begin{matrix} i_1 \\ i \end{matrix}\right]\right)\left(\frac{d^{n-k-i}}{dx^{n-k-i}} b_k(x)\right)\left(\frac{d^{i-i_1}}{dx^{i-i_1}} b_{i_1+1}(x)\right)\right)\right) \qquad (10.39)$$

To calculate the sum of this series we shall use (10.23):

$$\xi_{s+1}(n) = \sum_{k=1}^{m}\left(\sum_{i_0=0}^{n-k}\left[\begin{matrix} i_0+k \\ n \end{matrix}\right]\left(\frac{d^{n-i_0-k}}{dx^{n-i_0-k}} b_k(x)\right)\xi_s(i_0)\right)$$

By the substitution of the summation index, $i_0 = z - k$, this formula is reduced as follows:

$$\xi_{s+1}(n) = \sum_{k=1}^{m}\left(\sum_{z=k}^{n}\left[\begin{matrix} z \\ n \end{matrix}\right]\left(\frac{d^{n-z}}{dx^{n-z}} b_k(x)\right)\xi_s(z-k)\right)$$

Let's rerpesent it as follows: (10.40)

$$\xi_{s+1}(n) = \left(\sum_{k=1}^{m}\left(\sum_{z=0}^{n}\left[\begin{matrix} z \\ n \end{matrix}\right]\left(\frac{d^{n-z}}{dx^{n-z}} b_k(x)\right)\xi_s(z-k)\right)\right)$$
$$- \left(\sum_{k=1}^{m}\left(\sum_{z=0}^{k-1}\left[\begin{matrix} z \\ n \end{matrix}\right]\left(\frac{d^{n-z}}{dx^{n-z}} b_k(x)\right)\xi_s(z-k)\right)\right)$$

This representation is most appropriate for the application of Leibniz formula; therefore, further we shall use it as the basic one to define the desired values of summands.

Assuming here that $s = 0$, we get:

$$\xi_1(n) = \left(\sum_{k=1}^{m}\left(\sum_{z=0}^{n}\begin{bmatrix} z \\ n \end{bmatrix}\left(\frac{d^{n-z}}{dx^{n-z}} b_k(x)\right)\xi_0(z-k)\right)\right) - \left(\sum_{k=1}^{m}\left(\sum_{z=0}^{k-1}\begin{bmatrix} z \\ n \end{bmatrix}\left(\frac{d^{n-z}}{dx^{n-z}} b_k(x)\right)\xi_0(z-k)\right)\right) \qquad (10.41)$$

In (10.38) let's assume that $n = z - k$. Then it takes the following form:

$$\xi_0(z-k) = \frac{\partial^{z-k}}{\partial x^{z-k}}\left(\sum_{i=0}^{m-1}\begin{bmatrix} i \\ z-k \end{bmatrix}[b_{i+1}(x)]_i\right)$$

Substituting this value to the first summand in the right-hand side of (10.41), we get (10.42)

$$\xi_1(n) = \left(\sum_{k=1}^{m}\left(\sum_{z=0}^{n}\begin{bmatrix} z \\ n \end{bmatrix}\left(\frac{d^{n-z}}{dx^{n-z}} b_k(x)\right)\left(\frac{\partial^{z-k}}{\partial x^{z-k}}\left(\sum_{i=0}^{m-1}\begin{bmatrix} i \\ z-k \end{bmatrix}[b_{i+1}(x)]_i\right)\right)\right)\right) - \left(\sum_{k=1}^{m}\left(\sum_{z=0}^{k-1}\begin{bmatrix} z \\ n \end{bmatrix}\left(\frac{d^{n-z}}{dx^{n-z}} b_k(x)\right)\xi_0(z-k)\right)\right)$$

Let's calculate the value of the sum:

$$\sum_{z=0}^{n}\begin{bmatrix} i \\ z-k \end{bmatrix}\begin{bmatrix} z \\ n \end{bmatrix}\left(\frac{\partial^{z}}{\partial x^{z}}[b_{i+1}(x)]_{i+k}\right)\left(\frac{d^{n-z}}{dx^{n-z}} b_k(x)\right)$$

To do this, let's use the Proposition 1.4., that is, (1.11):

$$L[\mathrm{C}(z-k,i)] = \left(\frac{1}{i!}\right)\cdot\left(\sum_{s_0=0}^{i}\left(\sum_{s_1=0}^{s_0}\left(\sum_{s_2=0}^{s_1}\left(\sum_{s_3=0}^{s_2}\frac{s_3^{\,s_1}(-1)^{(s_3+s_2)}}{\Gamma(s_2-s_3+1)\,\Gamma(s_3+1)}\right)\mathrm{C}(n,s_2)\,s_2!\; r_{i-s_0+1}(i)\,\mathrm{C}(s_0,s_0-s_1)\,(-k)^{(s_0-s_1)}\left(\frac{\partial^{n-s_2}}{\partial x^{n-s_2}}\left(\left(\frac{\partial^{s_2}}{\partial x^{s_2}}u\right)v\right)\right)\right)\right)\right)$$

Then we shall get

$$\sum_{z=0}^{n}\left[\begin{matrix} i \\ z-k \end{matrix}\right]\left[\begin{matrix} z \\ n \end{matrix}\right]\left(\frac{\partial^{z}}{\partial x^{z}}[b_{i+1}(x)]_{i+k}\right)\left(\frac{d^{n-z}}{dx^{n-z}}b_k(x)\right)=\left(\frac{1}{i!}\right)\cdot\left(\sum_{s_0=0}^{i}\left(\sum_{s_1=0}^{s_0}\left(\sum_{s_2=0}^{s_1}\left(\sum_{s_3=0}^{s_2}\frac{s_3^{s_1}(-1)^{(s_3+s_2)}}{\Gamma(s_2-s_3+1)\Gamma(s_3+1)}\right)\mathrm{C}(n,s_2)\,s_2!\,r_{i-s_0+1}(i)\,\mathrm{C}(s_0,s_0-s_1)(-k)^{(s_0-s_1)}\left(\frac{\partial^{n-s_2}}{\partial x^{n-s_2}}\left(\left(\frac{\partial^{s_2}}{\partial x^{s_2}}[b_{i+1}(x)]_{i+k}\right)b_k(x)\right)\right)\right)\right)\right)$$

Substituting this value to (10.42), we obtain:

$$\xi_1(n)=\left(\sum_{k=1}^{m}\left(\sum_{i=0}^{m-1}\left(\left(\frac{1}{i!}\right)\cdot\left(\sum_{s_0=0}^{i}\left(\sum_{s_1=0}^{s_0}\left(\sum_{s_2=0}^{s_1}\left(\sum_{s_3=0}^{s_2}\frac{s_3^{s_1}(-1)^{(s_3+s_2)}}{\Gamma(s_2-s_3+1)\Gamma(s_3+1)}\right)\mathrm{C}(n,s_2)\,s_2!\,r_{i-s_0+1}(i)\,\mathrm{C}(s_0,s_0-s_1)(-k)^{(s_0-s_1)}\left(\frac{\partial^{n-s_2}}{\partial x^{n-s_2}}\left(\left(\frac{\partial^{s_2}}{\partial x^{s_2}}[b_{i+1}(x)]_{i+k}\right)b_k(x)\right)\right)\right)\right)\right)\right)\right)\right)\right)-\left(\sum_{k=1}^{m}\left(\sum_{z=0}^{k-1}\left[\begin{matrix} z \\ n \end{matrix}\right]\left(\frac{d^{n-z}}{dx^{n-z}}b_k(x)\right)\xi_0(z-k)\right)\right)$$

However, this formula can be further simplified. For this purpose, let’s use (1.16):

$$\sum_{i_0=0}^{m-1}L_{i_0,k_0}[\mathrm{C}(i-z,i_0)]=\sum_{i_1=0}^{m-1}\left(\sum_{i_0=0}^{m-1-i_1}(-1)^{i_0}\,\mathrm{C}(k_0-1+i_0,k_0-1)\,\mathrm{C}(n,i_1)\left(\frac{\partial^{n-i_1}}{\partial x^{n-i_1}}\left([b_{i_1+i_0}(x)]_{i_0+k_0}\,b_{k_0}(x)\right)\right)\right)$$

Upon transformations we obtain the final formula: (10.43)

$$(-1)^{i_0}\,\mathrm{C}(n,i_1)\left(\sum_{k_0=1}^{m}\mathrm{C}(k_0-1+i_0,k_0-1)\left(\frac{\partial^{n-i_1}}{\partial x^{n-i_1}}\left([b_{i_1+i_0+1}(x)]_{i_0+k_0}\,b_{k_0}(x)\right)\right)\right)-\left(\sum_{i_1=1}^{m}\left(\sum_{i_0=0}^{m-i_1}\mathrm{C}(n,i_0)\left(\frac{d^{n-i_0}}{dx^{n-i_0}}b_{i_0+i_1}(x)\right)\right)\xi_0(-i_1)\right)$$

**10.4.3. The $\xi_2(n)$ coefficient** is defined by (10.20)

$$\xi_2(n)=\sum_{k=1}^{m}\left(\sum_{i=0}^{-k+n}\left(\left[\begin{matrix} k+i \\ n \end{matrix}\right]\cdot\left(\left(\frac{d^{n-k-i}}{dx^{n-k-i}}b_k(x)\right)\left(\sum_{k_1=1}^{m}\left(\sum_{i_1=0}^{i-k_1}\left[\begin{matrix} i_1+k_1 \\ i \end{matrix}\right]\left(\frac{d^{i-i_1-k_1}}{dx^{i-i_1-k_1}}b_{k_1}(x)\right)\left(\sum_{i_2=0}^{m-1}\left[\begin{matrix} i_2 \\ i_1 \end{matrix}\right]\left(\frac{d^{i_1-i_2}}{dx^{i_1-i_2}}b_{1+i_2}(x)\right)\right)\right)\right)\right)\right)\right)$$

In order to calculate the sum of this series, let’s use the recurrence relation (10.40) where we assume $s = 1$ :

$$\xi_2(n) = \left( \sum_{k=1}^{m} \left( \sum_{z=0}^{n} \left[ \begin{matrix} z \\ n \end{matrix} \right] \left( \frac{d^{n-z}}{dx^{n-z}} b_k(x) \right) \xi_1(z-k) \right) \right)$$
$$- \left( \sum_{k=1}^{m} \left( \sum_{z=0}^{k-1} \left[ \begin{matrix} z \\ n \end{matrix} \right] \left( \frac{d^{n-z}}{dx^{n-z}} b_k(x) \right) \xi_1(z-k) \right) \right) \qquad (10.44)$$

Assuming in (10.43) $n = z - k$ we have:

$$\xi_1(z-k) = \left( \sum_{i_1=0}^{m-1} \left( \sum_{i_0=0}^{m-1-i_1} (-1)^{i_0} \mathrm{C}(z-k, i_1) \right.\right.$$
$$\left( \sum_{k_0=1}^{m} \mathrm{C}(k_0 - 1 + i_0, k_0 - 1) \left( \frac{\partial^{z-k-i_1}}{\partial x^{z-k-i_1}} \left( [b_{i_1+i_0+1}(x)]_{i_0+k_0} \; b_{k_0}(x) \right) \right) \right) \Bigg) \Bigg)$$
$$- \left( \sum_{i_1=1}^{m} \left( \sum_{i_0=0}^{m-i_1} \mathrm{C}(z-k, i_0) \left( \frac{d^{z-k-i_0}}{dx^{z-k-i_0}} b_{i_0+i_1}(x) \right) \right) \xi_0(-i_1) \right)$$

Substituting this value to (10.44), we have: (10.45)

$$\xi_2(n) = \left( \sum_{k=1}^{m} \left( \sum_{z=0}^{n} \left[ \begin{matrix} z \\ n \end{matrix} \right] \left( \frac{d^{n-z}}{dx^{n-z}} b_k(x) \right) \left( \left( \sum_{i_1=0}^{m-1} \left( \sum_{i_0=0}^{m-1-i_1} (-1)^{i_0} \mathrm{C}(z-k, i_1) \right.\right.\right.\right.\right.$$
$$\left( \sum_{k_0=1}^{m} \mathrm{C}(k_0 - 1 + i_0, k_0 - 1) \left( \frac{\partial^{z-k-i_1}}{\partial x^{z-k-i_1}} \left( [b_{i_1+i_0+1}(x)]_{i_0+k_0} \; b_{k_0}(x) \right) \right) \right) \Bigg) \Bigg)$$
$$- \left( \sum_{i_1=1}^{m} \left( \sum_{i_0=0}^{m-i_1} \mathrm{C}(z-k, i_0) \left( \frac{d^{z-k-i_0}}{dx^{z-k-i_0}} b_{i_0+i_1}(x) \right) \right) \xi_0(-i_1) \right) \Bigg) \Bigg) \Bigg)$$
$$- \left( \sum_{k=1}^{m} \left( \sum_{z=0}^{k-1} \left[ \begin{matrix} z \\ n \end{matrix} \right] \left( \frac{d^{n-z}}{dx^{n-z}} b_k(x) \right) \xi_1(z-k) \right) \right)$$

Hence, we need to define the sum of the partial summable sequences:

$$\sum_{z=0}^{n} \left[ \begin{matrix} z \\ n \end{matrix} \right] \mathrm{C}(z-k, i_1) \left( \frac{d^{n-z}}{dx^{n-z}} b_k(x) \right) \left( \frac{\partial^{z-k-i_1}}{\partial x^{z-k-i_1}} \left( [b_{i_1+i_0+1}(x)]_{i_0+k_0} \; b_{k_0}(x) \right) \right)$$

$$\sum_{z=0}^{n} \left[ \begin{matrix} z \\ n \end{matrix} \right] \mathrm{C}(z-k, i_0) \left( \frac{d^{n-z}}{dx^{n-z}} b_k(x) \right) \left( \frac{d^{z-k-i_0}}{dx^{z-k-i_0}} b_{i_0+i_1}(x) \right)$$

In accordance with (1.11), we have:

$$\sum_{z=0}^{n} \left[ \begin{matrix} z \\ n \end{matrix} \right] C(z-k, i_1) \left( \frac{d^{n-z}}{dx^{n-z}} b_k(x) \right) \left( \frac{\partial^z}{\partial x^z} \left[ [b_{i_1+i_0+1}(x)]_{i_0+k_0} \; b_{k_0}(x) \right]_{k+i_1} \right) = \left( \frac{1}{i_1!} \right) \cdot \left( \sum_{s_0=0}^{i_1} \left( \sum_{s_1=0}^{s_0} \left( \sum_{s_2=0}^{s_1} \left( \sum_{s_3=0}^{s_2} \frac{s_3^{\,s_1} (-1)^{(s_3+s_2)}}{\Gamma(s_2-s_3+1)\,\Gamma(s_3+1)} \right) C(n, s_2)\, s_2!\, r_{i_1-s_0+1}(i_1)\, C(s_0, s_0-s_1) \right.\right.\right.$$

$$\left.\left.\left. (-k)^{(s_0-s_1)} \left( \frac{\partial^{n-s_2}}{\partial x^{n-s_2}} \left( \left( \frac{\partial^{s_2}}{\partial x^{s_2}} \left[ [b_{i_1+i_0+1}(x)]_{i_0+k_0} \; b_{k_0}(x) \right]_{k+i_1} \right) b_k(x) \right) \right) \right) \right) \right)$$

$$\sum_{z=0}^{n} \left[ \begin{matrix} z \\ n \end{matrix} \right] C(z-k, i_0) \left( \frac{d^{n-z}}{dx^{n-z}} b_k(x) \right) \left( \frac{\partial^z}{\partial x^z} [b_{i_0+i_1}(x)]_{k+i_0} \right) = \left( \frac{1}{i_0!} \right) \cdot \left( \sum_{s_0=0}^{i_0} \left( \sum_{s_1=0}^{s_0} \left( \sum_{s_2=0}^{s_1} \right.\right.\right.$$

$$\left( \sum_{s_3=0}^{s_2} \frac{s_3^{\,s_1} (-1)^{(s_3+s_2)}}{\Gamma(s_2-s_3+1)\,\Gamma(s_3+1)} \right) C(n, s_2)\, s_2!\, r_{i_0-s_0+1}(i_0)\, C(s_0, s_0-s_1)\, (-k)^{(s_0-s_1)}$$

$$\left.\left.\left. \left( \frac{\partial^{n-s_2}}{\partial x^{n-s_2}} \left( \left( \frac{\partial^{s_2}}{\partial x^{s_2}} [b_{i_0+i_1}(x)]_{k+i_0} \right) b_k(x) \right) \right) \right) \right) \right)$$

Substituting these values to (10.45), upon transformations with consideration of (1.16), we eventually obtain: (10.46)

$$\xi_2(n) = \left( \sum_{i_2=0}^{m-1} \left( \sum_{i_1=0}^{m-1-i_2} \left( \sum_{i_0=0}^{m-1-i_1-i_2} (-1)^{(i_0+i_1)} C(n, i_2) \left( \sum_{k_1=1}^{m} \left( \sum_{k_0=1}^{m} C(k_0-1+i_0, k_0-1) \right.\right.\right.\right.\right.$$

$$\left.\left.\left.\left.\left. C(k_1-1+i_1, k_1-1) \left( \frac{\partial^{n-i_2}}{\partial x^{n-i_2}} \left( \left[ [b_{i_2+i_1+i_0+1}(x)]_{i_0+k_0} \; b_{k_0}(x) \right]_{i_1+k_1} b_{k_1}(x) \right) \right) \right) \right) \right) \right) - \left( \right.$$

$$\sum_{i_0=0}^{m-i_1-i_2} (-1)^{i_0} C(n, i_1) \left( \sum_{k_0=1}^{m} C(k_0-1+i_0, k_0-1) \left( \frac{\partial^{n-i_1}}{\partial x^{n-i_1}} \left( [b_{i_1+i_0+i_2}]_{i_0+k_0} \; b_{k_0}(x) \right) \right) \right)$$

$$\left.\left.\left. \right) \right) \xi_0(-i_2) \right) - \left( \sum_{i_1=1}^{m} \left( \sum_{i_0=0}^{m-i_1} C(n, i_0) \left( \frac{d^{n-i_0}}{dx^{n-i_0}} b_{i_0+i_1}(x) \right) \right) \xi_1(-i_1) \right)$$

Quite similary, we obtain: (13.47)

$$
\begin{aligned}
\xi_3(n) = {} & \left( \sum_{i_3=0}^{m-1} \left( \sum_{i_2=0}^{m-1-i_3} \left( \sum_{i_1=0}^{m-1-i_2-i_3} \left( \sum_{i_0=0}^{m-1-i_1-i_2-i_3} (-1)^{(i_0+i_1+i_2)} \mathrm{C}(n, i_3) \left( \sum_{k_2=1}^{m} \left( \sum_{k_1=1}^{m} \left( \sum_{k_0=1}^{m} \right.\right.\right.\right.\right.\right.\right. \\
& \mathrm{C}(k_0-1+i_0, k_0-1)\, \mathrm{C}(k_1-1+i_1, k_1-1)\, \mathrm{C}(k_2-1+i_2, k_2-1) \\
& \left.\left.\left.\left.\left.\left.\left. \left( \frac{\partial^{n-i_3}}{\partial x^{n-i_3}} \left( \left[ \left[ \left[ b_{i_3+i_2+i_1+i_0+1}(x) \right]_{i_0+k_0} b_{k_0}(x) \right]_{i_1+k_1} b_{k_1}(x) \right]_{i_2+k_2} b_{k_2}(x) \right) \right) \right) \right) \right) \right) \right) \right) \right) - \left( \right. \\
& \sum_{i_3=1}^{m} \left( \sum_{i_2=0}^{m-i_3} \left( \sum_{i_1=0}^{m-i_3-i_2} \left( \sum_{i_0=0}^{m-i_3-i_1-i_2} (-1)^{(i_0+i_1)} \mathrm{C}(n, i_2) \left( \sum_{k_1=1}^{m} \left( \sum_{k_0=1}^{m} \mathrm{C}(k_0-1+i_0, k_0-1) \right.\right.\right.\right.\right.\right. \\
& \left.\left.\left.\left.\left.\left. \mathrm{C}(k_1-1+i_1, k_1-1) \left( \frac{\partial^{n-i_2}}{\partial x^{n-i_2}} \left( \left[ \left[ b_{i_2+i_1+i_0+i_3}(x) \right]_{i_0+k_0} b_{k_0}(x) \right]_{i_1+k_1} b_{k_1}(x) \right) \right) \right) \right) \right) \right) \right) \right) \\
& \left. \xi_0(-i_3) \right) - \left( \sum_{i_2=1}^{m} \left( \sum_{i_1=0}^{m-i_2} \left( \right.\right.\right. \\
& \sum_{i_0=0}^{m-i_1-i_2} (-1)^{i_0} \mathrm{C}(n, i_1) \left( \sum_{k_0=1}^{m} \mathrm{C}(k_0-1+i_0, k_0-1) \left( \frac{\partial^{n-i_1}}{\partial x^{n-i_1}} \left( \left[ b_{i_1+i_0+i_2} \right]_{i_0+k_0} b_{k_0}(x) \right) \right) \right) \\
& \left.\left.\left. \right) \right) \xi_1(-i_2) \right) - \left( \sum_{i_1=1}^{m} \left( \sum_{i_0=0}^{m-i_1} \mathrm{C}(n, i_0) \left( \frac{d^{n-i_0}}{dx^{n-i_0}} b_{i_0+i_1}(x) \right) \right) \xi_2(-i_1) \right)
\end{aligned}
$$

and so on, in accordance with the program written in Maple 10.

## The program for printing $\xi_{m,s}(-n)$ coefficients

```
> restart;
```

```
> m:=m:s:=3: # Enter the order of the ODE - "m" and the number of the
coefficient - "s".
```

$> \; C0 := h \rightarrow \prod_{l=0}^{h-1} \mathrm{C}(k_l - 1 + i_l, k_l - 1)$

$> \; P0 := h \rightarrow \sum_{i_1=1}^{m} \left( \sum_{i_0=0}^{m-i_1} \mathrm{C}(-n, i_0) \left[ b_{i_0+i_1}(x) \right]_{n+i_0} \right) \xi_{m, h-1}(-i_1)$

```
> Ii:=h->add(i[j],j=0..h):
```

$> \; C1 := h \rightarrow (-1)^{\mathrm{Ii}(h-1)} \mathrm{C}(-n, i_h)$

```
> Iik:=h->i[h]+k[h]:
> B0(0):=b[1+Ii(s)](x):
```

> **for** $z$ **from** 0 **to** $s$ **do** $\mathrm{B0}(z+1) := [\mathrm{B0}(z)]_{\mathrm{Iik}(z)}\, b_{k_z}(x)$ **end do**; $z := 'z'$

```
> R0(0):=C0(s)*[B0(s)][n+i[s]]:
```

> **for** $z$ **from** 0 **to** $s$ **do** $\mathrm{R0}(z+1) := \sum_{k_z=1}^{m} \mathrm{R0}(z)$ **end do**; $z := 'z'$

```
> Ii1:=h->add(i[j],j=h..z):
> Ii2:=h->add(i[j],j=h..s):
> R1(0):=C1(s)*R0(s):
```

> **for** $z$ **from** 0 **to** $s$ **do** $\mathrm{R1}(z+1) := \sum_{i_z=0}^{m-1-\mathrm{Ii2}(z+1)} \mathrm{R1}(z)$ **end do**; $z := 'z'$

```
> R2(s):=R1(s+1)-P0(s):
> Ii2:=h->h:
> for z from 2 to s do:
> B1(0):=b[add(i[j],j=0..z)](x);
> for z0 from 0 to z-1 do B1(z0+1):=[B1(z0)][Iik(z0)]*b[k[z0]](x)
od:z0:='z0':
> R3(0):=Product(C(k[l]-1+i[l],k[l]-1),l = 0 .. z-2)*[B1(z-1)][n+i[z-1]];
> for z0 from 0 to z do R3(z0+1):=Sum(R3(z0),k[z0]=1..m) od;z0:='z0':
> R4(0):=(-1)^(add(i[j],j=0..z-2))*C(-n,i[z-1])*R3(z-1);
> for z0 from 0 to z do R4(z0+1):=Sum(R4(z0),i[z0]=0..m-Ii1(z0+1))
od;z0:='z0':
> R5[z]:=Sum(R4(z)*xi[m,s-z](-i[z]),i[z]=1..m);
> od:
> xi[m,s](-n)=R2(s)-add(R5[z1],z1=2..s);
>
##########################################################################
###############
```

In a general case, we can write down for any (10.48)

$$\xi_{m,s}(n) = \left( \sum_{i_{[s,0]}=0}^{[m-1,\, m-1-(s,s),\, m-1-(1,s)]} (-1)^{(0,s-1)} C(n, i_s) \left( \sum_{k_{[s-1,0]}=1}^{m} \left( \prod_{l=0}^{s-1} C(k_l - 1 + i_l, k_l - 1) \right) \left( \frac{\partial^{n-i_s}}{\partial x^{n-i_s}} \left( [b_{1+(0,s)}(x)]^{i_{[0,s-1]} + k_{[0,s-1]}} \, b_{k_{[0,s-1]}}(x) \right) \right) \right) \right)$$

$$- \left( \sum_{i_1=1}^{m} \left( \sum_{i_0=0}^{m-i_1} C(n, i_0) \left( \frac{d^{n-i_0}}{dx^{n-i_0}} b_{i_0+i_1}(x) \right) \right) \xi_{s-1}(-i_1) \right) - \left( \sum_{z=2}^{s} \left( \sum_{i_{[z,0]}=[1,0]}^{[m,\, m-(z,z),\, m-(z-1,z),\, m-(1,z)]} \left( \sum_{k_{[0,z-2]}=1}^{m} (-1)^{(0,z-2)} C(n, i_{z-1}) \right. \right. \right.$$

$$\left. \left. \left. \left( \prod_{l=0}^{z-2} C(k_l - 1 + i_l, k_l - 1) \right) \left( \frac{\partial^{n-i_{z-2+1}}}{\partial x^{n-i_{z-2+1}}} \left( [b_{[0,z]}(x)]^{i_{[0,z-2]} + k_{[0,z-2]}} \, b_{k_{[0,z-2]}}(x) \right) \right) \right) \right) \right)$$

Thus, we have solved the task of finding the functions, which fully define the desired solution.

## 11. The algorithm for linear ODEs. Examples.

**Summing up the above, taking into account the relevance of this topic, let's describe in detail the algorithm for general solutions of linear ODEs with variable coefficients.**

Let's consider a linear ODE of $m$ order:

$$\frac{\partial^m}{\partial x^m} Y = \sum_{p=1}^{m} a_p(x) \left( \frac{\partial^{m-p}}{\partial x^{m-p}} Y \right) \tag{11.1}$$

Here $Y = Y_m$ is the desired function, and $a_p(x)$ are the specified functions.

The general solution of ODE (11.1) is defined by the following formula:

(11.2)

where $C_i$ are arbitrary constants and partial solutions of $Y_{m,i}(x)$ are defined by the following formulas:

$$Y_{m,i}(x) = \frac{(-1)^i (1 - \alpha_m(-1))\, x^{(i-1)}}{(i-1)!} + \left( \sum_{k=2}^{i} \frac{(-1)^{(k+i)} \alpha_m(-k)\, x^{(i-k)}}{(i-k)!} \right) \tag{11.3}$$

Here

(11.4)

$\xi_{m,i}(n)$ are the desired functions defined by the following formulas

(11.5)

$$\xi_{m,0}(n) = \frac{\partial^n}{\partial x^n}\left( \sum_{i=0}^{m-1} \begin{bmatrix} i \\ n \end{bmatrix} [b_{i+1}(x)]_i \right)$$

$$(-1)^{i_0} C(n, i_1) \left( \sum_{k_0=1}^{m} C(k_0 - 1 + i_0, k_0 - 1) \left( \frac{\partial^{n-i_1}}{\partial x^{n-i_1}} \left( [b_{i_1+i_0+1}(x)]_{i_0+k_0} \; b_{k_0}(x) \right) \right) \right) ) ) )$$

$$- \left( \sum_{i_1=1}^{m} \left( \sum_{i_0=0}^{m-i_1} C(n, i_0) \left( \frac{d^{n-i_0}}{dx^{n-i_0}} b_{i_0+i_1}(x) \right) \right) \xi_0(-i_1) \right)$$

$$\xi_{m,2}(n) = \left( \sum_{i_2=0}^{m-1} \left( \sum_{i_1=0}^{m-1-i_2} \left( \sum_{i_0=0}^{m-1-i_1-i_2} (-1)^{(i_0+i_1)} C(n, i_2) \left( \sum_{k_1=1}^{m} \left( \sum_{k_0=1}^{m} C(k_0 - 1 + i_0, k_0 - 1) \right.\right.\right.\right.\right.$$

$$C(k_1 - 1 + i_1, k_1 - 1) \left( \frac{\partial^{n-i_2}}{\partial x^{n-i_2}} \left( \left[ [b_{i_2+i_1+i_0+1}(x)]_{i_0+k_0} \; b_{k_0}(x) \right]_{i_1+k_1} \; b_{k_1}(x) \right) \right) ) ) ) ) ) ) - \Bigg($$

$$\sum_{i_0=0}^{m-i_1-i_2} (-1)^{i_0} C(n, i_1) \left( \sum_{k_0=1}^{m} C(k_0 - 1 + i_0, k_0 - 1) \left( \frac{\partial^{n-i_1}}{\partial x^{n-i_1}} \left( [b_{i_1+i_0+i_2}]_{i_0+k_0} \; b_{k_0}(x) \right) \right) \right)$$

$$) ) \xi_0(-i_2) \Bigg) - \left( \sum_{i_1=1}^{m} \left( \sum_{i_0=0}^{m-i_1} C(n, i_0) \left( \frac{d^{n-i_0}}{dx^{n-i_0}} b_{i_0+i_1}(x) \right) \right) \xi_1(-i_1) \right)$$

In a general case we can write down for any : (11.6)

$$\xi_{m,s}(n) = \left( \sum_{i_{[s,0]}=0}^{[m-1,\, m-1-(s,s),\, m-1-(1,s)]} (-1)^{(0,s-1)} C(n, i_s) \left( \sum_{k_{[s-1,0]}=1}^{m} \right.\right.$$

$$\left( \prod_{l=0}^{s-1} C(k_l - 1 + i_l, k_l - 1) \right) \left( \frac{\partial^{n-i_s}}{\partial x^{n-i_s}} \left( [b_{1+(0,s)}(x)]_{i_{[0,s-1]}+k_{[0,s-1]}} \; b_{k_{[0,s-1]}}(x) \right) \right) ) ) ) )$$

$$- \left( \sum_{i_1=1}^{m} \left( \sum_{i_0=0}^{m-i_1} C(n, i_0) \left( \frac{d^{n-i_0}}{dx^{n-i_0}} b_{i_0+i_1}(x) \right) \right) \xi_{m,s-1}(-i_1) \right) - \left( \sum_{z=2}^{s} \right($$

$$\sum_{i_{[z,0]}=[1,0]}^{[m,\, m-(z,z),\, m-(z-1,z),\, m-(1,z)]} \left( \sum_{k_{[0,z-2]}=1}^{m} (-1)^{(0,z-2)} C(n, i_{z-1}) \right.$$

$$\left( \prod_{l=0}^{z-2} C(k_l - 1 + i_l, k_l - 1) \right) \left( \frac{\partial^{n-i_{z-2+1}}}{\partial x^{n-i_{z-2+1}}} \left( [b_{0,z}(x)]_{i_{[0,z-2]}+k_{[0,z-2]}} \; b_{k_{[0,z-2]}}(x) \right) \right) )$$

**Partial solutions of**

**are defined as follows**:

**<u>Step 1</u>: Calculation of** $b_k(x)$ , $k = 1, 2 .. m$ **functions**

$$b_k(x) = \sum_{i=1}^{k} (-1)^i \, C(m-i, m-k) \left( \frac{d^{k-i}}{dx^{k-i}} a_i(x) \right) \qquad k = 1, 2 \,..\, m \qquad (11.7)$$

**<u>Step 2:</u> Calculation of functions** in the following matrix:

$$\begin{bmatrix} \xi_{m,0}(-1) & \xi_{m,1}(-1) & \xi_{m,2}(-1) & [\,] & \xi_{m,N}(-1) \\ \xi_{m,0}(-2) & \xi_{m,1}(-2) & \xi_{m,2}(-2) & [\,] & \xi_{m,N}(-2) \\ [\,] & [\,] & [\,] & [\,] & [\,] \\ \xi_{m,0}(-m) & \xi_{m,1}(-m) & \xi_{m,2}(-m) & [\,] & \xi_{m,N}(-m) \end{bmatrix}$$

Here $N$ is a natural number that defines the number of serial members needed to achieve the required accuracy. In its limit (for the calculation of analytical solution) $N = \infty$ .
As soon as the functions are defined by (11.5), let's transform them in such form as to define the values of . To do this, let's use the basic definition of the operation of integration as an operation which is inverse to the operation of differentiation, and vice versa. According to this definition, the following equality holds true:

(11.8)

Thus, the formula for $\xi_{m,0}(-n)$ shall be as follows:

**Definition of the $\xi_{m,0}(-n)$ coefficient**

As soon as

$$\xi_{m,0}(n) = \frac{\partial^n}{\partial x^n} \left( \sum_{i=0}^{m-1} \begin{bmatrix} i \\ n \end{bmatrix} [b_{i+1}(x)]_i \right)$$

assuming that $n = -n$ , we get:

$$\xi_{m,0}(-n) = \frac{\partial^{-n}}{\partial x^{-n}} \left( \sum_{i=0}^{m-1} \begin{bmatrix} i \\ n \end{bmatrix} [b_{i+1}(x)]_i \right)$$

With consideration of (11.8), we finally obtain:

$$\xi_{m,0}(-n) = \sum_{i=0}^{m-1} \begin{bmatrix} i \\ -n \end{bmatrix} [b_{i+1}(x)]_{i+n} \qquad (11.9)$$

**Definition of the $\xi_1(-n)$ coefficient**

As soon as

$$(-1)^{i_0} C(n, i_1)\left(\sum_{k_0=1}^{m} C(k_0-1+i_0, k_0-1)\left(\frac{\partial^{n-i_1}}{\partial x^{n-i_1}}\left([b_{i_1+i_0+1}(x)]_{i_0+k_0}\, b_{k_0}(x)\right)\right)\right)\Bigg)\Bigg)\Bigg)$$

$$-\left(\sum_{i_1=1}^{m}\left(\sum_{i_0=0}^{m-i_1} C(n, i_0)\left(\frac{d^{n-i_0}}{dx^{n-i_0}} b_{i_0+i_1}(x)\right)\right)\xi_{m,0}(-i_1)\right)$$

assuming that $n = -n$, we get:

$$\xi_{m,1}(-n) = \left(\sum_{i_1=0}^{m-1}\left(\sum_{i_0=0}^{m-1-i_1}\right.\right.$$

$$(-1)^{i_0} C(-n, i_1)\left(\sum_{k_0=1}^{m} C(k_0-1+i_0, k_0-1)\left(\frac{\partial^{-n-i_1}}{\partial x^{-n-i_1}}\left([b_{i_1+i_0+1}(x)]_{i_0+k_0}\, b_{k_0}(x)\right)\right)\right)\Bigg)\Bigg)$$

$$-\left(\sum_{i_1=1}^{m}\left(\sum_{i_0=0}^{m-i_1} C(-n, i_0)\left(\frac{d^{-n-i_0}}{dx^{-n-i_0}} b_{i_0+i_1}(x)\right)\right)\xi_{m,0}(-i_1)\right)$$

With consideration of (11.8), the desired formula takes the following form: (11.10)

$$\sum_{i_0=0}^{m-1-i_1} (-1)^{i_0} C(-n, i_1)\left(\sum_{k_0=1}^{m} C(k_0-1+i_0, k_0-1)\left[[b_{i_1+i_0+1}(x)]_{i_0+k_0}\, b_{k_0}(x)\right]_{n+i_1}\right)\Bigg)$$

$$\Bigg)-\left(\sum_{i_1=1}^{m}\left(\sum_{i_0=0}^{m-i_1} C(-n, i_0)\,[b_{i_0+i_1}(x)]_{n+i_0}\right)\xi_{m,0}(-i_1)\right)$$

**Definition of the $\xi_{m,2}(-n)$ coefficient**

As soon as

$$\xi_{m,2}(n) = \left(\sum_{i_2=0}^{m-1}\left(\sum_{i_1=0}^{m-1-i_2}\left(\sum_{i_0=0}^{m-1-i_1-i_2} (-1)^{(i_0+i_1)} C(n, i_2)\left(\sum_{k_1=1}^{m}\left(\sum_{k_0=1}^{m} C(k_0-1+i_0, k_0-1)\right.\right.\right.\right.\right.$$

$$C(k_1-1+i_1, k_1-1)\left(\frac{\partial^{n-i_2}}{\partial x^{n-i_2}}\left(\left[[b_{i_2+i_1+i_0+1}(x)]_{i_0+k_0}\, b_{k_0}(x)\right]_{i_1+k_1} b_{k_1}(x)\right)\right)\Bigg)\Bigg)\Bigg)\Bigg)\Bigg)\Bigg) - \Bigg($$

$$\sum_{i_0=0}^{m-i_1-i_2} (-1)^{i_0} C(n, i_1)\left(\sum_{k_0=1}^{m} C(k_0-1+i_0, k_0-1)\left(\frac{\partial^{n-i_1}}{\partial x^{n-i_1}}\left([b_{i_1+i_0+i_2}]_{i_0+k_0}\, b_{k_0}(x)\right)\right)\right)\Bigg)$$

$$\Bigg)\Bigg)\xi_{m,0}(-i_2)\Bigg) - \left(\sum_{i_1=1}^{m}\left(\sum_{i_0=0}^{m-i_1} C(n, i_0)\left(\frac{d^{n-i_0}}{dx^{n-i_0}} b_{i_0+i_1}(x)\right)\right)\xi_{m,1}(-i_1)\right)$$

assuming that $n = -n$, we obtain:

$$\xi_{m,2}(-n) = \left(\sum_{i_2=0}^{m-1}\left(\sum_{i_1=0}^{m-1-i_2}\left(\sum_{i_0=0}^{m-1-i_1-i_2}(-1)^{(i_0+i_1)}\,\mathrm{C}(-n,i_2)\left(\sum_{k_1=1}^{m}\left(\sum_{k_0=1}^{m}\mathrm{C}(k_0-1+i_0,k_0-1)\,\mathrm{C}(k_1-1+i_1,k_1-1)\left(\frac{\partial^{-n-i_2}}{\partial x^{-n-i_2}}\left(\left[\left[[b_{i_2+i_1+i_0+1}(x)]\quad b_{k_0}(x)\right]_{i_0+k_0}\quad b_{k_1}(x)\right]_{i_1+k_1}\right)\right)\right)\right)\right)\right)\right) - \left(\sum_{i_2=1}^{m}\left(\sum_{i_1=0}^{m-i_2}\left(\sum_{i_0=0}^{m-i_1-i_2}(-1)^{i_0}\,\mathrm{C}(-n,i_1)\left(\sum_{k_0=1}^{m}\mathrm{C}(k_0-1+i_0,k_0-1)\left(\frac{\partial^{-n-i_1}}{\partial x^{-n-i_1}}\left(\left[b_{i_1+i_0+i_2}\right]\quad b_{k_0}(x)\right]_{i_0+k_0}\right)\right)\right)\right)\xi_{m,0}(-i_2)\right) - \left(\sum_{i_1=1}^{m}\left(\sum_{i_0=0}^{m-i_1}\mathrm{C}(-n,i_0)\left(\frac{d^{-n-i_0}}{dx^{-n-i_0}}b_{i_0+i_1}(x)\right)\right)\xi_{m,1}(-i_1)\right)$$

With consideration of (11.8), the desired formula takes the following form: (11.11)

$$\xi_{m,2}(-n) = \left(\sum_{i_2=0}^{m-1}\left(\sum_{i_1=0}^{m-1-i_2}\left(\sum_{i_0=0}^{m-1-i_1-i_2}(-1)^{(i_0+i_1)}\,\mathrm{C}(-n,i_2)\left(\sum_{k_1=1}^{m}\left(\sum_{k_0=1}^{m}\mathrm{C}(k_0-1+i_0,k_0-1)\,\mathrm{C}(k_1-1+i_1,k_1-1)\left[\left[\left[[b_{i_2+i_1+i_0+1}(x)]\quad b_{k_0}(x)\right]_{i_0+k_0}\quad b_{k_1}(x)\right]_{i_1+k_1}\right]_{n+i_2}\right)\right)\right)\right)\right) - \left(\sum_{i_2=1}^{m}\left(\sum_{i_1=0}^{m-i_2}\left(\sum_{i_0=0}^{m-i_1-i_2}(-1)^{i_0}\,\mathrm{C}(-n,i_1)\left(\sum_{k_0=1}^{m}\mathrm{C}(k_0-1+i_0,k_0-1)\left[\left[b_{i_1+i_0+i_2}\right]\quad b_{k_0}(x)\right]_{i_0+k_0}\right]_{n+i_1}\right)\right)\right)\xi_{m,0}(-i_2)\right) - \left(\sum_{i_1=1}^{m}\left(\sum_{i_0=0}^{m-i_1}\mathrm{C}(-n,i_0)\,[b_{i_0+i_1}(x)]_{n+i_0}\right)\xi_{m,1}(-i_1)\right)$$

Quite similarly, we obtain:

$$\xi_{m,3}(-n) = ( \sum_{i_3=0}^{m-1} ( \sum_{i_2=0}^{m-1-i_3} ( \sum_{i_1=0}^{m-1-i_2-i_3} ( \sum_{i_0=0}^{m-1-i_1-i_2-i_3} (-1)^{(i_0+i_1+i_2)} \mathrm{C}(-n, i_3) ( \sum_{k_2=1}^{m} ( \sum_{k_1=1}^{m} ($$

$$\sum_{k_0=1}^{m} ( \prod_{l=0}^{2} \mathrm{C}(k_l - 1 + i_l, k_l - 1) )$$

$$[[[[b_{1+i_0+i_1+i_2+i_3}(x)]_{i_0+k_0} b_{k_0}(x)]_{i_1+k_1} b_{k_1}(x)]_{i_2+k_2} b_{k_2}(x)]_{n+i_3} )))))))$$

$$- ( \sum_{i_1=1}^{m} ( \sum_{i_0=0}^{m-i_1} \mathrm{C}(-n, i_0) [b_{i_0+i_1}(x)]_{n+i_0} ) \xi_{m,2}(-i_1) ) - ( \sum_{i_2=1}^{m} ( \sum_{i_1=0}^{m-i_2} ( \sum_{i_0=0}^{m-i_1-i_2}$$

$$(-1)^{i_0} \mathrm{C}(-n, i_1) ( \sum_{k_0=1}^{m} ( \prod_{l=0}^{0} \mathrm{C}(k_l - 1 + i_l, k_l - 1) ) [[b_{i_0+i_1+i_2}(x)]_{i_0+k_0} b_{k_0}(x)]_{n+i_1} )))$$

$$\xi_{m,1}(-i_2) ) - ( \sum_{i_3=1}^{m} ( \sum_{i_2=0}^{m-i_3} ( \sum_{i_1=0}^{m-i_2-i_3} ( \sum_{i_0=0}^{m-i_1-i_2-i_3} (-1)^{(i_0+i_1)} \mathrm{C}(-n, i_2) ( \sum_{k_1=1}^{m}$$

$$( \sum_{k_0=1}^{m} ( \prod_{l=0}^{1} \mathrm{C}(k_l - 1 + i_l, k_l - 1) ) [[[b_{i_0+i_1+i_2+i_3}(x)]_{i_0+k_0} b_{k_0}(x)]_{i_1+k_1} b_{k_1}(x)]_{n+i_2} ))))$$

$$\xi_{m,4}(-n) = ( \sum_{i_4=0}^{m-1} ( \sum_{i_3=0}^{m-1-i_4} ( \sum_{i_2=0}^{m-1-i_3-i_4} ( \sum_{i_1=0}^{m-1-i_2-i_3-i_4} ( \sum_{i_0=0}^{m-1-i_1-i_2-i_3-i_4} (-1)^{(i_0+i_1+i_2+i_3)}$$

$$\mathrm{C}(-n, i_4) ( \sum_{k_3=1}^{m} ( \sum_{k_2=1}^{m} ( \sum_{k_1=1}^{m} ( \sum_{k_0=1}^{m} ( \prod_{l=0}^{3} \mathrm{C}(k_l - 1 + i_l, k_l - 1) )$$

$$\Bigg[\Bigg[\Bigg[\Bigg[\Big[b_{1+i_0+i_1+i_2+i_3+i_4}(x)\Big]_{i_0+k_0} b_{k_0}(x)\Bigg]_{i_1+k_1} b_{k_1}(x)\Bigg]_{i_2+k_2} b_{k_2}(x)\Bigg]_{i_3+k_3} b_{k_3}(x)\Bigg]_{n+i_4}\Bigg)\Bigg)$$

$$\Bigg)\Bigg)\Bigg)\Bigg)\Bigg)\Bigg)\Bigg)\Bigg) - \left(\sum_{i_1=1}^{m}\left(\sum_{i_0=0}^{m-i_1} \mathrm{C}(-n, i_0)\,[b_{i_0+i_1}(x)]_{n+i_0}\right)\xi_{m,3}(-i_1)\right) - \Bigg(\sum_{i_2=1}^{m}\Bigg(\sum_{i_1=0}^{m-i_2}\Bigg(\sum_{i_0=0}^{m-i_1-i_2}$$

$$(-1)^{i_0}\,\mathrm{C}(-n, i_1)\left(\sum_{k_0=1}^{m}\left(\prod_{l=0}^{0}\mathrm{C}(k_l-1+i_l, k_l-1)\right)\Big[[b_{i_0+i_1+i_2}(x)]_{i_0+k_0} b_{k_0}(x)\Big]_{n+i_1}\right)\Bigg)\Bigg)$$

$$\xi_{m,2}(-i_2)\Bigg) - \Bigg(\sum_{i_3=1}^{m}\Bigg(\sum_{i_2=0}^{m-i_3}\Bigg(\sum_{i_1=0}^{m-i_2-i_3}\Bigg(\sum_{i_0=0}^{m-i_1-i_2-i_3}(-1)^{(i_0+i_1)}\,\mathrm{C}(-n, i_2)\Bigg(\sum_{k_1=1}^{m}$$

$$\Bigg(\sum_{k_0=1}^{m}\left(\prod_{l=0}^{1}\mathrm{C}(k_l-1+i_l, k_l-1)\right)\Bigg[\Big[[b_{i_0+i_1+i_2+i_3}(x)]_{i_0+k_0} b_{k_0}(x)\Big]_{i_1+k_1} b_{k_1}(x)\Bigg]_{n+i_2}\Bigg)\Bigg)\Bigg)$$

$$\Bigg)\Bigg)\xi_{m,1}(-i_3)\Bigg) - \Bigg(\sum_{i_4=1}^{m}\Bigg(\sum_{i_3=0}^{m-i_4}\Bigg(\sum_{i_2=0}^{m-i_3-i_4}\Bigg(\sum_{i_1=0}^{m-i_2-i_3-i_4}\Bigg(\sum_{i_0=0}^{m-i_1-i_2-i_3-i_4}(-1)^{(i_0+i_1+i_2)}$$

$$\mathrm{C}(-n, i_3)\Bigg(\sum_{k_2=1}^{m}\Bigg(\sum_{k_1=1}^{m}\Bigg(\sum_{k_0=1}^{m}\left(\prod_{l=0}^{2}\mathrm{C}(k_l-1+i_l, k_l-1)\right)$$

$$\Bigg[\Bigg[\Bigg[[b_{i_0+i_1+i_2+i_3+i_4}(x)]_{i_0+k_0} b_{k_0}(x)\Bigg]_{i_1+k_1} b_{k_1}(x)\Bigg]_{i_2+k_2} b_{k_2}(x)\Bigg]_{n+i_3}\Bigg)\Bigg)\Bigg)\Bigg)\Bigg)\Bigg)\Bigg)\Bigg)\xi_{m,0}(-i_4)$$

and so on.

**Let's present the program for printing $\xi_{m,s}(-n)$ coefficients**

```
> restart;
```

```
> m:=m:s:=3: # Enter the order of the ODE - "m", and the number of the
coefficient - "s".
```

> $C0 := h \to \prod_{l=0}^{h-1} \mathrm{C}(k_l - 1 + i_l, k_l - 1)$

> $P0 := h \to \sum_{i_1=1}^{m} \left( \sum_{i_0=0}^{m-i_1} \mathrm{C}(-n, i_0) \left[ b_{i_0+i_1}(x) \right]_{n+i_0} \right) \xi_{m,h-1}(-i_1)$

```
> Ii:=h->add(i[j],j=0..h):
```

> $C1 := h \to (-1)^{\mathrm{Ii}(h-1)} \mathrm{C}(-n, i_h)$

```
> Iik:=h->i[h]+k[h]:
> B0(0):=b[1+Ii(s)](x):
```

> **for** $z$ **from** 0 **to** $s$ **do** $\mathrm{B0}(z+1) := [\mathrm{B0}(z)]_{\mathrm{Iik}(z)}\ b_{k_z}(x)$ **end do**; $z := 'z'$

```
> R0(0):=C0(s)*[B0(s)][n+i[s]]:
```

> **for** $z$ **from** 0 **to** $s$ **do** $\mathrm{R0}(z+1) := \sum_{k_z=1}^{m} \mathrm{R0}(z)$ **end do**; $z := 'z'$

```
> Ii1:=h->add(i[j],j=h..z):
> Ii2:=h->add(i[j],j=h..s):
> R1(0):=C1(s)*R0(s):
```

> **for** $z$ **from** 0 **to** $s$ **do** $\mathrm{R1}(z+1) := \sum_{i_z=0}^{m-1-\mathrm{Ii2}(z+1)} \mathrm{R1}(z)$ **end do**; $z := 'z'$

```
> R2(s):=R1(s+1)-P0(s):
> Ii2:=h->h:
> for z from 2 to s do:
> B1(0):=b[add(i[j],j=0..z)](x);
> for z0 from 0 to z-1 do B1(z0+1):=[B1(z0)][Iik(z0)]*b[k[z0]](x)
od:z0:='z0':
> R3(0):=Product(C(k[l]-1+i[l],k[l]-1),l = 0 .. z-2)*[B1(z-1)][n+i[z-1]];
> for z0 from 0 to z do R3(z0+1):=Sum(R3(z0),k[z0]=1..m) od;z0:='z0':
> R4(0):=(-1)^(add(i[j],j=0..z-2))*C(-n,i[z-1])*R3(z-1);
> for z0 from 0 to z do R4(z0+1):=Sum(R4(z0),i[z0]=0..m-Ii1(z0+1))
od;z0:='z0':
> R5[z]:=Sum(R4(z)*xi[m,s-z](-i[z]),i[z]=1..m);
```

```
> od:
> xi[m,s](-n)=R2(s)-add(R5[z1],z1=2..s);
> ###############################################################################################
```

**The definition of $\xi_{m,s}(-n)$ coefficient**

As soon as

$$\xi_{m,s}(n)=\left(\sum_{i_{[s,0]}=0}^{[m-1,m-1-(s,s),m-1-(1,s)]}(-1)^{(0,s-1)}C(n,i_s)\left(\sum_{k_{[s-1,0]}=1}^{[m]}\left(\prod_{l=0}^{s-1}C(k_l-1+i_l,k_l-1)\right)\left(\frac{\partial^{n-i_s}}{\partial x^{n-i_s}}\left([b_{1+(0,s)}(x)]_{i_{[0,s-1]}+k_{[0,s-1]}}\,b_{k_{[0,s-1]}}(x)\right)\right)\right)\right)$$
$$-\left(\sum_{i_1=1}^{m}\left(\sum_{i_0=0}^{m-i_1}C(n,i_0)\left(\frac{d^{n-i_0}}{dx^{n-i_0}}b_{i_0+i_1}(x)\right)\right)\xi_{m,s-1}(-i_1)\right)-\left(\sum_{z=2}^{s}\left(\sum_{i_{[z,0]}=[1,0]}^{[m,m-(z,z),m-(z-1,z),m-(1,z)]}\left(\sum_{k_{[0,z-2]}=1}^{[m]}(-1)^{(0,z-2)}C(n,i_{z-1})\right.\right.\right.$$
$$\left.\left.\left.\left(\prod_{l=0}^{z-2}C(k_l-1+i_l,k_l-1)\right)\left(\frac{\partial^{n-i_{z-2+1}}}{\partial x^{n-i_{z-2+1}}}\left([b_{0,z}(x)]_{i_{[0,z-2]}+k_{[0,z-2]}}\,b_{k_{[0,z-2]}}(x)\right)\right)\right)\right)\right)$$

assuming that $n=-n$ , we obtain:

$$\xi_{m,s}(-n)=\left(\sum_{i_{[s,0]}=0}^{[m-1,m-1+(-s,-s),m-1+(-1,-s)]}(-1)^{(0,s-1)}C(-n,i_s)\left(\sum_{k_{[s-1,0]}=1}^{[m]}\left(\prod_{l=0}^{s-1}C(k_l-1+i_l,k_l-1)\right)\left(\frac{\partial^{-n-i_s}}{\partial x^{-n-i_s}}\left([b_{1+(0,s)}(x)]_{i_{[0,s-1]}+k_{[0,s-1]}}\,b_{k_{[0,s-1]}}(x)\right)\right)\right)\right)$$
$$-\left(\sum_{i_1=1}^{m}\left(\sum_{i_0=0}^{m-i_1}C(-n,i_0)\left(\frac{d^{-n-i_0}}{dx^{-n-i_0}}b_{i_0+i_1}(x)\right)\right)\xi_{m,s-1}(-i_1)\right)-\left(\sum_{z=2}^{s}\left(\sum_{i_{[z,0]}=[1,0]}^{[m,m+(-z,-z),m+(-z+1,-z),m+(-1,-z)]}\left(\sum_{k_{[0,z-2]}=1}^{[m]}(-1)^{(0,z-2)}C(-n,i_{z-1})\right.\right.\right.$$
$$\left.\left.\left.\left(\prod_{l=0}^{z-2}C(k_l-1+i_l,k_l-1)\right)\left(\frac{\partial^{-n-i_{z-1}}}{\partial x^{-n-i_{z-1}}}\left([b_{0,z}(x)]_{i_{[0,z-2]}+k_{[0,z-2]}}\,b_{k_{[0,z-2]}}(x)\right)\right)\right)\right)\right)$$

With consideration of (14.8), the desired formula takes the following form: (11.12)

$$\xi_{m,s}(-n) = \left( \sum_{i_{[s,0]}=0}^{[m-1,\, m-1-(s,s),\, m-1-(1,s)]} (-1)^{(0,s-1)}\, C(-n, i_s) \left( \sum_{k_{[s-1,0]}=1}^{[m]} \left( \prod_{l=0}^{s-1} C(k_l - 1 + i_l, k_l - 1) \right) \left[ [b_{1+(0,s)}(x)]_{i_{[0,s-1]}+k_{[0,s-1]}} \; b_{k_{[0,s-1]}}(x) \right]_{n+i_s} \right) \right) - \left( \sum_{i_1=1}^{m} \left( \sum_{i_0=0}^{m-i_1} C(-n, i_0)\, [b_{i_0+i_1}(x)]_{n+i_0} \right) \xi_{m,s-1}(-i_1) \right) - \left( \sum_{z=2}^{s} \left( \sum_{i_{[z,0]}=[1,0]}^{[m,\, m-(z,z),\, m-(z+1,z),\, m-(1,z)]} \left( \sum_{k_{[0,z-2]}=1}^{[m]} (-1)^{(0,z-2)}\, C(-n, i_{z-1}) \left( \prod_{l=0}^{z-2} C(k_l - 1 + i_l, k_l - 1) \right) \left[ [b_{0,z}(x)]_{i_{[0,z-2]}+k_{[0,z-2]}} \; b_{k_{[0,z-2]}}(x) \right]_{n+i_{z-1}} \xi_{m,s-z}(-i_z) \right) \right) \right)$$

Thus, we have established the algorithm to define $\xi_{m,s}(-n)$ summands, which form, in a general case, an (m x N) matrix:

$$\begin{bmatrix} \xi_{m,0}(-1) & \xi_{m,1}(-1) & \xi_{m,2}(-1) & [\ ] & \xi_{m,N}(-1) \\ \xi_{m,0}(-2) & \xi_{m,1}(-2) & \xi_{m,2}(-2) & [\ ] & \xi_{m,N}(-2) \\ [\ ] & [\ ] & [\ ] & [\ ] & [\ ] \\ \xi_{m,0}(-m) & \xi_{m,1}(-m) & \xi_{m,2}(-m) & [\ ] & \xi_{m,N}(-m) \end{bmatrix}$$

The prinicipal task here is the finding of all values for the elements of this matrix.

**<u>Step 3:</u> Calculation of the function**

$$k = 1, 2 \,..\, m \tag{11.13}$$

In a general case, upon the definition of change pattern for the $\xi_{m,s}(-k)$ functions, we can calculate the following series:

$$k = 1, 2 \,..\, m \tag{11.14}$$

and in this case we obtain the exact solution of the original equation.

**<u>Step 4:</u> Calculation of the partial solutions:**

$$Y_{m,i}(x) = \frac{(-1)^i (1 - \alpha_m(-1))\, x^{(i-1)}}{(i-1)!} + \left( \sum_{k=2}^{i} \frac{(-1)^{(k+i)} \alpha_m(-k)\, x^{(i-k)}}{(i-k)!} \right) \qquad i = 1, 2 \,..\, m \tag{11.15}$$

**<u>Step 5:</u> Calculation of the relative absolute accuracy of the resulting partial solution, $Y_i(x)$, at a specified interval of change of the independent variable *x* for the ODE (11.1):**

$$\delta(Y_{m,i}(x)) = \left| \frac{\left(\frac{d^m}{dx^m} Y_{m,i}(x)\right) - \left(\sum_{p=1}^{m} a_p(x) \left(\frac{d^{m-p}}{dx^{m-p}} Y_{m,i}(x)\right)\right)}{a_m(x)\, Y_{m,i}(x)} \right| \qquad i = 1, 2 .. m \qquad (11.17)$$

In the case where it is less than the specified $\delta_z(Y_{m,i}(x))$, i.e. where the following equality holds true:

then the resulting partial solution $Y_{m,i}(x)$ is deemed established.

**<u>Step 6:</u> Formula for the general solution of the ODE (11.1):**

(11.18)

where $C_i$ are arbitrary constants.

##########################################################################################

## Below there is a Program (in MAPLE 10), which calculates, for the following type of m-order ODEs:

$$\frac{\partial^m}{\partial x^m} y = \sum_{i=1}^{m} a_i(x) \left(\frac{\partial^{m-i}}{\partial x^{m-i}} y\right)$$

**the values of all $\xi_{m,s}(-n)$, $s = 1, 2 .. N$ matrix coefficeices:**

$$\Xi_{m,n} = \begin{bmatrix} \xi_{m,0}(-1) & \xi_{m,1}(-1) & \xi_{m,2}(-1) & [\;] & \xi_{m,N}(-1) \\ \xi_{m,0}(-2) & \xi_{m,1}(-2) & \xi_{m,2}(-2) & [\;] & \xi_{m,N}(-2) \\ [\;] & [\;] & [\;] & [\;] & [\;] \\ \xi_{m,0}(-m) & \xi_{m,1}(-m) & \xi_{m,2}(-m) & [\;] & \xi_{m,N}(-m) \end{bmatrix}$$

**with the use of the following formulas:**

$$\xi_{m,0}(-n) = \sum_{i=0}^{m-1} \begin{bmatrix} i \\ -n \end{bmatrix} [b_{i+1}(x)]_{i+n}$$

$$\xi_{m,s}(-n) = \Bigg( \sum_{i_{[s,0]}=0}^{[m-1,\, m-1-(s,s),\, m-1-(1,s)]} (-1)^{(0,s-1)}\, \mathrm{C}(-n, i_s) \Bigg($$

$$\sum_{k_{[s-1,0]}=1}^{[m]} \left(\prod_{l=0}^{s-1} \mathrm{C}(k_l - 1 + i_l, k_l - 1)\right) \left[ [b_{1+(0,s)}(x)]_{i_{[0,s-1]}+k_{[0,s-1]}} \quad b_{k_{[0,s-1]}}(x) \right]_{n+i_s}$$

$$\Bigg)\Bigg) - \left(\sum_{i_1=1}^{m} \left(\sum_{i_0=0}^{m-i_1} \mathrm{C}(-n, i_0)\, [b_{i_0+i_1}(x)]_{n+i_0}\right) \xi_{m,s-1}(-i_1)\right) - \Bigg(\sum_{z=2}^{s} \Bigg($$

$$\sum_{i_{[z,0]}=[1,0]}^{[m,m-(z,z),m-(z+1,z),m-(1,z)]}\left(\sum_{k_{[0,z-2]}=1}^{[m]}(-1)^{(0,z-2)}\,C(-n,i_{z-1})\right.$$
$$\left.\left(\prod_{l=0}^{z-2}C(k_l-1+i_l,k_l-1)\right)\left[\,[b_{0,z}(x)]_{i_{[0,z-2]}+k_{[0,z-2]}}\;b_{k_{[0,z-2]}}(x)\right]_{n+i_{z-1}}\xi_{m,s-z}(-i_z)\right)\Bigg)\Bigg)$$

**- the values of**

$$\alpha_m(-n)=\sum_{s=0}^{N}\xi_{m,s}(-n),\ n=-1,-2\,..-m$$

**- partial solutions for the ODE:**

$$Y_{m,i}(x)=\frac{(-1)^i\,(1-\alpha_m(-1))\,x^{(i-1)}}{(i-1)!}+\left(\sum_{k=2}^{i}\frac{(-1)^{(k+i)}\,\alpha_m(-k)\,x^{(i-k)}}{(i-k)!}\right),\ i=1,2\,..\,m$$

**- absolute relative error of the obtained solution**

$$\Delta(Y_{m,i}(x)):=\left|\frac{\left(\frac{d^m}{dx^m}Y_{m,i}(x)\right)-\left(\sum_{k=1}^{m}a_k(x)\left(\frac{d^{m-k}}{dx^{m-k}}Y_{m,i}(x)\right)\right)}{a_m(x)\,Y_{m,i}(x)}\right|,\ i=1\,..\,m$$

**- Graphical change of this relative error** $\Delta(Y_i(x))$ **in the range of** $[X_{nach}=-10.01,X_{kon}=10.01]$

```
> restart;alias(y=y(x),C=binomial,Ex=expand,Q=Q(x)):

> X[nach]:=-10.01:X[kon]:=10.01:

> m:=2:# Enter the order of the ODE - "m"

> a[1](x):=0:a[2](x):=x*ln(x):a[3](x):=sqrt(x):a[4](x):=2*x:

> diff(y,x$m)=add(a[i](x)*diff(y,x$(m-i)),i=1..m-1)+a[m](x)*y;

> dsolve(%);

> N:=4:#Enter the number of the coefficients - from "s=1" to "s=N".

> dn := proc (n::integer, w) local k, Ds; global resd, x; option remember;
if n = 0 then resd := w elif 0 < n then resd := diff(w,`$`(x,n)) else
Ds := w; for k to -n do Ds := int(Ds,x) end do; resd := Ds end if;
RETURN(resd) end proc:

> for k0 to m do b[k0](x):=add((-1)^(i)*C(m-i,m-k0)*dn(k0-i,a[i](x)),i =
1 .. k0) od;

> diff(Q,x$m)=add(b[i](x)*diff(Q,x$(m-i)),i=1..m-1)+b[m](x)*Q;

> for n to m do

> for s to N do
```

```
> C0:=h->Product(C(k[l]-1+i[l],k[l]-1),l=0..h-1):

> P0 := proc (h) options operator, arrow; Sum(Sum(C(-n,i[0])*Dn(-(n+i[0]),b[i[0]+i[1]](x)),i[0] = 0 .. m-i[1])*xI[m,h-1](-i[1]),i[1] = 1 .. m) end proc:

> Ii:=h->add(i[j],j=0..h):

> Ii1:=h->add(i[j],j=h..z):

> Ii2:=h->add(i[j],j=h..s):

> C1 := proc (h) options operator, arrow; (-1)^(Ii(h-1))*C(-n,i[h]) end proc:

> Iik:=h->i[h]+k[h]:

> B0(0):=b[1+Ii(s)](x):

> for z from 0 to s do B0(z+1):=Dn(-Iik(z),B0(z))*b[k[z]](x) od:z:='z':

> R0(0):=C0(s)*Dn(-(n+i[s]),B0(s)):

> for z from 0 to s do R0(z+1):=Sum(R0(z),k[z]=1..m) od:z:='z':

> R1(0):=C1(s)*R0(s):

> for z from 0 to s do R1(z+1) := Sum(R1(z),i[z] = 0 .. m-1-Ii2(z+1)) end do; z := 'z':

> R2(s):=R1(s+1)-P0(s):

> Ii2:=h->h:

> for z from 2 to s do:

> B1(0):=b[add(i[j],j=0..z)](x);

> for z0 from 0 to z-1 do B1(z0+1):=Dn(-(Iik(z0)),B1(z0))*b[k[z0]](x) od:z0:='z0':

> R3(0):=Product(C(k[l]-1+i[l],k[l]-1),l = 0 .. z-2)*Dn(-(n+i[z-1]),B1(z-1));

> for z0 from 0 to z do R3(z0+1):=Sum(R3(z0),k[z0]=1..m) od;z0:='z0':

> R4(0):=(-1)^(add(i[j],j=0..z-2))*C(-n,i[z-1])*R3(z-1);

> for z0 from 0 to z do R4(z0+1) := Sum(R4(z0),i[z0] = 0 .. m-Ii1(z0+1)) end do; z0 := 'z0':

> R5[z]:=Sum(R4(z)*xI[m,s-z](-i[z]),i[z]=1..m);

> od:

> R(m,s,n):=(xI[m,s](-n)=R2(s)-add(R5[z1],z1=2..s));

> od:

> od:z:='z':k:='k':
```

```
> Xi[m,N]=matrix([seq([seq(xi[m,s0](-n0),s0=0..N)],n0=1..m)]);

> for s0 from 0 to N do

> if s0=0 then for n0 to m do xI[m,0](-n0):=value(add(C(-n0,i0)*dn(-(n0+i0),b[i0+1](x)),i0 = 0 .. m-1)) od else  for n0 to m do xI[m,s0](-n0):=expand(eval(subs(Sum=add,Product=product,Dn=dn,rhs(R(m,s0,n0))))) od fi;

> od:

> for s0 from 0 to N do

> if s0=0 then for n0 to m do print(xi[m,s0](-n0)=xI[m,0](-n0)) od else for n0 to m do print(xi[m,s0](-n0)=xI[m,s0](-n0)) od fi;

> od:

> for k to m do alpha[m](-k):=add(xI[m,s](-k),s=0..N) od;k:='k':

> for i to m do R[i] := Y(x) =collect(sort(combine(Ex(Ex(Ex((-1)^i*(1-alpha[m](-1))*x^(i-1)/(i-1)!+add((-1)^(k+i)*alpha[m](-k)*x^(i-k)/(i-k)!,k = 2 .. i))))),x),[seq(ln(x)^(N-i),i=0..N-1),x],factor) end do;

> for i to m do Delta(i) := expand(subs(R[i],diff(Y(x),`$`(x,m))-add(a[p](x)*dn(m-p,Y(x)),p = 1 .. m))/(a[m](x)*subs(R[i],Y(x)))); for s from X[nach] to X[kon] do he1(s) := abs(evalf(subs(x = s,Delta(i)))) end do; x := 'x'; with(plots); print(k); for s from X[nach] to X[kon] do if abs(he1(s)) < .1 then print(delta(s) = he1(s),x = s) else NULL end if end do; s := 's'; x := 'x'; print(pointplot({seq([n, he1(n)],n = X[nach] .. X[kon])},numpoints = 20,axes = BOXED,color = red,xtickmarks = 20,ytickmarks = 20,,scaling = UNCONSTRAINED,title = "ACCURACY OF CALCULATION ")) end do:
```

##########################################################################################

**Let's consider examples of linear ordinary differential equations of different orders, starting with the second order, to demonstrate its validity.**
**Example No. 1. Find a general solution of the following linear second-order ODE:**

(11.19)

**Solution**. **In accordance with (14.47) let's calculate the functions**:

$= -x$

$$b_2(x) = -\left(\frac{d}{dx} a_1(x)\right) + a_2(x) = -1 + \sqrt{x}$$

Let's calculate the summands: .

**The** $\xi_0(-n)$ **summand**. In accordance with (14.49), we have, for $m = 2$ :

$$\xi_0(-n) = \sum_{i=0}^{1} \left[\begin{matrix} i \\ -n \end{matrix}\right] [b_{i+1}(x)]_{i+n}$$

Hence we obtain the following formulas:

$$\xi_0(-1) := \int b_1(x)\,dx - \iint b_2(x)\,dx\,dx \qquad \xi_0(-2) := \iint b_1(x)\,dx\,dx - 2\iiint b_2(x)\,dx\,dx\,dx$$

Thus, the required values are:

$$\xi_0(-1) := \int -x\,dx - \iint -1 + \sqrt{x}\,dx\,dx = -\frac{4x^{\left(\frac{5}{2}\right)}}{15} \qquad \xi_0(-2) := \iint -x\,dx\,dx - 2\iiint -1 + \sqrt{x}\,dx\,dx\,dx =$$

**The $\xi_1(-n)$ summand**. In accordance with (14.49), we have, for $m = 2$ :

$$\left(\sum_{i_0=0}^{1-i_1} (-1)^{i_0}\, \mathrm{C}(-n, i_1) \left(\sum_{k_0=1}^{2} \mathrm{C}(k_0 - 1 + i_0, k_0 - 1) \left[ [b_{i_1+i_0+1}(x)]_{i_0+k_0}\, b_{k_0}(x) \right]_{n+i_1}\right)\right)$$
$$- \left(\sum_{i_1=1}^{2} \left(\sum_{i_0=0}^{2-i_1} \mathrm{C}(-n, i_0)\, [b_{i_0+i_1}(x)]_{n+i_0}\right) \xi_0(-i_1)\right)$$

Hence we obtain the following formulas:

$$\begin{aligned}
\xi_1(-1) = {} & -\frac{1\left(\int b_1(x)\,dx\right)^2}{2} - \iiint b_2(x)\,dx\,dx\,b_1(x)\,dx - \iiint b_2(x)\,dx\,b_1(x)\,dx\,dx \\
& + \iiint b_1(x)\,dx\,dx\,b_2(x)\,dx - 2\iiiint b_2(x)\,dx\,dx\,dx\,b_2(x)\,dx \\
& - \iiiint b_2(x)\,dx\,dx\,b_2(x)\,dx\,dx + 2\int b_1(x)\,dx \iint b_2(x)\,dx\,dx - \left(\iint b_2(x)\,dx\,dx\right)^2 \\
& - \iint b_1(x)\,dx\,dx \int b_2(x)\,dx + 2\iiint b_2(x)\,dx\,dx\,dx \int b_2(x)\,dx
\end{aligned}$$

$$\begin{aligned}
\xi_1(-2) := {} & \frac{1\int\left(\int b_1(x)\,dx\right)^2 dx}{2} - \iiiint b_2(x)\,dx\,dx\,b_1(x)\,dx\,dx \\
& - 2\iiiint b_2(x)\,dx\,b_1(x)\,dx\,dx\,dx + \iiiint b_1(x)\,dx\,dx\,b_2(x)\,dx\,dx \\
& - 2\iiiiint b_2(x)\,dx\,dx\,dx\,b_2(x)\,dx\,dx - 2\iiiiint b_2(x)\,dx\,dx\,b_2(x)\,dx\,dx\,dx \\
& - \iint b_1(x)\,dx\,dx \int b_1(x)\,dx + 2\iiint b_2(x)\,dx\,dx\,dx \int b_1(x)\,dx
\end{aligned}$$

Thus, the required values are:

$$\xi_1(-1) =$$

$$\xi_1(-2) := \frac{x^5}{40} - \frac{202\, x^{\left(\frac{11}{2}\right)}}{10395} - \frac{x^6}{105}$$

Quite similarly, we obtain:

$$\xi_2(-1) :=$$

$$-\frac{16\, x^{\left(\frac{13}{2}\right)}}{3003} - \frac{103\, x^7}{39690} - \frac{4\, x^{\left(\frac{15}{2}\right)}}{14625}, \xi_2(-2) = \frac{x^7}{336} - \frac{9049\, x^{\left(\frac{15}{2}\right)}}{4054050} - \frac{1109\, x^8}{623700} - \frac{124\, x^{\left(\frac{17}{2}\right)}}{580125}$$

$$\xi_3(-1) := -\frac{32\, x^{\left(\frac{17}{2}\right)}}{58905} - \frac{2\, x^{10}}{658125} - \frac{23974\, x^{\left(\frac{19}{2}\right)}}{416645775} - \frac{19049\, x^9}{58378320},$$

$$\xi_3(-2) = -\frac{2366443\, x^{\left(\frac{21}{2}\right)}}{53469541125} - \frac{86\, x^{11}}{34459425} - \frac{746737\, x^{10}}{3405402000} - \frac{2011\, x^{\left(\frac{19}{2}\right)}}{9189180} + \frac{1\, x^9}{3456}$$

$$\xi_4(-1) =$$

$$-\frac{64\, x^{\left(\frac{21}{2}\right)}}{1382535} - \frac{8\, x^{\left(\frac{25}{2}\right)}}{378421875} - \frac{383737\, x^{11}}{12129717600} - \frac{91589\, x^{12}}{137493105750} - \frac{30321563\, x^{\left(\frac{23}{2}\right)}}{4194263373300},$$

$$\xi_4(-2) := -\frac{26632\, x^{\left(\frac{27}{2}\right)}}{1486062703125} - \frac{54932745817\, x^{\left(\frac{25}{2}\right)}}{9961375511587500} + \frac{1\, x^{11}}{42240} - \frac{79843\, x^{12}}{3789283680}$$

$$-\frac{193447\, x^{13}}{356463607500} - \frac{12395081\, x^{\left(\frac{23}{2}\right)}}{674632838880}$$

$$\xi_5(-1) := -\frac{27446369\, x^{14}}{328139428617000} - \frac{1841832823\, x^{\left(\frac{27}{2}\right)}}{2708035264935000} - \frac{4\, x^{15}}{39734296875} - \frac{128\, x^{\left(\frac{25}{2}\right)}}{37855125}$$

$$-\frac{57826141\, x^{13}}{22878852796800} - \frac{29422126\, x^{\left(\frac{29}{2}\right)}}{6190283353629375}, \xi_5(-2) = -\frac{134675041993\, x^{\left(\frac{31}{2}\right)}}{3358228719343935937 5}$$

$$-\frac{3137018107\, x^{14}}{1870082229375360} - \frac{1218965329\, x^{\left(\frac{27}{2}\right)}}{910754332488000} - \frac{1689320593\, x^{15}}{24903438778968750}$$

$$-\frac{389\, x^{16}}{4458188109375} + \frac{1\, x^{13}}{599040} - \frac{17731412543069\, x^{\left(\frac{29}{2}\right)}}{343189309125212 55000}$$

$$\xi_6(-1) := -\frac{12765560791\, x^{\left(\frac{33}{2}\right)}}{21400023011543426250} - \frac{4641299551073\, x^{\left(\frac{31}{2}\right)}}{89085534027349860000} - \frac{16\, x^{\left(\frac{35}{2}\right)}}{45893112890625}$$
$$-\frac{7334349511\, x^{15}}{42501868849440000} - \frac{2190648114461\, x^{16}}{284018738586382800000} - \frac{484620533\, x^{17}}{21046963402339875000}$$
$$-\frac{256\, x^{\left(\frac{29}{2}\right)}}{1185622515}, \xi_6(-2) = -\frac{72126276704340979\, x^{\left(\frac{33}{2}\right)}}{1825629848822480680980000}$$
$$-\frac{31215984803\, x^{\left(\frac{31}{2}\right)}}{363896953291872000} - \frac{136275177383006\, x^{\left(\frac{35}{2}\right)}}{271512791958957220546875} + \frac{1\, x^{15}}{9676800}$$
$$-\frac{5009151926287\, x^{16}}{43759924167383424000} - \frac{44973727170989\, x^{17}}{7206975491629463550000} - \frac{2956\, x^{\left(\frac{37}{2}\right)}}{9622256002734375}$$
$$-\frac{8909310511\, x^{18}}{4477638292458581250000}$$

and so on.
From that we can establish the functional pattern of change of the summands $\xi_{2k+1}(-1)$ and $\xi_{2k}(-1)$.

$$\xi_{2s}(-1) = \left(\sum_{p=s}^{2s} b_{s,p}\, x^{\left(\frac{5}{2}+p+3s\right)}\right) + \left(\sum_{p=s}^{2s-1} c_{s,p}\, x^{(6s-p+2)}\right)$$
$$\xi_{2s+1}(-1) = \left(\sum_{p=s}^{2s} g_{s,p}\, x^{\left(\frac{5}{2}+p+3s+2\right)}\right) + \left(\sum_{p=s}^{2s} r_{s,p}\, x^{(6s+5-p)}\right) \qquad (11.20)$$

Here are the desired numerical coefficients.
Proceeding the same way, we can establish the functional pattern of change of the summands: $\xi_{2k+1}(-2)$ and $\xi_{2k}(-2)$.

$$\xi_{2s}(-2) = \left(\sum_{p=s}^{2s} \rho_{s,p}\, x^{\left(\frac{5}{2}+p+3s+1\right)}\right) + \left(\sum_{p=s}^{2s} \tau_{s,p}\, x^{(6s+3-p)}\right)$$
$$\xi_{2s+1}(-2) = \left(\sum_{p=s}^{2s} \chi_{s,p}\, x^{\left(\frac{5}{2}+\frac{2p}{2}+3s+3\right)}\right) + \left(\sum_{p=s}^{2s+1} \omega_{s,p}\, x^{(6s+6-p)}\right) \qquad (11.21)$$

Here are the desired numerical coefficients.
Thus, the required values of are defined by the following formulas:

$$\alpha(-k) = \sum_{s=0}^{\infty} \left(\xi_{2s}(-k) + \xi_{2s+1}(-1)\right) \qquad k = 1, 2$$

or, in expanded view :

$$\alpha(-1) = \sum_{s=0}^{\infty} \left(\left(\sum_{p=s}^{2s} b_{s,p}\, x^{\left(\frac{5}{2}+p+3s\right)}\right) + \left(\sum_{p=s}^{2s-1} c_{s,p}\, x^{(6s-p+2)}\right) + \left(\sum_{p=s}^{2s} g_{s,p}\, x^{\left(\frac{5}{2}+p+3s+2\right)}\right) + \left(\sum_{p=s}^{2s} r_{s,p}\, x^{(6s+5-p)}\right)\right) \qquad (11.22)$$

$$\alpha(-2)=\sum_{s=0}^{\infty}\left(\left(\sum_{p=s}^{2s}\rho_{s,p}\,x^{\left(\frac{5}{2}+p+3s+1\right)}\right)+\left(\sum_{p=s}^{2s}\tau_{s,p}\,x^{(6s+3-p)}\right)+\left(\sum_{p=s}^{2s}\chi_{s,p}\,x^{\left(\frac{5}{2}+p+3s+3\right)}\right)+\left(\sum_{p=s}^{2s+1}\omega_{s,p}\,x^{(6s+6-p)}\right)\right) \quad (11.23)$$

**The desired solutions are written down with the use of (11.3):**

$$Y_1(x)=-1+\left(\sum_{s=0}^{\infty}\left(\left(\sum_{p=s}^{2s}b_{s,p}\,x^{\left(\frac{5}{2}+p+3s\right)}\right)+\left(\sum_{p=s}^{2s-1}c_{s,p}\,x^{(6s-p+2)}\right)+\left(\sum_{p=s}^{2s}g_{s,p}\,x^{\left(\frac{5}{2}+p+3s+2\right)}\right)+\left(\sum_{p=s}^{2s}r_{s,p}\,x^{(6s+5-p)}\right)\right)\right) \quad (11.24)$$

$$Y_2(x)=\left(1-\left(\sum_{s=0}^{\infty}\left(\left(\sum_{p=s}^{2s}b_{s,p}\,x^{\left(\frac{5}{2}+p+3s\right)}\right)+\left(\sum_{p=s}^{2s-1}c_{s,p}\,x^{(6s-p+2)}\right)+\left(\sum_{p=s}^{2s}g_{s,p}\,x^{\left(\frac{5}{2}+p+3s+2\right)}\right)+\left(\sum_{p=s}^{2s}r_{s,p}\,x^{(6s+5-p)}\right)\right)\right)\right)x+\left(\sum_{s=0}^{\infty}\left(\left(\sum_{p=s}^{2s}\rho_{s,p}\,x^{\left(\frac{5}{2}+p+3s+1\right)}\right)+\left(\sum_{p=s}^{2s}\tau_{s,p}\,x^{(6s+3-p)}\right)+\left(\sum_{p=s}^{2s}\chi_{s,p}\,x^{\left(\frac{5}{2}+p+3s+3\right)}\right)+\left(\sum_{p=s}^{2s+1}\omega_{s,p}\,x^{(6s+6-p)}\right)\right)\right) \quad (11.25)$$

By way of successive substitution of the obtained solutions - - to the initial ODE (11.19), we establish (in a manner similar to that previously used for solving second-order ODE) recurrence relations for the calculation of numerical coefficients.

(11.26)

$$(p-1-6s)(p-2-6s)\,c_{s,p}+(p-1-6s)(p-2-6s)\,r_{s,p+3}+(p-6s)\,c_{s,p+2}+(-6s-5)\,r_{s,p+5}+p\,r_{s,p}=0$$

$$(2p+6s+5)(2p+6s+3)\,g_{s,p-2}+(2p+6s+5)(2p+6s+3)\,b_{s,p}-2\,(g_{s,p-4}+b_{s,p-2})(1+2p+6s)=0$$

$$(12+24s-4p)\,c_{s,p}+(-12+4p-24s)\,\tau_{s,p}+(12+24s-4p)\,r_{s,p+3}-4\,(6s+5-p)(6s+4-p)\,r_{s,p+1}+(-36p+80-48ps+144s^2+96s+4p^2)\,\tau_{s,p-2}-4\,c_{s,p-2}(6s+3-p)(6s+4-p)+4\,\omega_{s,p+1}(6s+5-p)(6s+4-p)+(-12+4p-24s)\,\omega_{s,p+3}+120\,s\,\tau_{s,p+2}=0$$

$$(22+12s+4p)g_{s,p}+(-12s-22-4p)\chi_{s,p}-(2p+15+6s)(13+2p+6s)b_{s,p+4}$$
$$+(2p+15+6s)(13+2p+6s)\rho_{s,p+4}+(22+12s+4p)b_{s,p+2}$$
$$+(2p+15+6s)(13+2p+6s)\chi_{s,p+2}-(2p+15+6s)(13+2p+6s)g_{s,p+2}$$
$$-2\rho_{s,p+2}(2p+6s+11)=0$$

$$\tau_{s,p-2}-r_{s,p+1}+\omega_{s,p+1}-c_{s,p-2}=0$$

$$\rho_{s,p}+\chi_{s,p-2}-g_{s,p-2}-b_{s,p}=0$$

Thus, the required general solution of the ODE (11.19) is defined by the following formula:

(11.27)

$C_i$ are arbitrary constants, and are defined by (11.24), (11.25). In this case, the numerical coefficients in these equations are found from the recurrence relations (11.26). **The task is solved**.

*Note:* ***Calculation of the relative error in accordance with (14.55) for each of the solutions demonstrates that in the [0,2] interval this relative error does not exceed 0.02, and the first seven terms of the series are enough for this. The desired partial solutions together with the relative error are as follows****:*

$$Y_1(x)=-1-\frac{4x^{\left(\frac{5}{2}\right)}}{15}-\frac{32x^{\left(\frac{17}{2}\right)}}{58905}-\frac{274740728791\,x^{\left(\frac{33}{2}\right)}}{21400023011543426250}-\frac{4641299551073\,x^{\left(\frac{31}{2}\right)}}{890855340273498600000}$$
$$-\frac{1x^5}{75}-\frac{4x^{20}}{4359845724609375}-\frac{64787428890700319\,x^{\left(\frac{35}{2}\right)}}{190851388019076534750000000}-\frac{64x^{\left(\frac{21}{2}\right)}}{1382535}$$
$$-\frac{27446369\,x^{14}}{328139428617000}-\frac{8x^{\left(\frac{9}{2}\right)}}{189}-\frac{16x^{\left(\frac{13}{2}\right)}}{3003}-\frac{103x^7}{39690}-\frac{1841832823\,x^{\left(\frac{27}{2}\right)}}{2708035264935000}$$
$$-\frac{917328514859\,x^{15}}{5312733606180000000}-\frac{4x^{\left(\frac{15}{2}\right)}}{14625}-\frac{2190648114461\,x^{16}}{284018738586382800000}$$
$$-\frac{779466376779397\,x^{17}}{754323168339861120000000}-\frac{1883704\,x^{\left(\frac{25}{2}\right)}}{553631203125}-\frac{383737\,x^{11}}{12129717600}$$
$$-\frac{91589\,x^{12}}{137493105750}-\frac{2x^{10}}{658125}-\frac{57826141\,x^{13}}{22878852796800}-\frac{23974\,x^{\left(\frac{19}{2}\right)}}{416645775}$$
$$-\frac{1366030126\,x^{\left(\frac{29}{2}\right)}}{6190283353629375}-\frac{47043047375059\,x^{\left(\frac{37}{2}\right)}}{86753832417230828850000}-\frac{30321563\,x^{\left(\frac{23}{2}\right)}}{4194263373300}$$
$$-\frac{2665769554417\,x^{19}}{922806209689076963250000}-\frac{19049\,x^9}{58378320}-\frac{856718372782751\,9x^{18}}{1493697149036575102620000}$$
$$-\frac{235781101\,x^{\left(\frac{39}{2}\right)}}{2920266172074657656250}$$

$$Y_2(x) = x + \frac{8\, x^{\left(\frac{17}{2}\right)}}{133875} + \frac{79268534468269\, x^{\left(\frac{33}{2}\right)}}{6295275340767174762000} + \frac{1811897241977\, x^{\left(\frac{31}{2}\right)}}{13843905831756000000} + \frac{1\, x^5}{40}$$

$$+ \frac{464403193394629\, x^{20}}{1169315090703242429821875000} + \frac{11118550640490221\, x^{\left(\frac{35}{2}\right)}}{1494780101967592440000000}$$

$$+ \frac{9641\, x^{\left(\frac{21}{2}\right)}}{725830875} + \frac{159944443\, x^{14}}{188166168590400} + \frac{1\, x^3}{6} + \frac{1\, x^7}{336} + \frac{4\, x^{\left(\frac{7}{2}\right)}}{35} + \frac{34\, x^{\left(\frac{11}{2}\right)}}{1485}$$

$$+ \frac{533041031\, x^{\left(\frac{27}{2}\right)}}{260513253000000} + \frac{21268766413501\, x^{15}}{178507849167648000000} + \frac{2\, x^6}{525} + \frac{163\, x^{\left(\frac{15}{2}\right)}}{52650}$$

$$+ \frac{11949630192209\, x^{16}}{205638741388080000000} + \frac{35284144904166467\, x^{17}}{4919961935619047116800000} + \frac{11961461\, x^{\left(\frac{25}{2}\right)}}{6975753159375}$$

$$+ \frac{45649\, x^{11}}{1884960000} + \frac{12110639\, x^{12}}{1146258313200} + \frac{273341\, x^{10}}{2554051500} + \frac{4417210841\, x^{13}}{2463876455040000}$$

$$+ \frac{271\, x^{\left(\frac{19}{2}\right)}}{835380} + \frac{80144783183\, x^{\left(\frac{29}{2}\right)}}{490270441607446500} + \frac{1189\, x^8}{1455300}$$

$$+ \frac{2028785616907679363\, x^{\left(\frac{37}{2}\right)}}{2470418656452573126300000000} + \frac{90119\, x^{\left(\frac{23}{2}\right)}}{3227908320}$$

$$+ \frac{506175625325751349\, x^{19}}{4600587219032651316069600000} + \frac{743661398\, x^{\left(\frac{41}{2}\right)}}{76326698392914891796875} + \frac{1\, x^9}{3456}$$

$$+ \frac{6519635286663831667\, x^{18}}{188095196545346494404000000} + \frac{1360265524133\, x^{\left(\frac{39}{2}\right)}}{156857701876641265500000}$$

$$+ \frac{1\, x^{21}}{10137159333984375}$$

$\delta(0.2) = 0.1963749884\ 10^{-17}, x = 0.2$
$\delta(1.2) = 0.5974263867\ 10^{-5}, x = 1.2$
$\delta(2.2) = 0.02477686822, x = 2.2$

$\delta(0.2) = 0.2507760209\ 10^{-17}, x = 0.2$
$\delta(1.2) = 0.3859974838\ 10^{-5}, x = 1.2$
$\delta(2.2) = 0.01710772920, x = 2.2$

The above-stated means that the calculation of summands is a very cumbersome task. In this case it is useful to utilize the following program developed to calculate partial solutions of linear ODEs of the following form:

$$\frac{d^m}{dx^m} \mathrm{Y}(x) = \sum_{p=1}^{m} a_p(x) \left( \frac{d^{m-p}}{dx^{m-p}} \mathrm{Y}(x) \right)$$

where $m$ is the order of the ODE, and $a_p(x)$ are the specified functions.
Partial solutions - - are calculated only with the relative accuracy defined by just seven terms -, which should be enough to determine the general analytical form of the desired solution.

**Example No. 2. Calculate partial solutions of the following ODE:**

(11.28)

**Solution: In accordance with (14.47), let's calculate the functions**

$$b_1(x) := 0$$

$$b_2(x) := 0$$

Let's calculate the summands using the above **program**:

$$\xi_0(-1) := -\frac{x^3 \ln(x)}{6} + \frac{11\, x^3}{36}$$

$$\xi_0(-2) := -\frac{x^4 \ln(x)}{8} + \frac{25\, x^4}{96}$$

$$\xi_0(-3) := -\frac{x^5 \ln(x)}{20} + \frac{137\, x^5}{1200}$$

$$\xi_1(-1) := -\frac{x^6 \ln(x)^2}{720} + \frac{23\, x^6 \ln(x)}{5400} - \frac{1477\, x^6}{648000}$$

$$\xi_1(-2) := -\frac{x^7 \ln(x)^2}{840} + \frac{2711\, x^7 \ln(x)}{705600} - \frac{311287\, x^7}{148176000}$$

$$\xi_1(-3) := -\frac{x^8 \ln(x)^2}{1920} + \frac{101\, x^8 \ln(x)}{57600} - \frac{1349\, x^8}{1382400}$$

$$\xi_2(-3) := -\frac{x^{11} \ln(x)^3}{887040} + \frac{11381\, x^{11} \ln(x)^2}{2276736000} - \frac{262430569\, x^{11} \ln(x)}{47333341440000} + \frac{4337898706921\, x^{11}}{2624160449433600000}$$

$$\xi_2(-2) := -\frac{x^{10} \ln(x)^3}{403200} + \frac{173\, x^{10} \ln(x)^2}{16128000} - \frac{1671821\, x^{10} \ln(x)}{142248960000} + \frac{1246810127\, x^{10}}{358467379200000}$$

$$\xi_2(-1) := -\frac{x^9 \ln(x)^3}{362880} + \frac{1177\, x^9 \ln(x)^2}{101606400} - \frac{534073\, x^9 \ln(x)}{42674688000} + \frac{78985223\, x^9}{21508042752000}$$

$$\xi_3(-3) := -\frac{x^{14} \ln(x)^4}{1117670400} + \frac{20291\, x^{14} \ln(x)^3}{4130909798400} - \frac{129227873\, x^{14} \ln(x)^2}{15904002723840000} + \frac{1482179880221\, x^{14} \ln(x)}{293905970336563200000} - \frac{11259303206003\, x^{14}}{1122186432194150400000}$$

$$\xi_3(-2) := -\frac{x^{13}\ln(x)^4}{518918400} + \frac{1297129\,x^{13}\ln(x)^3}{124664956416000} - \frac{42571146511\,x^{13}\ln(x)^2}{2495792427448320000} + \frac{631077790360889\,x^{13}\ln(x)}{5995891727701843968000} - \frac{2682230364401164069\,x^{13}}{12861187755920455311360000000}$$

$$\xi_3(-1) := -\frac{x^{12}\ln(x)^4}{479001600} + \frac{48977\,x^{12}\ln(x)^3}{4425974784000} - \frac{3400217\,x^{12}\ln(x)^2}{189333365760000} + \frac{2080291347473\,x^{12}\ln(x)}{188939552359219200000} - \frac{270831271860293\,x^{12}}{12470010455708467200000}$$

and so on.

Then the required solutions, which are formed by the first seven values: , are defined by:
(11.29)

$$Y_1(x) = -\frac{\ln(x)^4 x^{12}}{479001600} + \left(-\frac{x^9}{362880} + \frac{48977\,x^{12}}{4425974784000}\right)\ln(x)^3 + \left(-\frac{x^6}{720} + \frac{1177\,x^9}{101606400} - \frac{3400217\,x^{12}}{189333365760000}\right)\ln(x)^2 + \left(\frac{23\,x^6}{5400} + \frac{2080291347473\,x^{12}}{188939552359219200000} - \frac{534073\,x^9}{42674688000} - \frac{x^3}{6}\right)\ln(x) + \frac{11\,x^3}{36} - \frac{270831271860293\,x^{12}}{12470010455708467200000} - 1 + \frac{78985223\,x^9}{21508042752000} - \frac{1477\,x^6}{648000}$$

$$Y_2(x) = \frac{x^{13}\ln(x)^4}{6227020800} + \left(\frac{x^{10}}{3628800} - \frac{494339\,x^{13}}{747989738496000}\right)\ln(x)^3 + \left(\frac{2250513983\,x^{13}}{2495792427448320000} + \frac{1\,x^7}{5040} - \frac{871\,x^{10}}{1016064000}\right)\ln(x)^2 + \left(-\frac{883\,x^7}{2116800} + \frac{325267\,x^{10}}{426746880000} - \frac{374028577528049\,x^{13}}{770900364990237081600000} + \frac{1\,x^4}{24}\right)\ln(x) - \frac{13\,x^4}{288} + \frac{16316172670946377\,x^{13}}{188980718046178118860800000} - \frac{29832967\,x^{10}}{153628876800000} + \frac{79361\,x^7}{444528000} + x$$

$$Y_3(x) = -\frac{x^{14}\ln(x)^4}{87178291200} + \left(\frac{5233\,x^{14}}{130898204236800} - \frac{x^{11}}{39916800}\right)\ln(x)^3 + \left(-\frac{111280619\,x^{14}}{232940626561843200} - \frac{1\,x^8}{40320} + \frac{197\,x^{11}}{3073593600}\right)\ln(x)^2 + \left(\frac{347\,x^8}{8467200} - \frac{1\,x^5}{120} - \frac{2784139\,x^{11}}{56800009728000} + \frac{21765605030527\,x^{14}}{944352947113040424960000}\right)\ln(x) - \frac{1\,x^2}{2} + \frac{47\,x^5}{7200} - \frac{175478272585183\,x^{14}}{46874246283974552002560000} + \frac{622561721\,x^{11}}{56232009630720000} - \frac{20921\,x^8}{1422489600}$$

From here we can easily write down the following functional representations for the required solutions:

$$Y_1(x) = \sum_{i=0}^{N}\left(\sum_{k=3i}^{3N} r_{k,i}\,x^k\right)\ln(x)^i \qquad (11.30)$$

$$Y_2(x) = \sum_{i=0}^{N} \left( \sum_{k=3i}^{3N} \ln(x)^i g_{k,i} x^{(k+1)} \right) \tag{11.31}$$

$$Y_3(x) = \sum_{i=0}^{N} \left( \sum_{k=3i}^{3N} q_{k,i} x^{(k+2)} \ln(x)^i \right) \tag{11.32}$$

where are numerical coefficients, $N$ is a natural number tending to infinity.
Substituting the obtained functional values to the initial equation (11.28), we obtain the following recurrence equations for the unknown coefficients:

$$k(k-1)(k-2) r_{k,i} + (3k^2 - 6k + 2)(i+1) r_{k,i+1} + (i+3)(i+2)(i+1) r_{k,i+3} + 3(i+2)(i+1)(k-1) r_{k,i+2} - r_{k-3,i-1} = 0$$

$$3ik(i-1) g_{k,i} + (3k^2 - 1)(i-1) g_{k,i-1} + k(k-1)(k+1) g_{k,i-2} + i(i-1)(1+i) g_{k,1+i} - g_{k-3,i-3} = 0$$

$$i(2 + 3k^2 + 6k) q_{k,i} + 3i(k+1)(i+1) q_{k,i+1} + i(i+2)(i+1) q_{k,i+2} + k^3 q_{k,i-1} + 3k^2 q_{k,i-1} + 2k q_{k,i-1} - q_{k-3,i-2} = 0$$

**The task is solved. Note: In the range of change of $x$ value from -10 to 10 the relative accuracy of the solutions calculated successively for $Y_1(x)$, $Y_2(x)$, $Y_3(x)$ is determined by the following graphs:**

$Y_1(x)$

$\delta(-3.2) = 0.01361087762, x = -3.2$
$\delta(-2.2) = 0.0005106432648, x = -2.2$
$\delta(-1.2) = 0.1507409925\ 10^{-5}, x = -1.2$
$\delta(-0.2) = 0.2911197268\ 10^{-14}, x = -0.2$
$\delta(0.8) = 0.4705009829\ 10^{-9}, x = 0.8$
$\delta(1.8) = 0.534002739\ 10^{-6}, x = 1.8$
$\delta(2.8) = 0.0000146381539, x = 2.8$
$\delta(3.8) = 0.000301726876, x = 3.8$
$\delta(4.8) = 0.0078713024, x = 4.8$
$\delta(5.8) = 0.0468762405, x = 5.8$

$Y_2(x)$

$\delta(-4.2) = 0.05292621898, x = -4.2$
$\delta(-3.2) = 0.004850226320, x = -3.2$
$\delta(-2.2) = 0.0001295937217, x = -2.2$
$\delta(-1.2) = 0.1565109996\ 10^{-6}, x = -1.2$
$\delta(-0.2) = 0.1872934287\ 10^{-15}, x = -0.2$
$\delta(0.8) = 0.1747630201\ 10^{-10}, x = 0.8$
$\delta(1.8) = 0.309035150\ 10^{-8}, x = 1.8$
$\delta(2.8) = 0.45413596\ 10^{-6}, x = 2.8$
$\delta(3.8) = 0.000101169727, x = 3.8$
$\delta(4.8) = 0.00155474534, x = 4.8$
$\delta(5.8) = 0.0061725691, x = 5.8$
$\delta(6.8) = 0.0069026771, x = 6.8$
$\delta(7.8) = 0.017802343, x = 7.8$
$\delta(8.8) = 0.086334797, x = 8.8$

$Y_3(x)$

$\delta(-4.2) = 0.02433135682, x = -4.2$
$\delta(-3.2) = 0.001728030305, x = -3.2$
$\delta(-2.2) = 0.00002730275997, x = -2.2$
$\delta(-1.2) = 0.2218154088\ 10^{-7}, x = -1.2$
$\delta(-0.2) = 0.2424658015\ 10^{-16}, x = -0.2$
$\delta(0.8) = 0.1627211974\ 10^{-11}, x = 0.8$
$\delta(1.8) = 0.113622098\ 10^{-9}, x = 1.8$
$\delta(2.8) = 0.37596498\ 10^{-7}, x = 2.8$
$\delta(3.8) = 0.673065013\ 10^{-5}, x = 3.8$
$\delta(4.8) = 0.0000309136272, x = 4.8$
$\delta(5.8) = 0.001597850951, x = 5.8$
$\delta(6.8) = 0.01173322970, x = 6.8$
$\delta(7.8) = 0.04392281537, x = 7.8$

As we can see, even with only first four summands we obtain a sufficiently high accuracy (i.e. low relative error) for the established partial solutions in the range of ***x*** from -3 to 5.

**Example No. 3. Find partial solutions of the following ODE:**

$$\left(\frac{d^5}{dx^5}\mathrm{Y}(x)\right) - x\,\mathbf{e}^{(-x)}\,\mathrm{Y}(x) \qquad (11.33)$$

**Solution. In accordance with (14.47), let's calculate the functions.**

Let's calculate the summands with the use of the above **program**:

$\xi_0(-2) := -5\,(6 + x)\,\mathbf{e}^{(-x)}$
$\xi_0(-3) := 15\,(7 + x)\,\mathbf{e}^{(-x)}$
$\xi_0(-4) := -35\,(8 + x)\,\mathbf{e}^{(-x)}$
$\xi_0(-5) := 70\,(9 + x)\,\mathbf{e}^{(-x)}$
$\xi_1(-1) := -\frac{1}{32}\,\mathbf{e}^{(-2x)}\,(20 + 10x + x^2)$
$\xi_1(-2) := \frac{15}{128}\,\mathbf{e}^{(-2x)}\,(47 + 22x + 2x^2)$
$\xi_1(-3) := -\frac{5}{256}\,\mathbf{e}^{(-2x)}\,(1359 + 600x + 50x^2)$
$\xi_1(-4) := \frac{15}{256}\,\mathbf{e}^{(-2x)}\,(1582 + 663x + 51x^2)$
$\xi_1(-5) := -\frac{15}{2}\,\mathbf{e}^{(-2x)}\,(35 + 14x + x^2)$
$\xi_2(-5) := \frac{5}{15116544}\,\mathbf{e}^{(-3x)}\,(22808765 + 3413217x^2 + 179643x^3 + 17066085x)$
$\xi_2(-4) := -\frac{5}{5038848}\,\mathbf{e}^{(-3x)}\,(2359890 + 1815173x + 376758x^2 + 20931x^3)$

$$\xi_2(-3) := \frac{5}{1679616} e^{(-3x)} (194659 + 154347 x + 33354 x^2 + 1962 x^3)$$

$$\xi_2(-2) := -\frac{5}{839808} e^{(-3x)} (17113 + 14031 x + 3168 x^2 + 198 x^3)$$

$$\xi_2(-1) := \frac{1}{69984} e^{(-3x)} (670 + 570 x + 135 x^2 + 9 x^3)$$

$$\xi_3(-5) := -\frac{5}{7739670528} e^{(-4x)} (70027060 x + 24906960 x^2 + 3320928 x^3 + 138372 x^4 + 61902395)$$

$$\xi_3(-4) := \frac{5}{660451885056} e^{(-4x)} (1502290155 + 1732241320 x + 630997648 x^2 + 86768880 x^3 + 3772560 x^4)$$

$$\xi_3(-3) := -\frac{5}{55037657088} e^{(-4x)} (28113835 + 33098272 x + 12371700 x^2 + 1758240 x^3 + 79920 x^4)$$

$$\xi_3(-2) := \frac{5}{13759414272} e^{(-4x)} (3600 x^4 + 1109483 + 1336090 x + 513540 x^2 + 75600 x^3)$$

$$\xi_3(-1) := -\frac{1}{2293235712} e^{(-4x)} (77065 + 95120 x + 37680 x^2 + 5760 x^3 + 288 x^4)$$

**and so on.**

Then the required solutions, which are formed by the first seven values - , are defined by:

(11.34)

$$Y_1(x) = \left( -\frac{785 x^2}{47775744} - \frac{5945 x}{143327232} - \frac{77065}{2293235712} - \frac{1 x^4}{7962624} - \frac{5 x^3}{1990656} \right) e^{(-4x)} + \left( \frac{1 x^3}{7776} + \frac{95 x}{11664} + \frac{5 x^2}{2592} + \frac{335}{34992} \right) e^{(-3x)} + \left( -\frac{5}{8} - \frac{5 x}{16} - \frac{1 x^2}{32} \right) e^{(-2x)} - 1$$

$$Y_2(x) = \left( \frac{4195 x^3}{95551488} + \frac{261535 x^2}{1146617856} + \frac{1 x^5}{7962624} + \frac{365 x^4}{95551488} + \frac{892855 x}{1719926784} + \frac{5547415}{13759414272} \right) e^{(-4x)} + \left( -\frac{145 x^3}{46656} - \frac{1 x^4}{7776} - \frac{26065 x}{279936} - \frac{35 x^2}{1296} - \frac{85565}{839808} \right) e^{(-3x)} + \left( \frac{205 x}{64} + \frac{35 x^2}{64} + \frac{1 x^3}{32} + \frac{705}{128} \right) e^{(-2x)} + (-10 x - 30 - x^2) e^{(-x)} + x$$

$$Y_3(x) = \left( -\frac{16415\,x^4}{382205952} - \frac{1\,x^6}{15925248} - \frac{420905\,x^3}{1146617856} - \frac{15640085\,x}{4586471424} - \frac{245\,x^5}{95551488} - \frac{11188135\,x^2}{6879707136} - \frac{140569175}{55037657088} \right) \mathbf{e}^{(-4x)} + \left( \frac{25\,x^4}{11664} + \frac{942865\,x}{1679616} + \frac{895\,x^3}{31104} + \frac{973295}{1679616} + \frac{1\,x^5}{15552} + \frac{6565\,x^2}{34992} \right) \mathbf{e}^{(-3x)} + \left( -\frac{6795}{256} - \frac{2205\,x}{128} - \frac{25\,x^3}{64} - \frac{495\,x^2}{128} - \frac{1\,x^4}{64} \right) \mathbf{e}^{(-2x)} + \left( 105 + 45\,x + \frac{15\,x^2}{2} + \frac{1\,x^3}{2} \right) \mathbf{e}^{(-x)} - \frac{1\,x^2}{2}$$

$$Y_4(x) = \left( \frac{3969745\,x^4}{13759414272} + \frac{775565\,x^3}{382205952} + \frac{659252815\,x^2}{82556485632} + \frac{27215\,x^5}{1146617856} + \frac{2587009175\,x}{165112971264} + \frac{1\,x^7}{47775744} + \frac{2503816925}{220150628352} + \frac{205\,x^6}{191102976} \right) \mathbf{e}^{(-4x)} + \left( -\frac{4655\,x^4}{279936} - \frac{655525}{279936} - \frac{137245\,x^3}{839808} - \frac{5997875\,x}{2519424} - \frac{1\,x^6}{46656} - \frac{85\,x^5}{93312} - \frac{742615\,x^2}{839808} \right) \mathbf{e}^{(-3x)} + \left( \frac{65\,x^4}{384} + \frac{1\,x^5}{192} + \frac{455\,x^3}{192} + \frac{11865}{128} + \frac{2235\,x^2}{128} + \frac{4185\,x}{64} \right) \mathbf{e}^{(-2x)} + \left( -\frac{10\,x^3}{3} - \frac{1\,x^4}{6} - 280 - 30\,x^2 - 140\,x \right) \mathbf{e}^{(-x)} + \frac{1\,x^3}{6}$$

$$Y_5(x) = \left( -\frac{112168989125\,x}{1981355655168} - \frac{185\,x^7}{573308928} - \frac{701155625\,x^3}{82556485632} - \frac{229600625\,x^4}{165112971264} - \frac{1943605\,x^5}{13759414272} - \frac{10065795625\,x^2}{330225942528} - \frac{309511975}{7739670528} - \frac{1\,x^8}{191102976} - \frac{20395\,x^6}{2293235712} \right) \mathbf{e}^{(-4x)} + \left( \frac{155\,x^6}{559872} + \frac{142325\,x^4}{1679616} + \frac{40242925\,x}{5038848} + \frac{1\,x^7}{186624} + \frac{1195\,x^5}{186624} + \frac{6852725\,x^3}{10077696} + \frac{133535\,x^2}{41472} + \frac{114043825}{15116544} \right) \mathbf{e}^{(-3x)} + \left( -\frac{25305\,x}{128} - \frac{625\,x^3}{64} - \frac{30525\,x^2}{512} - \frac{5\,x^5}{96} - \frac{1\,x^6}{768} - \frac{525}{2} - \frac{725\,x^4}{768} \right) \mathbf{e}^{(-2x)} + \left( \frac{25\,x^3}{2} + \frac{25\,x^4}{24} + \frac{1\,x^5}{24} + 350\,x + \frac{175\,x^2}{2} + 630 \right) \mathbf{e}^{(-x)} - \frac{1\,x^4}{24}$$

From here we can easily write down the following functional representations for the required solutions:

$$Y_{1+s}(x) = \sum_{i=0}^{N} \left( \sum_{k=0}^{i} r_{k,i}\, x^{(k+s)} \right) \mathbf{e}^{(-i\,x)},\ s = 0\,..\,4 \qquad (11.35)$$

where , $t_{k,i}, p_{k,i}$ are numerical coefficients, $N$ is a natural number tending to infinity. Substituting the obtained functional values to the initial equation (11.33) we obtain the following recurrence equations for the unknown numerical coefficients:

$$-k\,(k-1)\,(k-2)\,(k-3)\,(1+k)\,r_{1+k,i} + r_{k-4,i}\,i^5 - 5\,i^4\,(k-3)\,r_{k-3,i} + r_{k-5,i-1} - 10\,i^2\,(k-2)\,(k-3)\,(k-1)\,r_{k-1,i} + 5\,i\,(k-2)\,(k-3)\,k\,(k-1)\,r_{k,i} + 10\,i^3\,(k-2)\,(k-3)\,r_{k-2,i} = 0$$

$$-k(k-1)(k-2)(2+k)(1+k)\,g_{1+k,i}+i^5 g_{k-4,i}-5\,i^4(k-2)\,g_{k-3,i}+g_{k-5,i-1}$$
$$-10\,k\,i^2(k-1)(k-2)\,g_{k-1,i}+5\,k\,i(k-1)(k-2)(1+k)\,g_{k,i}$$
$$+10\,i^3 g_{k-2,i}(k-1)(k-2)=0$$

$$-k(k-1)(3+k)(k+2)(1+k)\,q_{1+k,i}+q_{k-4,i}\,i^5-5\,i^4(k-1)\,q_{k-3,i}+q_{k-5,i-1}$$
$$-10\,i^2 k(k-1)(1+k)\,q_{k-1,i}+5\,i\,k(k-1)(k+2)(1+k)\,q_{k,i}$$
$$+10\,i^3 k\,q_{k-2,i}(k-1)=0$$

$$5\,i(4+k)(3+k)(k+2)(1+k)\,t_{1+k,i}+i^5 t_{k-3,i}+10\,i^3(k+2)(1+k)\,t_{k-1,i}$$
$$-10\,i^2(3+k)(k+2)(1+k)\,t_{k,i}-(5+k)(4+k)(3+k)(k+2)(1+k)\,t_{k+2,i}$$
$$-5\,i^4(1+k)\,t_{k-2,i}+t_{k-4,i-1}=0$$

$$-(5+k)(4+k)(3+k)(k+2)(1+k)\,p_{1+k,i}+p_{k-4,i}\,i^5-5\,i^4(1+k)\,p_{k-3,i}+p_{k-5,i-1}$$
$$-10\,i^2(3+k)(k+2)(1+k)\,p_{k-1,i}+5\,i(4+k)(3+k)(k+2)(1+k)\,p_{k,i}$$
$$+10\,i^3(k+2)(1+k)\,p_{k-2,i}=0$$

**The task is solved.**
**Example No. 4. Establish the functional structure of the partial solutions of following ODE:**

$$\left(\frac{d^7}{dx^7}\mathrm{Y}(x)\right)+\left(\frac{d^4}{dx^4}\mathrm{Y}(x)\right)-x\left(\frac{d}{dx}\mathrm{Y}(x)\right)-\sin(x)\,\mathrm{Y}(x)=0 \tag{11.36}$$

**Solution. In accordance with (14.47), let's calculate the functions**

$$b_3(x):=1$$

$$b_6(x):=x$$

Let's calculate the $\xi_k(-n),\ n=1..7,\ k=0,1,2..N$ summands using the above **program**:

$$\xi_0(-2):=\frac{x^8}{40320}+\frac{x^4}{8}-7\sin(x)$$
$$\xi_0(-3):=\frac{x^9}{51840}+\frac{x^5}{20}+28\cos(x)$$
$$\xi_0(-4):=\frac{x^{10}}{129600}+\frac{x^6}{72}+84\sin(x)$$
$$\xi_0(-5):=\frac{x^{11}}{475200}+\frac{x^7}{336}-210\cos(x)$$
$$\xi_0(-6):=\frac{x^{12}}{2280960}+\frac{x^8}{1920}-462\sin(x)$$

$$\xi_0(-7) := \frac{x^{13}}{13478400} + \frac{x^9}{12960} + 924\cos(x)$$

$$\xi_1(-1) := \frac{x^{10}}{1209600} + \frac{x^6}{320} - \frac{545\,x^4}{384} + \frac{10625\,x^2}{128} + \left(\frac{1}{6}x^3 - 27\,x\right)\cos(x) + \left(-\frac{7\,x^2}{2} + 76\right)\sin(x) - \frac{1}{256}\cos(2\,x) - \frac{141569}{256}$$

$$\xi_1(-2) := \frac{x^{15}}{163459296000} + \frac{x^{11}}{1425600} + \frac{53\,x^7}{20160} - \frac{679\,x^5}{640} + \frac{5383\,x^3}{128} + \frac{102123\,x}{256} - \frac{21}{512}\sin(2\,x) + \left(\frac{1}{40320}x^8 - \frac{7}{180}x^6 + \frac{71}{8}x^4 - 476\,x^2 + 2901\right)\cos(x) + \left(-\frac{1}{720}x^7 + \frac{7}{10}x^5 - \frac{238}{3}x^3 + 1729\,x\right)\sin(x)$$

$$\xi_1(-3) := \frac{x^{16}}{190207180800} + \frac{29\,x^{12}}{95800320} + \frac{13\,x^8}{11520} - \frac{2359\,x^6}{5760} + \frac{7735\,x^4}{768} + \frac{66185\,x^2}{256} + \frac{119}{512}\cos(2\,x) + \frac{1215095}{512} + \left(-\frac{1}{960}x^8 + \frac{7}{15}x^6 - \frac{315}{8}x^4 + 10338\right)\sin(x) + \left(\frac{1}{51840}x^9 - \frac{1}{36}x^7 + \frac{53}{10}x^5 - \frac{476}{3}x^3 - 3094\,x\right)\cos(x)$$

$$\xi_1(-4) := \frac{29\,x^{17}}{12703122432000} + \frac{x^{13}}{11321856} + \frac{17\,x^9}{51840} - \frac{103\,x^7}{960} + \frac{609\,x^5}{640} + \frac{14175\,x^3}{128} - \frac{42525\,x}{256} + \left(\frac{1}{129600}x^{10} - \frac{1}{96}x^8 + \frac{127}{72}x^6 - \frac{119}{3}x^4 - 27240\right)\cos(x) + \left(-\frac{7}{17280}x^9 + \frac{1}{6}x^7 - \frac{707}{60}x^5 - 5460\,x\right)\sin(x) + \frac{483}{512}\sin(2\,x)$$

$$\xi_1(-5) := \frac{61\,x^{18}}{91462481510400} + \frac{31\,x^{14}}{1585059840} + \frac{x^{10}}{13824} - \frac{979\,x^8}{46080} - \frac{4907\,x^6}{23040} + \frac{108395\,x^4}{3072} - \frac{108395\,x^2}{1024} + \left(-\frac{7}{64800}x^{10} + \frac{1}{24}x^8 - \frac{1883}{720}x^6 - 65163\right)\sin(x) + \left(\frac{1}{475200}x^{11} - \frac{7}{2592}x^9 + \frac{47}{112}x^7 - \frac{119}{15}x^5 + 9828\,x\right)\cos(x) - \frac{12564629}{2048} - \frac{6293}{2048}\cos(2\,x)$$

$$\xi_1(-6) := \frac{x^{19}}{6788231049600} + \frac{x^{15}}{285768000} + \frac{49\,x^{11}}{3801600} - \frac{31\,x^9}{9216} - \frac{179\,x^7}{1536} + \frac{45689\,x^5}{5120} - \frac{45689\,x^3}{1024} + \frac{137067\,x}{2048} + \left(-\frac{7}{316800}x^{11} + \frac{7}{864}x^9 - \frac{7}{15}x^7 + 18564\,x\right)\sin(x) - \frac{34965}{4096}\sin(2\,x) + \left(\frac{1}{2280960}x^{12} - \frac{7}{12960}x^{10} + \frac{151}{1920}x^8 - \frac{119}{90}x^6 + 151795\right)\cos(x)$$

$$\xi_1(-7) := \frac{67\,x^{20}}{2555569336320000} + \frac{11\,x^{16}}{20901888000} + \frac{x^{12}}{518400} - \frac{19\,x^{10}}{43200} - \frac{29\,x^8}{960} + \frac{28\,x^6}{15} - 14\,x^4 + 42\,x^2 + \left(-\frac{7}{1900800}x^{12} + \frac{7}{5400}x^{10} - \frac{403}{5760}x^8 + 352716\right)\sin(x) + \left(\frac{1}{13478400}x^{13} - \frac{7}{79200}x^{11} + \frac{317}{25920}x^9 - \frac{17}{90}x^7 - 37128\,x\right)\cos(x) + 13545 + 21\cos(2\,x)$$

$$\xi_2(-7) := \frac{79\,x^{27}}{146391189406285824000000} + \frac{61\,x^{23}}{988113623777280000} + \frac{661\,x^{19}}{511113867264000} - \frac{17359\,x^{17}}{9238634496000} - \frac{12211\,x^{15}}{747242496000} + \frac{4639199\,x^{13}}{28466380800} - \frac{1862999\,x^{11}}{116785152}$$

$$- \frac{17558041\, x^{9}}{11796480} + \frac{535356211\, x^{7}}{2949120} - \frac{1181402425\, x^{5}}{393216} + \frac{6556604957\, x^{3}}{393216}$$
$$- \frac{735245621789\, x}{262144} + \left( \frac{67}{2555569336320000} x^{20} - \frac{59}{418784256000} x^{18} \right.$$
$$+ \frac{3071}{20901888000} x^{16} - \frac{5791}{94348800} x^{14} + \frac{67241}{5702400} x^{12} - \frac{53953}{51840} x^{10} + \frac{358751}{11520} x^{8}$$
$$\left. + \frac{81641}{5760} x^{6} - \frac{81641}{768} x^{4} + \frac{761258729}{256} x^{2} - \frac{3886171440379}{8192} \right) \cos(x) + \left( \right.$$
$$- \frac{7}{2652300288000} x^{19} + \frac{11}{2115072000} x^{17} - \frac{6269}{1886976000} x^{15} + \frac{69323}{74131200} x^{13}$$
$$- \frac{18673}{152064} x^{11} + \frac{2842027}{414720} x^{9} - \frac{1267061}{17920} x^{7} - \frac{215649}{5120} x^{5} + \frac{215649}{1024} x^{3}$$
$$\left. - \frac{152236064547}{2048} x \right) \sin(x)$$
$$+ \left( - \frac{7}{324403200} x^{12} + \frac{161}{11059200} x^{10} - \frac{5023}{3932160} x^{8} + \frac{6734357987}{524288} \right) \sin(2\, x)$$
$$+ \frac{53845547}{4353564672} \cos(3\, x) +$$
$$\left( \frac{1}{3450470400} x^{13} - \frac{119}{162201600} x^{11} + \frac{18911}{106168320} x^{9} - \frac{1553}{368640} x^{7} - 1239\, x \right) \cos(2\, x)$$
$$\xi_2(-6) := \frac{2579\, x^{26}}{8729252405337784320000000} + \frac{157\, x^{22}}{386653157130240000} + \frac{1933\, x^{18}}{228656203776000}$$
$$- \frac{1769\, x^{16}}{135862272000} + \frac{404123\, x^{14}}{697426329600} + \frac{93521\, x^{12}}{99532800} - \frac{31632469\, x^{10}}{265420800} - \frac{6055653\, x^{8}}{1310720}$$
$$+ \frac{2464827253\, x^{6}}{2949120} - \frac{4449630143\, x^{4}}{393216} + \frac{4924507217\, x^{2}}{131072} - \frac{27524075747}{262144} +$$
$$\left( \frac{1}{583925760} x^{12} - \frac{119}{26542080} x^{10} + \frac{1501}{1310720} x^{8} - \frac{10871}{368640} x^{6} + \frac{1529605917}{262144} \right) \cos(2\, x)$$
$$+ \left( - \frac{1}{67005480960} x^{18} + \frac{13}{435456000} x^{16} - \frac{2449}{125798400} x^{14} + \frac{8027}{1425600} x^{12} - \frac{827}{1080} x^{10} \right.$$
$$\left. + \frac{258421}{5760} x^{8} - \frac{143641}{288} x^{6} - \frac{34951}{192} x^{4} - \frac{95112185}{64} x^{2} + \frac{1607459972267}{8192} \right) \sin(x)$$
$$+ \left( - \frac{7}{54067200} x^{11} + \frac{161}{1769472} x^{9} - \frac{2791}{327680} x^{7} + \frac{87357879}{131072} x \right) \sin(2\, x)$$
$$- \frac{6462295}{1451188224} \sin(3\, x) + \left( \frac{1}{67882310496000} x^{19} - \frac{7}{8724672000} x^{17} + \frac{3163}{3714984000} x^{15} \right.$$
$$- \frac{6001}{16473600} x^{13} + \frac{823063}{11404800} x^{11} - \frac{2751451}{414720} x^{9} + \frac{6796817}{32256} x^{7} + \frac{1127819}{15360} x^{5}$$
$$\left. - \frac{1127819}{3072} x^{3} - \frac{69616294517}{2048} x \right) \cos(x)$$
$$\xi_2(-5) := \frac{4799\, x^{25}}{369314524841213952000000} + \frac{83\, x^{21}}{37103080734720000} + \frac{47\, x^{17}}{1016249794560}$$
$$- \frac{112561\, x^{15}}{1494484992000} + \frac{198007\, x^{13}}{28466380800} + \frac{22660189\, x^{11}}{5109350400} - \frac{26925581\, x^{9}}{37158912} + \frac{34423\, x^{7}}{5160960}$$
$$+ \frac{148708147\, x^{5}}{49152} - \frac{139347971\, x^{3}}{4096} + \frac{27306124909\, x}{32768}$$
$$+ \left( \frac{1}{121651200} x^{11} - \frac{119}{5308416} x^{9} + \frac{6305}{1032192} x^{7} - \frac{10871}{61440} x^{5} + \frac{23591351}{65536} x \right) \cos(2\, x)$$

$$
-\frac{339227}{241864704}\cos(3x)+\left(-\frac{19}{279189504000}x^{17}+\frac{7}{50544000}x^{15}-\frac{499}{5391360}x^{13}+\frac{4931}{178200}x^{11}-\frac{407567}{103680}x^{9}+\frac{1978693}{8064}x^{7}-\frac{11568721}{3840}x^{5}-\frac{481775}{768}x^{3}+\frac{7681818095}{512}x\right)\sin(x)+\left(-\frac{7}{11059200}x^{10}+\frac{23}{49152}x^{8}-\frac{2345}{49152}x^{6}-\frac{336800989}{131072}\right)\sin(2x)+\left(\frac{61}{91462481510400}x^{18}-\frac{53}{14370048000}x^{16}+\frac{487}{121927680}x^{14}-\frac{1003}{570240}x^{12}+\frac{375211}{1036800}x^{10}-\frac{3766711}{107520}x^{8}+\frac{27766757}{23040}x^{6}+\frac{323561}{1024}x^{4}-\frac{747030459}{1024}x^{2}+\frac{311360082703}{4096}\right)\cos(x)
$$

$$
\xi_2(-4) := \frac{4453\,x^{24}}{1034080669555399065600000}+\frac{43\,x^{20}}{4336723722240000}+\frac{19\,x^{16}}{93405312000}-\frac{5807\,x^{14}}{16605388800}+\frac{732037\,x^{12}}{15328051200}+\frac{375539\,x^{10}}{23224320}-\frac{14248211\,x^{8}}{4128768}+\frac{75462193\,x^{6}}{737280}+\frac{248520959\,x^{4}}{32768}-\frac{2509967113\,x^{2}}{32768}+\frac{62459811677}{65536}+\left(\frac{29}{12703122432000}x^{17}-\frac{1}{77837760}x^{15}+\frac{1775}{124540416}x^{13}-\frac{371}{57024}x^{11}+\frac{455}{324}x^{9}-\frac{2932567}{20160}x^{7}+\frac{2153027}{384}x^{5}+\frac{416353}{384}x^{3}+\frac{1623287711}{256}x\right)\cos(x)+\left(\frac{1}{33177600}x^{10}-\frac{17}{196608}x^{8}+\frac{18919}{737280}x^{6}-\frac{10871}{12288}x^{4}-\frac{69003427}{65536}\right)\cos(2x)+\left(-\frac{1}{4257792000}x^{16}+\frac{23}{47174400}x^{14}-\frac{851}{2534400}x^{12}+\frac{27127}{259200}x^{10}-\frac{140881}{8960}x^{8}+\frac{2053619}{1920}x^{6}-\frac{11636873}{768}x^{4}+\frac{89064521}{256}x^{2}-\frac{13721905749}{512}\right)\sin(x)+\left(-\frac{7}{2949120}x^{9}+\frac{23}{12288}x^{7}-\frac{35189}{163840}x^{5}-\frac{3128811}{16384}x\right)\sin(2x)+\frac{3745}{10077696}\sin(3x)
$$

$$
\xi_2(-3) := \frac{103\,x^{23}}{107716736412020736000}+\frac{269\,x^{19}}{8109673360588800}+\frac{353\,x^{15}}{523069747200}-\frac{8711\,x^{13}}{7116595200}+\frac{34927\,x^{11}}{159667200}+\frac{381349\,x^{9}}{9289728}-\frac{62139773\,x^{7}}{5160960}+\frac{107733893\,x^{5}}{163840}+\frac{742719539\,x^{3}}{98304}-\frac{20718002027\,x}{65536}+\left(-\frac{1}{163840}x^{8}+\frac{161}{30720}x^{6}-\frac{11739}{16384}x^{4}+\frac{48945257}{131072}\right)\sin(2x)+\left(\frac{1}{190207180800}x^{16}-\frac{47}{1556755200}x^{14}+\frac{661}{19160064}x^{12}-\frac{1427}{86400}x^{10}+\frac{102121}{26880}x^{8}-\frac{276279}{640}x^{6}+\frac{5092197}{256}x^{4}+\frac{42931441}{256}x^{2}-\frac{2078104777}{256}\right)\cos(x)+\left(\frac{1}{13271040}x^{9}-\frac{17}{73728}x^{7}+\frac{6309}{81920}x^{5}-\frac{10871}{3072}x^{3}-\frac{1596917}{16384}x\right)\cos(2x)+\left(-\frac{17}{31135104000}x^{15}+\frac{43}{37065600}x^{13}-\frac{1579}{1900800}x^{11}+\frac{7093}{25920}x^{9}-\frac{297701}{6720}x^{7}+\frac{6549359}{1920}x^{5}-\frac{23393327}{384}x^{3}-\frac{654217297}{256}x\right)\sin(x)+\frac{49}{629856}\cos(3x)
$$

$$\xi_2(-2) := \frac{x^{22}}{936667273148006400} + \frac{53\,x^{18}}{711374856192000} + \frac{131\,x^{14}}{87178291200} - \frac{689\,x^{12}}{239500800} + \frac{73621\,x^{10}}{116121600} + \frac{49321\,x^{8}}{860160} - \frac{6638971\,x^{6}}{245760} + \frac{208428941\,x^{4}}{98304} - \frac{942091753\,x^{2}}{32768} + \left(\frac{1}{163459296000}x^{15} - \frac{1}{27799200}x^{13} + \frac{61}{1425600}x^{11} - \frac{8017}{362880}x^{9} + \frac{37523}{6720}x^{7} - \frac{1390901}{1920}x^{5} + \frac{6005015}{128}x^{3} - \frac{247804805}{256}x\right)\cos(x) - \frac{77}{6718464}\sin(3\,x) - \frac{37077149193}{65536} + \left(\frac{1}{10321920}x^{8} - \frac{119}{368640}x^{6} + \frac{6317}{49152}x^{4} - \frac{10871}{1024}x^{2} + \frac{6199799}{65536}\right)\cos(2\,x) + \left(-\frac{1}{1556755200}x^{14} + \frac{1}{712800}x^{12} - \frac{553}{518400}x^{10} + \frac{3859}{10080}x^{8} - \frac{49643}{720}x^{6} + \frac{205965}{32}x^{4} - \frac{8355723}{32}x^{2} + \frac{1796605139}{1024}\right)\sin(x) + \left(-\frac{1}{122880}x^{7} + \frac{161}{20480}x^{5} - \frac{11767}{8192}x^{3} + \frac{1539447}{32768}x\right)\sin(2\,x)$$

$$\xi_2(-1) := \frac{x^{17}}{11856247603200} + \frac{x^{13}}{593049600} - \frac{109\,x^{11}}{31933440} + \frac{10441\,x^{9}}{11612160} + \frac{1153\,x^{7}}{80640} - \frac{573171\,x^{5}}{20480} + \frac{17790277\,x^{3}}{6144} - \frac{674393563\,x}{8192} + \left(\frac{1}{1536}x^{3} - \frac{3743}{16384}x\right)\cos(2\,x) + \left(-\frac{21\,x^{2}}{1024} + \frac{29417}{32768}\right)\sin(2\,x) + \left(\frac{1}{1209600}x^{10} - \frac{1}{480}x^{8} + \frac{281}{320}x^{6} - \frac{15323}{128}x^{4} + \frac{696913}{128}x^{2} - \frac{21750561}{512}\right)\cos(x) - \frac{1}{1119744}\cos(3\,x) + \left(-\frac{1}{17280}x^{9} + \frac{1}{20}x^{7} - \frac{1869}{160}x^{5} + \frac{89959}{96}x^{3} - \frac{1363255}{64}x\right)\sin(x)$$

Thus

$$\alpha(-n) = \sum_{i=0}^{2} \xi_i(-n),\ n = 1\,..\,7$$

and partial solutions for these values (in this case) are as follows:

$$Y_1(x) = \frac{1\,x^{17}}{11856247603200} + \frac{1\,x^{13}}{593049600} - \frac{109\,x^{11}}{31933440} + \frac{1\,x^{10}}{1209600} + \frac{10441\,x^{9}}{11612160} + \frac{1153\,x^{7}}{80640} + \frac{1\,x^{6}}{320} - \frac{573171\,x^{5}}{20480} - \frac{545\,x^{4}}{384} + \frac{17791301\,x^{3}}{6144} + \frac{10625\,x^{2}}{128} - \frac{674393563\,x}{8192} - \frac{141825}{256} - \frac{1\cos(3\,x)}{1119744} + \left(-\frac{1\,x^{9}}{17280} + \frac{1\,x^{7}}{20} - \frac{1869\,x^{5}}{160} + \frac{89959\,x^{3}}{96} - \frac{7\,x^{2}}{2} - \frac{1363255\,x}{64} + 76\right)\sin(x)$$

$$+\left(-\frac{21x^2}{1024}+\frac{29417}{32768}\right)\sin(2x)+$$
$$\left(\frac{1x^{10}}{1209600}-\frac{1x^8}{480}+\frac{281x^6}{320}-\frac{15323x^4}{128}+\frac{1x^3}{6}+\frac{696913x^2}{128}-27x-\frac{21751073}{512}\right)$$
$$\cos(x)+\left(\frac{1x^3}{1536}-\frac{3743x}{16384}-\frac{1}{256}\right)\cos(2x)$$

$$Y_2(x)=\frac{1x^{22}}{9366672731480064000}-\frac{1x^{18}}{101624979456000}+\frac{1x^{15}}{163459296000}-\frac{1x^{14}}{5448643200}$$
$$+\frac{257x^{12}}{479001600}-\frac{1x^{11}}{7983360}-\frac{3421x^{10}}{12902400}+\frac{111131x^8}{2580480}-\frac{1x^7}{2016}+\frac{239081x^6}{245760}+\frac{43x^5}{120}$$
$$-\frac{25406529x^4}{32768}-\frac{2621x^3}{64}+\frac{1755482499x^2}{32768}+\frac{60987x}{64}+\frac{1\cos(3x)x}{1119744}$$
$$-\frac{37077149193}{65536}+\left(-\frac{1x^{14}}{1556755200}+\frac{1x^{12}}{712800}-\frac{523x^{10}}{518400}+\frac{671x^8}{2016}-\frac{1x^7}{720}\right.$$
$$\left.-\frac{16493x^6}{288}+\frac{7x^5}{10}+\frac{16498x^4}{3}-\frac{455x^3}{6}-\frac{15348191x^2}{64}+1653x+\frac{1796597971}{1024}\right)$$
$$\sin(x)+\left(-\frac{1x^7}{122880}+\frac{161x^5}{20480}-\frac{11599x^3}{8192}+\frac{755015x}{16384}-\frac{21}{512}\right)\sin(2x)-\frac{77\sin(3x)}{6718464}$$
$$+\left(\frac{1x^{15}}{163459296000}-\frac{1x^{13}}{27799200}+\frac{67x^{11}}{1596672}-\frac{7261x^9}{362880}+\frac{1x^8}{40320}+\frac{15811x^7}{3360}-\frac{7x^6}{180}\right.$$
$$\left.-\frac{36283x^5}{60}+\frac{209x^4}{24}+\frac{2654051x^3}{64}-449x^2-\frac{473858537x}{512}+2901\right)\cos(x)$$
$$+\left(\frac{1x^8}{10321920}-\frac{119x^6}{368640}+\frac{2095x^4}{16384}-\frac{170193x^2}{16384}+\frac{1x}{256}+\frac{6199799}{65536}\right)\cos(2x)$$

$$Y_3(x)=-\frac{1x^{23}}{897639470100172800000}+\frac{1x^{19}}{1192599023616000}-\frac{1x^{16}}{1162377216000}$$
$$+\frac{1x^{15}}{65383718400}-\frac{179x^{13}}{3321077760}+\frac{1x^{12}}{68428800}+\frac{137x^{11}}{3991680}-\frac{60683x^9}{6635520}+\frac{1x^8}{16128}$$
$$+\frac{9881x^7}{10080}-\frac{671x^6}{11520}-\frac{732785x^5}{49152}+\frac{457x^4}{48}-\frac{119341645x^3}{24576}-\frac{213701x^2}{512}$$
$$+\frac{77\sin(3x)x}{6718464}+\frac{8179573583x}{32768}+\frac{1215095}{512}+\left(\frac{1x^{15}}{10378368000}-\frac{1x^{13}}{4118400}\right.$$
$$+\frac{1181x^{11}}{5702400}-\frac{3055x^9}{36288}+\frac{1x^8}{2880}+\frac{189577x^7}{10080}-\frac{7x^6}{30}-\frac{1636317x^5}{640}+\frac{917x^4}{24}$$
$$\left.+\frac{4549099x^3}{24}-1691x^2-\frac{4413467159x}{1024}+10338\right)\sin(x)$$
$$+\left(\frac{1x^8}{491520}-\frac{161x^6}{61440}+\frac{11627x^4}{16384}-\frac{3049477x^2}{65536}+\frac{21x}{512}+\frac{48945257}{131072}\right)\sin(2x)+\Bigg($$
$$-\frac{1x^{16}}{1162377216000}+\frac{1x^{14}}{172972800}-\frac{49x^{12}}{6220800}+\frac{2057x^{10}}{453600}-\frac{1x^9}{181440}-\frac{5167x^8}{3840}+\frac{1x^7}{90}$$

$$+\frac{894283\,x^6}{3840}-\frac{419\,x^5}{120}-\frac{777615\,x^4}{32}+\frac{1823\,x^3}{6}+\frac{1141193911\,x^2}{1024}-5995\,x$$
$$-\frac{2078097609}{256}\Bigg)\cos(x)+$$
$$\left(-\frac{1\,x^9}{46448640}+\frac{17\,x^7}{184320}-\frac{6289\,x^5}{122880}+\frac{684515\,x^3}{98304}-\frac{1\,x^2}{512}-\frac{12587467\,x}{65536}+\frac{119}{512}\right)$$
$$\cos(2\,x)+\left(-\frac{1\,x^2}{2239488}+\frac{49}{629856}\right)\cos(3\,x)$$

$$Y_4(x)=\frac{1\,x^{24}}{1216565493594587136000}-\frac{1\,x^{20}}{16550353797120000}+\frac{1\,x^{17}}{11856247603200}$$
$$-\frac{1\,x^{16}}{871782912000}+\frac{1\,x^{14}}{207567360}-\frac{1\,x^{13}}{691891200}-\frac{2809\,x^{12}}{729907200}+\frac{10211\,x^{10}}{7257600}-\frac{1\,x^9}{145152}$$
$$-\frac{5226017\,x^8}{20643840}+\frac{1\,x^7}{120}+\frac{65846557\,x^6}{2949120}-\frac{2467\,x^5}{1280}-\frac{41004777\,x^4}{65536}+144\,x^3$$
$$-\frac{5681013591\,x^2}{131072}-\frac{1300145\,x}{512}+$$
$$\left(-\frac{1\,x^9}{2949120}+\frac{23\,x^7}{40960}-\frac{34909\,x^5}{163840}+\frac{1147231\,x^3}{49152}-\frac{21\,x^2}{1024}-\frac{73975745\,x}{131072}+\frac{483}{512}\right)$$
$$\sin(2\,x)+\Bigg(-\frac{1\,x^{16}}{99632332800}+\frac{1\,x^{14}}{34594560}-\frac{73\,x^{12}}{2534400}+\frac{947\,x^{10}}{67200}-\frac{1\,x^9}{17280}$$
$$-\frac{318529\,x^8}{80640}+\frac{1\,x^7}{20}+\frac{103751\,x^6}{144}-\frac{1379\,x^5}{120}-\frac{62392385\,x^4}{768}+\frac{5111\,x^3}{6}$$
$$+\frac{7742852515\,x^2}{2048}-15798\,x-\frac{13721862741}{512}\Bigg)\sin(x)$$
$$+\left(\frac{1\,x^3}{6718464}-\frac{49\,x}{629856}\right)\cos(3\,x)+\Bigg(\frac{1\,x^{17}}{11856247603200}-\frac{1\,x^{15}}{1556755200}$$
$$+\frac{233\,x^{13}}{230630400}-\frac{9167\,x^{11}}{13305600}+\frac{1\,x^{10}}{1209600}+\frac{181957\,x^9}{725760}-\frac{1\,x^8}{480}-\frac{56489\,x^7}{1008}+\frac{629\,x^6}{720}$$
$$+\frac{2115865\,x^5}{256}-\frac{229\,x^4}{2}-\frac{1976924225\,x^3}{3072}+\frac{9089\,x^2}{2}+\frac{462673165\,x}{32}-27240\Bigg)$$
$$\cos(x)+\Bigg(\frac{1\,x^{10}}{309657600}-\frac{17\,x^8}{983040}+\frac{18871\,x^6}{1474560}-\frac{257161\,x^4}{98304}+\frac{1\,x^3}{1536}+\frac{18975135\,x^2}{131072}$$
$$-\frac{119\,x}{512}-\frac{69003427}{65536}\Bigg)\cos(2\,x)+\frac{62459811677}{65536}$$
$$+\left(-\frac{77\,x^2}{13436928}+\frac{3745}{10077696}\right)\sin(3\,x)$$

$$Y_5(x)=-\frac{1\,x^{25}}{19584861165821952000000}+\frac{1\,x^{21}}{258035061473280000}-\frac{1\,x^{18}}{145508493312000}$$
$$+\frac{1\,x^{17}}{12703122432000}-\frac{1\,x^{15}}{2594592000}+\frac{1\,x^{14}}{7925299200}+\frac{983\,x^{13}}{2656862208}-\frac{435899\,x^{11}}{2554675200}$$
$$+\frac{1\,x^{10}}{1451520}+\frac{984673\,x^9}{23224320}-\frac{1\,x^8}{960}-\frac{8121157\,x^7}{1290240}+\frac{14783\,x^6}{46080}+\frac{38041521\,x^5}{65536}$$
$$-\frac{219707\,x^4}{6144}-\frac{694554729\,x^3}{32768}+\frac{9975\,x^2}{8}-\frac{7847561859\,x}{65536}$$

$$+\left(\frac{77\,x^3}{40310784}-\frac{3745\,x}{10077696}\right)\sin(3\,x)-\frac{12564629}{2048}+\Bigg(\frac{1\,x^{10}}{22118400}-\frac{23\,x^8}{245760}$$

$$+\frac{11641\,x^6}{245760}-\frac{6128371\,x^4}{786432}+\frac{7\,x^3}{1024}+\frac{99006233\,x^2}{262144}-\frac{483\,x}{512}-\frac{336800989}{131072}\Bigg)\sin(2\,x)$$

$$+\Bigg(-\frac{1\,x^{11}}{2554675200}+\frac{17\,x^9}{6635520}-\frac{699\,x^7}{286720}+\frac{457591\,x^5}{655360}-\frac{1\,x^4}{6144}-\frac{25362803\,x^3}{393216}$$

$$+\frac{119\,x^2}{1024}+\frac{46297389\,x}{32768}-\frac{6293}{2048}\Bigg)\cos(2\,x)+\Bigg(\frac{1\,x^{17}}{1154829312000}-\frac{1\,x^{15}}{353808000}$$

$$+\frac{263\,x^{13}}{80870400}-\frac{4699\,x^{11}}{2494800}+\frac{1\,x^{10}}{129600}+\frac{469321\,x^9}{725760}-\frac{1\,x^8}{120}-\frac{255943\,x^7}{1680}+\frac{1841\,x^6}{720}$$

$$+\frac{186710663\,x^5}{7680}-285\,x^4-\frac{11788608239\,x^3}{6144}+10629\,x^2+\frac{5350920209\,x}{128}-65163\Bigg)$$

$$\sin(x)+\Bigg(-\frac{1\,x^{18}}{145508493312000}+\frac{1\,x^{16}}{16982784000}-\frac{839\,x^{14}}{7925299200}+\frac{10099\,x^{12}}{119750400}$$

$$-\frac{1\,x^{11}}{9979200}-\frac{89261\,x^{10}}{2419200}+\frac{1\,x^9}{3240}+\frac{34739\,x^8}{3360}-\frac{839\,x^7}{5040}-\frac{94379201\,x^6}{46080}+\frac{3673\,x^5}{120}$$

$$+\frac{993867129\,x^4}{4096}-\frac{4061\,x^3}{2}-\frac{11396376521\,x^2}{1024}+37068\,x+\frac{311359222543}{4096}\Bigg)\cos(x)$$

$$+\left(-\frac{1\,x^4}{26873856}+\frac{49\,x^2}{1259712}-\frac{339227}{241864704}\right)\cos(3\,x)$$

$$Y_6(x)=\frac{1\,x^{26}}{353764439584741785600000}-\frac{1\,x^{22}}{4407845991284736000}+\frac{1\,x^{19}}{2027418340147200}$$

$$-\frac{1\,x^{18}}{200074178304000}+\frac{1\,x^{16}}{35582976000}-\frac{1\,x^{15}}{100590336000}-\frac{461\,x^{14}}{14529715200}$$

$$+\frac{1987\,x^{12}}{111476736}-\frac{1\,x^{11}}{15966720}-\frac{1174213\,x^{10}}{206438400}+\frac{1\,x^9}{8640}+\frac{3021211\,x^8}{2580480}-\frac{11\,x^7}{240}$$

$$-\frac{663305989\,x^6}{3932160}+\frac{18327\,x^5}{2560}+\frac{7088916165\,x^4}{524288}-\frac{1282127\,x^3}{3072}-\frac{20920090371\,x^2}{65536}$$

$$+6202\,x+\Bigg(-\frac{1\,x^{11}}{194641920}+\frac{23\,x^9}{1769472}-\frac{8317\,x^7}{983040}+\frac{3833909\,x^5}{1966080}-\frac{7\,x^4}{4096}$$

$$-\frac{124036721\,x^3}{786432}+\frac{483\,x^2}{1024}+\frac{106039717\,x}{32768}-\frac{34965}{4096}\Bigg)\sin(2\,x)$$

$$+\left(-\frac{77\,x^4}{161243136}+\frac{3745\,x^2}{20155392}-\frac{6462295}{1451188224}\right)\sin(3\,x)-\frac{27524075747}{262144}+\Bigg($$

$$\frac{1\,x^{12}}{24524881920}-\frac{17\,x^{10}}{53084160}+\frac{31457\,x^8}{82575360}-\frac{858451\,x^6}{5898240}+\frac{1\,x^5}{30720}+\frac{31750471\,x^4}{1572864}$$

$$-\frac{119\,x^3}{3072}-\frac{116186129\,x^2}{131072}+\frac{6293\,x}{2048}+\frac{1529605917}{262144}\Bigg)\cos(2\,x)$$

$$+\left(\frac{1\,x^5}{134369280}-\frac{49\,x^3}{3779136}+\frac{339227\,x}{241864704}\right)\cos(3\,x)+\Bigg(\frac{1\,x^{19}}{2027418340147200}$$

$$-\frac{1\,x^{17}}{211718707200}+\frac{12587\,x^{15}}{1307674368000}-\frac{55141\,x^{13}}{6227020800}+\frac{1\,x^{12}}{95800320}+\frac{2711\,x^{11}}{591360}$$

$$-\frac{1\,x^{10}}{25920}-\frac{570947\,x^9}{362880}+\frac{1049\,x^8}{40320}+\frac{507949\,x^7}{1260}-\frac{2299\,x^6}{360}-\frac{1385049587\,x^5}{20480}$$

$$
+\frac{15277\,x^{4}}{24}+\frac{8067942521\,x^{3}}{1536}-23448\,x^{2}-\frac{450591811577\,x}{4096}+151795\Bigg)\cos(x)+\Bigg(
$$

$$
-\frac{1\,x^{18}}{15243746918400}+\frac{1\,x^{16}}{4151347200}-\frac{3947\,x^{14}}{12454041600}+\frac{3209\,x^{12}}{14968800}-\frac{1\,x^{11}}{1140480}
$$

$$
-\frac{641953\,x^{10}}{7257600}+\frac{1\,x^{9}}{864}+\frac{1054031\,x^{8}}{40320}-\frac{329\,x^{7}}{720}-\frac{673337\,x^{6}}{120}+\frac{8569\,x^{5}}{120}
$$

$$
+\frac{16550114803\,x^{4}}{24576}-4453\,x^{3}-\frac{30607293891\,x^{2}}{1024}+83727\,x+\frac{1607456187563}{8192}\Bigg)
$$

$$
\sin(x)
$$

$$
Y_{7}(x)=-\frac{1\,x^{27}}{6980044519498943692800000}+\frac{53\,x^{23}}{4308669456480829440000}
$$

$$
-\frac{1\,x^{20}}{31191051386880000}+\frac{1\,x^{19}}{3379030566912000}-\frac{1\,x^{17}}{529296768000}
$$

$$
+\frac{1\,x^{16}}{1394852659200}+\frac{269\,x^{15}}{108972864000}-\frac{262681\,x^{13}}{159411732480}+\frac{1\,x^{12}}{191600640}
$$

$$
+\frac{662197\,x^{11}}{1021870080}-\frac{1\,x^{10}}{86400}-\frac{64160111\,x^{9}}{371589120}+\frac{11\,x^{8}}{1920}+\frac{2543261\,x^{7}}{73728}-\frac{73301\,x^{6}}{61440}
$$

$$
-\frac{17261837111\,x^{5}}{3932160}+\frac{1281161\,x^{4}}{12288}+\frac{31053340361\,x^{3}}{131072}-\frac{12666731\,x^{2}}{4096}
$$

$$
-\frac{353860773021\,x}{131072}+13545+\Bigg(\frac{1\,x^{12}}{1946419200}-\frac{23\,x^{10}}{14745600}+\frac{4991\,x^{8}}{3932160}-\frac{1841453\,x^{6}}{4718592}
$$

$$
+\frac{7\,x^{5}}{20480}+\frac{149067209\,x^{4}}{3145728}-\frac{161\,x^{3}}{1024}-\frac{511516747\,x^{2}}{262144}+\frac{34965\,x}{4096}+\frac{6734357987}{524288}\Bigg)
$$

$$
\sin(2\,x)+\Bigg(\frac{77\,x^{5}}{806215680}-\frac{3745\,x^{3}}{60466176}+\frac{6462295\,x}{1451188224}\Bigg)\sin(3\,x)
$$

$$
+\Bigg(-\frac{1\,x^{6}}{806215680}+\frac{49\,x^{4}}{15116544}-\frac{339227\,x^{2}}{483729408}+\frac{53845547}{4353564672}\Bigg)\cos(3\,x)+\Bigg(
$$

$$
-\frac{1\,x^{20}}{31191051386880000}+\frac{1\,x^{18}}{2931489792000}-\frac{1091\,x^{16}}{1394852659200}+\frac{61\,x^{14}}{74131200}
$$

$$
-\frac{1\,x^{13}}{1037836800}-\frac{476563\,x^{12}}{958003200}+\frac{1\,x^{11}}{237600}+\frac{372341\,x^{10}}{1814400}-\frac{1259\,x^{9}}{362880}-\frac{177545\,x^{8}}{2688}
$$

$$
+\frac{263\,x^{7}}{240}+\frac{1099514983\,x^{6}}{73728}-\frac{18371\,x^{5}}{120}-\frac{1351342165\,x^{4}}{768}+9454\,x^{3}
$$

$$
+\frac{614184679939\,x^{2}}{8192}-188923\,x-\frac{3886163870971}{8192}\Bigg)\cos(x)+\Bigg(\frac{1\,x^{19}}{222793224192000}
$$

$$
-\frac{1\,x^{17}}{54286848000}+\frac{1711\,x^{15}}{62270208000}-\frac{1\,x^{13}}{46592}+\frac{1\,x^{12}}{11404800}+\frac{33457\,x^{11}}{3193344}-\frac{1\,x^{10}}{7200}
$$

$$
-\frac{1377361\,x^{9}}{362880}+\frac{79\,x^{8}}{1152}+\frac{664721\,x^{7}}{630}-\frac{5149\,x^{6}}{360}-\frac{22027454527\,x^{5}}{122880}+\frac{5363\,x^{4}}{4}
$$

$$
+\frac{13777782951\,x^{3}}{1024}-\frac{102291\,x^{2}}{2}-\frac{2216400445751\,x}{8192}+352716\Bigg)\sin(x)+\Bigg(
$$

$$
-\frac{1\,x^{13}}{265686220800}+\frac{17\,x^{11}}{486604800}-\frac{3775\,x^{9}}{74317824}+\frac{294433\,x^{7}}{11796480}-\frac{1\,x^{6}}{184320}-\frac{12712713\,x^{5}}{2621440}
$$

$$
+\frac{119\,x^{4}}{12288}+\frac{17472185\,x^{3}}{49152}-\frac{6293\,x^{2}}{4096}-\frac{1854402333\,x}{262144}+21\Bigg)\cos(2\,x)
$$

Analyzing the obtained solutions, we can easy establish a common structural functional form of the desired partial solutions:

$$Y_{1+q}(x) = \left(\sum_{k=0}^{s} \sin(kx)\left(\sum_{i=0}^{9s-1-(13-q)k} r_{i,k}^{\,q}\, x^i\right)\right) + \left(\sum_{k=0}^{s} \cos(kx)\left(\sum_{i=0}^{9s-(13-q)k} g_{i,k}^{\,q}\, x^i\right)\right)$$
$$+ \left(\sum_{i=0}^{9s-(6-q)} m_i^{\,q}\, x^i\right), q = 0 .. 6$$

where $r_{i,k}^{\,q}$ are the numerical coefficients. They are determined by substituting these partial solutions to the initial ODE (17.8). **The task is solved.**

**Example No. 5. Determine partial solutions - - of the following ODE:**

$$\left(\frac{d^{11}}{dx^{11}} \mathrm{Y}(x)\right) - \sqrt{x^3}\, \mathrm{Y}(x) = 0 \tag{11.37}$$

**which, in the interval [-10.2, 9.8] of the independent variable *x* give the following relative error:**

$$\delta(Y_i(x)) = \left| \frac{\left(\frac{d^{11}}{dx^{11}} Y_i(x)\right) - \sqrt{x^3}\, Y_i(x)}{\sqrt{x^3}\, Y_i(x)} \right| \tag{11.38}$$

**not exceeding** $10^{(-4)}$**.**

**Solution. In accordance with (14.47), let's calculate the following functions:**

Let's calculate the $\xi_k(-n), n = 1 .. 11, k = 0, 1, 2 .. N$ summands using the above **program**:

$$\xi_0(-1) := -\frac{2048\,(x^3)^{\left(\frac{1}{2}\right)} x^{11}}{2635284526875}$$

$$\xi_0(-2) := -\frac{4096}{6468425656875} \sqrt{x^3}\, x^{12}$$

$$\xi_0(-3) := -\frac{16384}{62528114683125} \sqrt{x^3}\, x^{13}$$

$$\xi_0(-4) := -\frac{32768}{447316512733125} \sqrt{x^3}\, x^{14}$$

$$\xi_0(-5) := -\frac{32768}{2108777845741875} \sqrt{x^3}\, x^{15}$$

$$\xi_0(-6) := -\frac{65536}{24602408200321875} \sqrt{x^3}\, x^{16}$$

$$\xi_0(-7) := -\frac{1048576}{2730867310235728125} \sqrt{x^3}\, x^{17}$$

$$\xi_0(-8) := -\frac{2097152}{4385451621731492812 5}\sqrt{x^3}\, x^{18}$$

$$\xi_0(-9) := -\frac{1048576}{1997816849899902281 25}\sqrt{x^3}\, x^{19}$$

$$\xi_0(-10) := -\frac{2097152}{40692374784803272781 25}\sqrt{x^3}\, x^{20}$$

$$\xi_0(-11) := -\frac{8388608}{18311568653161472751 5625}\sqrt{x^3}\, x^{21}$$

$$\xi_1(-1) := -\frac{x^{25}}{2289469627775865691406 25}$$

$$\xi_1(-2) := -\frac{x^{26}}{2563329631098337664062 50}$$

$$\xi_1(-3) := -\frac{3338\, x^{27}}{1906732746092498471412890625}$$

$$\xi_1(-4) := -\frac{124379\, x^{28}}{23643486051546981045519843750 0}$$

$$\xi_1(-5) := -\frac{2565359\, x^{29}}{215493487155528198672023718750 00}$$

$$\xi_1(-6) := -\frac{4055003\, x^{30}}{187441080543565947957322406250000}$$

$$\xi_1(-7) := -\frac{4594343419\, x^{31}}{13974669759925559249958171997968750 00}$$

$$\xi_1(-8) := -\frac{2714906861\, x^{32}}{63189811088359050521549995121250000 00}$$

$$\xi_1(-9) := -\frac{185459608753\, x^{33}}{376345211724142639843084076206350000000 0}$$

$$\xi_1(-10) := -\frac{30984049607633\, x^{34}}{614948075957249073503599380521175900000000 0}$$

$$\xi_1(-11) := -\frac{2147483648\, x^{35}}{461950177251539267956429823107854492187 5}$$

$$\xi_2(-11) := -\frac{150875956300234034777784}{102364108729952345800723452758801377763200732255751373539752807}$$

$$\sqrt{x^3}\, x^{46}$$

$$\xi_2(-10) := -\frac{118274989849305936579 32}{754261853799648863794804388749062783518321185042378541871862792}$$

$$\sqrt{x^3}\, x^{45}$$

$$\xi_2(-9) := -$$

$$\frac{1239477992094923984338 4}{824237134559756322776517924189398440603916873877267873963500976562 5}$$

$$\sqrt{x^3}\, x^{44}$$

$$\xi_2(-8) := -$$

$$\frac{283652225813002982 72}{220915876322636377050795476866630512088961906694523686401367187 5}$$

$$\sqrt{x^3}\, x^{43}$$

$$\xi_2(-7) := -$$

$$\frac{373584153557782016}{388778447569280208699112115779910201320766858945283325195312 5}\sqrt{x^3}\, x^{42}$$

$\xi_2(-6) :=$

$$-\frac{1996872725777792}{3233684188462681991204071151102080059888730149862060546875}\sqrt{x^3}\,x^{41}$$

$$\xi_2(-5) := -\frac{23547989251328}{7106998216401498881767189343081494637118088241455078125}\sqrt{x^3}\,x^{40}$$

$$\xi_2(-4) := -\frac{547018961408}{38384286170843159727035089929397353669520210693359375}\sqrt{x^3}\,x^{39}$$

$$\xi_2(-3) := -\frac{14789208064}{321015535478975051480580202515867951835891845703125}\sqrt{x^3}\,x^{38}$$

$$\xi_2(-2) := -\frac{91136}{915818749351611881333265441970620906005859375}\sqrt{x^3}\,x^{37}$$

$$\xi_2(-1) := -\frac{2048}{19009668789960197781032371811355989501953125}\sqrt{x^3}\,x^{36}$$

**In this case:**

$$\alpha(-1) := -\frac{2048\sqrt{x^3}\,x^{36}}{19009668789960197781032371811355989501953125} - \frac{1\,x^{25}}{228946962777586569140625} - \frac{2048\sqrt{x^3}\,x^{11}}{2635284526875}$$

$$\alpha(-2) := -\frac{91136}{915818749351611881333265441970620906005859375}\sqrt{x^3}\,x^{37} - \frac{1}{256332963109833766406250}x^{26} - \frac{4096}{6468425656875}\sqrt{x^3}\,x^{12}$$

$$\alpha(-3) := -\frac{14789208064}{321015535478975051480580202515867951835891845703125}\sqrt{x^3}\,x^{38} - \frac{3338}{1906732746092498471412890625}x^{27} - \frac{16384}{62528114683125}\sqrt{x^3}\,x^{13}$$

$$\alpha(-4) := -\frac{547018961408}{38384286170843159727035089929397353669520210693359375}\sqrt{x^3}\,x^{39} - \frac{124379}{236434860515469810455198437500}x^{28} - \frac{32768}{447316512733125}\sqrt{x^3}\,x^{14}$$

$$\alpha(-5) := -\frac{23547989251328}{7106998216401498881767189343081494637118088241455078125}\sqrt{x^3}\,x^{40} - \frac{2565359}{21549348715552819867202371875000}x^{29} - \frac{32768}{2108777845741875}\sqrt{x^3}\,x^{15}$$

$\alpha(-6) :=$

$$-\frac{1996872725777792}{3233684188462681991204071151102080059888730149862060546875}\sqrt{x^3}\,x^{41} - \frac{4055003}{187441080543565947957322406250000}x^{30} - \frac{65536}{24602408200321875}\sqrt{x^3}\,x^{16}$$

$\alpha(-7) := -$

$$\frac{373584153557782016}{3887784475692802086991121157799102013207668589452833251953125}\sqrt{x^3}\,x^{42} - \frac{4594343419}{1397466975992555924995817199796875000}x^{31} - \frac{1048576}{2730867310235728125}\sqrt{x^3}\,x^{17}$$

$$\alpha(-8) := -\frac{2836522258130029827 2}{2209158763226363770507954768666305120889619066945236864013671875}\sqrt{x^3}\,x^{43} - \frac{2714906861}{631898110883590505215499951212 5000000}x^{32} - \frac{2097152}{4385451621731492812 5}\sqrt{x^3}\,x^{18}$$

$$\alpha(-9) := -\frac{12394779920949239843384}{824237134559756322776517924189398440603916873877267873963500976 5625}\sqrt{x^3}\,x^{44} - \frac{185459608753}{3763452117241426398430840762063500000000}x^{33} - \frac{1048576}{1997816849899902281 25}\sqrt{x^3}\,x^{19}$$

$$\alpha(-10) := -\frac{11827498984930593657932}{754261853799648863794804388749062783518321185042378541871862 7929}\sqrt{x^3}\,x^{45} - \frac{30984049607633}{61494807595724907350359938052117590000000000}x^{34} - \frac{2097152}{40692374784803272781 25}\sqrt{x^3}\,x^{20}$$

$$\alpha(-11) := -\frac{15087595630023403477784}{1023641087299523458007234527588013777632007322557513735397528076}\sqrt{x^3}\,x^{46} - \frac{2147483648}{4619501772515392679564298231078544921875}x^{35} - \frac{8388608}{1831156865316147275156 25}\sqrt{x^3}\,x^{21}$$

**Then the required solutions are defined by the following formulas:** (17.12)

$$Y_1(x) = -\frac{2048\sqrt{x^3}\,x^{36}}{1900966878996019778103237181135598950195312 5} - \frac{1\,x^{25}}{2289469627775865691406 25} - \frac{2048\sqrt{x^3}\,x^{11}}{2635284526875} - 1$$

$$Y_2(x) = \frac{1024\,(x^3)^{\left(\frac{1}{2}\right)}\,x^{37}}{1245513499118192158613241001080044432167968 75} + \frac{1\,x^{26}}{2142943571598210287156250} + \frac{2048\,(x^3)^{\left(\frac{1}{2}\right)}\,x^{12}}{14230536445125} + x$$

$$Y_3(x) = -\frac{8192\,(x^3)^{\left(\frac{1}{2}\right)}\,x^{38}}{1930920514159248429958377157990183168937695312 5} - \frac{4\,x^{27}}{119851772611528475345953125} - \frac{1024\,(x^3)^{\left(\frac{1}{2}\right)}\,x^{13}}{58955079558375} - \frac{1\,x^2}{2}$$

$$Y_4(x) = \frac{2048\,(x^3)^{\left(\frac{1}{2}\right)}\,x^{39}}{112794250174068161028176787324670819409830078125}$$

$$+ \frac{1\, x^{28}}{5099575422882682186288 59375} + \frac{1024\, (x^3)^{\left(\frac{1}{2}\right)} x^{14}}{609202488769875} + \frac{1\, x^3}{6}$$

$$Y_5(x) = - \frac{1024\, (x^3)^{\left(\frac{1}{2}\right)} x^{40}}{148358338890154133111924331743433148808448 6328125}$$

$$- \frac{1\, x^{29}}{985917915090651889349 1281250} - \frac{256\, (x^3)^{\left(\frac{1}{2}\right)} x^{15}}{1827607466309625} - \frac{1\, x^4}{24}$$

$$Y_6(x) = \frac{512\, (x^3)^{\left(\frac{1}{2}\right)} x^{41}}{21272703788230602757276599524615077005260009765625}$$

$$+ \frac{1\, x^{30}}{2095574515880940250640859 37500} + \frac{256\, (x^3)^{\left(\frac{1}{2}\right)} x^{16}}{24602408200321875} + \frac{1\, x^5}{120}$$

$$Y_7(x) = - \frac{1024\, (x^3)^{\left(\frac{1}{2}\right)} x^{42}}{1306327285122303765258999747668808565709471455078125}$$

$$- \frac{1\, x^{31}}{480724793943087693497013 1406250} - \frac{128\, (x^3)^{\left(\frac{1}{2}\right)} x^{17}}{182057820682381875} - \frac{1\, x^6}{720}$$

$$Y_8(x) = \frac{32\, (x^3)^{\left(\frac{1}{2}\right)} x^{43}}{1326971614474543859881178848158745838916469248046875}$$

$$+ \frac{1\, x^{32}}{117636184870779106173386745000000} + \frac{128\, (x^3)^{\left(\frac{1}{2}\right)} x^{18}}{2923634414487661875} + \frac{1\, x^7}{5040}$$

$$Y_9(x) = - \frac{8\, (x^3)^{\left(\frac{1}{2}\right)} x^{44}}{1132936177026206668470635762581299937070331067382 8125}$$

$$- \frac{1\, x^{33}}{30461580503380694861740146600000 00} - \frac{16\, (x^3)^{\left(\frac{1}{2}\right)} x^{19}}{63088953154733756 25} - \frac{1\, x^8}{40320}$$

$$Y_{10}(x) = \frac{4\, (x^3)^{\left(\frac{1}{2}\right)} x^{45}}{2021355965080201415693444957945005168975216086621 09375}$$

$$+ \frac{1\, x^{34}}{8298415657629051407676539316000000 0} + \frac{16\, (x^3)^{\left(\frac{1}{2}\right)} x^{20}}{116263927956580779375}$$

$$+ \frac{1\, x^9}{362880}$$

$$Y_{11}(x) = -\frac{16\,(x^3)^{\left(\frac{1}{2}\right)}\,x^{46}}{300223730668180600380317353417233867786298101013183593 75} - \frac{1\,x^{35}}{236775446753002825490770822875000000 0} - \frac{8\,(x^3)^{\left(\frac{1}{2}\right)}\,x^{21}}{113736451261872501562 5} - \frac{1\,x^{10}}{3628800}$$

**Direct calculation of the relative error for the obtained partial solutions (17.12) demonstrates that all of them satisfy the established condition (11.38). To save space, let's present the calculated values of relative error for the first partial solution $Y_1(x)$:**
**The first partial solution**

$\delta(-10.2) = 0.2273782166\ 10^{-5}, x = -10.2$
$\delta(-9.2) = 0.1724022946\ 10^{-6}, x = -9.2$
$\delta(-8.2) = 0.9709498590\ 10^{-8}, x = -8.2$
$\delta(-7.2) = 0.3758685307\ 10^{-9}, x = -7.2$
$\delta(-6.2) = 0.8833793128\ 10^{-11}, x = -6.2$
$\delta(-5.2) = 0.6271561989\ 10^{-13}, x = -5.2$
$\delta(-4.2) = 0.2533403403\ 10^{-16}, x = -4.2$
$\delta(-3.2) = 0.9451027597\ 10^{-21}, x = -3.2$
$\delta(-2.2) = 0.7468027920\ 10^{-27}, x = -2.2$
$\delta(-1.2) = 0.1003811413\ 10^{-36}, x = -1.2$
$\delta(-0.2) = 0.6621863738\ 10^{-66}, x = -0.2$
$\delta(0.8) = 0.2501669389\ 10^{-43}, x = 0.8$
$\delta(1.8) = 0.4027854924\ 10^{-30}, x = 1.8$
$\delta(2.8) = 0.6319030462\ 10^{-23}, x = 2.8$
$\delta(3.8) = 0.5865651511\ 10^{-18}, x = 3.8$
$\delta(4.8) = 0.3022483817\ 10^{-14}, x = 4.8$
$\delta(5.8) = 0.1233801170\ 10^{-11}, x = 5.8$
$\delta(6.8) = 0.8572904007\ 10^{-10}, x = 6.8$
$\delta(7.8) = 0.2753898252\ 10^{-8}, x = 7.8$
$\delta(8.8) = 0.5642758111\ 10^{-7}, x = 8.8$
$\delta(9.8) = 0.8247551309\ 10^{-6}, x = 9.8$

**The task is solved.**
**Example No. 6. Define partial solutions -- for the following ODE:**

$$\left(\frac{d^{13}}{dx^{13}}\mathrm{Y}(x)\right) - \mathbf{e}^{(4x)}\,\mathrm{Y}(x) = 0 \qquad (11.39)$$

**which, in the interval [ $-\infty$ , 5] of the independent variable *x* give the following relative error:**

$$\delta(Y_i(x)) = \left|\frac{\left(\frac{d^{13}}{dx^{13}}Y_i(x)\right) - \mathbf{e}^{(4x)}\,Y_i(x)}{\mathbf{e}^{(4x)}\,Y_i(x)}\right| \qquad (11.40)$$

**not exceeding** $10^{(-4)}$ **.**

**Solution. In accordance with (14.47), let's calculate the following function:**

Let's calculate the $\xi_k(-n), n = 1 .. 13, k = 0, 1, 2 .. N$ summands using the above **program**:

$$\xi_0(-1) := -\frac{1}{67108864}\mathbf{e}^{(4x)}$$
$$\xi_0(-2) := -\frac{13}{268435456}\mathbf{e}^{(4x)}$$
$$\xi_0(-3) := -\frac{91}{1073741824}\mathbf{e}^{(4x)}$$
$$\xi_0(-4) := -\frac{455}{4294967296}\mathbf{e}^{(4x)}$$
$$\xi_0(-5) := -\frac{455}{4294967296}\mathbf{e}^{(4x)}$$
$$\xi_0(-6) := -\frac{1547}{17179869184}\mathbf{e}^{(4x)}$$
$$\xi_0(-7) := -\frac{4641}{68719476736}\mathbf{e}^{(4x)}$$
$$\xi_0(-8) := -\frac{12597}{274877906944}\mathbf{e}^{(4x)}$$
$$\xi_0(-9) := -\frac{62985}{2199023255552}\mathbf{e}^{(4x)}$$
$$\xi_0(-10) := -\frac{146965}{8796093022208}\mathbf{e}^{(4x)}$$
$$\xi_0(-11) := -\frac{323323}{35184372088832}\mathbf{e}^{(4x)}$$
$$\xi_0(-12) := -\frac{676039}{140737488355328}\mathbf{e}^{(4x)}$$
$$\xi_0(-13) := -\frac{676039}{281474976710656}\mathbf{e}^{(4x)}$$
$$\xi_1(-1) := -\frac{1}{36893488147419103232}\mathbf{e}^{(8x)}$$
$$\xi_1(-2) := -\frac{39}{295147905179352825856}\mathbf{e}^{(8x)}$$
$$\xi_1(-3) := -\frac{793}{2361183241434822606848}\mathbf{e}^{(8x)}$$
$$\xi_1(-4) := -\frac{11193}{18889465931478580854784}\mathbf{e}^{(8x)}$$
$$\xi_1(-5) := -\frac{61607}{75557863725914323419136}\mathbf{e}^{(8x)}$$
$$\xi_1(-6) := -\frac{35217}{37778931862957161709568}\mathbf{e}^{(8x)}$$
$$\xi_1(-7) := -\frac{2227589}{2417851639229258349412352}\mathbf{e}^{(8x)}$$

$$\xi_1(-8) := -\frac{15643563}{193428131138340667952988160}\, \mathbf{e}^{(8\,x)}$$

$$\xi_1(-9) := -\frac{199012021}{3094850098213450687247810560}\, \mathbf{e}^{(8\,x)}$$

$$\xi_1(-10) := -\frac{1163645223}{24758800785707605497982484480}\, \mathbf{e}^{(8\,x)}$$

$$\xi_1(-11) := -\frac{6326776729}{198070406285660843983859875840}\, \mathbf{e}^{(8\,x)}$$

$$\xi_1(-12) := -\frac{32281129881}{1584563250285286751870879006720}\, \mathbf{e}^{(8\,x)}$$

$$\xi_1(-13) := -\frac{4641}{37778931862957161709568}\, \mathbf{e}^{(8\,x)}$$

In this case, the desired partial solutions take the following form:

$$Y_1 := \mathrm{Y}(x) = -1 - \frac{1\,\mathbf{e}^{(4\,x)}}{67108864} - \frac{1\,\mathbf{e}^{(8\,x)}}{36893488147419103232}$$

$$Y_2(x) = x + \left(\frac{1\,x}{67108864} - \frac{13}{268435456}\right)\mathbf{e}^{(4\,x)} + \left(\frac{1\,x}{36893488147419103232} - \frac{39}{295147905179352825856}\right)\mathbf{e}^{(8\,x)}$$

$$Y_3(x) = -\frac{1\,x^2}{2} + \left(-\frac{1\,x^2}{134217728} + \frac{13\,x}{268435456} - \frac{91}{1073741824}\right)\mathbf{e}^{(4\,x)} + \left(-\frac{1\,x^2}{73786976294838206464} + \frac{39\,x}{295147905179352825856} - \frac{793}{2361183241434822606848}\right)\mathbf{e}^{?}$$

$$Y_4(x) = \frac{1\,x^3}{6} + \left(\frac{1\,x^3}{402653184} - \frac{13\,x^2}{536870912} + \frac{91\,x}{1073741824} - \frac{455}{4294967296}\right)\mathbf{e}^{(4\,x)} + \left(\frac{1\,x^3}{221360928884514619392} - \frac{39\,x^2}{590295810358705651712} + \frac{793\,x}{2361183241434822606848} - \frac{11193}{18889465931478580854784}\right)\mathbf{e}^{(8\,x)}$$

$$Y_5(x) = -\frac{1\,x^4}{24} + \left(-\frac{1\,x^4}{1610612736} + \frac{13\,x^3}{1610612736} - \frac{91\,x^2}{2147483648} + \frac{455\,x}{4294967296} - \frac{455}{4294967296}\right)\mathbf{e}^{(4\,x)} + \left(-\frac{1\,x^4}{885443715538058477568} + \frac{13\,x^3}{590295810358705651712} - \frac{793\,x^2}{4722366482869645213696} + \frac{11193\,x}{18889465931478580854784} - \frac{61607}{75557863725914323419136}\right)\mathbf{e}^{(8\,x)}$$

$$Y_6(x) = \frac{1\,x^5}{120} + \left(\frac{1\,x^5}{8053063680} - \frac{13\,x^4}{6442450944} + \frac{91\,x^3}{6442450944} - \frac{455\,x^2}{8589934592} + \frac{455\,x}{4294967296} - \frac{1547}{17179869184}\right) \mathbf{e}^{(4\,x)} + \left(\frac{1\,x^5}{4427218577690292387840} - \frac{13\,x^4}{2361183241434822606848} + \frac{793\,x^3}{14167099448608935641088} - \frac{11193\,x^2}{37778931862957161709568} + \frac{61607\,x}{75557863725914323419136} - \frac{35217}{37778931862957161709568}\right) \mathbf{e}^{(8\,x)}$$

$$Y_7(x) = -\frac{1\,x^6}{720} + \left(-\frac{1\,x^6}{48318382080} + \frac{13\,x^5}{32212254720} - \frac{91\,x^4}{25769803776} + \frac{455\,x^3}{25769803776} - \frac{455\,x^2}{8589934592} + \frac{1547\,x}{17179869184} - \frac{4641}{68719476736}\right) \mathbf{e}^{(4\,x)} + \left(-\frac{1\,x^6}{26563311466141754327040} + \frac{13\,x^5}{11805916207174113034240} - \frac{793\,x^4}{56668397794435742564352} + \frac{3731\,x^3}{37778931862957161709568} - \frac{61607\,x^2}{151115727451828646838272} + \frac{35217\,x}{37778931862957161709568} - \frac{2227589}{2417851639229258349412352}\right) \mathbf{e}^{(8\,x)}$$

$$Y_8(x) = \frac{1\,x^7}{5040} + \left(\frac{1\,x^7}{338228674560} - \frac{13\,x^6}{193273528320} + \frac{91\,x^5}{128849018880} - \frac{455\,x^4}{103079215104} + \frac{455\,x^3}{25769803776} - \frac{1547\,x^2}{34359738368} + \frac{4641\,x}{68719476736} - \frac{12597}{274877906944}\right) \mathbf{e}^{(4\,x)} + \left(\frac{1\,x^7}{185943180262992280289280} - \frac{13\,x^6}{70835497243044678205440} + \frac{793\,x^5}{283341988972178712821760} - \frac{3731\,x^4}{151115727451828646838272} + \frac{61607\,x^3}{453347182355485940514816} - \frac{35217\,x^2}{75557863725914323419136} + \frac{2227589\,x}{2417851639229258349412352} - \frac{15643563}{19342813113834066795298816}\right) \mathbf{e}^{(8\,x)}$$

$$Y_9(x) = -\frac{1\,x^8}{40320} + \left(-\frac{1\,x^8}{2705829396480} + \frac{13\,x^7}{1352914698240} - \frac{91\,x^6}{773094113280} + \frac{91\,x^5}{103079215104} - \frac{455\,x^4}{103079215104} + \frac{1547\,x^3}{103079215104} - \frac{4641\,x^2}{137438953472} + \frac{12597\,x}{274877906944} - \frac{62985}{2199023255552}\right) \mathbf{e}^{(4\,x)} + \left(-\frac{1\,x^8}{1487545442103938242314240} + \frac{13\,x^7}{495848480701312747438080} - \frac{793\,x^6}{1700051933833072276930560} + \frac{3731\,x^5}{755578637259143234191360} - \frac{61607\,x^4}{1813388729421943762059264} + \frac{11739\,x^3}{75557863725914323419136} - \frac{2227589\,x^2}{4835703278458516698824704}\right.$$

$$+ \frac{15643563\,x}{19342813113834066795298816} - \frac{199012021}{309485009821345068724781056}\Bigg) e^{(8\,x)}$$

$$Y_{10}(x) = \frac{1\,x^9}{362880} + \Bigg( \frac{1\,x^9}{24352464568320} - \frac{13\,x^8}{10823317585920} + \frac{13\,x^7}{773094113280}$$

$$- \frac{91\,x^6}{618475290624} + \frac{91\,x^5}{103079215104} - \frac{1547\,x^4}{412316860416} + \frac{1547\,x^3}{137438953472}$$

$$- \frac{12597\,x^2}{549755813888} + \frac{62985\,x}{2199023255552} - \frac{146965}{8796093022208}\Bigg) e^{(4\,x)} + \Bigg($$

$$\frac{1\,x^9}{13387908978935444180828160} - \frac{13\,x^8}{3966787845610501979504640}$$

$$+ \frac{793\,x^7}{11900363536831505938513920} - \frac{3731\,x^6}{4533471823554859405148160}$$

$$+ \frac{61607\,x^5}{9066943647109718810296320} - \frac{11739\,x^4}{302231454903657293676544}$$

$$+ \frac{2227589\,x^3}{14507109835375550096474112} - \frac{15643563\,x^2}{38685626227668133590597632}$$

$$+ \frac{199012021\,x}{309485009821345068724781056} - \frac{1163645223}{2475880078570760549798248448}\Bigg) e^{(8\,x)}$$

$$Y_{11}(x) = -\frac{1\,x^{10}}{3628800} + \Bigg( -\frac{1\,x^{10}}{243524645683200} + \frac{13\,x^9}{97409858273280} - \frac{13\,x^8}{6184752906240}$$

$$+ \frac{13\,x^7}{618475290624} - \frac{91\,x^6}{618475290624} + \frac{1547\,x^5}{2061584302080} - \frac{1547\,x^4}{549755813888}$$

$$+ \frac{4199\,x^3}{549755813888} - \frac{62985\,x^2}{4398046511104} + \frac{146965\,x}{8796093022208} - \frac{323323}{35184372088832}\Bigg) e^{(4\,x)}$$

$$+ \Bigg( -\frac{1\,x^{10}}{133879089789354441808281600} + \frac{13\,x^9}{35701090610494517815541760}$$

$$- \frac{793\,x^8}{95202908294652047508111360} + \frac{533\,x^7}{4533471823554859405148160}$$

$$- \frac{61607\,x^6}{54401661882658312861777920} + \frac{11739\,x^5}{1511157274518286468382720}$$

$$- \frac{2227589\,x^4}{58028439341502200385896448} + \frac{5214521\,x^3}{38685626227668133590597632}$$

$$- \frac{199012021\,x^2}{618970019642690137449562112} + \frac{1163645223\,x}{2475880078570760549798248448}$$

$$- \frac{6326776729}{19807040628566084398385987584}\Bigg) e^{(8\,x)}$$

$$Y_{12}(x) = \frac{1\,x^{11}}{39916800} + \Bigg( \frac{1\,x^{11}}{2678771102515200} - \frac{13\,x^{10}}{974098582732800} + \frac{13\,x^9}{55662776156160}$$

$$- \frac{13\,x^8}{4947802324992} + \frac{13\,x^7}{618475290624} - \frac{1547\,x^6}{12369505812480} + \frac{1547\,x^5}{2748779069440}$$

$$- \frac{4199\,x^4}{2199023255552} + \frac{20995\,x^3}{4398046511104} - \frac{146965\,x^2}{17592186044416} + \frac{323323\,x}{35184372088832}$$

$$- \frac{676039}{140737488355328}\Bigg) e^{(4\,x)} + \Bigg( \frac{1\,x^{11}}{1472669987682898859891097600}$$

$$
-\frac{13\,x^{10}}{357010906104945178155417600}+\frac{793\,x^{9}}{856826174651868427573002240}
-\frac{533\,x^{8}}{36267774588438875241185280}+\frac{8801\,x^{7}}{54401661882658312861777920}
-\frac{3913\,x^{6}}{3022314549036572936765440}+\frac{2227589\,x^{5}}{290142196707511001929482240}
-\frac{5214521\,x^{4}}{154742504910672534362390528}+\frac{199012021\,x^{3}}{1856910058928070412348686336}
-\frac{1163645223\,x^{2}}{4951760157141521099596496896}+\frac{6326776729\,x}{19807040628566084398385987584}
-\frac{32281129881}{158456325028528675187087900672}\Bigg)\mathrm{e}^{(8\,x)}
$$

$$
Y_{13}(x)=-\frac{1\,x^{12}}{479001600}+\Bigg(-\frac{1\,x^{12}}{32145253230182400}+\frac{13\,x^{11}}{10715084410060800}
-\frac{13\,x^{10}}{556627761561600}+\frac{13\,x^{9}}{44530220924928}-\frac{13\,x^{8}}{4947802324992}+\frac{221\,x^{7}}{12369505812480}
-\frac{1547\,x^{6}}{16492674416640}+\frac{4199\,x^{5}}{10995116277760}-\frac{20995\,x^{4}}{17592186044416}+\frac{146965\,x^{3}}{52776558133248}
-\frac{323323\,x^{2}}{70368744177664}+\frac{676039\,x}{140737488355328}-\frac{676039}{281474976710656}\Bigg)\mathrm{e}^{(4\,x)}+\Bigg(
-\frac{1\,x^{12}}{17672039852194786318693171200}+\frac{13\,x^{11}}{3927119967154396959709593600}
-\frac{793\,x^{10}}{8568261746518684275730022400}+\frac{533\,x^{9}}{326409971295949877170667520}
-\frac{8801\,x^{8}}{435213295061266502894223360}+\frac{559\,x^{7}}{3022314549036572936765440}
-\frac{2227589\,x^{6}}{1740853180245066011576893440}+\frac{5214521\,x^{5}}{773712524553362671811952640}
-\frac{199012021\,x^{4}}{7427640235712281649394745344}+\frac{387881741\,x^{3}}{4951760157141521099596496896}
-\frac{6326776729\,x^{2}}{39614081257132168796771975168}+\frac{32281129881\,x}{158456325028528675187087900672}
-\frac{4641}{37778931862957161709568}\Bigg)\mathrm{e}^{(8\,x)}
$$

Direct calculation of the relative error for the obtained partial solutions demonstrates that all of them satisfy the established condition (11.40). To save space, let's present the calculated values of relative error for the first partial solution $Y_1(x)$:

$\delta(-10.2)=0.9876897758\ 10^{-55},\ x=-10.2$

$\delta(-9.2)=0.2944261721\ 10^{-51},\ x=-9.2$

$\delta(-8.2)=0.8776720523\ 10^{-48},\ x=-8.2$

$\delta(-7.2)=0.2616303513\ 10^{-44},\ x=-7.2$

$\delta(-6.2)=0.7799090840\ 10^{-41},\ x=-6.2$

$\delta(-5.2)=0.2324876217\ 10^{-37},\ x=-5.2$

$\delta(-4.2)=0.6930358324\ 10^{-34},\ x=-4.2$

$\delta(-3.2)=0.2065910701\ 10^{-30},\ x=-3.2$

$\delta(-2.2)=0.6158393020\ 10^{-27},\ x=-2.2$

$\delta(-1.2)=0.1835791080\ 10^{-23},\ x=-1.2$

$\delta(-0.2)=0.5472416047\ 10^{-20},\ x=-0.2$

$\delta(0.8) = 0.1631303644\ 10^{-16}, x = 0.8$
$\delta(1.8) = 0.4862752358\ 10^{-13}, x = 1.8$
$\delta(2.8) = 0.1448017033\ 10^{-9}, x = 2.8$
$\delta(3.8) = 0.4078519634\ 10^{-6}, x = 3.8$
$\delta(4.8) = 0.0003031082538, x = 4.8$
$\delta(5.8) = 0.02107510366, x = 5.8$

**The task is solved.**

**Conclusions: The question is: "How can we optimally use the transformation of variables to facilitate the solution of ODEs?" The answer to this question is clear: transform the variables to remove the derivatives of order. Extensive experience with linear ordinary differential equations demonstrates that the absence of this derivative considerably speeds up the definition of the desired solution with the specified accuracy. For example: let's take the following linear ODE:**

$$\left(\frac{d^2}{dx^2} \mathrm{Y}(x)\right) - x\left(\frac{d}{dx} \mathrm{Y}(x)\right) - x^2\,\mathrm{Y}(x) = 0$$

**In the [-10.2, 10.2] interval of the independent variable *x*, using only the first five summands -** $\xi_i(-k), k = 1, 2, i = 0, 1, 2, 3, 4$ **, we shall calculate partial solutions - - and their relative errors -** $\delta(x_i)$ **at every point of the specified interval.**

$$Y_1(x) = -\frac{1\,x^{20}}{8089804800} - \frac{5267\,x^{18}}{555761606400} - \frac{54287\,x^{16}}{186810624000} - \frac{199\,x^{14}}{43243200} - \frac{571\,x^{12}}{11975040} - \frac{41\,x^{10}}{113400} - \frac{3\,x^8}{1120} - \frac{1\,x^6}{90} - \frac{1\,x^4}{12} - 1$$

$\delta(-1.2) = 0.0008710744155, x = -1.2$
$\delta(-0.2) = 0.1002411666\ 10^{-10}, x = -0.2$
$\delta(0.8) = 0.00001304234665, x = 0.8$
$\delta(1.8) = 0.04092510073, x = 1.8$

$$Y_2(x) = \frac{1\,x^{21}}{25662873600} + \frac{182177\,x^{19}}{49691625984000} + \frac{41257\,x^{17}}{296406190080} + \frac{9979\,x^{15}}{3592512000} + \frac{9781\,x^{13}}{283046400} + \frac{41\,x^{11}}{134400} + \frac{119\,x^9}{51840} + \frac{13\,x^7}{1008} + \frac{3\,x^5}{40} + \frac{1\,x^3}{6} + x$$

$\delta(-1.2) = 0.001779916264, x = -1.2$
$\delta(-0.2) = 0.6797179871\ 10^{-9}, x = -0.2$
$\delta(0.8) = 0.00005683642779, x = 0.8$
$\delta(1.8) = 0.04755209583, x = 1.8$

**Upon the substitution:**

**we can reduce this ODE to the following form:**

$$\frac{d^2}{dx^2} Z(x) = \frac{1\, Z(x)\,(-2 + 5x^2)}{4}$$

**Using the same program, let's calculate the required partial solutions - - again with the use of only the first five summands:**

$$Z_1(x) = -\frac{125\,x^{20}}{331358404608} + \frac{1007525\,x^{18}}{182111963185152} - \frac{369071\,x^{16}}{2678117105664} + \frac{150511\,x^{14}}{111588212736} - \frac{1769\,x^{12}}{68124672} + \frac{21601\,x^{10}}{116121600} - \frac{1721\,x^8}{645120} + \frac{71\,x^6}{5760} - \frac{11\,x^4}{96} + \frac{1\,x^2}{4} - 1$$

$\delta(-3.2) = 0.0084184799, x = -3.2$
$\delta(-2.2) = 0.000244861968, x = -2.2$
$\delta(-1.2) = 0.9886333\ 10^{-10}, x = -1.2$
$\delta(-0.2) = 0.7689107508\ 10^{-15}, x = -0.2$
$\delta(0.8) = 0.1511631525\ 10^{-10}, x = 0.8$
$\delta(1.8) = 0.578832713\ 10^{-5}, x = 1.8$
$\delta(2.8) = 0.00245727593, x = 2.8$

$$Z_2(x) = \frac{125\,x^{21}}{1051151302656} - \frac{3593285\,x^{19}}{3460127300517888} + \frac{655873\,x^{17}}{15175996932096} - \frac{456691\,x^{15}}{1673823191040} + \frac{24709\,x^{13}}{2656862208} - \frac{4891\,x^{11}}{116121600} + \frac{6641\,x^9}{5806080} - \frac{131\,x^7}{40320} + \frac{31\,x^5}{480} - \frac{1\,x^3}{12} + x$$

$\delta(-3.2) = 0.02330519447, x = -3.2$
$\delta(-2.2) = 0.0000225601438, x = -2.2$
$\delta(-1.2) = 0.32505420\ 10^{-10}, x = -1.2$
$\delta(-0.2) = 0.6828845697\ 10^{-16}, x = -0.2$
$\delta(0.8) = 0.619307551\ 10^{-12}, x = 0.8$
$\delta(1.8) = 0.24388445\ 10^{-6}, x = 1.8$
$\delta(2.8) = 0.00331525789, x = 2.8$

**As we can see, in the first case the relative error is less than 0.1 only in the [-1.2, 1.8] interval, while in the latter case it is in the [-3.2, 2.8] interval, and the transformed equation gives a lower relative error for the desired solutions, that confirms the conclusion. Therefore, it provides a higher accuracy of the desired solution with fewer calculations.**

## 12. Non-homogeneous ODEs.

Conventional methods to obtain partial solutions of a non-homogeneous ODE:

$$\left(\frac{d^m}{dx^m} Y(x)\right) + \left(\sum_{p=1}^{m} b_p(x)\left(\frac{d^{m-p}}{dx^{m-p}} Y(x)\right)\right) = F(x) \qquad (12.1)$$

are not simple or perfect because, in order to find a partial solution of a non-homogenous ODE, in a general case we need to know the Wronskian, i.e. all partial solutions of a homogeneous ODE. Partial solution can also be defined through special functions, such as Green function, or by other methods [V.F. Zaitsev, A.D. Polyanin. Reference Book. Ordinary Differential Equations. Publishing House “Physics & Mathematical Literature”, Moscow, 2001. 576 pages].
Essentially, to write down the formula for a partial solution of an ODE of *m* order it is enough to know $m-1$ partial solutions; so, we do not need to calculate (albeit with the use of a well-known formula) another solution.
Above we have proved that a non-homogeneous second-order ODE (5.28) has a partial solution defined by (5.29), i.e. it defines the desired partial solution based on only one partial solution of a homogeneous ODE in symbolic form.
Let’s present a general approach to the obtainment of formulas for the partial solution of a linear non-homogeneous ODE, taking a third-order ODE as an example. For this purpose we shall prove the following **theorem**:

**Assume that are independent partial solutions of the following homogeneous third-order ODE:**

$$\left(\frac{d^3}{dx^3}\mathrm{Y}(x)\right)+b_1(x)\left(\frac{d^2}{dx^2}\mathrm{Y}(x)\right)+b_2(x)\left(\frac{d}{dx}\mathrm{Y}(x)\right)+b_3(x)\,\mathrm{Y}(x)=0 \tag{12.2}$$

**that is, the following equalities hold true:**

$$\left(\frac{d^3}{dx^3}y_i(x)\right)+b_1(x)\left(\frac{d^2}{dx^2}y_i(x)\right)+b_2(x)\left(\frac{d}{dx}y_i(x)\right)+b_3(x)\,y_i(x)=0 \qquad i=1,2$$

**then the partial solution of the non-homogenous ODE**

$$\left(\frac{d^3}{dx^3}\mathrm{Y}(x)\right)+b_1(x)\left(\frac{d^2}{dx^2}\mathrm{Y}(x)\right)+b_2(x)\left(\frac{d}{dx}\mathrm{Y}(x)\right)+b_3(x)\,\mathrm{Y}(x)=\mathrm{F}(x) \tag{12.3}$$

**is defined by the following formula:**

$$\mathrm{Y}(x):=y_2(x)\int\left(\frac{d}{dx}\left(\frac{y_1(x)}{y_2(x)}\right)\right)\int\left(\frac{\mathbf{e}^{\left(-\int b_1(x)\,dx\right)}}{\left(\frac{d}{dx}\left(\frac{y_1(x)}{y_2(x)}\right)\right)^2 y_2(x)^3}\right)\cdot\left(\int \mathbf{e}^{\left(\int b_1(x)\,dx\right)}\left(\frac{d}{dx}\left(\frac{y_1(x)}{y_2(x)}\right)\right)y_2(x)^2\,\mathrm{F}(x)\,dx\right)dx\,dx \tag{12.4}$$

**Proof: Let’s obtain partial solution of the non-homogenous ODE (15.3) in the following form:**

$$y:=\mathrm{B}(x)\int \mathrm{Z}(x)\int \mathrm{V}(x)\int \mathrm{m}(x)\,\mathrm{F}(x)\,dx\,dx\,dx \tag{12.5}$$

where are the desired functions.
Substituting (15.5) to the non-homogeneous ODE (15.3), we obtain, upon transformations

$$\left(b_1(x)\left(\frac{d^2}{dx^2}B(x)\right)+b_2(x)\left(\frac{d}{dx}B(x)\right)+\left(\frac{d^3}{dx^3}B(x)\right)+b_3(x)\,B(x)\right)$$

$$\int Z(x)\int V(x)\int m(x)\,F(x)\,dx\,dx\,dx+\left(3\left(\frac{d^2}{dx^2}B(x)\right)Z(x)+B(x)\left(\frac{d^2}{dx^2}Z(x)\right)\right.$$

$$+2\,b_1(x)\left(\frac{d}{dx}B(x)\right)Z(x)+b_1(x)\,B(x)\left(\frac{d}{dx}Z(x)\right)+3\left(\frac{d}{dx}B(x)\right)\left(\frac{d}{dx}Z(x)\right)$$

$$\left.+B(x)\,Z(x)\,b_2(x)\right)\int V(x)\int m(x)\,F(x)\,dx\,dx+\left(B(x)\,Z(x)\left(\frac{d}{dx}V(x)\right)\right.$$

$$\left.+3\left(\frac{d}{dx}B(x)\right)Z(x)\,V(x)+B(x)\,Z(x)\,b_1(x)\,V(x)+2\left(\frac{d}{dx}Z(x)\right)B(x)\,V(x)\right)$$

$$\int m(x)\,F(x)\,dx+B(x)\,Z(x)\,V(x)\,m(x)\,F(x)=F(x)$$

This implies that the following equalities shall hold true:

$$\left(\frac{d^3}{dx^3}B(x)\right)+b_1(x)\left(\frac{d^2}{dx^2}B(x)\right)+b_2(x)\left(\frac{d}{dx}B(x)\right)+b_3(x)\,B(x)=0 \tag{12.6}$$

$$3\left(\frac{d^2}{dx^2}B(x)\right)Z(x)+B(x)\left(\frac{d^2}{dx^2}Z(x)\right)+2\,b_1(x)\left(\frac{d}{dx}B(x)\right)Z(x)+b_1(x)\,B(x)\left(\frac{d}{dx}Z(x)\right)$$
$$+3\left(\frac{d}{dx}B(x)\right)\left(\frac{d}{dx}Z(x)\right)+B(x)\,Z(x)\,b_2(x)=0$$

(12.7)

$$B(x)\,Z(x)\left(\frac{d}{dx}V(x)\right)+3\left(\frac{d}{dx}B(x)\right)Z(x)\,V(x)+B(x)\,Z(x)\,b_1(x)\,V(x)$$
$$+2\left(\frac{d}{dx}Z(x)\right)B(x)\,V(x)=0 \tag{12.8}$$

$$B(x)\,Z(x)\,V(x)\,m(x)=1 \tag{12.9}$$

Obviously, the $B(x)$ function is the solution of the homogeneous ODE (12.2), i.e. . In this case, equality (12.6) of this system is satisfied identically. The second equality, i.e. (12.7), is satisfied if we assume:

(12.10)

where $y_2(x)$ is the second partial solution of ODE (12.2). In this case the equality (12.7) is satisfied identically, and from (12.8) we get:

$$V(x):=\frac{C_1}{y_1(x)^3\left(\frac{d}{dx}\left(\frac{y_2(x)}{y_1(x)}\right)\right)^2\mathbf{e}^{\left(\int b_1(x)\,dx\right)}} \tag{12.11}$$

$C_1$ is an arbitrary constant.

According to the equalities (12.10) and (12.11), from (12.9) we obtain:

$$\mathrm{m}(x) := \frac{e^{\left(\int b_1(x)\,dx\right)}\left(\left(\frac{d}{dx}y_2(x)\right)y_1(x) - y_2(x)\left(\frac{d}{dx}y_1(x)\right)\right)}{C_1}$$

Substituting the specified values of the $\mathrm{B}(x), \mathrm{Z}(x), \mathrm{V}(x), \mathrm{m}(x), \mathrm{F}(x)$ functions to (12.5), we obtain (12.4).
**The theorem is proved.**
Proceeding the same way, we can get new symbolic representations for the partial solution of a non-homogeneous ODE, principally for the equations of any order.
Proceeding the same way, we can prove the following
**Theorem: A non-homogenous ODE of fourth order:**

$$\left(\frac{\partial^4}{\partial x^4}y\right) + b_1(x)\left(\frac{\partial^3}{\partial x^3}y\right) + b_2(x)\left(\frac{\partial^2}{\partial x^2}y\right) + b_3(x)\left(\frac{\partial}{\partial x}y\right) + b_4(x)\,y = \mathrm{F}(x)$$

**has the following partial solution:**

$$y = y_3\int\left(\frac{\partial}{\partial x}\left(\frac{y_2}{y_3}\right)\right)\int\left(\frac{\partial}{\partial x}\left(\frac{y_1}{y_2}\right)\right)\int\frac{e^{\left(-\int b_1\,dx\right)}\int e^{\left(\int b_1\,dx\right)}y_3^{\,3}\left(\frac{\partial}{\partial x}\left(\frac{y_1}{y_2}\right)\right)\left(\frac{\partial}{\partial x}\left(\frac{y_2}{y_3}\right)\right)^2\mathrm{F}(x)\,dx}{\left(\frac{\partial}{\partial x}\left(\frac{y_1}{y_2}\right)\right)^2\left(\frac{\partial}{\partial x}\left(\frac{y_2}{y_3}\right)\right)^3 y_3^{\,4}}\,dx\,dx\,dx$$

**where $y_1, y_2, y_3$ are partial solutions of the following homogeneous ODE:**

$$\left(\frac{\partial^4}{\partial x^4}y\right) + b_1(x)\left(\frac{\partial^3}{\partial x^3}y\right) + b_2(x)\left(\frac{\partial^2}{\partial x^2}y\right) + b_3(x)\left(\frac{\partial}{\partial x}y\right) + b_4(x)\,y = 0$$

**Theorem: A non-homogenous ODE of fifth order:**

$$\left(\frac{\partial^5}{\partial x^5}y\right) + b_1(x)\left(\frac{\partial^4}{\partial x^4}y\right) + b_2(x)\left(\frac{\partial^3}{\partial x^3}y\right) + b_3(x)\left(\frac{\partial^2}{\partial x^2}y\right) + b_4(x)\left(\frac{\partial}{\partial x}y\right) + b_5(x)\,y = \mathrm{F}(x)$$

**has the following partial solution:**

$$y = y_4 \int \left(\frac{\partial}{\partial x}\left(\frac{y_3}{y_4}\right)\right) \int \left(\frac{\partial}{\partial x}\left(\frac{y_2}{y_3}\right)\right) \int \left(\frac{\partial}{\partial x}\left(\frac{y_1}{y_2}\right)\right) \int \frac{\mathbf{e}^{\left(-\int b_1\,dx\right)} \int \mathbf{e}^{\left(\int b_1\,dx\right)} y_4^{\,4} \left(\frac{\partial}{\partial x}\left(\frac{y_1}{y_2}\right)\right)\left(\frac{\partial}{\partial x}\left(\frac{y_2}{y_3}\right)\right)^2 \left(\frac{\partial}{\partial x}\left(\frac{y_3}{y_4}\right)\right)^3 \mathrm{F}(x)\,dx}{\left(\frac{\partial}{\partial x}\left(\frac{y_1}{y_2}\right)\right)^2 \left(\frac{\partial}{\partial x}\left(\frac{y_2}{y_3}\right)\right)^3 \left(\frac{\partial}{\partial x}\left(\frac{y_3}{y_4}\right)\right)^4 y_4^{\,5}}\,dx\,dx \quad dx\,dx$$

Obviously, using the above formulas it is easy to write down a similar general formula for a non-homogeneous ODE (12.1):

**Proposition: Partial solution of the following non-homogeneous linear ODE:**

$$\left(\frac{d^m}{dx^m}\mathrm{Y}(x)\right) + \left(\sum_{s=1}^{m} b_s(x)\left(\frac{d^{m-s}}{dx^{m-s}}\mathrm{Y}(x)\right)\right) = \mathrm{F}(x)$$

**is defined by the following formula:**

$$\mathrm{Y}(x) = y_{m-1} \int \left(\frac{\partial}{\partial x}\left(\frac{y_{m-2}}{y_{m-1}}\right)\right) \int \left(\frac{\partial}{\partial x}\left(\frac{y_{m-3}}{y_{m-2}}\right)\right) \int \left(\frac{\partial}{\partial x}\left(\frac{y_{m-4}}{y_{m-3}}\right)\right) [\;] \int \left(\frac{\partial}{\partial x}\left(\frac{y_2}{y_3}\right)\right) \int \left(\frac{\partial}{\partial x}\left(\frac{y_1}{y_2}\right)\right) \int \frac{\mathbf{e}^{\left(-\int b_1(x)\,dx\right)} \int \mathbf{e}^{\left(\int b_1(x)\,dx\right)} \left(\prod_{i=0}^{m-3}\left(\frac{\partial}{\partial x}\left(\frac{y_{i+1}}{y_{i+2}}\right)\right)^{(1+i)}\right) y_{m-1}^{\,(m-1)}\,\mathrm{F}(x)\,dx}{\left(\prod_{i=0}^{m-3}\left(\frac{\partial}{\partial x}\left(\frac{y_{i+1}}{y_{i+2}}\right)\right)^{(2+i)}\right) y_{m-1}^{\,(m-0)}}\,dx\,d$$

**where  are partial solutions of the following homogeneous linear ODE:**

$$\left(\frac{\partial^m}{\partial x^m}y\right)+\left(\sum_{s=1}^{m} b_s(x)\left(\frac{\partial^{m-s}}{\partial x^{m-s}}y\right)\right)=0$$

As we can see, unlike the conventional approach, to obtain a partial solution of a non-homogeneous linear ODE we require only $m-1$ partial solutions, though in this case we have to take an $m$ -fold multiple integral.

$G_i(a)$